\documentclass[12pt]{amsart}
\usepackage{amssymb}
\usepackage[noadjust]{cite}
\usepackage{booktabs}
\usepackage{url}
\usepackage{hyphenat}
\usepackage{mathtools}
\usepackage{cancel} 
\usepackage{ifpdf}
\usepackage[T1]{fontenc}
\usepackage[utf8]{inputenc}
\usepackage{mdwtab}
\usepackage{enumitem}
\usepackage{xr-hyper}
\ifpdf
  \usepackage[pdftex]{graphicx}
  \usepackage[pdftex,margin=1in]{geometry}
  \usepackage[bookmarks=true, bookmarksopen=true,%
    bookmarksdepth=3,bookmarksopenlevel=2,%
    colorlinks=true,%
    linkcolor=blue,%
    citecolor=blue,%
    filecolor=blue,%
    menucolor=blue,%
    urlcolor=blue]{hyperref}
  \hypersetup{pdftitle={Computing HF\textasciicircum\ by factoring mapping classes}}
  \hypersetup{pdfauthor={Robert Lipshitz, Peter S. Ozsváth, and Dylan
      P. Thurston}}
\else
  \usepackage[dvips]{graphicx}
  \usepackage[dvips,margin=1in]{geometry}
  \usepackage[draft]{hyperref}
\fi
\usepackage{color}
\usepackage{tikz}

\usetikzlibrary{matrix,arrows,positioning}
\tikzstyle{every picture}=[> = latex']
\tikzset{cdlabel/.style={execute at begin node=$\scriptstyle,execute at end node=$}}
\tikzset{algarrow/.style={->, thick}}
\tikzset{blgarrow/.style={->, thick}}
\tikzset{clgarrow/.style={->, thick}}
\tikzset{tensoralgarrow/.style={double, double equal sign distance, -implies}}
\tikzset{tensorblgarrow/.style={->, thin, double}}
\tikzset{tensorclgarrow/.style={->, thin, double}}
\tikzset{modarrow/.style={->, dashed}}
\tikzset{othmodarrow/.style={->, thick}}
\tikzset{Amodar/.style={->, dashed}}
\tikzset{Dmodar/.style={->, dashed}}
\tikzset{boxed/.style={shape=rectangle,draw}}

\newcommand{\ZZ}{\mathbb Z}

\newcommand{\FF}{\mathbb F}
\newcommand{\NN}{\mathbb N}

\newcommand{\bD}{\mathbb{D}}

\newcommand{\co}{\colon}

\newcommand{\OneHalf}{{\textstyle\frac{1}{2}}}
\newcommand{\OneQuart}{{\textstyle\frac{1}{4}}}

\newcommand{\bdy}{\partial}

\newcommand{\lbracket}{[}
\newcommand{\rbracket}{]}

\newcommand{\spinc}{\mathfrak s}
\DeclareMathOperator{\divis}{div}

\DeclareMathOperator{\Hom}{Hom}
\DeclareMathOperator{\Ext}{Ext}

\DeclareMathOperator{\spin}{spin}
\newcommand{\SpinC}{\spin^c}

\DeclareMathOperator{\ind}{ind}

\DeclareMathOperator{\inv}{inv}

\DeclareMathOperator{\gr}{gr}

\DeclareMathOperator{\supp}{supp}

\DeclareMathOperator{\Aut}{Aut}
\DeclareMathOperator{\Out}{Out}
\DeclareMathOperator{\Coeff}{Coeff}
\DeclareMathOperator{\Stab}{Stab}
\DeclareMathOperator{\Isom}{Isom}

\DeclareMathOperator{\Sp}{Sp}

\theoremstyle{plain}
\newtheorem{introthm}{Theorem}
\numberwithin{equation}{section}
\newtheorem{proposition}[equation]{Proposition}
\newtheorem{citethm}[equation]{Theorem}
\newtheorem{lemma}[equation]{Lemma}
\newtheorem{corollary}[equation]{Corollary}

\newtheorem{theorem}[equation]{Theorem}

\theoremstyle{definition}
\newtheorem{definition}[equation]{Definition}

\theoremstyle{remark}

\newtheorem{remark}[equation]{Remark}
\newtheorem{convention}[equation]{Convention}

\newcommand{\HF}{\mathit{HF}}
\newcommand{\HFa}{\widehat {\HF}}

\newcommand{\CFa}{\widehat {\mathit{CF}}}

\newcommand{\x}{\mathbf x}
\newcommand{\y}{\mathbf y}

\newcommand\HH{\mathit{HH}}
\newcommand\Hochschild\HH

\newcommand{\MCG}{\mathit{MCG}}

\newcommand{\Ainf}{\mathcal A_\infty}

\newcommand{\Alg}{\mathcal{A}}

\newcommand{\Idem}{\mathcal{I}}

\newcommand\Gen{\mathfrak{S}}

\renewcommand{\S}{\Gen}

\newcommand{\Blg}{\mathcal{B}}

\newcommand{\alphas}{{\boldsymbol{\alpha}}}
\newcommand{\betas}{{\boldsymbol{\beta}}}

\newcommand{\bSigma}{\widebar{\Sigma}}
\newcommand{\balpha}{\widebar{\alpha}}
\newcommand{\balphas}{\widebar{\alphas}}
\newcommand{\cM}{\mathcal{M}}
\newcommand{\Mod}{\cM}

\newcommand{\CFD}{\mathit{CFD}}
\newcommand{\CFDD}{\mathit{CFDD}}
\newcommand{\CFA}{\mathit{CFA}}

\newcommand{\CFDA}{\mathit{CFDA}}
\newcommand{\CFDAa}{\widehat{\CFDA}}
\newcommand{\CFAA}{\mathit{CFAA}}
\newcommand{\CFAAa}{\widehat{\CFAA}}
\newcommand{\CFDa}{\widehat{\CFD}}

\newcommand{\CFDDa}{\widehat{\CFDD}}

\newcommand{\CFAa}{\widehat{\CFA}}

\newcommand{\cZ}{\mathcal{Z}}
\newcommand{\PtdMatchCirc}{\cZ}
\newcommand{\PMC}{\PtdMatchCirc}

\newcommand{\CircPts}{{\mathbf{a}}}

\newcommand\DGA{A}

\newcommand{\dg}{\textit{dg} }
\newcommand{\DD}{\textit{DD}}

\newcommand{\DA}{\textit{DA}}

\newcommand{\AAm}{\textit{AA}} 

\newcommand\Id{\mathbb{I}}
\newcommand\Ground{\mathbf k}

\newcommand\DT{\boxtimes}

\newcommand\Tensor{\mathcal T}
\newcommand\Zmod[1]{\mathbb{Z}/{#1}\mathbb{Z}}
\newcommand{\Field}{\FF_2}

\newcommand{\Heegaard}{\mathcal{H}}
\newcommand{\HD}{\Heegaard}

\newcommand{\HB}{\mathsf{H}}

\renewcommand{\th}{^\text{th}}

\newcommand{\bigGroup}{G'}
\newcommand{\smallGroup}{G}

\newcommand{\GmodTwo}{\smallGroup_{\ZZ/2}}
\newcommand{\gpElt}{\gamma}
\newcommand{\othgpElt}{\gamma'}
\newcommand{\othh}{h'}
\newcommand{\bigGrSet}{S'}

\newcommand{\grb}{\gr'}

\newcommand{\op}{{\mathrm{op}}}
\newcommand\PunctF{F^\circ}
\newcommand{\drY}{Y_{dr}}

\newcommand{\Region}{\mathcal{R}}

\newcommand{\Hyph}{\text{-}}

\DeclareMathOperator{\Mor}{Mor}

\newcommand{\sos}[3]{\mathbin{{}_{#1}\mathord#2_{#3}}}
\newcommand{\lsub}[2]{{}_{#1}#2}
\newcommand{\lsup}[2]{{}^{#1}\mskip-.6\thinmuskip#2}

\newcommand{\std}{{\mathrm{std}}}

\newcommand\Domain{Q}
\newcommand\Chain{K}
\newcommand\OthChain{L}
\newcommand\DDmod{\widehat{\mathcal{DD}}}
\newcommand\Dmod{\widehat{\mathcal D}}

\newcommand\SetS{\mathbf{s}}
\newcommand\SetT{\mathbf{t}}
\newcommand\PuncF{F^{\circ}}

\makeatletter
\newcommand\honestalg[3]{\bigl\lbracket
\begin{smallmatrix} #1\@ifempty{#3}{}{&#3} \\ #2 \end{smallmatrix}
\bigr\rbracket}

\makeatother

\newcommand{\glsit}[1]{\index{#1}}

\newread\testin

\graphicspath{{draws/}{mpdraws/}}
\makeatletter
\def\input@path{{}{draws/}}
\makeatother

\def\mathcenter#1{%
  \vcenter{\hbox{$#1$}}%
}

\DeclareRobustCommand{\widebar}[1]{\overline{#1}{}}

\makeatletter
\newcommand\mi@kern[1]{%
  \settowidth\@tempdima{$\mi@obj^{#1}$}
  \kern-\@tempdima
  #1
  \settowidth\@tempdima{$\mi@obj$}
  \kern\@tempdima
}

\newtoks\mi@toksp
\newtoks\mi@toksb
\DeclareRobustCommand{\manyindices}[5]{
  \def\mi@obj{#5}
  \mi@toksp\expandafter{\mi@kern{#2}}
  \mi@toksb\expandafter{\mi@kern{#1}}
  \@mathmeasure4\textstyle{#5_{#1}^{#2}}
  \@mathmeasure6\textstyle{#5_{#3}^{#4}}
  \dimen0-\wd6 \advance\dimen0\wd4
  \@mathmeasure8\textstyle{\hphantom{{}_{#1}^{#2}}#5^{\the\mi@toksp#4}_{\the\mi@toksb#3}}
  \hbox to \dimen0{}{\kern-\dimen0\box8}
}
\makeatother 

\usepackage{ifpdf}
\ifpdf
  \let\textalt\texorpdfstring
\else
  \newcommand{\textalt}[2]{#1}
\fi

\begin{document}
\title[Computing \smash{$\HFa$} by factoring mapping classes]
{Computing {$\HFa$} by factoring mapping classes}

\author[Lipshitz]{Robert Lipshitz}
\thanks{RL was supported by NSF grant DMS-0905796, the Mathematical
  Sciences Research Institute, and a Sloan Research
  Fellowship.}
\address{Department of Mathematics, Columbia University\\
  New York, NY 10027}
\email{lipshitz@math.columbia.edu}

\author[Ozsv\'ath]{Peter~S.~Ozsv\'ath}
\thanks{PSO was supported by NSF grant DMS-0505811, the Mathematical
  Sciences Research Institute, and a Clay
  Senior Scholar Fellowship.}
\address {Department of Mathematics, Princeton University\\ 
Princeton, NJ 08544}
\email {petero@math.princeton.edu}

\author[Thurston]{Dylan~P.~Thurston}
\thanks {DPT was supported by NSF
  grant DMS-1008049, the Mathematical
  Sciences Research Institute, and a Sloan Research Fellowship.}
\address{Department of Mathematics,
         Indiana University\\
         Bloomington, IN 47405}
\email{dpthurst@indiana.edu}

\begin{abstract}
  Bordered Heegaard Floer homology is an invariant for three-manifolds
  with boundary. In particular, this invariant associates to a handle decomposition
  of a surface~$F$ a differential graded algebra, and to an
  arc slide between two handle decompositions, a bimodule over the two
  algebras. In this paper, we describe these bimodules for arc slides
  explicitly, and
  then use them to give a combinatorial description of $\HFa$ of a
  closed three-manifold, as well as the bordered Floer homology of any
  $3$-manifold with boundary.
\end{abstract} 

\maketitle 

\newcounter{bean}

\tableofcontents
\section{Introduction}

Heegaard Floer homology is an invariant of three-manifolds defined
using Heegaard diagrams and holomorphic
disks~\cite{OS04:HolomorphicDisks}. This invariant is the homotopy
type of a chain complex over a polynomial algebra in a formal variable
$U$. The present paper will focus on $\HFa(Y)$ (with coefficients in
$\Field$), which is the homology of the $U=0$ specialization.  This
variant is simpler to work with, but it still encodes
interesting information about the underlying three-manifold $Y$ (for
instance, the Thurston norm~\cite{OS04:ThurstonNorm}).
Although the definition of $\HFa$ involves holomorphic disks, the
work of Sarkar and Wang~\cite{SarkarWang07:ComputingHFhat} allows one to calculate
$\HFa$ explicitly from a Heegaard diagram for $Y$ satisfying certain
properties. (See also~\cite{OzsvathStipsiczSzabo:Nice}.)

Bordered Heegaard Floer homology~\cite{LOT1} is an invariant for
three-manifolds with parameterized boundary. A pairing theorem
from~\cite{LOT1} allows one to reconstruct the invariant $\HFa(Y)$ of a
closed three-manifold $Y$ which is decomposed along a separating
surface $F$ in terms of the bordered invariants of the pieces.

In this paper we use bordered Floer theory to give another algorithm
to compute $\HFa(Y)$ (with $\Field$-coefficients). This algorithm is logically independent
of~\cite{SarkarWang07:ComputingHFhat}, and quite natural (both aesthetically
and mathematically, see~\cite{LOTCobordisms}). It is also practical
for computer implementation; some computer computations are described
in Section~\ref{sec:AModules}.
A Heegaard decomposition of $Y$ is determined by an automorphism
$\phi$ of the
Heegaard surface, and the complex we describe here is
associated to a suitable factorization of $\phi$. To explain
this in slightly more detail, we recall some of the basics
of the bordered Heegaard Floer homology package.

We represent oriented surfaces by {\em pointed
  matched circles}, which are essentially handle decompositions with a little extra
structure; see~\cite{LOT1} or the review in
Section~\ref{subsec:AlgPMC} below. To a pointed matched circle $\PMC$,
bordered Heegaard Floer theory associates a differential graded algebra
$\Alg(\PMC)$. 
A {\em $\PMC$-bordered $3$-manifold}
\index{$3$-manifold!$\PMC$-bordered}%
is three-manifold $Y_0$ equipped with an orientation-preserving
identification of its boundary $\bdy Y_0$ with the surface associated
to the pointed
matched circle $\PMC$. To a $\PMC$-bordered $3$-manifold, bordered
Heegaard Floer theory associates modules over the algebra
$\Alg(\PMC)$. Specifically, if $Y_1$ is a $\PMC$-bordered
three-manifold, there is an associated module
$\CFAa(Y_1)$, which is a right $\Ainf$-module over
$\Alg(\PMC)$. Similarly, if $Y_2$ is a $(-\PMC)$-bordered three-manifold, 
we obtain a different module $\CFDa(Y_2)$ which is
a left differential  module over $\Alg(\PMC)$. A key property of the bordered
invariants~\cite{LOT1} states that the Heegaard Floer complex
$\CFa(Y)$ of the closed three-manifold obtained by gluing $Y_1$ and
$Y_2$ along the above identifications is calculated by the (derived)
tensor product of $\CFAa(Y_1)$ and $\CFDa(Y_2)$; we call results of
this sort ``pairing theorems''.
\index{pairing theorem}%

In~\cite{LOT2} we construct bimodules which can
be used to change the boundary parameterizations. These are special
cases of bimodules associated to three-manifolds with two boundary
components.  In that paper, a duality property is also established,
which allows us to formulate all aspects of the theory purely in terms
of $\CFDa$.  This has two advantages: $\CFDa$ involves fewer (and simpler)
holomorphic curves than $\CFAa$ does, and $\CFDa$ is always an
honest differential module, rather than a more general $\Ainf$ module
(like $\CFAa$). This point of view is further pursued
in~\cite{LOTHomPair}.  Suppose that a closed
three-manifold $Y$ is separated into two pieces $Y_1$ and $Y_2$ along a
separating surface $F$. It is shown in~\cite{LOTHomPair} that
$\HFa(Y)$ is obtained as the homology of the chain
complex of maps from the bimodule associated to the identity cobordism
of $F$ into $\CFDa(Y_1)\otimes_{\Field} \CFDa(Y_2)$; this result is restated as
Theorem~\ref{thm:PairingTheorem} below.  (It is also shown
in~\cite{LOTHomPair} that the Heegaard Floer invariant for
$Y=Y_1\cup_\bdy Y_2$ is identified via a ``Hom pairing theorem'' with
$\Ext_{\Alg(\PMC)}(\CFDa(-Y_1),\CFDa(Y_2))$, the homology of the chain
complex of homomorphisms from $\CFDa(-Y_1)$ to $\CFDa(Y_2)$. We will
usually work with the other formulation in order to avoid orientation
reversals.)

A Heegaard decomposition of a closed three-manifold $Y$ is a
decomposition of $Y$ along a surface into two particular simple
pieces: handlebodies. The invariants of handlebodies with suitable boundary parametrizations are easy to calculate. Thus, the key step to
calculating $\HFa$ from a Heegaard splitting is calculating the bimodule
associated to a surface automorphism~$\psi$ to allow us to match up
the two boundary parametrizations. We approach this problem as
follows.  Any diffeomorphism $\psi$ between two surfaces
associated to two (possibly different) pointed matched circles can be
factored into a sequence of elementary pieces, called {\em
  arc-slides}.  Then Theorem~\ref{thm:PairingTheorem} allows us to
compute the bimodule associated to
$\psi$ as a suitable composition of bimodules associated
to arc-slides.

Thus, the primary task of the present paper is to calculate the type
\DD\  bimodule associated to an arc-slide. These calculations
turn out to follow from simple geometric constraints coming
from the Heegaard diagram, combined with algebraic constraints imposed
from the relation that $\partial^2=0$.

For most of this paper, we work entirely within the context of
``type $D$'' invariants. This allows us to avoid much of the algebraic
complication of $\Ainf$-modules which are built into the bordered
theory: the bordered invariants we consider are simply
differential modules over a differential algebra. (In practice, it may
be convenient to take homology to make our
complexes smaller; the cost of doing this is to work with $\Ainf$
modules. We return to this point in Section~\ref{sec:AModules}.)

We now turn to an explicit description of all the ingredients
in our calculation of $\HFa(Y)$.

\subsection{Algebras for pointed matched circles}
\label{subsec:AlgPMC}

\glsit{$\PMC$}\glsit{$Z$}\glsit{$\mathbf{a}$}\glsit{$[4k]$}\glsit{$z$}\glsit{$M=M_\PMC$}%
As defined in~\cite{LOT1}, \emph{a pointed matched circle} $\PMC$
\index{pointed matched circle}%
is an oriented circle $Z$, equipped
with a basepoint $z\in Z$, and $4k$ basepoints ${\mathbf a}$, which we
label in order
$1,\dots,4k$, and which are matched in pairs. This matching is encoded
in a two-to-one function $M=M_\PMC\co 1,\dots,4k \to 1,\dots,2k$.
We sometimes denote $\mathbf{a}$ by $[4k]$, and
write $[4k]/M$ for the range of the matching $M$.  This matching is
further required to satisfy the following combinatorial property: surgering
out the $2k$ pairs of matched points (thought of as embedded
zero-spheres in $Z$) results in a  one-manifold which is connected.

A pointed matched circle specifies a surface $\PuncF(\PMC)$ with a
single boundary component by filling in $Z$ with a disk and then
attaching $1$-handles along the pairs of matched points in
$\mathbf{a}$. Capping off the boundary component with a disk, we
obtain a compact, oriented surface $F(\PMC)$.
\index{pointed matched circle!surface associated to}%
\glsit{$\PunctF(\PMC)$}\glsit{$F(\PMC)$}%

Given $\PMC$, there is an associated differential graded algebra
$\Alg(\PMC)$, introduced in~\cite{LOT1}. We briefly recall the
construction here, and describe some more of its properties in
Section~\ref{subsec:AlgebraPMC}; for more, see~\cite{LOT1}.

\glsit{$\Alg(\PMC)$}%
\index{pointed matched circle!algebra associated to}%
\index{algebra associated to a pointed matched circle}%
The algebra $\Alg(\PMC)$ is generated as an $\Field$-vector space by strands
diagrams $\upsilon$, defined as follows. First, cut the circle
along $z$ so that it is an interval $I_Z$. Strands diagrams are collections of
non-decreasing, linear paths in $[0,1]\times I_Z$, each of which starts
and ends in one of the distinguished points in the matched circle. The
strands are required to satisfy the following properties:
\index{strands diagram}%
\begin{itemize}
  \item Horizontal strands come in matching pairs: i.e.,
    if $i\neq j$ and $M(i)=M(j)$, then 
    $\upsilon$ contains a horizontal strand starting at $i$ if
    and only if it contains the corresponding strand at $j$.
    \index{horizontal strand}%
    \index{strand!horizontal}%
  \item If there is a non-horizontal strand in $\upsilon$ starting at
    some point $i$, then there is no other strand in $\upsilon$
    starting at either $j$ with $M(j)=M(i)$.
  \item If there is a  non-horizontal strand in $\upsilon$ ending at
    some point $i$, then there is no other strand in $\upsilon$
    ending at either $j$ with $M(j)=M(i)$.
\end{itemize}
When we think of a strands diagram as an element
of the algebra, we call the corresponding element a {\em basic
  generator}.
\index{basic generator}%

\begin{figure}
  \centering
  \includegraphics{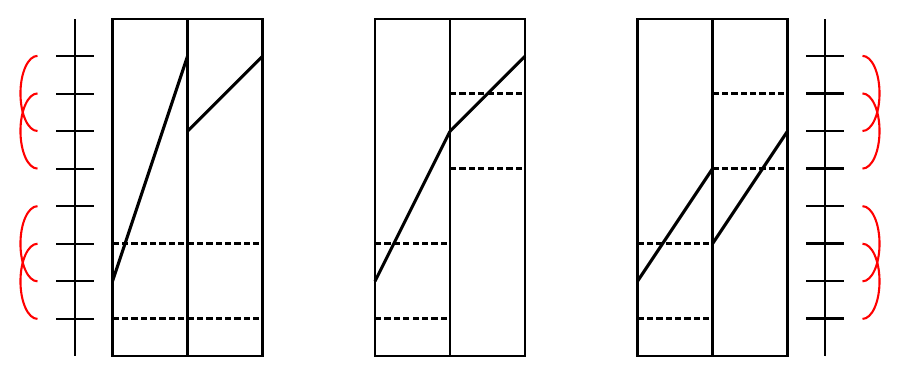}
  \caption{\textbf{Vanishing products on $\Alg$.} The five cases in which the multiplication on $\Alg$ vanishes are illustrated. The product on the left vanishes for both of the first two reasons; the product in the center for both the third and fourth reasons, and the product on the right for the last reason. Horizontal strands are drawn dashed, to illustrate that they have weight $1/2$. All pictures are in $\Alg(\PMC,0)$ where $\PMC$ represents a surface of genus $2$.}
  \label{fig:VanishingProds}
\end{figure}

Multiplication in this algebra is defined in terms of basic
generators.  Given two strands diagrams, their product is gotten by
concatenating the diagrams, when possible, throwing out one of any
given pair of horizontal strands if necessary, and then homotoping straight
the piecewise linear juxtaposed paths (while fixing endpoints), and
declaring the result to be $0$ if this homotopy decreases the total
number of crossings. More precisely, suppose
$\upsilon$ and $\tau$ are two strands diagrams.  We declare the
product to be zero if any of the following conditions is satisfied:
\begin{enumerate}
  \item $\upsilon$ contains a non-horizontal strand whose terminal
    point is not the initial point of any strand in $\tau$ or
  \item $\tau$ contains a non-horizontal strand whose 
    initial point is not the terminal point of any strand in $\upsilon$ or
  \item $\upsilon$ contains a horizontal strand which has the property
    that both it and its matching strand have terminal points which
    are not initial points of strands in $\tau$ or
  \item $\tau$ contains a horizontal strand which has the property
    that both it and its matching strand have initial points which
    are not terminal points of strands in $\upsilon$ or
  \item the concatenation $\upsilon*\tau$ of $\upsilon$ and $\tau$
    contains two piecewise 
    linear strands which cross each other twice.
\end{enumerate}
(See Figure~\ref{fig:VanishingProds} for an illustration.)
Otherwise, we take the resulting diagram $\upsilon*\tau$, remove any
horizontal strands which do not go all the way across, and then pull
all piecewise linear strands straight (fixing the endpoints). 

The differential of a strands diagram is a sum of terms, one for each
crossing.  The term corresponding to a crossing $c$ is gotten by
forming the upward resolution at $c$ (i.e., if two strands meet at $c$,
we replace them by a nearby approximation by two non-decreasing paths which 
do not cross at $c$, in such a manner that the two initial points and two
terminal points are the same). If this resolved diagram has a double-crossing,
we set it equal to zero. Otherwise, once again, the corresponding term is gotten
by pulling the strands straight and dropping any horizontal strand
whose mate is no longer present. We denote the differential on these
algebras by $d$; differentials on modules will usually be denoted
$\bdy$.
\glsit{$d$}%

The product and differential endow $\Alg(\PMC)$ with the structure of
a differential algebra; i.e., $d(a\cdot b)=(da)\cdot b + a\cdot (db)$.

\index{weight of strands diagram}%
\index{strands diagram!weight}%
A strands diagram $\upsilon$ has a total weight, which is gotten by
counting each non-horizontal strand with weight $1$, each horizontal
strand with weight $1/2$, and then subtracting the genus $k$. Let
$\Alg(\PMC,i)\subset\Alg(\PMC)$ be the subalgebra generated by weight
$i$ strands diagrams. This, of course, is a differential subalgebra.

\glsit{$\SetS$}\glsit{$I(\SetS)$}%
Note also that for each subset $\SetS$ of $[4k]/M$, there is a
corresponding idempotent $I(\SetS)$, consisting of the collection of
horizontal strands $[0,1]\times M^{-1}(\SetS)$. These are the
minimal idempotents of $\Alg(\PMC)$, with respect to the partial order
$I_1\leq I_2$ if $I_1I_2=I_1$.

A strands diagram $\upsilon$ also has an underlying one-chain in
$H_1(Z,\mathbf{a})$, which we
denote in this paper by $\supp(\upsilon)$ and called the \emph{support
of $\upsilon$}. At any position $q$
between two consecutive marked points $p_i$ and $p_{i+1}$ in $Z$, the
local multiplicity of $\supp(\upsilon)$ is the intersection number of
$\upsilon$ with $[0,1]\times q$.

\index{chord}%
A {\em chord} is an interval $[i,j]$ connecting two elements in
$[4k]$. A chord $\xi$ determines an algebra element $a(\xi)$, which is
represented by the sum of all strands diagrams in which the strand
from $i$ to $j$ is the only non-horizontal strand. We denote
the set of chords by ${\mathcal C}(\PMC)$, or simply ${\mathcal C}$.
\glsit{$a(\xi)$}\glsit{$\mathcal{C}=\mathcal{C}(\PMC)$}%

\subsection{The identity Type \DD\  bimodule}
\label{sec:intro:DDId}
Before introducing the bimodules for arc-slides, we first describe a
simpler bimodule $\DDmod(\Id_{\PMC})$ associated, in a suitable sense, to the identity
map.  Motivation for calculating
this invariant comes from its prominent role in one version of the
pairing theorem, quoted as Theorem~\ref{thm:PairingTheorem}, below.

\begin{definition}
  \label{def:ComplementaryId}
  Let $\SetS,\SetT\subset [4k]/M_{\PMC}$ be subsets
  with the property that $\SetS$ and $\SetT$ form a partition 
  of $[4k]/M_{\PMC}$. Then we say that the
  corresponding idempotents
  $I(\SetS)$ and $I(\SetT)$ are {\em complementary idempotents}.
  \index{idempotents!complementary}%
  \index{complementary idempotents}%
\end{definition}

Our bimodules have the following special form:

\begin{definition}
  \label{def:DDstruct}
  Let $\Alg$ and $\Blg$ be two \dg algebras.  A {\em \DD\
    bimodule} over $\Alg$ and $\Blg$ is a \dg bimodule $M$ which, as
  a bimodule, splits into summands isomorphic to $\Alg i \otimes j\Blg$
  for various choices of idempotent $i$ and $j$ 
  (but the differential need not respect this splitting). (See also
  Section~\ref{subsubsec:TypeDStructures}.)
  \index{bimodule!type \DD}%
  \index{DD bimodule}%
\end{definition}

\begin{definition}
  \label{def:DDmodId}
  The module $\DDmod(\Id_{\PMC})$ is generated by all pairs of
  \glsit{$\DDmod(\Id_{\PMC})$}%
  complementary idempotents.  This means that its elements are of the
  form $r i\otimes i' s$, where $r,s\in\Alg(\PMC)$,
  and $i$ and $i'$ are complementary idempotents. The differential on
  $\DDmod(\Id_{\PMC})$ is determined by the Leibniz rule and the fact that
  \[
  \bdy(i\otimes i')=\sum_{\xi\in\mathcal{C}} i a(\xi)\otimes a(\xi)i'.
  \]
  (Here and later, the symbol $\otimes$ denotes tensor product over
  $\Field$, unless otherwise specified.)
  \glsit{$\otimes$}%
\end{definition}
In particular, the differential on $\DDmod(\Id_{\PMC})$ is determined
by an element 
\[A=\sum_{\xi\in {\mathcal C}} a(\xi)\otimes a(\xi)\in
\Alg(\PMC)\otimes\Alg(\PMC).\]
\glsit{$A$}%
If we let $*$ denote the action of $\Alg(\PMC)\otimes\Alg(\PMC)$ on
itself by multiplication on the outside then
the fact that $\DDmod(\Id_{\PMC})$, as defined above, is a
chain complex is equivalent to the fact that 
$dA+A*A=0$. (See Proposition~\ref{prop:DDsquaredZero} for an algebraic verification that $\DDmod(\Id_{\PMC})$ is,
indeed, a chain complex.)
The relevance of $\DDmod(\Id_{\PMC})$ to bordered Floer homology
arises from the following:

\begin{introthm}
  \label{thm:DDforIdentity}
  The bimodule $\DDmod(\Id_{\PMC})$ is canonically homotopy equivalent to the type
  \DD\ bimodule of the identity map defined using pseudoholomorphic curves
  in~\cite{LOT2}.
\end{introthm}
More precisely, in the notation of~\cite{LOT2},
$$\DDmod(\Id_{\PMC})\cong (\Alg(\PMC)\otimes\Alg(-\PMC))\DT \CFDDa(\Id_{\PMC}),$$
where we identify right actions by $\Alg(\PMC)$ with left actions by
$\Alg(-\PMC)=\Alg(\PMC)^{\op}$.

\subsection{\DD\ bimodules for arc-slides}
\label{sec:dd-arc-slide-intro}
We turn now to bimodules for arc-slides. Before doing this, we recall
briefly the notion of an arc-slide, and introduce some notation.

\begin{figure}
\begin{center}
\input{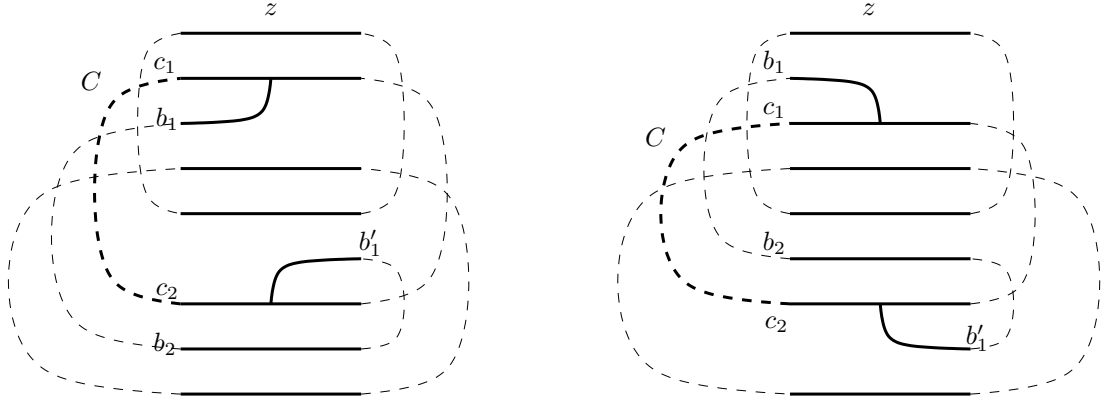}
\end{center}
\caption {{\bf Arc-slides.}
\label{fig:ArcslideMatching}
Two examples of arc-slides connecting pointed matched circles
for genus $2$ surfaces. In both cases, the foot $b_1$ is sliding
over the matched pair $C=\{c_1,c_2\}$ (indicated by the darker dotted
matching) at $c_1$.}
\end{figure}

Let $\PMC$ be a pointed matched circle, and fix two matched pairs
$C=\{c_1,c_2\}$ and $B=\{b_1,b_2\}$. Suppose moreover that $b_1$ and
$c_1$ are adjacent, in the sense that there is an arc $\sigma$ connecting $b_1$
and $c_1$ which does not contain the basepoint $z$ or any other
point $p_i\in\mathbf{a}$.
\glsit{$B$}\glsit{$C$}\glsit{$b_1$, $b_2$}\glsit{$c_1$, $c_2$}%
Then we can form a new pointed matched circle $\PMC'$ which agrees
everywhere with $\PMC$, except that $b_1$ is replaced by a new
distinguished point $b_1'$, which now is adjacent to $c_2$ and $b_1'$ is positioned so that
the orientation on the arc from $b_1$ to $c_1$ is opposite to the
orientation of the arc from $b_1'$ to~$c_2$. In this case, we say that
$\PMC'$ and $\PMC$ differ by an {\em arc-slide $m$ of $b_1$ over $C$
  at $c_1$} or, more succinctly, an \emph{arc-slide of $b_1$ over
  $c_1$}, and write $m\co \PMC
\rightarrow \PMC'$. Let $\sigma'$ denote the
arc connecting $c_2$ and $b_1'$. See Figure~\ref{fig:ArcslideMatching} for two
examples.
\index{arc-slide}%

\begin{figure}
\begin{center}
\input{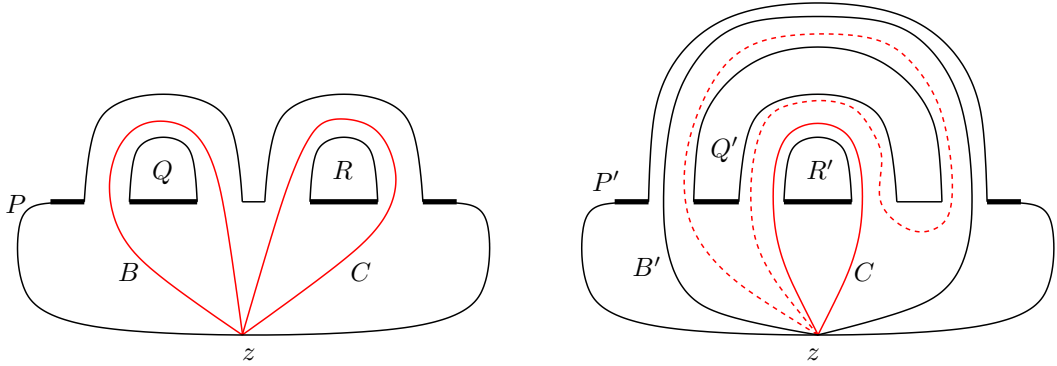}
\end{center}
\caption {{\bf The arc-slide diffeomorphism.}
\label{fig:ArcSlide}
The local case of an arc-slide diffeomorphism. Left: a
pair of pants with boundary components labeled $P$, $Q$, and
$R$, and two distinguished curves $B$ and $C$. Right:
another pair of pants with boundary components $P'$, $Q'$, $R'$ and
distinguished curves $B'$ and $C$. The arc-slide diffeomorphism
carries $B$ to the dotted curve on the right, the curve labeled
$C$ on the left to the curve labeled $C$ on the right, and boundary
components $P$, $Q$, and $R$ to $P'$, $Q'$, and $R'$ respectively.
This diffeomorphism can be extended to a diffeomorphism between
surfaces associated to pointed matched circles: in such a surface
there are further handles attached along the four dark
intervals; however, our diffeomorphism carries the four dark intervals
on the left to the four dark intervals on the right and hence extends
to a diffeomorphism as stated. (This is only one of several possible
configurations of $B$ and $C$: they could also be nested or linked.)}
\end{figure}

Note that if $\PMC$ and $\PMC'$ differ by an
arc-slide, then there is a canonical diffeomorphism from
$\PuncF(\PMC)$ to $\PuncF(\PMC')$; see Figure~\ref{fig:ArcSlide}.
\index{arc-slide!diffeomorphism associated to}%
\index{diffeomorphism associated to arc-slide}%
We will denote this diffeomorphism $\PuncF(m)$.
\glsit{$\PunctF(m)$}%

Let $m$ be an arc-slide taking the pointed matched circle $\PMC$ to
the pointed matched circle $\PMC'$.  Our next goal is to describe an
$\Alg(\PMC)\Hyph\Alg(\PMC')$-bimodule associated to $m$, which
we denote $\lsub{\Alg(\PMC)}\DDmod(m\co\PMC\rightarrow
\PMC')_{\Alg(\PMC')}$, or just $\DDmod(m)$.

To describe the generators of $\DDmod(m\co\PMC\rightarrow\PMC')$ we
need two extensions of the notion of complementary idempotents to the
case of arc-slides.

\begin{definition}
  \label{def:ComplementarySlide}
  Let $\SetS\subset [4k]/M_{\PMC}$ and $\SetT\subset [4k]/M'_{\PMC'}$
  be subsets with the property that $\SetS$ and $\SetT$ form a
  partition of $[4k]/{M_{\PMC'}}$ (where we have suppressed the identification
  between the matched pairs $[4k]/M_{\PMC}$ and $[4k]/M_{\PMC'}$).
  We say that the corresponding idempotents $I(\SetS)$ and $I(\SetT)$
  in $\Alg(\PMC)$ and $\Alg(\PMC')$ are {\em complementary idempotents}.
  An idempotent in $\Alg(\PMC)\otimes \Alg(\PMC')$ of the form
  $i\otimes i'$, where $i$ and $i'$ are complementary idempotents,
  is also called an {\em idempotent of type $X$}.
  \index{idempotents!complementary}%
  \index{idempotent!of type $X$}%
  \index{complementary idempotents}%
\end{definition}

In a similar vein, we have the following:

\begin{definition}
  \label{def:SubComplementary}
  Two elementary idempotents $i$ of $\Alg(\PMC)$ and $i'$ of $\Alg(\PMC')$
  are \emph{sub\hyp complementary idempotents} if $i=I(\SetS)$ and
  $i'=I(\SetT)$ where $\SetS\cap \SetT$ consists of the matched
  pair of the feet of $C$, while $\SetS\cup\SetT$ contains all the
  matched pairs, except for the pair of feet of $B$.  
  An idempotent in $\Alg(\PMC)\otimes\Alg(\PMC')$ of the
  form $i\otimes i'$
  where $i$ and $i'$ are sub-complementary idempotents is also called an 
  {\em idempotent of type $Y$}. Two elementary idempotents $i$ of $\Alg(\PMC)$ and
  $i'$ of $\Alg(\PMC')$ are said to be {\em near-complementary} if they are either
  complementary or sub-complementary.
  \index{idempotents!sub-complementary}%
  \index{idempotent!of type $Y$}%
  \index{sub-complementary idempotents}%
  \index{idempotents!near-complementary}%
  \index{near-complementary idempotents}%
\end{definition}

Given a chord $\xi$ for $\PMC$, let $a(\xi)$ be the
algebra element in $\Alg(\PMC)$ associated to $\xi$. Similarly, given
a chord $\xi$ for $\PMC'$, let $a'(\xi)$ be the algebra element in
$\Alg(\PMC')$ associated to the chord $\xi$.

\begin{definition}\label{def:RestrictedSupport}
  The {\em restricted support} $\supp_R(a)$ of a basic generator $a\in
  \Alg(\PMC)$ is the image of $\supp(a)$ under the map
  $H_1(Z,\CircPts)\to H_1(Z,\CircPts\setminus b_1)$ gotten by
  contracting $\sigma$ to a point. In other words, the restricted
  support of $a$ is the collection of local multiplicities of the
  associated one-chain at all the regions except $\sigma$. Similarly,
  if $a\in\Alg(\PMC')$, then the restricted support $\supp_R(a)$ of
  $a$ is the image of $\supp(a)$ under the map $H_1(Z,\CircPts')\to
  H_1(Z,\CircPts'\setminus b_1')$
  gotten by contracting $\sigma'$ to a
  point.
  \glsit{$\supp_R$}%
  \index{support!restricted}%
  \index{restricted!support}%

  A \emph{short near-chord} is a non-zero algebra element of the form
  $(i\cdot a \cdot j)\otimes (j'\cdot b'\cdot i')$ with the following
  four properties:
  \begin{enumerate}
    \item the pairs $(i\otimes i')$ and $(j\otimes j')$ are near-compementary idempotents;\label{item:near-comp}
    \item $\supp_R(a)=\supp_R(b')$;
    \item the support of at least one of $a$ or $b'$ is non-zero; and
    \item the lengths of the (unrestricted) support of $a$ and the (unrestricted) support of $b'$ are both no greater than $1$.
  \end{enumerate}
  (In a particular degenerate case, we also allow one more kind of
  short near-chord; see Definition~\ref{def:short-near-chord}.)%
\end{definition}

The above definition of short near-chords includes elements of the form
$(i\cdot a(\sigma)\cdot j) \otimes j'$ where both $(i,j')$ and $(j,j')$
are near-complementary idempotents; and also $i\otimes (j'\cdot
a'(\sigma')\cdot i')$ where both $(i,i')$ and $(i,j')$ are
near-complementary. (Note that an element $a\otimes b$ with
$\supp(a)=\sigma$ and $\supp(b)=\sigma'$ does not satisfy Property~(\ref{item:near-comp}) of Definition~\ref{def:RestrictedSupport}.)

\begin{definition}
  \label{def:Arc-SlideBimodule}
  Let $m\co \PMC\to \PMC'$ be an arc-slide.
  Let ${N}$ be any type \DD\  bimodule over $\Alg(\PMC)$ and $\Alg(\PMC')$.
  Suppose ${N}$ satisfies the following properties:
  \begin{enumerate}[label=(AS-\arabic*),ref=AS-\arabic*]
    \item 
      \label{AS:Generators}
      As an $\Alg(\PMC)\Hyph\Alg(\PMC')$-bimodule, $N$ has the form
      \[
      N=\bigoplus_{i\otimes i'\text{
          near-complementary}}\Alg(\PMC)i\otimes i'\Alg(\PMC').
      \]
    \item 
      \label{AS:NearDiagonalSubalgebra}
        For each generator $I=i\otimes i'$ of $N$ the differential of
      $I$ has the form
      \begin{equation}
        \bdy(I)=\sum_{J=j\otimes j'}\sum_k (i \cdot \nu_k \cdot j)\otimes 
        (j'\cdot \nu'_k\cdot i')\label{eq:nuk}
      \end{equation}
      where the $\nu_k$ and $\nu'_k$ are strand diagrams with the same
      restricted support (and $k$ ranges over some index set)
      and $J$ runs through the generators of $N$. (For the $\nu_k$ as
      in Formula~\eqref{eq:nuk}, we will say that
      the differential on $N$ \emph{contains}
      \index{contains}\index{differential contains}%
      $(i\cdot \nu_k\cdot j)\otimes(j'\cdot \nu_k'\cdot i')$.)
    \item 
      \label{AS:Graded}
      ${N}$ is graded (see Section~\ref{subsec:AlgebraPMC} below)
      by a $\lambda$-free grading set~$S$ (see
      Definition~\ref{def:lambdaFree}).
    \item \label{AS:ShortChords}
      All short near-chords appear in the differential; i.e. 
      given a generator $I=i\otimes i'$ of $N$,
      the differential of $I$ contains all short near-chords 
      of the form
      $(i \cdot \nu \cdot j)\otimes (j'\cdot \nu'
      \cdot i')$.
  \end{enumerate}
  Then we say that $N$ is an {\em arc-slide bimodule} for $m$.
  \index{arc-slide bimodule}%
  \index{bimodule!arc-slide b.|see{arc-slide bimodule}}%
\end{definition}

\glsit{$\#$}\glsit{$q$}%
\index{connected sum!of pointed matched circles}%
Let $\PMC$ and $\PMC_0$ be pointed matched circles. We can form their
connected sum $\PMC\#\PMC_0$.
Given any idempotent $I(\SetS_0)$ of $\Alg(\PMC_0)$, we have a quotient map
$$q\co \Alg(\PMC\#\PMC_0)\rightarrow\Alg(\PMC).$$
(For more on this, see Subsection~\ref{subsec:QuotientMap}.)

\begin{definition}
  \label{def:StableModule}
  \index{stable!arc-slide bimodule|see{arc-slide bimodule, stable}}%
  \index{arc-slide bimodule!stable}%
  We say that an arc-slide bimodule $N$ for $m\co \PMC\rightarrow \PMC'$
  is {\em stable} if for any other $\PMC_0$ and idempotent in $\Alg(\PMC_0)$, 
  and either choice of connect sum $\PMC\#\PMC_0$,
  there is an arc-slide bimodule $M$ for 
  $m_0\co \PMC\#\PMC_0\rightarrow \PMC'\#\PMC_0$
  with the property that 
  $N=q_*(M)$, where $q_*$ denotes induction of bimodules, i.e.,
  $q_*(M)=\Alg(\PMC)\otimes_{\Alg(\PMC\#\PMC_0)}M\otimes_{\Alg(\PMC'\#\PMC_0)}\Alg(\PMC')$.
\index{induction}%
\end{definition}

\begin{remark}
  \label{rmk:Stability}
  In fact, stability is much weaker than it might appear from the
  above definition.  From the proof of
  Proposition~\ref{prop:UniqueArc-Slide}, one can see that $N$ is
  stable if there exists some pointed matched circle $\PMC_0$
  of genus greater than one and a
  single associated idempotent $I$ in $\Alg(\PMC_0)$ with weight zero, 
  so that for both choices of connected sum, there are arc-slide bimodules $M$
  and $N$ as in Definition~\ref{def:StableModule} with $N=q_*(M)$.
\end{remark}

The following is proved in Section~\ref{sec:Arc-Slides}. 

\begin{proposition}
  \label{prop:UniqueArc-Slide}
  Let $m\co \PMC\to \PMC'$ be an arc-slide.
  Then, up to isomorphism, there is a
  unique stable type \DD\  arc-slide bimodule for~$m$
  (as defined in
  Definitions~\ref{def:Arc-SlideBimodule} and~\ref{def:StableModule}).
  \index{arc-slide bimodule!stable!is unique}%
\end{proposition}
The proof is constructive: after making some explicit choices, the coefficients
in the differential of the arc-slide bimodule are uniquely determined.

\begin{definition}
  \glsit{$\DDmod(m\co\PMC\rightarrow \PMC')$}%
  Let $\DDmod(m\co\PMC\rightarrow \PMC')$ be the arc-slide
  bimodule for $m$.
\end{definition}

In~\cite{LOT2}, it is shown that for any mapping class $\phi\co
\PuncF(\PMC)\rightarrow \PuncF(\PMC')$ which fixes the boundary,
there is an associated type \DD\ bimodule $\CFDDa(\phi)$. Given an
arc-slide $m\co \PMC\to\PMC'$, let
$\CFDDa(\PunctF(m))$ denote this construction, applied to the canonical
diffeomorphism $\PunctF(m)\co \PuncF(\PMC) \to \PuncF(\PMC')$ specified by $m$.

\begin{introthm}
  \label{thm:DDforArc-Slides}
  The bimodule $\DDmod(m\co\PMC\rightarrow\PMC')$ is canonically homotopy
  equivalent to the type \DD\ bimodule $\CFDDa(\PunctF(m))$ associated in~\cite{LOT2} to the
  arc-slide diffeomorphism from $\PuncF(\PMC)$ to $\PuncF(\PMC')$.
\end{introthm}
More precisely, in the notation of~\cite{LOT2}, if $\DDmod(m)$ is an arc-slide bimodule 
given in Proposition~\ref{prop:UniqueArc-Slide},  then (up to
homotopy) there is a canonical homotopy equivalence
$$\DDmod(m)\simeq (\Alg(\PMC)\otimes\Alg(-\PMC'))\DT\CFDDa(\PunctF(m)).$$
where we identify right actions by $\Alg(\PMC')$ with left actions by
$\Alg(-\PMC')=\Alg(\PMC')^{\op}$. 

\subsection{Modules associated to a handlebody}
\label{subsec:Handlebodies}

We now describe the modules associated to a handlebody. First we
consider the case of a handlebody with a standard framing, and then we
show how the arc-slide bimodules can be used to change the framing.

\subsubsection{The \textalt{$0$}{0}-framed handlebody}\label{sec:0-fr-hb}
We start by fixing some notation.
Let $\PMC^1$ denote the unique genus $1$ pointed matched
circle. $\PMC^1$ consists of an oriented
circle $Z$ equipped with a basepoint $z$ and two pairs $\{a,a'\}$ and
$\{b,b'\}$ of matched points. As we travel along $Z$ in the positive
direction starting at $z$ we encounter the points $a,b,a',b'$ in
that order. Note that the pair $\{a,a'\}$ specifies a simple closed
curve on $F(\PMC^1)$, as does the pair $\{b,b'\}$.
\glsit{$\PMC^1$}%

\begin{figure}
  \includegraphics{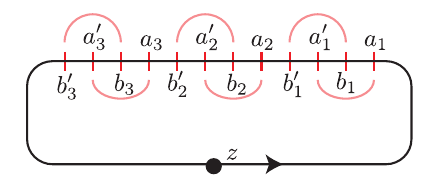}
\caption{{\bf Split pointed matched circle.}
\label{fig:StandardPMC}
The genus $3$ case is illustrated.}
\end{figure}

\glsit{$\PMC^g_0$}%
Let $\PMC^g_0=\#^g\PMC^1$ be the \emph{split pointed matched circle}
\index{split pointed matched circle}\index{pointed matched circle!split}%
describing a
surface of genus $g$, which is obtained by taking the connect sum
of $g$ copies of $\PMC^1$. The circle $\PMC^g_0$ has $4g$ marked points,
which we label in order $a_1,b_1,a_1',b_1',a_2,\dots, b_g'$, as well
as a basepoint $z$. See Figure~\ref{fig:StandardPMC}.

\begin{figure}
\begin{center}
\input{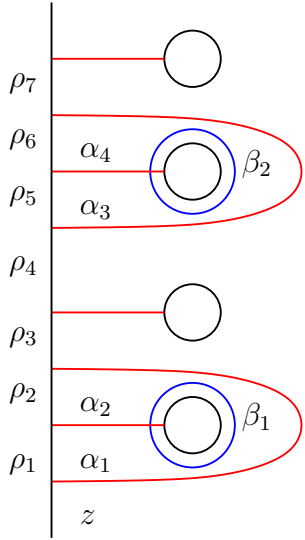}
\end{center}
\caption {{\bf Heegaard diagram for the $0$-framed genus two handlebody.}
\label{fig:GenusTwoBorderedDiagram}}
\end{figure}

The \emph{$0$-framed solid torus}
\index{solid torus!$0$-framed}\index{$0$-framed!solid torus}%
\glsit{$\HB^1=(H^1,\phi^1_0)$}%
$\HB^1=(H^1,\phi^1_0)$ is the solid torus with boundary
$-{F(\PMC^1)}$ in which $\{a,a'\}$ bounds a disk; let $\phi^1_0$
denote the preferred diffeomorphism $-{F(\PMC^1)}\to
\bdy H^1$.
The \emph{$0$-framed handlebody of genus $g$}
\index{$0$-framed!handlebody of genus $g$}\index{handlebody of genus $g$!$0$-framed}%
\glsit{$\HB^g=(H^g,\phi^g_0)$}%
$\HB^g=(H^g,\phi^g_0)$ is a
boundary connect sum of $g$ copies of $\HB^1$. Our conventions are
illustrated by the bordered Heegaard diagrams in Figure~\ref{fig:GenusTwoBorderedDiagram}.

\glsit{$\Dmod(\HB^g)$}%
We next give a combinatorial model $\Dmod(\HB^g)$ for the type $D$ module
$\CFDa(\HB^g)$ associated to $\HB^g$.
Let $\mathbf{s}=\{a_i,a_i'\}_{i=1}^g$. The module $\Dmod(\HB^g)$ is generated over the
algebra by a single idempotent $I=I(\mathbf{s})$,
and equipped with the differential determined by 
$$\bdy(I)=\sum_{i=1}^g a(\xi_i)I,$$ 
where $\xi_i$ is the arc in $\PMC^g$ connecting $a_i$ and $a'_i$.

A straightforward calculation (see Section~\ref{sec:CompleteProof}) shows:
\begin{proposition}\label{prop:CFD-of-0-fr-hdlbdy}
  The module $\Dmod(\HB^g)$ is homotopy equivalent to the module
  $\CFDa(\HB^g)=\CFDa(H^g,\phi^g_0)$ as defined (via holomorphic
  curves) in~\cite{LOT1}.
\end{proposition}

\subsubsection{Handlebodies with arbitrary framings}
Before turning to handlebodies with arbitrary framings, we pause for an
algebraic interlude. Let $M$ and $N$ be type \DD\ bimodules over $A$
and $B$. Define $\Mor(M,N)$, the \emph{chain complex of bimodule
  morphisms from $M$ to $N$},
\index{bimodule!morphisms!chain complex of|see{morphisms, chain complex of}}%
\index{morphisms!chain complex of}%
\glsit{$\Mor(M,N)$}%
to be the space of bimodule maps from $M$
to $N$, equipped with a differential given by
\[
\bdy(f) = f\circ \bdy_M + \bdy_N\circ f.
\]
(Under technical assumptions on $A$ and $B$ satisfied by the
algebras in bordered Floer theory, the homology of $\Mor(M,N)$ is the
\emph{Hochschild cohomology} $\HH^*(M,N)$ of $M$ with
$N$. See~\cite{LOTHomPair} for a little further discussion.)

Now, let $(H^g,\phi_g^0\circ \psi)$ be a handlebody with
arbitrary framing. Here, $\psi\co {-F{(\PMC)}}\to -{F(\PMC_0^g)}$ for some
genus $g$ pointed matched circle $\PMC$.  Fix a factorization
$\psi=\psi_1\circ\cdots\circ \psi_n$ of $\psi$ into
arc-slides. Let $\psi_i\co -{F(\PMC_{i})}\to -{F(\PMC_{i-1})}$. Here,
$\PMC_0=\PMC_0^g$ and $\PMC_n=\PMC$.

\begin{figure}
  \centering
  \input{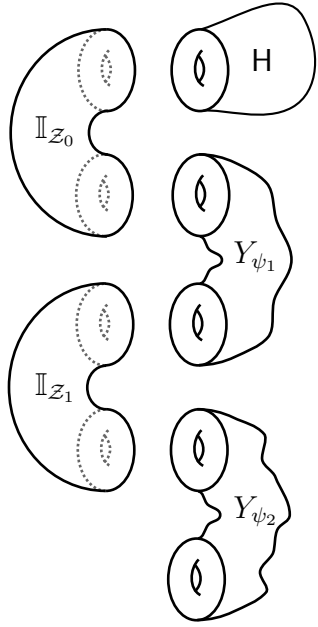}
  \caption{\textbf{Changing framing by gluing mapping cylinders.} This
    is an illustration of Formula~\eqref{eq:ChangeFraming}. The pieces
    labeled $Y_{\psi_1}$ and $Y_{\psi_2}$ represent the mapping
    cylinders of $\psi_1$ and $\psi_2$ (see
    Section~\ref{subsec:MCGroupoid}), and the pieces labeled
    $\Id_{\PMC_i}$ represent copies of $F(\PMC_i)\times[0,1]$. The fact that the $\Id$ pieces face right indicates they have been reflected, i.e., dualized.}
  \label{fig:glue-cylinders}
\end{figure}

As discussed in Section~\ref{sec:dd-arc-slide-intro},
associated to each $\psi_i$ is a bimodule $\lsub{\Alg(\PMC_{i})}\DDmod(\psi_{i})_{\Alg(\PMC_{i-1})}$. Define
\begin{equation}\label{eq:ChangeFraming}
\Dmod(H^g,\phi^0_g\circ \psi)=\Mor(\DDmod(\Id_{\PMC_{n-1}})\otimes
\cdots \otimes \DDmod(\Id_{\PMC_0}) ,\DDmod(\psi_n)\otimes
\DDmod(\psi_{n-1})\otimes\cdots\otimes\DDmod(\psi_1)\otimes \Dmod(\HB)),
\end{equation}
the chain complex of morphisms of
$\Alg(\PMC_{n-1})\otimes\cdots\otimes \Alg(\PMC_0)$-bimodules. This
complex retains a left action by $\Alg(\PMC)$, from the left action on
$\DDmod(\psi_n)$. (This is illustrated schematically in Figure~\ref{fig:glue-cylinders}.)

\begin{introthm}
  \label{thm:Handlebodies}
  The module $\Dmod(H^g,\phi^0_g\circ \psi)$ is homotopy
  equivalent to the module $\CFDa(H^g,\phi^0_g\circ \psi)$ as defined
  (via holomorphic curves) in~\cite{LOT1}.
\end{introthm}
\emph{A priori}, the module $\Dmod(H^g,\phi^0_g\circ \psi)$ depends
not just on $\psi$ but also the factorization into
arc-slides. Theorem~\ref{thm:Handlebodies} implies that, up to
homotopy equivalence, $\Dmod(H^g,\phi^0_g\circ \psi)$ is independent
of the factorization. In fact, this homotopy equivalence is canonical
up to homotopy, a fact that will be used (and proved) in
\cite{LOTCobordisms}.

\subsection{Assembling the pieces: calculating \textalt{$\HFa$}{HF\textasciicircum} from a Heegaard splitting}

Let $Y$ be a closed, oriented three-manifold presented by a Heegaard
splitting $Y=H_1\cup_{\Sigma} H_2$, where $H_1$ and $H_2$ are handlebodies, with
$\partial H_1=-\Sigma$ and $\partial H_2=\Sigma$.  
Thinking of both $H_1$ and $H_2$
as a standard bordered handlebody $\HB_0$, we
can think of the gluing map identifying the two boundaries as a map
$\psi\co F(\PMC_0^g)\to -{F(\PMC_0^g)}$.

Using $\psi$ we get a module $\lsub{\Alg(-\PMC_0^g)}\Dmod(H^g,\phi^0_g\circ
\psi)$, which is a left module over
$\Alg(-\PMC_0^g)$. Using the identification $\Alg(-\PMC_0^g)=
\Alg(\PMC_0^g)^\op$, we can view $\Dmod(H^g,\phi^0_g\circ
\psi)_{\Alg(\PMC_0^g)}$ as a right module over $\Alg(\PMC_0^g)$. So,
$\lsub{\Alg(\PMC_0^g)}\Dmod(H^g,\phi^0_g) \otimes
\Dmod(H^g,\phi^0_g\circ \psi)_{\Alg(\PMC_0^g)}$ becomes an
$\Alg(\PMC_0^g)$-bimodule.

\begin{introthm}
  \label{thm:HFa}
  The chain complex $\CFa(Y)$, as defined
  in~\cite{OS04:HolomorphicDisks} via holomorphic curves, is homotopy
  equivalent to
  \begin{equation*}
    \Mor(\lsub{\Alg(\PMC_0^g)}\DDmod(\Id_{\PMC_0^g})_{\Alg(\PMC_0^g)},
    \lsub{\Alg(\PMC_0^g)}\Dmod(H^g,\phi^0_g)
    \otimes \Dmod(H^g,\phi^0_g\circ \psi)_{\Alg(\PMC_0^g)}),
  \end{equation*}
  the chain complex of bimodule morphisms from
  $\DDmod(\Id_{\PMC_0^g})$ to $\Dmod(H^g,\phi^0_g ) \otimes
  \Dmod(H^g,\phi^0_g\circ\psi)$.
\end{introthm}
(Compare Theorem~\ref{thm:PairingTheorem} in Section~\ref{sec:pairing-theorem}.)

To keep the exposition simple, we have suppressed relative gradings
and $\SpinC$-structures from the introduction.  However, this
information can be extracted in a natural way from the tensor
products, once the gradings on the constituent modules have been
calculated, i.e., once we have graded analogues of
Theorem~\ref{thm:DDforArc-Slides} and Theorem~\ref{thm:Handlebodies}.
We return to a graded analogue of Theorem~\ref{thm:HFa} in Section~\ref{sec:CompleteProof}.

\subsection{More computations: open books, bordered invariants}
Theorem~\ref{thm:DDforArc-Slides} can be combined with slight variants
on Proposition~\ref{prop:CFD-of-0-fr-hdlbdy} for other kinds of
computations as well. 

As a first example, suppose one is given an open book decomposition of
a $3$-manifold~$Y$, with connected binding. Let $F$ denote the fiber of
the open book and $\psi\co F\to F$ the monodromy. Let
$\overline{\psi}\co F\cup_\bdy -F\to F\cup_\bdy -F$ be the result of
extending $\psi$ by the identity map of $-F$. Fix a pointed matched
circle $\PMC$ representing $F$. There is a particular bordered
handlebody, which we call the self-gluing handlebody $\HB_{sg}$, so
that $\bdy\HB_{sg}=\PMC\cup(-\PMC)$ and 
\[
Y= \HB_{sg}\cup_{\overline{\psi}}\HB_{sg};
\]
see Definition~\ref{def:sghb}. Factoring $\psi$ into arc-slides gives
a formula for $\CFa(Y)$, analogous to Theorem~\ref{thm:HFa}, in terms
of the module $\CFDa(\HB_{sg})$ and the bimodules for the arc-slides. The
invariant $\CFDa(\HB_{sg})$ is computed in Theorem~\ref{thm:sghb},
completing this algorithm for computing $\CFa(Y)$ from an open book
decomposition.  See Section~\ref{sec:ok-computer} for a
(computer-assisted) example.

As a second example, one can also compute the bordered invariants of
arbitrary bordered $3$-manifolds $Y$. The new ingredient here is a
computation of the invariants of elementary cobordisms,
Proposition~\ref{prop:CalculateElementaryCobordism}. (The answer is an
amalgam of $\CFDa(\HB^1)$ and $\CFDDa(\Id)$.) From there, one
decomposes an arbitrary bordered $3$-manifold as a composition of
elementary cobordisms with standard framings and mapping cylinders of
arc-slides, and again uses an analogue of the formula from
Theorem~\ref{thm:HFa} to obtain $\CFDa(Y)$. See
Section~\ref{sec:ElementaryCobordisms} for more details.

One can obtain $\CFAa(Y)$ from $\CFDa(Y)$ via the duality
theorem~\cite[Proposition~\ref*{LOT2:prop:DDAA-duality}]{LOT2}; see
also Section~\ref{sec:computing-type-a}. This gives a
finite-dimensional but fairly large for model for $\CFAa$. One can obtain
a smaller model by using homological perturbation theory; this is
reviewed, with an example, in
Sections~\ref{sec:hpt}--\ref{sec:AA-torus-id}.

\subsection{Organization}
In Section~\ref{sec:Background},
we give some of the background on bordered Floer theory needed for
this paper; for further
details the reader is referred to~\cite{LOT1} and~\cite{LOT2}.  In
Section~\ref{sec:DDforIdentity}, we calculate the \DD\ bimodule for
the identity map, verifying Theorem~\ref{thm:DDforIdentity}.  The
proof follows from inspecting the relevant Heegaard diagram and
applying the relations which are forced by $\partial^2=0$. In
Section~\ref{sec:Arc-Slides}, we calculate the \DD\ bimodules for
arc-slides. This uses similar reasoning to the proof of
Theorem~\ref{thm:DDforIdentity}. In Section~\ref{sec:GenusOne}, we explain the specialization of these to the genus one case, for concreteness.
In Section~\ref{sec:mcg-grading}
we compute gradings on the arc-slide bimodules (needed for a suitably graded
analogue of Theorem~\ref{thm:HFa}).
Note that we do not at present have a conceptual description of the bimodules for arbitrary
surface diffeomorphisms (rather, they have to be factored into arc-slides, and the corresponding
bimodules have to be composed); however we do give an intrinsic description of its corresponding
grading set. This is done in Proposition~\ref{prop:mu-is-ell-1}.
In Section~\ref{sec:CompleteProof},
we compute the invariant for handlebodies with the preferred framing
$\phi_0$, which is quite easy, and assemble the ingredients to prove
the main result, Theorem~\ref{thm:HFa} (as well as a graded version).

The computations of the arc-slide bimodules also lead quickly to
a description of the bordered Heegaard Floer invariants for arbitrary bordered
three-manifolds. The main ingredient beyond what we
have explained so far is the invariant associated to an elementary
cobordism that adds or removes a handle. This is an easy
generalization of the calculations for handlebodies, and is discussed
in Section~\ref{sec:ElementaryCobordisms}.

The point of view of $\Ainf$-modules, which we have
otherwise avoided in this paper, allows one to trade generators for
complexity of the differential, and is useful in practice. This is
discussed in Section~\ref{sec:AModules}, along with some examples.

\subsection{Further remarks}

Theorem~\ref{thm:HFa} gives a purely combinatorial description of
$\HFa(Y)$, with coefficients in $\Zmod{2}$, in terms of a mapping
class of a corresponding Heegaard splitting.  We point out again
that this calculation is independent
of the methods of~\cite{SarkarWang07:ComputingHFhat}.
Indeed, the methods of this paper are based on general properties of
bordered invariants, together with some
very crude input coming from the Heegaard diagrams (see especially
Theorems~\ref{thm:DDforIdentity} and~\ref{thm:DDforArc-Slides}).  The
particular form of the bimodules is then forced by algebraic
considerations (notably $\bdy^2=0$).

The \DD\ bimodule for the identity (as described in
Theorem~\ref{thm:DDforIdentity}) was also calculated
in~\cite[Theorem~\ref*{HomPair:thm:PreciseDD}]{LOTHomPair}, by
different methods. The proof of Theorem~\ref{thm:DDforIdentity} is
included here (despite its redundancy with results of~\cite{LOTHomPair}),
since it is a model for the more complicated
Theorem~\ref{thm:DDforArc-Slides}.

In the present paper, we have calculated the $\HFa$ variant of
Heegaard Floer homology for closed three-manifolds. This is also a key
component in the combinatorial description of the invariant for
cobordisms, which will be given in a future paper~\cite{LOTCobordisms}.

\subsection{Acknowledgements}
We thank Bohua Zhan for helpful comments on a draft of this paper. We also thank the referee for a thorough reading of the paper and many suggestions and corrections.


\section{Preliminaries}
\label{sec:Background}
In this section we will review most of the background on bordered
Floer theory needed later in the paper.  In
Section~\ref{subsec:MCGroupoid}, we recall the mapping class group (or
rather, groupoid) relevant to our considerations.  In
Section~\ref{subsec:AlgebraPMC}, we amplify the remarks in the
introduction regarding the algebras $\Alg(\PMC)$ associated to pointed
matched circles.  In Section~\ref{subsec:Bordered}, we review the
basics of the type $D$ modules associated to $3$-manifolds with one
boundary component. We also introduce the notion of the {\em
  coefficient algebra} of a type $D$ structure, which is used later
in the calculation of arc-slide bimodules. In
Section~\ref{subsec:Bimodules}, we review the case of type \DD\
modules for three-manifolds with two boundary components, and
introduce their coefficient
algebras.
In Section~\ref{sec:pairing-theorem}, we turn
to the versions of the pairing theorem that will be used in this paper. For
more details on any of these topics, the reader is referred
to~\cite{LOT1},~\cite{LOT2} and~\cite{LOTHomPair}.

This section does not discuss the type $A$ module associated to a
bordered $3$-manifold with one boundary component, nor the type \DA\
or \AAm\ modules associated to a bordered $3$-manifold with two
boundary components. By~\cite[Proposition~\ref*{LOT2:prop:DDAA-duality}]{LOT2} (or any of
several results from~\cite{LOTHomPair}), these invariants can be
recovered from the type $D$ and \DD\ invariants. We use these results
to circumvent explicitly using type $A$ modules in most of the paper,
though we return to them in
Section~\ref{sec:AModules}.

\subsection{The (strongly-based) mapping class groupoid}
\label{subsec:MCGroupoid}

As discussed in the introduction, the main work in this paper consists
in computing the bimodules associated to arc-slides. Since arc-slides
connect different pointed matched circles, they correspond to maps
between different (though homeomorphic) surfaces. To put this
phenomenon in a more general context, we recall some basic properties
of a certain mapping class groupoid.

\glsit{$k$}%
\glsit{$\PMC$}%
\glsit{$Z$}\glsit{$\mathbf{a}$}\glsit{$M$}\glsit{$z$}%
\index{pointed matched circle}%
Fix an integer $k$.  Let $\PMC=(Z,\mathbf{a},M,z)$ be a pointed
matched circle on $4k$ points. We can associate to $\PMC$ a surface
$\PunctF(\PMC)$ as follows. Let $D$ be a disk with boundary
$Z$. Attach a $2$-dimensional $1$-handle to $\bdy D$ along each pair
of matched points in $\mathbf{a}$. The result is a surface
$\PunctF(\PMC)$ with one boundary component, and a basepoint $z$ on
that boundary component. Let $F(\PMC)$ denote the result of filling
the boundary component of $\PunctF(\PMC)$ with a disk; we call the
disk $F(\PMC)\setminus \PunctF(\PMC)$ the \emph{preferred disk}
\index{preferred disk}\index{disk, preferred}%
in $F(\PMC)$; the basepoint $z$ lies on the boundary of the preferred disk.
\glsit{$\PunctF(\PMC)$}\glsit{$F(\PMC)$}%

(The construction of $F(\PMC)$ given here agrees with \cite{LOT1} and
\cite{LOT2}, and differs superficially
from the construction in \cite{LOTHomPair}.)

\index{mapping class, strongly based}%
\index{strongly based!mapping class}%
\glsit{$\MCG_0(\PMC,\PMC_2)$}%
Given pointed matched circles $\PMC_1$ and $\PMC_2$, the 
set of {\em strongly-based mapping classes from $\PMC_1$ to $\PMC_2$}, denoted $\MCG_0(\PMC,\PMC_2)$,
is the set of orientation-preserving, isotopy class of
homeomorphisms $\phi\co\PunctF(\PMC_1)\to
\PunctF(\PMC_2)$ carrying $z_1$ to $z_2$, where $z_i\in
\partial \PunctF(\PMC_i)$ is the basepoint; 
\[
\MCG_0(\PMC_1,\PMC_2)=\{\phi\co\PunctF(\PMC_1)\stackrel{\cong}{\to}
\PunctF(\PMC_2)\mid \phi(z_1)=z_2\}/\text{isotopy}.
\] 
(The subscript
$0$ on the mapping class group indicates that maps respect the boundary and the basepoint.)
In the case where $\PMC_1=\PMC_2$, this set naturally forms a group,
which we call the {\em strongly-based mapping class group}.
\index{mapping class group, strongly based}%
\index{strongly based!mapping class group}%

\index{mapping class groupoid!genus $k$}%
\index{genus $k$ mapping class groupoid}%
\glsit{$\MCG_0(k)$}%
More generally, the {\em{strongly-based genus $k$ mapping class groupoid}}
$\MCG_0(k)$ is the category whose objects are pointed matched circles
with $4k$ points and with morphism set between $\PMC_1$ and $\PMC_2$
given by $\MCG_0(\PMC_1,\PMC_2)$. 

\glsit{$\PunctF(m)$}%
Recall that when $\PMC$ and $\PMC'$ differ
by an arc-slide, there is a canonical
strongly-based diffeomorphism $\PunctF(m)\co F(\PMC)\rightarrow F(\PMC')$,
as pictured in Figure~\ref{fig:ArcSlide}.

Any morphism in the mapping class groupoid can be factored as a
product of arc-slides; see, for
example,~\cite[Theorem~5.3]{Bene08:ChordDiagrams}.  One proof:
consider Morse functions $f$ on $\PunctF(\PMC)$ and $f'$ on
$\PunctF(\PMC')$ inducing the pointed matched circles. Let $\phi\co
\PunctF(\PMC)\to \PunctF(\PMC')$ be an orientation-preserving
diffeomorphism. The Morse functions $f$ and $\phi^*(f')=f'\circ \phi$
can be connected by a generic one-parameter family of Morse functions
$f_t$. The finitely many times $t$ for which $f_t$ has a flow-line
from between two index $1$ critical points give the sequence of
arc-slides connecting $\PMC$ and $\PMC'$.

\begin{figure}
  \centering
  \includegraphics[scale=.5]{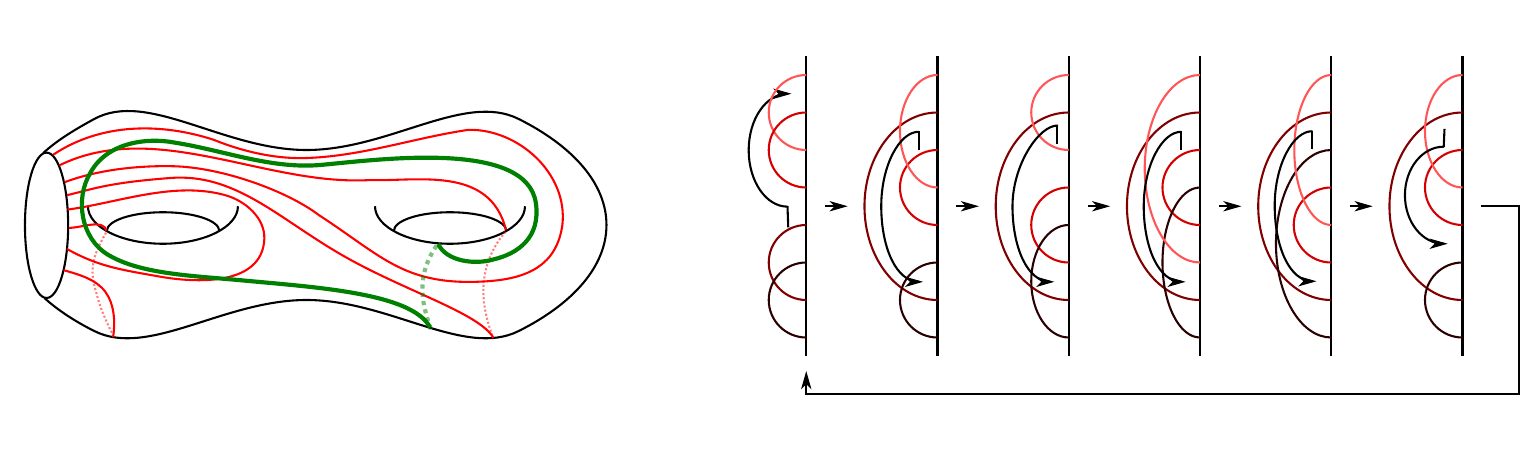}
  \caption{\textbf{Factoring a Dehn twist into arc-slides.} Left: a
    genus $2$ surface specified by a pointed matched circle, and a
    curve $\gamma$ (drawn in thick green) in it. Right: a sequence of arc-slides whose
    composition is a Dehn twist around $\gamma$.}
  \label{fig:factor-Dehn}
\end{figure}

For instance, any Dehn twist can be factored as a product of
arc-slides. The key point to doing this in practice is the following:
\begin{lemma}\label{lemma:factor-Dehn}
  Let $\PMC=(Z,\mathbf{a},M,z)$ be a pointed matched circle and
  $\{b,b'\}\subset\mathbf{a}$ a matched pair in $\PMC$. Consider the
  sequence of arc-slides where one slides each of the points in
  $\mathbf{a}$ between $b$ and $b'$ over $\{b,b'\}$ once, in
  turn.
  This product of arc-slides
  is a factorization of the Dehn twist around the curve in
  $\PunctF(\PMC)$ specified by $\{b,b'\}$.
\end{lemma}
(See Figure~\ref{fig:factor-Dehn} and compare \cite[Lemma 8.3]{AndersenBenePenner09:MCGroupoid}.)
\begin{proof}
  The proof is left to the reader.
\end{proof}

In particular, for a genus $1$ pointed matched circle, arc-slides are
Dehn twists. For an illustration of the factorization of a more
interesting Dehn twist in the genus $2$ case, see
Figure~\ref{fig:factor-Dehn}.

\subsubsection{Strongly bordered \textalt{$3$}{3}-manifolds and mapping cylinders}
It will be convenient to think of strongly-based diffeomorphisms in
terms of their mapping cylinders. Given a strongly-based
diffeomorphism $\psi\co \PunctF(\PMC_1)\to \PunctF(\PMC_2)$, we can
extend $\psi$ by the identity map on $\bD^2$ to a diffeomorphism
$\psi\co F(\PMC_1)\to F(\PMC_2)$. Consider the $3$-manifold
$[0,1]\times F(\PMC_2)$. This manifold is equipped with
orientation-preserving identifications
\begin{align*}
  \psi&\co {{-F(\PMC_1)}}\to \{0\}\times F(\PMC_2)\subset [0,1]\times
  F(\PMC_2)\\
  \Id&\co F(\PMC_2)\to \{1\}\times F(\PMC_2)\subset [0,1]\times F(\PMC_2).
\end{align*}
It is also equipped with a cylinder $[0,1]\times (F(\PMC_2)\setminus
\PunctF(\PMC_2))$, which is essentially the same data as a framed arc
$\gamma$ connecting the two boundary components of $[0,1]\times
F(\PMC_2)$. We call the data $([0,1]\times F(\PMC_2),\psi, \Id,
\gamma)$ the \emph{strongly bordered $3$-manifold associated to
  $\psi$} or the \emph{mapping cylinder of $\psi$};
\index{strongly bordered $3$-manifold associated to diffeomorphism}%
\index{mapping cylinder of diffeomorphism}%
compare~\cite[Construction~\ref*{LOT2:construct:MCG-to-Bordered}]{LOT2}. Let $Y_\psi$ denote
the mapping cylinder of $\psi$.
\glsit{$([0,1]\times F(\PMC_2),\psi, \Id, \gamma)$}%
\glsit{$Y_\psi$}%

Observe that for $\psi_{12}\co F(\PMC_1)\to F(\PMC_2)$ and
$\psi_{23}\co F(\PMC_2)\to F(\PMC_3)$, 
$Y_{\psi_{23}\circ \psi_{12}}$ is orientation-preserving homeomorphic
to $Y_{\psi_{12}}\cup_{F(\PMC_2)}Y_{\psi_{23}}$.

\index{strongly bordered $3$-manifold with two boundary components}%
More generally, a \emph{strongly bordered $3$-manifold with two
  boundary components} consists of a $3$-manifold $Y$ with two
boundary components $\bdy_LY$ and $\bdy_RY$, diffeomorphisms
$\phi_L\co F(\PMC_L)\to \bdy_LY$ and $\phi_R\co F(\PMC_R)\to \bdy_RY$
for some pointed matched circles $\PMC_L$ and $\PMC_R$, and a framed
arc in $Y$ connecting the basepoints $z$ in $F(\PMC_L)$ and
$F(\PMC_R)$, and so that the framing points into the preferred disk of
$F(\PMC_L)$ and $F(\PMC_R)$ at the two boundary components.
\glsit{$\bdy_LY$, $\bdy_RY$}\glsit{$\phi_L$, $\phi_R$}\glsit{$\PMC_L$, $\PMC_R$}%

\subsection{More on the algebra associated to a pointed matched circle}
\label{subsec:AlgebraPMC}
\newcommand\Lgen{g_L}
\newcommand\Rgen{g_R}

Fix a pointed matched circle $\PMC$, as in
Section~\ref{subsec:AlgPMC}, with basepoint $z$.
Let ${\mathbf a}\subset \PMC$ denote the set of points which are matched.

Each strands diagram $\upsilon$ has an associated one-chain $\supp(\upsilon)$, which is
an element of $H_1(\PMC\setminus\{z\},{\mathbf a})$. This is gotten by
projecting the strands diagram, thought of as a one-chain in
$[0,1]\times(\PMC\setminus\{z\})$, onto $\PMC$.
(The one-chain $\supp(\upsilon)$ is denoted $[\upsilon]$
in~\cite{LOT1}; we have chosen to change
notation here in order to avoid a conflict with the
standard notation for a closed interval.)
\glsit{$\supp(\upsilon)$}%

Recall that if $\upsilon$ is a strands diagram then we call its associated algebra
element a {\em basic generator}
\index{basic generator}%
for the algebra $\Alg(\PMC)$; we
usually do not distinguish between the strands diagram and its
associated algebra element, writing, for example $\supp(a)$ when $a$
is a basic generator.

In general, for a set $\boldsymbol{\xi} = \{\xi_1, \dots, \xi_k\}$ of
chords on $\PMC$ with endpoints on $\CircPts$, there is an algebra
element $a(\boldsymbol{\xi})$, in which the moving strands correspond
to the $\xi_i$ and we sum over all valid ways of adding horizontal
strands. (If some $\xi_i, \xi_j$, $i\neq j$, have their initial (respectively
terminal) endpoints the
same or matched, we define $a(\boldsymbol{\xi})=0$: in this case there
are no valid
ways of adding horizontal strands.) We will also abuse notation
slightly, and write
$a(X)$ for $X$ a subset of $\PMC$ with boundary only at points in
$\CircPts$: this means $a(\boldsymbol{\xi})$, where $\boldsymbol{\xi}$
is the set of connected components of~$X$.  (Each connected component
is an interval, of course.)
\glsit{$a(\boldsymbol{\xi})$}\glsit{$a(X)$}%

\index{homogeneous algebra element}%
An element of the algebra is called {\em homogeneous} if it can be
written as a sum of basic generators so that each basic generator in
the sum
\begin{itemize}
\item has the same associated one-chain,
\item has the same initial (and hence, in view of the previous
  condition, terminal) idempotent, and in particular has the same
  weight, and
\item has the same number of crossings.
\end{itemize}

\subsubsection{The opposite algebra}\label{subsubsec:opposite-algebra}
Suppose that $\PMC=(Z,\mathbf{a},M,z)$ is a pointed matched
circle. Let $-\PMC$ denote its reverse, i.e., the pointed matched
circle obtained by reversing the orientation on $Z$.  There is an
obvious orientation-reversing map $r\co \PMC\rightarrow -\PMC$,
and hence an identification between chords for $\PMC$ and chords for
$-\PMC$. It is easy to see that this map $r$ induces an isomorphism
$\Alg(\PMC)^\op\cong \Alg(-\PMC)$, where $\Alg(\PMC)^{\op}$ denotes the
opposite algebra to $\Alg(\PMC)$.
\glsit{$-\PMC$}\glsit{$r$}\glsit{$\Alg(\PMC)^{\op}$}%

In particular, left $\Alg(-\PMC)$-modules correspond to left
$\Alg(\PMC)^\op$-modules, and hence to right $\Alg(\PMC)$-modules.

\subsubsection{Gradings}
\label{sec:gradings-alg}

\glsit{$G(\PMC)$}\glsit{$\lambda$}\glsit{$\gr$}%
\index{gradings!on $\Alg(\PMC)$}%
The algebra $\Alg(\PMC)$ is graded in the following sense.
There is a group $G(\PMC)$, equipped with a distinguished central element $\lambda$,
and a function $\gr$ from basic generators of $\Alg(\PMC)$ to $G(\PMC)$,
with the following properties:
\begin{itemize}
  \item If $\upsilon$ and $\tau$ are basic generators of $\Alg(\PMC)$, and
    $\upsilon\cdot \tau\neq 0$,
    then $\gr(\upsilon\cdot \tau)=\gr(\upsilon)\cdot \gr(\tau)$.
  \item If $\tau$ appears with non-zero multiplicity in $d\upsilon$ 
    then $\gr(\upsilon)=\lambda\cdot \gr(\tau)$.
\end{itemize}

In fact, there are two choices of grading group for $\Alg(\PMC)$. The smaller
one, which is more natural from the point of view the pairing theorem,
is a Heisenberg group on the first homology of the underlying surface.
\index{Heisenberg group}\index{grading refinement data}%
Gradings in this smaller set depend on a further universal choice of grading
refinement data, as in~\cite[Section~\ref*{LOT2:sec:Gradings}]{LOT2}
(see also Section~\ref{sec:gr-of-gens}), although different
choices of refinement data lead to canonically equivalent module
categories. However, we will generally work with the big grading group $G'(\PMC)$
in this paper.

\glsit{$G'(\PMC)$}%
More precisely, the big grading group $G'(\PMC)$ is a $\ZZ$ central extension of
$H_1(Z \setminus\{z\}, \mathbf{a})$, realized explicitly as pairs
$(j,\alpha)\in \frac{1}{2}\ZZ\times H_1(Z\setminus \{z\},\mathbf{a})$
subject to a congruence condition $j\equiv \epsilon(\alpha)$ $\pmod
1$, for the function $\epsilon\co
H_1(Z\setminus\{z\},\mathbf{a})\to \OneHalf\ZZ/\ZZ$ 
given by $1/4$ the number of parity changes in the support of $a$;
see~\cite[Section~\ref*{LOT2:sec:Gradings}]{LOT2}. The multiplication is given by
\[
(j_1,\alpha_1)\cdot(j_2,\alpha_2)=(j_1+j_2+m(\alpha_2,\bdy\alpha_1),\alpha_1+\alpha_2)
\]
where $m(\alpha,x)$ is the local multiplicity of $\alpha$ at $x$;
$m(\alpha_2,\bdy \alpha_1)$
is a $\OneHalf\ZZ$-valued extension of the intersection form on $H_1(F(\PMC))\subset
H_1(Z\setminus \{z\},\mathbf{a})$. The distinguished central element
is $\lambda=(1,0)$.
The $G'(\PMC)$-grading of a strands diagram is given by
\[
\grb(a) \coloneqq (\iota(a), \supp(a))
\]
\glsit{$\grb(a)$}%
where $\supp(a) \in H_1(\PMC\setminus\{z\},\mathbf{a})$ is as defined
above, and $\iota(a)$ records the number of crossings plus a
correction term:
\glsit{$\iota(a)$}%
\[
  \iota(a) \coloneqq \inv(a) - m(\supp(a), \SetS),
\]
where $I(\SetS)$ is the initial idempotent of $a$.
See~\cite[Section~\ref*{LOT1:sec:gradings-algebra}]{LOT1} for further details.

Homogeneous algebra elements (as defined earlier) live
in a single
grading. For a sum of basic generators with the same left and right
idempotents, the converse is true: homogeneity with respect to the
grading (for either grading group) implies homogeneity as defined above.

\begin{lemma}\label{lem:NegGradings}
  If a basic generator~$a$ is not an idempotent, then $\iota(a) \le
  -k/2$, where $k$ is the number of intervals in a minimal
  expression of $\supp(a)$ as a sum of intervals.
\end{lemma}

\begin{proof}
  This is essentially
  \cite[Lemma~\ref*{LOT2:lem:NegativeGradingsOnAlgebra}]{LOT2}.  The
  argument there shows that $\iota(a) \le -k'/2$, where $k'$ is the
  number of moving strands in $a$; but if $k$ is as given in the
  statement, then $k' \ge k$.
\end{proof}

\subsubsection{The quotient map}
\label{subsec:QuotientMap}

Recall that if $\PMC$ and $\PMC_0$ are pointed matched circles then we can
form their connected sum $\PMC\#\PMC_0$.
\index{connected sum}%
Note that there are two natural choices of
where to put the basepoint in $\PMC\#\PMC_0$.

Given any idempotent $I_0=I(\SetS_0)$ for $\PMC_0$, we have a quotient map
$$q\co \Alg(\PMC\#\PMC_0)\rightarrow\Alg(\PMC)$$
\glsit{$q$}%
defined as follows.
The idempotents for $\PMC\#\PMC_0$ have the form
$I(\SetS\amalg \SetT)$, where $\SetS$ (respectively $\SetT$) is a
subset of the matched pairs in $\PMC$ (respectively $\PMC_0$). The quotient map $q$ is determined by
its action on the idempotents:
$$q(I(\SetS\amalg \SetT))=\left\{\begin{array}{ll}
    I(\SetS) & {\text{if $\SetT=\SetS_0$}} \\
    0 & {\text{otherwise}}\end{array}\right.
$$
and also the property that
$q(a)=0$ unless $\supp(a)\subset \PMC\subset \PMC\#\PMC_0$.

\glsit{$Q$}%
The map $q$  can be promoted to a map 
$$Q\co \Alg(\PMC\#\PMC_0)\otimes\Alg(-\PMC'\#\PMC_0)
\rightarrow \Alg(\PMC)\otimes\Alg(-\PMC').$$

The map $q$ is used in the definition of stability for arc-slide
bimodules. The notation is somewhat lacking, since we have not
specified how we have taken the connect sum of $\PMC$ and $\PMC_0$
(i.e., in which of the two possible regions in the connect sum we have
placed the basepoint). This information is not important, however,
since stability uses both possible choices.

\subsection{Bordered invariants of \textalt{$3$}{3}-manifolds with connected
  boundary}
\label{subsec:Bordered}

We recall the basics of the bordered Heegaard Floer invariant $\CFDa(Y)$ for a
bordered three-manifold. 

\subsubsection{Bordered Heegaard diagrams and \textalt{$\CFDa$}{CFD\textasciicircum}}
\begin{definition}
\label{def:BorderedDiagram}
A \emph{bordered Heegaard diagram}, is a quadruple
$\HD=(\widebar{\Sigma},\widebar{\alphas},\betas,z)$ consisting of
\begin{itemize}
\item a compact, oriented surface $\widebar{\Sigma}$ with one boundary component, of
  some genus~$g$;
\item a $g$-tuple of pairwise-disjoint circles
  $\betas=\{\beta_1,\dots,\beta_g\}$ in the interior of~$\Sigma$;
\item a $(g+k)$-tuple of pairwise-disjoint curves
  $\widebar{\alphas}$ in
  $\widebar{\Sigma}$, consisting of $g-k$ circles
  $\alphas^c=(\alpha_1^c,\dots,\alpha_{g-k}^c)$ in the
  interior of $\widebar{\Sigma}$ and $2k$ arcs
  $\balphas^a=(\balpha_1^a,\dots,\balpha_{2k}^a)$
  in $\widebar{\Sigma}$ with boundary on $\partial\widebar{\Sigma}$ (and transverse to
  $\partial{\widebar{\Sigma}}$); and
\item a point $z$ in $(\partial\widebar{\Sigma})\setminus(\alphas\cap\partial\widebar{\Sigma})$,
\end{itemize}
such that $\betas\pitchfork\balphas$ and
$\bSigma\setminus \balphas$ and $\bSigma
\setminus \betas$ are connected.
\end{definition}
(As in~\cite{LOT1}, we let $\Sigma$ denote the interior of
$\overline{\Sigma}$ and $\alphas=\overline{\alphas}\cap\Sigma$; and
will often blur the distinction between
$(\overline{\Sigma},\overline{\alphas},\betas,z)$ and
$(\Sigma,\alphas,\betas,z)$.)

The boundary of a  Heegaard diagram
$\HD=(\Sigma,\overline{\alphas},\betas,z)$ is naturally a pointed, matched
circle as follows. The boundary
$(\partial\Sigma)$ inherits its basepoint from
$z\in \partial \Sigma$, and the points
$(\overline{\alphas}\cap\partial{\overline \Sigma})$, can be paired off
according to which arc they belong to. 
\index{pointed matched circle!boundary of bordered Heegaard diagram}%

A bordered Heegaard diagram specifies an oriented
three-manifold with boundary $Y$, along with an orientation preserving
diffeomorphism $\phi\co F(\PMC)\rightarrow \partial Y$; i.e., a $\PMC$-bordered
three-manifold.

We briefly recall the construction of the type $D$ module associated
to a bordered Heegaard diagram $\HD$. Let $\PMC$ be the matched circle
appearing on the boundary of $\HD$. The type $D$ module associated to $\HD$ is
a left module over $\Alg(-\PMC)$, where $-\PMC$ is the reverse of
$\PMC$.

Let $\Gen(\HD)$ be the set of subsets
$\x\subset \alphas \cap \betas$ with the following properties:
\begin{itemize}
\item $\x$ contains exactly one element on each $\beta$ circle,
\item $\x$ contains exactly one element on each $\alpha$ circle, and
\item $\x$ contains at most one element on each $\alpha$ arc.
\end{itemize}
Let $X(\HD)$ be the $\Field$--vector space spanned by $\S(\HD)$.  For
$\x \in \S(\HD)$, let $o(\x) \subset [2k]$ be the set of $\alpha$-arcs
occupied by~$\x$.  Define $I_D(\x)$ to be $I([2k]\setminus o(\x))$;
that is, the idempotent corresponding to the complement of
$o(\x)$.  We can now define an action of the subalgebra of
idempotents~$\Idem$ inside~$\Alg(-\PMC)$ on $X(\HD)$ via
\index{idempotents!subalgebra of}%
\glsit{$\Idem$}%
\begin{equation}\label{eq:def-idem-DMod}
  I(\SetS) \cdot \x =
  \begin{cases}
    \x&I(\SetS) = I_D(\x)\\
    0&\text{otherwise,}
  \end{cases}
\end{equation}
where $\SetS$ is a $k$-element subset of $[2k]$. As a module, let
$\CFDa(\HD)=\Alg(-\PMC)\otimes_{\Idem}X(\HD)$.
\glsit{$\CFDa(\HD)$}\index{type $D$ module!associated to bordered $3$-manifold}%

Fix generators $\x,\y\in\Gen(\HD)$. Two-chains in $\Sigma$ which
connect $\x$ and $\y$ in a suitable sense can be organized into
{\em homology classes}, denoted $\pi_2(\x,\y)$; we say elements of
$\pi_2(\x,\y)$ \emph{connect $\x$ to $\y$}.
\index{connect $\x$ to $\y$}%
\index{homology classes}\glsit{$\pi_2(\x,\y)$}%
(To justify the terminology ``homology class,''
note that the difference between any two elements of $\pi_2(\x,\y)$ 
can be thought of as a two-dimensional homology class in $Y$;
the notation is justified by its interpretation in terms of the symmetric
product, see~\cite{OS04:HolomorphicDisks}.) Given a homology class
$B\in\pi_2(\x,\y)$ and
asymptotics specified by a vector $\vec{\rho}$, there is an associated
moduli space of holomorphic curves $\Mod^B(\x,\y,\vec{\rho})$.
Counting points in this moduli space
gives rise to an algebra element
\index{holomorphic curves, moduli space of}\index{moduli space of holomorphic curves}%
\begin{equation}
  \label{eq:Differential}
n^B_{\x,\y}=\sum_{\{\vec\rho\,\mid\,\ind(B,\vec{\rho}) = 1\}}\#(\Mod^B(\x, \y; \vec{\rho}))
    a(-\vec{\rho})\in\Alg(-\PMC_L).
\end{equation}
\index{holomorphic curves, moduli space of}%

The algebra elements $n_{\x,\y}^B$ can be assembled to define an operator
$$\delta^1\co X(\HD)\rightarrow \Alg(\PMC)\otimes X(\HD)$$
by
\begin{equation}\label{eq:delta-on-CFD}
  \delta^1(\x) \coloneqq \sum_{\y\in\S(\HD)}\sum_{B\in\pi_2(\x,\y)}n^B_{\x,\y} \otimes \y.
\end{equation}
\glsit{$\delta^1$}%

Let $\CFDa(\HD)$ denote the space $\Alg(-\PMC)\otimes
X(\HD)$. We endow $\CFDa(\HD)$ with a differential $\partial$
induced from the above map $\delta^1$, and the differential on the
algebra $\Alg(-\PMC)$, via the Leibniz rule:
$$\partial(a\otimes \x) = (da)\otimes \x + a\cdot \delta^1(\x).$$
\glsit{$\bdy$}%
The differential, together with the obvious left action on
$\Alg(-\PMC)$, gives $\CFDa(\HD)$ the structure of a left differential
module over $\Alg(-\PMC)$. (The proof involves studying one-parameter
families of holomorphic curves; see~\cite[Section~\ref*{LOT1:sec:typed-d-sq-zero}]{LOT1}.)

The differential on $\CFDa(\HD)$ has the following key property.
Suppose that $B\in\pi_2(\x,\y)$ gives a non-zero contribution of
$a\otimes \y$ to $\partial \x$. Then, $\supp(a)$ is calculated by the
local multiplicities of $B$ at $\partial{\overline \Sigma}$.

Up to homotopy equivalence, the module $\CFDa(\HD)$ depends only on
the bordered $3$-manifold specified by $\HD$. Thus, given a bordered
$3$-manifold $Y$ we will write $\CFDa(Y)$ to denote the homotopy type
of $\CFDa(\HD)$ for any bordered Heegaard diagram~$\HD$
representing~$Y$.

\subsubsection{Type \textalt{$D$}{D} structures}\label{subsubsec:TypeDStructures}

The special structure of $\CFDa(\HD)$ can be formalized in the following:

\begin{definition}
  \index{type $D$ structure}\glsit{$\lsup{\Alg}N$}\glsit{$\delta^1$}%
  Let $\Alg$ be a \dg algebra over a ground ring $\Ground$. A left
  \emph{type $D$ structure} over $\Alg$ is a left $\Ground$-module
  $\lsup{\Alg}N$, together with a degree $0$ map $\delta^1 \colon N
  \to \DGA[1] \otimes N$, satisfying the structural equation
  \begin{equation}
    \label{eq:TypeDStructEq}
    (\mu_2\otimes \Id_N) \circ (\Id_{\Alg}\otimes \delta^1)\circ \delta^1 + (\mu_1\otimes \Id_N)\circ \delta^1 = 0. 
  \end{equation}
\end{definition}
(Here, $\mu_1$ and $\mu_2$ denote the differential and multiplication
on $\Alg$, respectively; the notation is drawn from the theory of $\Ainf$-algebras.)

Given a type $D$ structure as above, we can form the associated module denoted
$\lsub{\Alg}N$ or $\Alg\DT N$,
\glsit{$\lsub{\Alg}N$}\glsit{$\Alg\DT N$}%
whose generators are $a\otimes \x$ with $a\in\Alg$ and $\x\in N$, algebra action by
$$a\cdot (b\otimes \x)=(a\cdot b)\otimes \x,$$
and differential given by
$$\partial(a\otimes \x)=(da)\otimes \x + a\cdot (\delta^1 \x).$$
(Here, $\otimes$ denotes tensor product over $\Ground$, not $\Field$. The structural equation for $\delta^1$,
Equation~\eqref{eq:TypeDStructEq}, is equivalent to the condition that $\partial^2=0$.)

The bordered invariant $\CFDa(\HD)$ is naturally a type $D$ structure over $\Alg(-\PMC)$.

The notion of a type $D$ structure has
an obvious analogue for bimodules: a \emph{type \DD\ structure} over
\dg algebras $\Alg$ and $\Blg$ is just a type $D$ structure over
$\Alg\otimes\Blg^\op$. 
\index{type \DD\ structure}%

\subsubsection{Some particular holomorphic curves}

We have not explained here precisely which curves contribute to
$\Mod^B(\x,\y,\vec{\rho})$. Rather than reviewing the general case, we
restrict our discussion to the main examples we  will need in this
paper. (See~\cite{LOT1} for further details.)

\begin{definition}
\label{def:Polygon}
Suppose that $P$ is a connected component of ${\overline
  \Sigma}\setminus (\alphas\cup\betas)$, and $P$ does not contain $z$.
Suppose moreover that $P$ is a $2n$-gon. Each side of $P$ is one of
three kinds:
\begin{enumerate}[label=(P-\arabic*),ref=P-\arabic*]
  \item \label{type:Beta} an arc contained in some $\beta_i$
  \item \label{type:Alpha} an arc contained in some $\alpha_i$ 
    (which might be of the form
    $\alpha_i^c$, $\alpha_i^{a}$)
  \item \label{type:Boundary} 
    an arc contained in $\partial \Sigma$.
\end{enumerate}
Traversing the boundary of $P$ with its induced orientation, one
alternates between meeting sides of type (\ref{type:Beta}) and
sequences of sides of types~(\ref{type:Alpha}) and~(\ref{type:Boundary}).
Suppose that $P$ has only one side of type (\ref{type:Boundary}).
We call such a component a {\em fundamental polygon}.
\index{fundamental polygon}%
\index{polygon, fundamental}%
\end{definition}

The intersection points of arcs of Type~(\ref{type:Beta}) and those
of Type~(\ref{type:Alpha}), which we call {\em corners},
\index{corner}%
can be
partitioned into two types: those which lie at the initial point of
the arc of Type~(\ref{type:Alpha}) (with its induced
orientation from $\bdy P$), and those which lie at the terminal point of the arc of
Type~(\ref{type:Alpha}). Let $\x_0$ denote the set of corners of the
first type, and let $\y_0$ denote the set of corners of the second
type.   Let $\x$ and $\y$ be two generators with the property that
$\x\setminus (\x\cap\y)=\x_0$ and $\y\setminus (\x\cap\y)=\y_0$.
Then, $P$ determines a homology class $B\in\pi_2(\x,\y)$.
\begin{definition}
  \label{def:PolygonConnects}
  In the above situation, we say that the fundamental polygon $P$ {\em
    connects $\x$ to $\y$ in the homology class $B$}, or simply that
  \emph{$P$ connects $\x$ to $\y$}.
  \index{connects $\x$ to $\y$}%
\end{definition}

\begin{lemma}\label{lem:polygons-represent}
  Suppose $P$ is a fundamental polygon that connects $\x$
  and $\y$ in the homology class $B$. Let $\xi$ be the chord in
  $\partial {\overline\Sigma}$ which lies on $\partial P$. Then
  $n_{\x,\y}^B=I(\x)\cdot a(-\xi)\cdot I(\y)$.
\end{lemma}
\begin{proof}
  This is an easy consequence of the definitions
  (see~\cite[Section~\ref*{LOT1:chap:type-d-mod}]{LOT1})) and a little
  complex analysis (cf.~\cite[Section
  9.5]{Rasmussen03:Knots}).
\end{proof}

\subsubsection{Gradings and the coefficient algebra}
\label{sec:coeff-algebra}
\index{$G$-set graded module}%
\index{set-graded module}%
In addition to the group-valued grading on the algebra as described in
Section~\ref{sec:gradings-alg}, the modules
$\CFDa(\HD)$ are $G$-set graded modules. This means that
there is a set $\bigGrSet(\HD)$ with a left
$G'(-\PMC)$-action and a grading function $\grb\co \Gen(\HD)\to
\bigGrSet(\HD)$ satisfying the following compatibility conditions:
if $m\in \CFDa(\HD)$ is
an $\bigGrSet(\HD)$-homogeneous element, and $a$ is a $G'(-\PMC)$-homogeneous algebra element with
$a\cdot m\neq 0$, then $a\cdot m$ is $\bigGrSet(\HD)$-homogeneous and
$\grb(a\cdot m)=\grb(a)\cdot \grb(m)$ (where
$\cdot$ means the left translation of $G'(-\PMC)$ on
$\bigGrSet(\HD)$); and if $m$ is an $\bigGrSet(\HD)$-homogeneous
element then $\bdy m$ is also $\bigGrSet(\HD)$-homogeneous 
and $\grb(\bdy m)=\lambda^{-1}\cdot \grb(m)$.
\glsit{$\bigGrSet(\HD)$}\glsit{$\grb$}%

To be more explicit about these gradings in the case of $\CFDa(\HD)$
for a Heegaard diagram~$\HD$ with boundary~$\PMC$, for any
generators $\x$, $\y$ and any $B \in \pi_2(\x, \y)$, define
\[
g'(B) \coloneqq (-e(B) - n_\x(B) - n_\y(B), \bdy^\bdy(B)) \in
  G'(\PMC),
\]
cf.~\cite[Equation~\ref*{LOT1:eq:Largeg}]{LOT1}.
Here $\bdy^\bdy(B)$ is $\bdy(B) \cap \bdy\Sigma$, the portion of $\bdy
B$ on the boundary of the Heegaard diagram.  To relate this to the
grading on the algebra $\Alg(-\PMC)$, let $R \co G'(\PMC) \to
G'(-\PMC)$ be the map induced by the orientation-reversing
map~$r\co \PMC \to -\PMC$ via the formula
$R(k,\alpha) = (k, r_*(\alpha))$.  
\glsit{$R$}%
The map $R$ is a group anti-homomorphism.
(Note that if $\alpha$ has positive multiplicities, then $r_*(\alpha)$
has negative multiplicities. In particular,
$\grb(a(-\rho))=(-1/2,-r_*(\supp(\rho)))$ while
$R(\grb(a(\rho)))=(-1/2,r_*(\supp(\rho)))$.)

The function $g'(B)$ satisfies the
crucial property that if $\vec{\rho}$ is any set of asymptotics
compatible with~$B$ and $a(-\vec{\rho}) \ne 0$ then
\[
R(g'(B)) \grb(a(-\vec{\rho})) = \lambda^{-\ind(B,\vec{\rho})};
\]
see~\cite[Lemma~\ref*{LOT1:lem:index-vs-gB}]{LOT1}.
\glsit{$S'(\HD)$}%
The grading set $S'(\HD)$ is therefore chosen in a suitable way so
that $\gr'(\x)$ and $\gr'(\y)$ are in the same $G'(-\PMC)$-orbit if
and only if
there is a domain $B \in \pi_2(\x, \y)$, and if there is such a $B$,
then
\begin{equation}\label{eq:grading-fund}
R(g'(B))\gr'(\x) = \gr'(\y)
\end{equation}
\cite[Equation~\ref*{LOT1:eq:BD-fund}]{LOT1}; this guarantees that the
grading on $\CFDa(\HD)$ is compatible with the grading on
$\Alg(-\PMC)$.
Explicitly, after fixing a base generator $\x_0$ for each
$\SpinC$-structure $\spinc$ we can set 
\[
S'(\HD,\spinc)=G'(-\PMC)/\langle R(g'(P))\mid P\in\pi_2(\x_0,\x_0)\rangle.
\] 
We then let $S'(\HD)=\amalg_{\spinc\in\SpinC(\HD)}S'(\HD,\spinc)$.
See~\cite[Chapter~\ref*{LOT1:chap:gradings}]{LOT1} for more details.

The $G'(\PMC)$-sets for the mapping cylinders will have the following
convenient property:
\begin{definition}
  \label{def:lambdaFree}
  A $G$-set $S$ is said to be \emph{$\lambda$-free} if for any $s\in
  S$ and $n\in\ZZ$, $\lambda^n\cdot s \neq s$.
\end{definition}

In this paper, we will also use an alternate way of
thinking about gradings, which we define in slightly greater
generality than we use in this paper.

\begin{definition}\label{def:based-alg}
  \index{based algebra}\index{algebra!based}%
  A \emph{based algebra} is an algebra over $\Field$ with a
  distinguished finite set of \emph{basic idempotents}, which are
  primitive, pairwise-orthogonal idempotents whose sum is the
  identity.
\end{definition}

A based algebra can also be thought of as a \dg category with a finite
number of elements.  The algebra $\Alg(\PMC)$ is a based algebra with
basic idempotents the idempotents $I(\SetS)$.

\begin{definition}\label{def:coeff-alg}
  \index{coefficient algebra}%
  \glsit{$\Coeff(M)$}%
  For a type $D$ structure $\lsup{A}M$ over a based algebra $A$, where $A$
  is graded by~$G$ and $M$ is graded by a $G$-set~$S$, the
  \emph{coefficient algebra} $\Coeff(M)$ of $M$ is the differential
  algebra spanned by triples $(\x, a, \y)$ 
  with $a\in A$ and $\x,\y \in M$
  so that
  \begin{itemize}
  \item $a = I \cdot a \cdot J$, $I \x = \x$ and $J \y = \y$ for some basic idempotents  $I$ and $J$, and
  \item there is a $k\in \ZZ$ so that $\lambda^k \gr(\x) =
    \gr(a) \gr(\y)$,
  \end{itemize}
  modulo the relations that the triples are linear in each factor:
  \begin{equation}\label{eq:coeff-alg-tens}
  \begin{split}
    (\x+\x',a,\y)&=(\x,a,\y)+(\x',a,\y) \\
    (\x,a+a',\y)&=(\x,a,\y)+(\x,a',\y) \\ 
    (\x,a,\y+\y')&=(\x,a,\y)+(\x,a,\y').
  \end{split}
  \end{equation}
  The differential is $\bdy(\x,a,\y) = (\x, \bdy a, \y)$, and
  the product is given by
  \[
  (\x_1, a_1, \y_1) \cdot (\x_2, a_2, \y_2) =
    \begin{cases}
      (\x_1, a_1 \cdot a_2, \y_2) & \y_1 = \x_2\\
      0  & \text{otherwise.}
    \end{cases}
  \]
\end{definition}

Thus each generator of $\lsup{A}M$ gives an idempotent of the coefficient algebra. For $A=\Alg(\PMC)$ and
$M=\lsup{A}\CFDa(\HD)$, this means that there is an idempotent of the
coefficient algebra corresponding to each element of $\Gen(\HD)$. The
elements of $\Coeff(\lsup{A}M)$ record algebra coefficients whose
gradings do not prevent them from
appearing in the differential,  as we see in the next lemma.

\begin{lemma}\label{lem:grading-coeff-implies}
  \index{coefficient algebra!grading on}%
  If $\lsup{A}M$ is a type $D$ structure over $A$, where $A$ is graded
  by~$G$ and $M$ is graded by~$S$, 
  let $n=\gcd\{m\in\NN\mid \lambda^ms=s\text{ for some }s\in S\}$ (or
  $0$ if $S$ is $\lambda$-free). Then
  the coefficient algebra $\Coeff(M)$ has a canonical grading
  $\gr(\x,a,\y) \in \ZZ/n\ZZ$, characterized as follows. By
  definition, there is some
  $k\in\ZZ$ such that $\lambda^k\gr(\x)=\gr(a)\gr(\y)$. Then
  $\gr(\x,a,\y)\equiv k\pmod{n}$. 

  With this grading, if $a\cdot \y$ appears in $\bdy
  \x$ then $(\x, a, \y)$ has grading~$-1$.
\end{lemma}

(Note that the divisibility of $\lambda$ in its action on~$S$ is
constant on each $G$-orbit, since $\lambda$ is central in $G$.)

\begin{proof}
  Since, by definition of $\Coeff(M)$, there is a $k$ so that
  $\lambda^k \gr(\x) = \gr(a) \gr(\y)$, this defines $\gr(\x, a, \y)$
  as an element of the cyclic subgroup of $G$ generated by $\lambda$,
  up to indeterminacy given by the divisibility of $\lambda$ in its
  action on~$\gr(\x)$.  By assumption, $n$ divides this divisibility,
  so we get a well-defined element of $\ZZ/n\ZZ$, as claimed.  It is
  elementary to check that this is a grading.
  The last statement follows from the
  assumption that $M$ is a graded differential module: $\gr(\bdy \x) =
  \lambda^{-1} \gr(\x)$.
\end{proof}

Thus, for a module~$M$ graded by a $\lambda$-free $G$-set,
$\Coeff(M)$ is $\ZZ$-graded.

If $\lsup{A}M$ has at most one generator per idempotent,
then we can view $\Coeff(M)$ as a subalgebra of~$A$.

Recall that there are two different gradings on our algebras
$\Alg(\PMC)$, one by $G(\PMC)$ and one by $G'(\PMC)$. These induce the
same grading on the coefficient algebra of any type $D$ structure over
$\Alg(\PMC)$:
\begin{lemma}\label{lem:we-dont-care}
  Let $\lsup{\Alg(\PMC)}M$ be a $G'(\PMC)$-set graded type $D$ structure over
  $\Alg(\PMC)$. Fix some collection of grading refinement data $\Xi$
  for $\PMC$. Let $\Coeff(M)$ denote the coefficient algebra of $M$ as
  a $G$-set graded module and $\Coeff'(M)$ the coefficient algebra of
  $M$ as a $G'$-set graded module. Then $\Coeff(M)=\Coeff'(M)$, and
  this identification respects the gradings on the two sides.
\end{lemma}
\begin{proof}
  Recall that for a generator $\x$ of $M$ with $\x=i\cdot \x$ for some
  minimal idempotent $i$, $\gr(\x)=\Xi(i)\cdot \grb(\x)$; and if
  $a\in\Alg(\PMC)$ is such that $j\cdot a\cdot i=a$ for minimal
  idempotents $i$ and $j$ then $\gr(a)=\Xi(j)\grb(a)\Xi(i)^{-1}$. The
  result follows.
\end{proof}

We now compute the grading on the coefficient algebra for a
Heegaard diagram more explicitly.  Loosely speaking, it is the Maslov
component of the grading on the algebra plus a correction term.

\begin{lemma}\label{lem:grading-coeff}
  Suppose $(\x,a,\y)\in\Coeff(\CFDa(\HD))$. Then there is a
  $B\in\pi_2(\x,\y)$ so that $-r_*(\bdy^\bdy(B)) = \supp(a)$. Moreover,
  for any such $B$,
  \[
  \gr(\x, a, \y) = \iota(a) - e(B) - n_{\x}(B) - n_{\y}(B). 
  \]
\end{lemma}
(The map $r_*$ appears in this lemma because $a\in\Alg(\PMC)$ where
$\PMC=-\bdy \HD$.)
\begin{proof}
  By definition of $\Coeff(\CFDa(\HD))$, $\lambda^k \grb(\x) = \grb(a)
  \grb(\y)$.  In particular, $\grb(\x)$ and $\grb(\y)$ are in the same
  $G'(-\PMC)$ orbit, so there is a domain $B$ connecting them.
  By definition of the grading on the coefficient algebra, we have
  \begin{align*}
  \lambda^{\gr(\x, a, \y)} \grb(\x) &= \grb(a) \grb(\y)\\
    &= \grb(a) R(g'(B)) \grb(\x).
  \end{align*}
  \glsit{$\gr(\x,a,\y)$}%
  By assumption, the homological components of $\grb(a)$ and $R(g'(B))$
  cancel each other, and give no correction to the Maslov component of
  the grading.  The Maslov components sum to the stated total.
\end{proof}

\begin{remark}
  The coefficient algebra is not invariant under homotopy equivalences
  of modules, as can be seen by comparing the coefficient algebra of
  an acyclic but non-zero module with that of the zero module (with no
  generators).
\end{remark}

\begin{remark}
  One can give a more abstract definition of the coefficient algebra
  as follows. Let $I$ denotes the subring of $A$ generated by the
  (distinguished) orthogonal idempotents and let $M^*$ denote the dual of
  $M$ over $I$. Then Formula~\eqref{eq:coeff-alg-tens} is equivalent
  to saying that the coefficient algebra is a subring of the tensor
  product $M^*\otimes_I A\otimes_I M$. The multiplication is induced
  by the obvious pairing $M\otimes M^*\to I$ and the multiplication on
  $A$.

  We can identify $M^*\otimes_I A\otimes_I M$ with the space
  $\Mor^A(\lsup{A}M,\lsup{A}M)$ of type $D$ structure morphisms (as
  in~\cite[Section~\ref*{LOT2:sec:cat-type-d-str}]{LOT2});
  multiplication corresponds to composition. The space
  $\Mor^A(\lsup{A}M,\lsup{A}M)$ is graded by a $\ZZ$-set. The
  coefficient algebra $\Coeff(M)$ is the subring of $\Mor^A(\lsup{A}M,\lsup{A}M)$
  generated by elements whose gradings lie in the same $\ZZ$-orbit as
  the identity map. Note, however, that the differential we have
  specified on the coefficient algebra is not induced by the
  differential on $\Mor^A(\lsup{A}M,\lsup{A}M)$.%
\end{remark}

\subsection{Bordered invariants of manifolds with two boundary components}
\label{subsec:Bimodules}
The ideas from Section~\ref{subsec:Bordered} were
extended to $3$-manifolds with two boundary components in
\cite{LOT2}. This extension
takes the form of bimodules of various types; we will focus on the
type \DD\ bimodules. The most important case for us is the case of
mapping cylinders of diffeomorphisms, though in
Section~\ref{sec:ElementaryCobordisms} we will also use elementary
cobordisms.

\subsubsection{Arced bordered Heegaard diagrams and \textalt{$\CFDDa$}{CFDD\textasciicircum}}
As explained in~\cite{LOT2}, bordered Heegaard Floer homology admits a
fairly straightforward generalization to the case of several boundary
components:

\begin{definition}
  \label{def:ArcedBordered}
  \index{Heegaard diagram!arced, bordered}%
  \index{arced bordered Heegaard diagram|see{Heegaard diagram, arced bordered}}%
  An \emph{arced bordered Heegaard diagram with two
    boundary components} (or just \emph{arced bordered Heegaard
    diagram}) is a quadruple
  \glsit{$\HD=(\overline{\Sigma},\overline{\alphas},\betas,{\mathbf z})$}%
  $\HD=(\overline{\Sigma},\overline{\alphas},\betas,{\mathbf z})$ where
  \begin{itemize}
  \item ${\overline\Sigma}$ is a
    compact surface of some genus $g$ with two boundary components,
    $\bdy_L\overline\Sigma$ and $\bdy_R\overline\Sigma$;
  \item $\betas$ is a $g$-tuple of pairwise-disjoint curves in the
    interior $\Sigma$ of $\overline{\Sigma}$;
  \item 
    \[
    \overline{\alphas}=\{\overbrace{\overline{\alpha}_1^{a,L},\dots,\overline{\alpha}_{2\Lgen}^{a,L}}^{\overline{\alphas}^{a,L}},\overbrace{\overline{\alpha}_1^{a,R},\dots,\overline{\alpha}_{2\Rgen}^{a,R}}^{\overline{\alphas}^{a,R}},\overbrace{\alpha_1^c,\dots,\alpha_{g-\Lgen-\Rgen}^c}^{\alphas^c}\}
    \]
    is a collection of pairwise disjoint embedded arcs with boundary
    on $\bdy_L\overline{\Sigma}$ (the $\overline{\alpha}_i^{a,L}$), arcs with
    boundary on $\bdy_R\overline{\Sigma}$ (the $\overline{\alpha}_i^{a,R}$), and
    circles (the $\alpha_i^c$) in the interior $\Sigma$ of
    $\overline{\Sigma}$; and
  \item ${\mathbf z}$ is a path in $\overline{\Sigma}\setminus(\overline\alphas\cup\betas)$ between
    $\bdy_L\overline{\Sigma}$ and $\bdy_R\overline{\Sigma}$.
  \end{itemize}
  These are required to satisfy:
  \begin{itemize}
  \item $\overline{\Sigma}\setminus\overline{\alphas}$ and
    $\overline{\Sigma}\setminus\betas$ are connected and
  \item $\overline{\alphas}$ intersects $\betas$ transversely.
  \end{itemize}
\end{definition}
\glsit{$\bdy_L\Sigma$, $\bdy_R\Sigma$}\glsit{${\overline{\alphas}^{a,L}}$,%
  ${\overline{\alphas}^{a,R}}$}

In this case, there are two pointed matched circles,
\begin{align*}
\PMC_L&=(\partial_L{\overline \Sigma},{\overline\alphas}^{a,L}\cap
{\partial_L{\overline\Sigma}},{\mathbf z}\cap\partial_L{\overline\Sigma}) \\
\PMC_R&=(\partial_R{\overline \Sigma},{\overline\alphas}^{a,R}\cap
{\partial_R{\overline\Sigma}},{\mathbf z}\cap\partial_R{\overline\Sigma})
\end{align*}
\index{pointed matched circle!boundary of bordered Heegaard diagram}%
\glsit{$\PMC_L$, $\PMC_R$}%

An arced bordered Heegaard diagram specifies a compact, oriented
three-manifold $Y$ with two boundary components, $\partial Y
= \partial_L Y \amalg \partial_R Y$, along with identifications
\[
  \phi_L \co F(\PMC_L) \rightarrow \partial Y \qquad\qquad
  \phi_R \co F(\PMC_R) \rightarrow \partial Y.
\] 
The data also specifies a framed arc connecting the two boundary
components of $Y$, as explained in~\cite[Section~\ref*{LOT2:sec:Diagrams}]{LOT2}, and hence
specifies $Y$ as a strongly bordered $3$-manifold.
\index{bordered $3$-manifold!specified by Heegaard diagram}%
\index{Heegaard diagram!arced, bordered!specifies strongly bordered $3$-manifold}%

To a arced bordered Heegaard diagram with two boundary components $\HD$
we associate a left $\Alg(-\PMC_L)\otimes \Alg(-\PMC_R)$ module
$\CFDDa(\HD)$, where
$\PMC_L$ and $\PMC_R$ are the pointed matched circles appearing
on the boundary of the Heegaard diagram at
$\partial_L{\overline\Sigma}$ and $\partial_R{\overline\Sigma}$
respectively.

\glsit{$\CFDDa(\HD)$}\glsit{$\Gen(\HD)$}\glsit{$o_L(\x)$, $o_R(\x)$}%
  \glsit{$I_{D,L}(\x)$, $I_{D,R}(\x)$}%
The module $\CFDDa(\HD)$ has generating set $\Gen(\HD)$ defined
exactly as in the one boundary component case.
If $\x$ is a generator, let $o_L(\x)$ (respectively $o_R(\x)$) denote the set of 
$\alpha^{L}$-arcs (respectively $\alpha^R$-arcs) occupied by $\x$.
Let $I_{D,L}(\x)$ (respectively $I_{D,R}(\x)$)
denote the idempotent in $\Alg(-\PMC_L)$
(respectively $\Alg(-\PMC_R)$)
corresponding to the complement of $o_L(\x)$ (respectively $o_R(\x)$).

Given a sequence $\vec{\rho}$ of chords in $\PMC_L\amalg \PMC_R$ let
$\vec{\rho}_L$ (respectively $\vec{\rho}_R$) denote the subsequence of
$\vec{\rho}$ consisting of the chords lying in $\PMC_L$ (respectively
$\PMC_R$). Modify Equation~\eqref{eq:Differential} as follows:
\begin{equation}
  \label{eq:DifferentialBimodule}
  n^B_{\x,\y}=\sum_{\{\vec\rho\,\mid\,\ind(B,\vec{\rho}) = 1\}}\#(\Mod^B(\x, \y; \vec{\rho}))
  a(-\vec{\rho_L})\otimes a(-\vec{\rho_R})\in\Alg(-\PMC_L)\otimes \Alg(-\PMC_R).
\end{equation}
Exactly as in Equation~\eqref{eq:Differential}, this determines a map
$$\delta^1\co X(\HD)\rightarrow \Alg(-\PMC_L)\otimes
\Alg(-\PMC_R)\otimes X(\HD),$$ 
\glsit{$\delta^1$}%
which can be used to build a differential on the 
space 
\[
\CFDDa(\HD)=\Alg(-\PMC_L)\otimes_{\Idem(-\PMC_L)}\Alg(-\PMC_R)\otimes_{\Idem(-\PMC_R)}
X(\HD).
\]

It is proved
in~\cite[Theorem~\ref*{LOT2:thm:gradedInvarianceOfBimodules}]{LOT2}
that the homotopy type of
$\CFDDa(\HD)$ depends only on the strongly bordered $3$-manifold
represented by $\HD$. So, if $Y$ is a strongly bordered $3$-manifold
then we will often write $\CFDDa(Y)$ to denote the module $\CFDDa(\HD)$
for some arced bordered Heegaard diagram $\HD$ representing $Y$.

Let $Y$ be a strongly bordered $3$-manifold with boundary
parameterized by $F(-\PMC_L)$ and $F(\PMC_R)$; i.e., a $(-\PMC_L)$-$\PMC_R$-bordered three-manifold.  Using the
identification $\Alg(-\PMC_R)=\Alg(\PMC_R)^{\op}$, we can view the
$\Alg(\PMC_L)\otimes\Alg(-\PMC_R)$-module $\CFDDa(Y)$ as an
$\Alg(\PMC_L)$-$\Alg(\PMC_R)$\hyp bimodule. When it is important to
indicate in which way we are viewing $\CFDDa(Y)$, we will write either
$\lsub{\Alg(\PMC_L),\Alg(-\PMC_R)}\CFDDa(Y)$ (for the bimodule with
two left actions) or $\lsub{\Alg(\PMC_L)}\CFDDa(Y)_{\Alg(\PMC_R)}$
(for the bimodule with one left and one right action).
\glsit{$\lsub{A,B}M$, $\lsub{A}M_B$}%

As a special case, we obtain bimodules associated to strongly-based
diffeomorphisms:
\begin{definition}
  Suppose $\psi\co F(-\PMC_1)\to F(-\PMC_2)$ is a strongly-based
  diffeomorphism. Let $Y_\psi$ denote the mapping cylinder of
  $\psi$. Then define
  \[
  \lsub{\Alg(\PMC_2)}\CFDDa(\psi)_{\Alg(\PMC_1)}=\lsub{\Alg(\PMC_2)}\CFDDa(Y_\psi)_{\Alg(\PMC_1)}.
  \]
\end{definition}
\glsit{$\lsub{\Alg(\PMC_2)}\CFDDa(Y)_{\Alg(\PMC_1)}$}%

\subsubsection{Polygons in diagrams with two boundary components}
\index{polygon, fundamental}%
\index{fundamental polygon}%
Again, fundamental polygons (in the sense of Definition~\ref{def:Polygon})
contribute to the differential. 
In the present case, chords on the boundary can be of two types:
chords contained in
$\partial_L{\overline\Sigma}$, and chords
contained in  $\partial_R{\overline\Sigma}$. We assume that, for our
polygon, there is at most one edge of each type, $\xi_L$ and $\xi_R$.
The associated algebra element is
$a_L(-\xi_L)\otimes
a_R(-\xi_R)$ if both $\xi_L$ and $\xi_R$ are present, 
$a(-\xi_L)\otimes 1$ or $1\otimes
a_R(-\xi_R)$ if only $\xi_L$ or $\xi_R$ is present, or $1\otimes 1$ if
neither is.

For a concrete example, the reader is invited to look ahead to the
bordered diagram displayed in Figure~\ref{fig:Genus2Flow}. There two
generators are indicated: $\x$ (which is indicated by black circles) and
$\y$ (which is indicated by white ones). There is a shaded
octagon from $\x$ to $\y$, which goes out to ${\partial
  {\overline\Sigma}}$ in two chords, denoted $\rho_5$ and $\sigma_3$.
This shows that $\partial \x$ contains a term $(\rho_5\otimes
\sigma_3)\otimes \y$.

This notion of polygons is still a little too restrictive for our
purpose.  In some cases, we will need to consider polygonal regions
which are obtained as unions of closures of components in ${\overline
  \Sigma}\setminus\alphas\cup\betas$. In order for such more general
polygons to contribute, we must have that $P$ is a union of $R_i$
which meet along edges which do not contain any component of $\x$ or~$\y$. 

\subsubsection{Gradings and the coefficient algebra for bimodules}
\label{sec:coeff-bimod}

As was the case for modules (Section~\ref{sec:coeff-algebra}), 
the bimodules $\CFDDa(\HD)$ are set-graded.  Suppose that $\HD$ has boundary
$\PMC_L \cup \PMC_R$.  Then there is a set $S'(\HD)$
with commuting left actions of $G'(-\PMC_L)$ and $G'(-\PMC_R)$
(in which the two actions of $\lambda$ agree), and the generators of
$\CFDDa(\HD)$ have gradings in $S'(\HD)$ which are compatible with
the differential and algebra actions in the natural
sense. If we extend $R$ to a map $R\co G'(\PMC_L)\times_\ZZ
G'(\PMC_R)\to G'(-\PMC_L)\times_\ZZ G'(-\PMC_R)$
(applying the map $R$ from Section~\ref{sec:coeff-algebra} to
both factors),
Equation~\eqref{eq:grading-fund} remains true.
\glsit{$R$}%
See~\cite[Section~\ref*{LOT2:sec:cf-gradings}]{LOT2} for further details.

(In Sections~\ref{sec:Arc-Slides} and~\ref{sec:mcg-grading}, we will
also use the analogous extension of $r_*$ from
Section~\ref{sec:coeff-algebra} to the disconnected case, gotten by
applying $r_*$ to each component. That is, $r_*\co
H_1(Z_L,\mathbf{a}_L)\times H_1(Z_R,\mathbf{a}_R)\to
H_1(-Z_L,\mathbf{a}_L)\times H_1(-Z_R,\mathbf{a}_R)$ is induced by the
orientation-reversing map $r\co (Z_L\amalg Z_r)\to -(Z_L\amalg Z_r)$.)
\glsit{$r_*$}\glsit{$r$}%

The coefficient algebra of a type \DD\ structure $\lsup{\Alg(\PMC)}M^{\Alg(\PMC')}$ is
defined as in Definition~\ref{def:coeff-alg}: $\Coeff(M)$ is generated
by triples $(\x, a_1 \otimes a_2, \y)$ where
\begin{itemize}
\item if $a_1 = I(\SetS_1)\cdot a_1\cdot I(\SetT_1)$ and $a_2 =
  I(\SetS_2) \cdot a_2 \cdot I(\SetT_2)$, then $\x = I(\SetS_1) I(\SetS_2)
  \x$ and $\y = I(\SetT_1) I(\SetT_2) \y$, and
\item there is a $k\in \ZZ$ so that $\lambda^k \grb(\x) = \grb(a_1)
   \grb(a_2) \grb(\y)$,
\end{itemize}
modulo the obvious analogue of Formula~\ref{eq:coeff-alg-tens}.
The differential is $\bdy(\x, a_1 \otimes a_2, \y) = (\x, \bdy (a_1
\otimes a_2), \y)$, and the product is
  \[
  (\x_1, a_1 \otimes b_1, \y_1) \cdot (\x_2, a_2 \otimes b_2, \y_2) =
    \begin{cases}
      (\x_1, (a_1 \cdot a_2) \otimes (b_1 \cdot b_2), \y_2) & \y_1 = \x_2\\
      0  & \text{otherwise,}
    \end{cases}
  \]
just as before.  

The rest of the theory carries through, as follows:
\begin{itemize}
\item The idempotents of the coefficient algebra correspond to the
  generating set of the type \DD\ structure $M$. In particular, if $M$
  has at most one generator per idempotent then we can view
  $\Coeff(M)$ as a subalgebra of $\Alg(\PMC)\otimes\Alg(\PMC')$.
\item As in Lemma~\ref{lem:grading-coeff-implies}, the coefficient
  algebra is graded by $\ZZ/n$, where $n$ is the divisibility of the
  kernel of the action of $\ZZ=\langle\lambda\rangle$ on the grading
  set $S$ of $M$. This action is characterized by $\gr(\x,a_1\otimes
  a_2,\y)\equiv k$ where
  \begin{equation}
    \lambda^k\gr(\x)=(\gr(a_1)\times\gr(a_2))\gr(\y).\label{eq:DD-coeff-gr}
  \end{equation}
  In particular, if $S$ is $\lambda$-free then $\Coeff(M)$ is
  $\ZZ$-graded.
\item The following analogue of Lemma~\ref{lem:grading-coeff} holds:
  \begin{lemma}\label{lem:grading-coeff-DD}
    Suppose $(\x,a_1\otimes a_2,\y)\in\Coeff(\CFDDa(\HD))$. Then there is a
    $B\in\pi_2(\x,\y)$ so that $-r_*(\bdy^\bdy(B)) = \supp(a_1)\amalg
    \supp(a_2)$. Moreover, for any such $B$,
    \[
    \gr(\x, a_1\otimes a_2, \y) = \iota(a_1)+\iota(a_2) - e(B) - n_{\x}(B) - n_{\y}(B). 
    \]
  \end{lemma}
\end{itemize}

\subsection{A pairing theorem}
\label{sec:pairing-theorem}

In this paper, we will employ a particular version of the
pairing theorem for reconstructing Heegaard Floer homology from the 
bordered invariants. Before stating it, suppose that $(Y_1,\phi_1)$ and $(Y_2,\phi_2)$ are
bordered three-manifolds with boundaries parameterized by $F(\PMC)$
and $F(-\PMC)$; i.e., we have homeomorphisms $\phi_1\co
F(-\PMC)\rightarrow \partial Y_1$ and $\phi_2\co
F(\PMC)\rightarrow \partial Y_2$.  In this case, the bordered
invariant $\CFDa(Y_1)$ is a left module over $\Alg(\PMC)$, while the
bordered invariant for $\CFDa(Y_2)$ is a left module over
$\Alg(-\PMC)\cong \Alg(\PMC)^{\op}$ and hence can be viewed as a
right module over $\Alg(\PMC)$. In particular, $\CFDa(Y_1)\otimes \CFDa(Y_2)$
can be viewed as an $\Alg(\PMC)\Hyph\Alg(\PMC)$ bimodule. 

\begin{citethm}
  \label{thm:PairingTheorem}
  Let $(Y_1,\phi_1)$ and $(Y_2,\phi_2)$ be bordered three-manifolds with
  parameterizations $\phi_1\co {-F(\PMC)}\to \partial Y_1$
  and $\phi_2\co F(\PMC)\to \partial Y_2$, and let 
  $$Y=Y_1 \sos{\partial Y_1}{\cup}{\partial Y_2} Y_2.$$ 
\index{morphisms, chain complex of}%
\index{chain complex of morphisms}%
\index{pairing theorem}%
\glsit{$\Mor$}%
  Then the chain
  complex $\CFa(Y)$ calculating $\HFa(Y)$ with coefficients in $\Zmod{2}$ is homotopy equivalent to
  \[
  \Mor(\lsub{\Alg(\PMC)}\CFDDa(\Id_{\PMC})_{\Alg(\PMC)},\lsub{\Alg(\PMC)}\CFDa(Y_1)\otimes \CFDa(Y_2)_{\Alg(\PMC)}),
  \]
  that is, the chain
  complex of $\Alg(\PMC)$-bimodule maps from
  $\lsub{\Alg(\PMC)}\CFDDa(\Id_{\PMC})_{\Alg(\PMC)}$ to $\lsub{\Alg(\PMC)}\CFDa(Y_1)\otimes \CFDa(Y_2)_{\Alg(\PMC)}$.
\end{citethm}
\begin{proof}
  This is a special case of~\cite[Corollary~\ref*{HomPair:cor:bimod-hom-pair}]{LOTHomPair}, where one
  boundary component of $Y_1$ and $Y_2$ is empty. (Note that the
  boundary Dehn twist appearing
  in~\cite[Corollary~\ref*{HomPair:cor:bimod-hom-pair}]{LOTHomPair}
  acts
  trivially on the invariant of a $3$-manifold with a single boundary
  component.)
\end{proof}

There is an analogue when $Y_2$ is a $3$-manifold with two boundary
components:
\begin{citethm}
  \label{thm:PairingTheoremModBimod}\index{pairing theorem}
  Let $(Y_1,\phi_1\co {-F(\PMC_1)}\to Y_1)$ be a bordered $3$-manifold
  with one boundary component and $(Y_2,\phi_2\co F(\PMC_1)\to
  \bdy_LY_2,\phi_3\co {-F(\PMC_2)}\to \bdy_RY_2)$ be a strongly bordered
  $3$-manifold with two boundary components. Let
  $$Y=Y_1 \sos{\partial Y_1}{\cup}{\partial_L Y_2} Y_2,$$
  a bordered $3$-manifold with boundary parameterized by $\phi_3\co
  F(\PMC_2)\to \bdy Y$.  Then $\CFDa(Y,\phi_3)$ is homotopy equivalent,
  as a differential $\Alg(\PMC_2)$-module, to
  \[
  \Mor(\lsub{\Alg(\PMC_1)}\CFDDa(\Id_{\PMC_1})_{\Alg(\PMC_1)},\lsub{\Alg(\PMC_1)}\CFDa(Y_1)\otimes \lsub{\Alg(\PMC_2)}\CFDDa(Y_2)_{\Alg(\PMC_1)}),
  \]
  the chain complex of $\Alg(\PMC)$-bimodule maps from
  $\lsub{\Alg(\PMC_1)}\CFDDa(\Id_{\PMC_1})_{\Alg(\PMC_1)}$ to
  $\lsub{\Alg(\PMC_1)}\CFDa(Y_1)\otimes
  \lsub{\Alg(\PMC_2)}\CFDDa(Y_2)_{\Alg(\PMC_1)}$.
\end{citethm}
\begin{proof}
  Again, this is a special case of~\cite[Corollary~\ref*{HomPair:cor:bimod-hom-pair}]{LOTHomPair}.
\end{proof}

In particular, for mapping classes we have:
\begin{corollary}
  \label{cor:PairingTheoremReparam}\index{pairing theorem}
  Let $(Y_1,\phi_1\co {-F(\PMC_1)}\to Y_1)$ be a bordered $3$-manifold
  with one boundary component and $\psi\co {-\PunctF(\PMC_2)}\to
  -\PunctF(\PMC_1)$ be a strongly-based diffeomorphism. Then
  $\CFDa(Y,\phi_1\circ \psi)$ is homotopy equivalent, as a differential
  $\Alg(\PMC_2)$-module, to 
  \[
  \Mor(\lsub{\Alg(\PMC_1)}\CFDDa(\Id_{\PMC_1})_{\Alg(\PMC_1)},\lsub{\Alg(\PMC_1)}\CFDa(Y_1)\otimes
  \lsub{\Alg(\PMC_2)}\CFDDa(\psi)_{\Alg(\PMC_1)}).
  \]
\end{corollary}

\begin{remark}
  The obvious analogue of Theorem~\ref{thm:PairingTheoremModBimod}
  when both $Y_1$ and $Y_2$ are strongly bordered $3$-manifolds with
  two boundary components is false. Rather, the chain complex of
  bimodule morphisms picks up an extra boundary Dehn twist;
  see~\cite[Corollary~\ref*{HomPair:cor:bimod-hom-pair}]{LOTHomPair} for more details.
\end{remark}

The isomorphisms in Theorems~\ref{thm:PairingTheorem}
and~\ref{thm:PairingTheoremModBimod} are graded isomorphisms in the
following sense. In Theorem \ref{thm:PairingTheorem}, $\CFDa(Y_1)$
(respectively $\CFDa(Y_2)$) is graded by a set $S_1$ (respectively
$S_2$) with a left (respectively right) action of $G(\PMC)$, the
(small) grading group associated to $\PMC$. The space of bimodule
homomorphisms is then graded by the set
\[
S_2\times_{G}G\times_G S_1 = S_2\times_G S_1=S_2\times S_1/[(xg,y)\sim (x,gy)].
\]
The center $\ZZ$ of $G$ acts on $S_2\times_G S_1$. As a $\ZZ$-set,
$S_2\times_G S_1$ decomposes into orbits
\[
S_2\times_G S_1=\coprod_i \ZZ/n_i.
\]
Each orbit corresponds to a $\SpinC$-structure on $Y=Y_1 \sos{\partial
  Y_1}{\cup}{\partial Y_2} Y_2$. If $\ZZ/n_i$ corresponds to the
$\SpinC$-structure $\spinc$ then $n_i=\divis(c_1(\spinc))$, and the
relative $\ZZ/n_i$-grading from $S_2\times_G S_1$ corresponds to the
relative $\ZZ/\divis(c_1(\spinc))$-grading in Heegaard Floer homology.
\index{pairing theorem!grading in}%

The story for Theorem~\ref{thm:PairingTheoremModBimod} is the same,
except that $S_1\times_G S_2$ retains a left action by $G(\PMC_2)$;
and the isomorphism of Theorem~\ref{thm:PairingTheoremModBimod} covers
an isomorphism of $G(\PMC_2)$-sets (in all orbits where the modules
are nontrivial).

See~\cite[Section~\ref*{LOT2:subsubsec:graded-pairing}]{LOT2} for further discussion in a closely
related context.


\section{The type \DD\ bimodule for the identity map}
\label{sec:DDforIdentity}

Recall that Theorem~\ref{thm:DDforIdentity} provides a model for the
type \DD\ bimodule for the identity map. The aim of the present
section is to prove that theorem. First we set up some notation.

In the standard Heegaard diagram for the identity map
(Figure~\ref{fig:Genus2Identity}), the generators are in one-to-one
correspondence with pairs of complementary idempotents, and each
domain has the same multiplicities on the left and right of the
diagram. It follows that the type \DD\ bimodule for the identity
diagram has a special form. To make this precise, think of
$\CFDDa(\Id_\PMC)$ as a left-left bimodule (compare
Section~\ref{subsubsec:opposite-algebra}). Then $\CFDDa(\Id_\PMC)$ is
induced from a module over a preferred subalgebra of
$\Alg(\PMC)\otimes\Alg(-\PMC)$:
\begin{definition}
  \index{diagonal subalgebra}%
  The {\em diagonal subalgebra} of $\Alg(\PMC)\otimes \Alg(-\PMC)$ is the
  algebra generated by
  elements of the form $(j\cdot a\cdot i)\otimes (j_o\cdot b\cdot i_o)$, where
  \begin{itemize}
  \item The support of $a$ is identified with
    the corresponding one for $b$; i.e., in the notation of 
    Section~\ref{subsec:AlgebraPMC},
    $r_*(\supp(a))=-\supp(b)$.
  \item The elements $i\in\Alg(\PMC)$ and $i_o\in\Alg(-\PMC)$ 
    are complementary idempotents.
  \item  The elements $j\in\Alg(\PMC)$ and $j_o\in\Alg(-\PMC)$ 
    are complementary idempotents.
  \end{itemize}
\end{definition}
Note that in view of the first condition above, the second two
conditions are redundant with one another. 

Some definitions for $\Alg(\PMC)$ extend in obvious ways to the
diagonal subalgebra; for instance, a \emph{basic generator} of the
diagonal subalgebra is an element $a\otimes a'$ of the diagonal
subalgebra so that $a$ and $a'$ are basic generators of $\Alg(\PMC)$
and $\Alg(-\PMC)$.
\index{basic generator!of diagonal subalgebra}%

\glsit{$\DDmod(\Id)$}%
We next rephrase the module $\DDmod(\Id)$ from the introduction as a
left-left module.
Letting ${\mathcal C}$ denote the set of connected chords for $\PMC$,
we have a map
$${\widetilde a}\colon {\mathcal C} \rightarrow \Alg(\PMC)\otimes \Alg(-\PMC)$$
(where here $\otimes$ is taken over $\Field$)
defined by
\glsit{$\widetilde{a}(\xi)$, $a(\xi)$, $a_o(\xi)$}%
$${\widetilde a}(\xi) = a(\xi)\otimes a_o(\xi),$$
where here $a(\xi)$ denotes the algebra element in $\Alg(\PMC)$
associated to $\xi$, and $a_o(\xi)$ denotes the algebra element in
$\Alg(-\PMC)$ specified by the chord $r(\xi)$.  Non-zero elements of
the form $I\cdot{\widetilde a}(\xi)\cdot J$, where $\xi$ is a chord
and each of $I$ and $J$ is a pair of complementary idempotents
(Definition~\ref{def:ComplementaryId}) in $\Alg(\PMC)\otimes
\Alg(-\PMC)$, are called {\em chord-like}. (The orientation-reversing
map $r$ is discussed in Section~\ref{subsubsec:opposite-algebra}.)

Fix a chord diagram $\PMC$, and consider the left-left
$\Alg(\PMC)$-$\Alg(-\PMC)$ bimodule $\DDmod(\Id_{\PMC})$ defined in 
Definition~\ref{def:DDmodId}.  In our present notation, the element
$A\in\Alg(\PMC)\otimes\Alg(-\PMC)$ which determines the differential can 
be written as
$$A=\sum_{\xi\in{\mathcal C}} {\widetilde a}(\xi).$$
\glsit{$A$}%
More explicitly, write a typical element of $\DDmod(\Id_{\PMC})$ as $c
\cdot I$, where $c$ is an element of $\Alg(\PMC)\otimes\Alg(-\PMC)$,
and $I=i\otimes i_o$ is a pair of complementary idempotents. The differential on
$\DDmod(\Id_{\PMC})$ is given by
$$\partial (c\cdot I) = 
(dc)\cdot I + 
\sum_{\xi\in{\mathcal C}} c\cdot {\widetilde a}(\xi)\cdot I.$$

\begin{lemma}
  \label{lem:MeetAlongBoundary}
  Let $\xi$ and $\eta$ be chords, and $I=(i,i_o)$ and $J=(j,j_o)$
  be pairs of complementary idempotents.  If $\xi$ and $\eta$ share an
  endpoint and $J\cdot {\widetilde a}(\xi)\cdot {\widetilde
    a}(\eta)\cdot I$ is non-zero, then there is a unique
  nontrivial factorization of $J\cdot {\widetilde a}(\xi)\cdot
  {\widetilde a}(\eta)\cdot I$ into homogeneous elements with
  non-trivial support in the diagonal subalgebra. Moreover,
  $J\cdot{\widetilde a}(\xi)\cdot {\widetilde a}(\eta)\cdot I$ appears
  in the differential $\partial(J\cdot {\widetilde a}(\xi\cup\eta)\cdot I)$.
\end{lemma}

\begin{proof}
  Let $a$ be a basic generator of $\Alg(\PMC)$ (i.e., one 
  represented by a strands diagram). We say that $a$ has a {\em
    break at $p$}
\index{break at $p$}%
if $p$ is the initial point of some strand in $a$
  and also the terminal point of some strand in $a$.
  Similarly, let $x\otimes x'$ be an element
  of the diagonal subalgebra, where $x$ and $x'$ 
  are basic generators. We say that $x\otimes x'$
  has a break at $p$ if either $x$ has a break at $p$ or
  $x'$ has a break at $r(p)$.
  \index{break at $p$}%

  Let ${\widetilde x}=x\otimes x'$ and ${\widetilde y}=y\otimes y'$ be
  a pair of basic generators in the diagonal subalgebra, with non-zero
  product.  Suppose that there is some position $p$ in the boundary of
  the support of both $x$ and $y$. We claim then that ${\widetilde
    x}\cdot {\widetilde y}$ has a break at $p$. There are two cases:
  either $p$ is an initial endpoint of the support of $x$, or it is a
  terminal endpoint of the support of $x$. If $p$ is an initial endpoint in the support of
  $x$ then it must also be a terminal endpoint of the support of $y$;
  thus, $x\cdot y$ has a break at $p$. Symmetrically, if $p$ is a
  terminal endpoint of the support of $x$, then it must be an initial
  endpoint in the support of $x'$, and hence $x'\cdot y'$ has a break
  at $r(p)$.  

  \begin{figure}
    \begin{center}
      \input{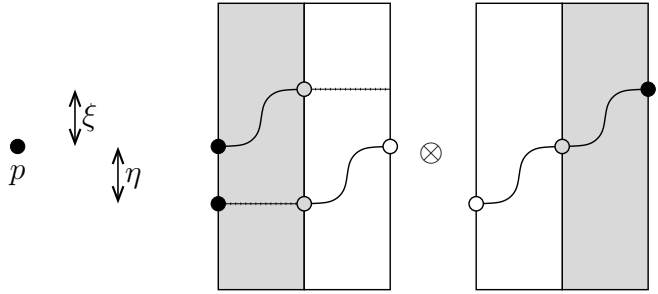}
    \end{center}
    \caption {{\bf Break in product.}
      \label{fig:BreakInProduct}
      The algebra element ${\widetilde a}(\xi)$ is contained in the two
      shaded boxes; the algebra element ${\widetilde a}(\eta)$ is contained
      in the two unshaded boxes. The chords $\xi$ and $\eta$ share a
      boundary point $p$.  The product $\widetilde{a}(\xi)\cdot \widetilde{a}(\eta)$ is gotten from the
      illustrated juxtaposition, and has a break at $p$. Note that the juxtaposition appears in
      $d{\widetilde a}(\xi\cup\eta)$.}
  \end{figure}
  
  Suppose now that $\xi$ and $\eta$ are chords which share an endpoint.
  Then, $t=J\cdot {\widetilde a}(\xi)\cdot{\widetilde
    a}(\eta)\cdot I$ has a unique break. Now, consider a
  factorization of $t$ into ${\widetilde x}\cdot {\widetilde y}$
  in the diagonal subalgebra. As in the previous paragraph, $t$ must have a
  break at any point $q$ where the support of ${\widetilde x}$ and
  ${\widetilde y}$ meet; since the product has a unique break,
  there must be a single such point $q$, and so $q$ agrees with $p$. From
  this, it is straightforward to see that the factorization coincides
  with the initial one, i.e., ${\widetilde x}=J \cdot {\widetilde a}(\xi)$
  and ${\widetilde y}={\widetilde a}(\eta)\cdot I$.

  To see that $J\cdot\widetilde{a}(\xi)\cdot\widetilde{a}(\eta)\cdot
  I$ appears in the differential of $J\cdot
  \widetilde{a}(\xi\cup\eta)\cdot I$, 
  suppose without loss of generality that the terminal point $p$ of
  $\xi$ coincides with the initial point of $\eta$, so that $j\cdot
  a(\xi)\cdot a(\eta)\cdot i$ has a break. In this case, 
  $a(\xi)\cdot a(\eta)$ appears in the
  differential $d a(\xi\cup\eta)$, and 
  ${\widetilde a}(\xi)\cdot {\widetilde a}(\eta)$ appears in
  $d{\widetilde a}(\xi\cup\eta)$.
\end{proof}

\begin{lemma}
  \label{lem:CancelInPairs}
  Let $I=(i,i_o)$ and $J=(j,j_o)$ be pairs of complementary 
  idempotents and $\xi$, $\eta$ be chords.
  If $\xi$ and $\eta$ do not share an endpoint then
  $$J\cdot {\widetilde a}(\xi)\cdot {\widetilde a}(\eta)\cdot I=
  J\cdot {\widetilde a}(\eta)\cdot {\widetilde a}(\xi)\cdot I.$$
\end{lemma}

\begin{proof}
  There are three cases on the endpoints of $\xi$ which are handled
  differently: the boundaries can be linked, an endpoint of $\xi$ can
  be matched with an endpoint of $\eta$, and the
  case where neither of the above two phenomena occur.  

  \begin{figure}
    \begin{center}
      \input{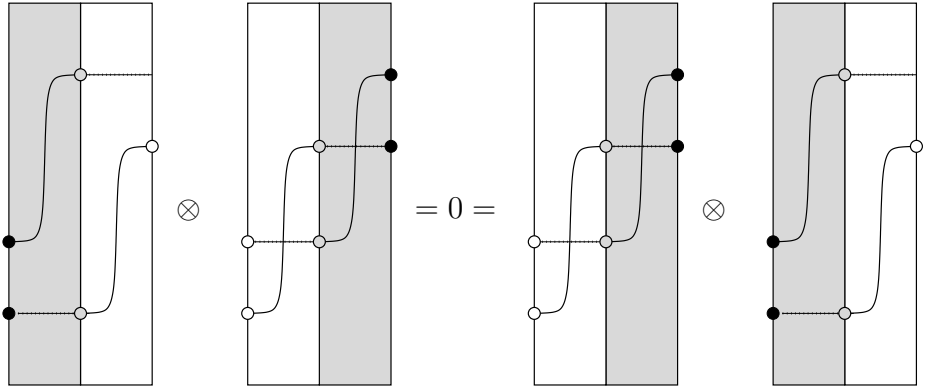}
    \end{center}
    \caption {{\bf Boundaries of $\xi$ and $\eta$ are linked.}
      \label{fig:CommuteInDiagAlg}
      We have illustrated here chords $\xi$ and $\eta$ whose boundaries are
      linked.  Both products ${\widetilde a}(\xi)\cdot {\widetilde
        a}(\eta)=0={\widetilde a}(\eta)\cdot {\widetilde a}(\xi)$, since one
      of the two sides always has a double-crossing in it.}
  \end{figure}
  
  In the first case, where the boundaries of $\xi$ and $\eta$ are linked
  (i.e., exactly one of the endpoints of $\xi$ is contained in the
  interior of $\eta$), we claim that ${\widetilde a}(\xi)\cdot
  {\widetilde a}(\eta)=0$.  To see this, write $x=i\cdot a(\xi)$,
  $x'=i_o\cdot a_o(\xi)$, $y=j \cdot a(\eta)$ $y_o=j_o\cdot
  a(\eta)$. Observe that all of these are basic generators. Next, note
  that exactly one of the juxtapositions $x * y$ or $x_o * y_o$
  contains a double-crossing. Thus, ${\widetilde a}(\xi)\cdot {\widetilde a}(\eta)=0$.
  By the same reasoning, ${\widetilde
    a}(\eta)\cdot {\widetilde a}(\xi)=0$.

  Consider next the second case, where some endpoint of $\xi$ is matched
  with some endpoint of $\eta$. If the initial boundary of $\xi$ is
  matched with the initial point of $\eta$ then
  $a(\xi)\cdot a(\eta)=a(\eta)\cdot a(\xi)=0$, so once again both
  terms vanish. The same reasoning applies if the terminal
  points are matched. Finally, suppose that
  the terminal
  point $p$ of $\xi$ is matched with the initial point $p'$ of
  $\eta$. Thus, $a(\xi)\cdot a(\eta)=0$. However, it is
  not guaranteed that $a(\eta)\cdot a(\xi)=0$. But for $j\cdot
  a(\eta)\cdot a(\xi)\cdot i$ to be nonzero, both
  $i$ and $j$ must contain the matched pair $\{p,p'\}$, and hence
  neither $i_o$ nor $j_o$ contains the
  matched pair $\{p,p'\}$. It follows that $j_o\cdot
  a_o(\eta)\cdot a_o(\xi)\cdot i_o=0$, and hence ${\widetilde
    a}(\eta)\cdot {\widetilde a}(\xi)=0$.

  In the third case, where none of the endpoints of $\xi$ are
  matched to endpoints of $\eta$, and their boundaries are unlinked,
  we have that $a(\xi)\cdot a(\eta)=a(\eta)\cdot a(\xi)$ and
  $a_o(\xi)\cdot a_o(\eta)=a_o(\eta)\cdot a_o(\xi)$, so the stated
  equality holds.
\end{proof}

It follows from the identification of $\DDmod(\Id_{\PMC})$ with
$\CFDDa(\Id_{\PMC})$ that $\partial^2=0$ on
$\DDmod(\Id_{\PMC})$. Although it is not strictly necessary for this
paper, it is not hard to give a combinatorial proof of this fact:
\begin{proposition}
  \label{prop:DDsquaredZero}
  The endomorphism $\partial$ on $\DDmod(\Id_{\PMC})$ is a
  differential, i.e., $\partial^2=0$.
\end{proposition}

\begin{proof}
  Let $I=(i,i_o)$.
  The $J=(j,j_o)$ coefficient of $\partial^2 I$ is given by
  \begin{align*}\sum_{\xi_1,\xi_2\in{\mathcal C}} &
    (j\cdot a(\xi_2)\cdot a(\xi_1)\cdot i) \otimes
    (j_o\cdot a_o(\xi_2)\cdot a_o(\xi_1)\cdot i_o) \\
    & + d \sum_{\xi\in{\mathcal C}} (j\cdot a(\xi)\cdot i) \otimes
    (j_o\cdot a_o(\xi)\cdot i_o).
  \end{align*}
  By Lemma~\ref{lem:CancelInPairs}, many of the terms in the first sum
  cancel in pairs (and some are individually zero), leaving only those
  terms where the $\xi_1$ and $\xi_2$ are two chords which share an
  endpoint. Indeed, according to Lemma~\ref{lem:MeetAlongBoundary}
  all such terms are left over in the sum. The last statement in
  Lemma~\ref{lem:MeetAlongBoundary} implies that all of these terms
  occur in the second sum; moreover, that these
  terms exactly cancel with the second sum.
\end{proof}

We investigate some of the algebraic properties of the diagonal subalgebras.
These will be useful when considering gradings.
It is not true that any homogeneous element in $\Alg(\PMC)$
can be factored as a product of chord-like elements. The corresponding
fact is, however, true for the diagonal subalgebra:

\begin{lemma}
  \label{lem:FactorDiagonalSubalgebra}
  Any homogeneous element of the diagonal subalgebra can be factored
  as a product of chord-like elements.   
\end{lemma}

\begin{proof}
  Fix a basic element of the diagonal subalgebra $\xi \otimes \xi'$
  with non-trivial support.  It suffices to 
  factor off a chord-like element; i.e., exhibit a
  factorization $\xi\otimes \xi' = (\eta\otimes \eta')\cdot
  (\zeta\otimes \zeta')$, where $\eta\otimes\eta'$ is also in the
  diagonal subalgebra and $\zeta\otimes\zeta'$ is
  chord-like. (Since this expresses $\xi\otimes\xi'$ as a product of
  an element which has strictly smaller support with a chord-like
  element, induction on the total support then gives the factorization
  claimed in the lemma.)

  To factor off the chord-like element, we proceed as follows. Let $s$
  be a moving strand in $\xi$ and $t$ be a moving strand in
  $\xi'$. For the strand $s$,  the initial point $s^-$ and the terminal point $s^+$
  are points in $\PMC$ with $s^- < s^+$; thinking of the initial and
  terminal points of $t$ ($t^-$ and $t^+$ respectively) as points in
  $\PMC$, as well, we have that $t^+ < t^-$. Now, consider pairs of
  such strands $(s,t)$ with the property that $t^+ < s^+$; and choose
  among these a pair of strands for which the distance $s^+ - t^+$ is
  minimized. We call such a strand pair $(s,t)$ {\em minimal}.

  For a minimal strand pair $(s,t)$, we claim that $s^- \leq t^+ < s^+
  \leq t^-$. This follows from the condition on the support of
  elements of the diagonal subalgebra. Specifically, if the stated
  inequalities do not hold, then either $t^-<s^+$ or $t^+<s^-$.  The
  case where $t^-<s^+$ can be divided into two subcases, depending on
  whether or not the support of $\xi'$ jumps at $t^-$, i.e.,
  $\bdy \supp(\xi')$ has a non-zero coefficient at~$t^-$.
  (Recall that $\supp(\xi) = \supp(\xi')$.) If the support of $\xi'$
  jumps at $t^-$, there is a different strand $u$ in $\xi$
  ending at $t^-$. But in this case $(u,t)$ satisfies $t^+<u^+<s^+$,
  contradicting minimality of $(s,t)$. If the support of $\xi'$ does
  not jump at~$t^-$, there must be a different strand $v$ in
  $\xi'$ which ends in~$t^-$. But in this case $(s,v)$ satisfies
  $t^+<v^+<s^+$, again contradicting minimality of $(s,t)$.  The case where
  $t^+<s^-$ is excluded similarly.   

  For a minimal strand pair $(s,t)$, let $\zeta$ be the algebra
  element in $\Alg(\PMC)$ associated to a single moving strand from
  $t^+$ to $s^+$, and whose terminal idempotent coincides with that
  for~$\xi$. Since $t^+$ appears in the terminal idempotent of $\xi'$,
  it does not appear in the terminal idempotent of $\xi$, and hence we
  can consider the algebra element $\eta$ which consists of all the
  moving strands in $\xi$, except that the strand $s$ is terminated at
  $t^+$ instead of at $s^-$. Now, $\eta\cdot \zeta=\xi$,
  unless $\eta\cdot \zeta=0$ because of the introduction of a
  double-crossing.  But $\eta\cdot \zeta$ cannot introduce a double-crossing,
  for that would mean that there is some strand $u$ in
  $\xi$ with $t^+<u^+<s^+$, contradicting the minimality of the
  strand pair $(s,t)$. Similarly, we can factor $\xi'=\eta'\cdot
  \zeta'$, where $\zeta'$ is the algebra element consisting of a
  single moving strand from $s^+$ to $t^+$ and whose terminal
  idempotent coincides with that for $\xi'$. So, $\zeta\otimes
  \zeta'$ is the desired chord-like element. 
\end{proof}

\begin{figure}
  \begin{center}
    \input{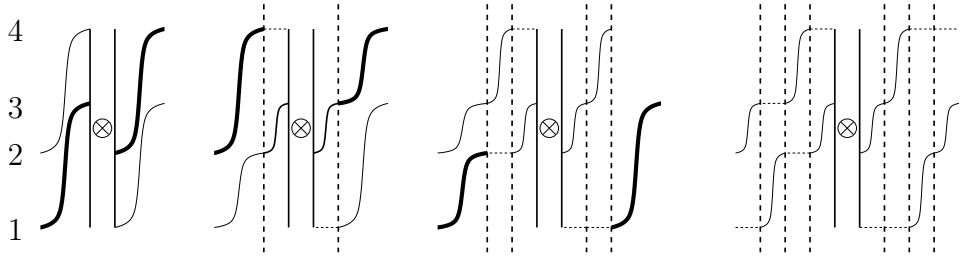}
  \end{center}
  \caption {{\bf Factoring into chord-like elements.}
    \label{fig:FactorAlgorithm}
    We start with the element consisting of two strands, one from $1$
    to $3$ and another from $2$ to $4$ on the $\PMC$ side, and a
    similar element on the $\PMC'$ side, as illustrated on the left.
    We then apply successively the algorithm from the proof of
    Lemma~\ref{lem:FactorDiagonalSubalgebra}, to factor off chord-like elements.
  At each stage, a minimal strand pair is illustrated with darker strands.}
\end{figure}

The above proof in fact gives an algorithm for factoring any given element of the diagonal subalgebra as a product of
chord-like elements. For an illustration of this algorithm, see Figure~\ref{fig:FactorAlgorithm}.

\begin{definition}
  \label{def:SpecialLengthThree}
  A chord $\xi$ is said to be a {\em special length three chord} if
  the following three conditions hold:
  \begin{itemize}
  \item $\xi$ has length three 
  \item  the terminal point $p$ of $\xi$ is matched
    with some other point $p'$ in the interior of $\xi$, and 
  \item the initial point $q$ of $\xi$ is matched with another point
    $q'$ in the interior of $\xi$.
  \end{itemize}
\end{definition}

\begin{figure}
    \begin{center}
      \input{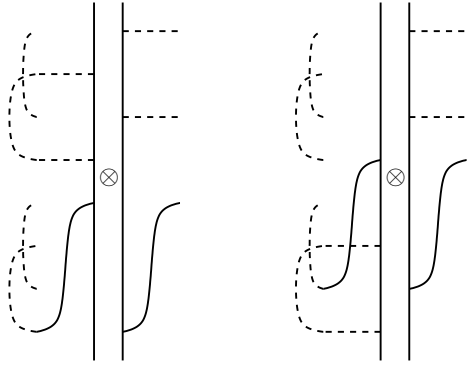}
    \end{center}
    \caption {{\bf Length three chords.}
      \label{fig:LengthThreeChord}
      Left: a special length three chord. This chord is a cycle. Right:
      a length three chord which is not special. This chord is not a cycle.}
  \end{figure}

See Figure~\ref{fig:LengthThreeChord} for an illustration.

\begin{lemma}\label{lem:special-3-special}
  If $\xi$ is a chord which is not length $1$ or special length $3$,
  and $I\widetilde{a}(\xi)J\neq 0$, then
  $d(I\widetilde{a}(\xi)J)\neq 0$.
\end{lemma}
\begin{proof}
  By hypothesis, there must be some position $p$ in the interior of
  support of $\xi$ which is not matched with either an initial or
  terminal point of $\xi$. The initial idempotent $I=i\otimes
  i_o$ must contain this position either in $i$ or in $i_o$ (since $i$
  and $i_o$ are complementary idempotents). Thus the differential of
  either $a(\xi)$ or $a_o(\xi)$ must contain a term corresponding to
  the resolution at $p$.
\end{proof}

Not every pointed matched circle has special length three chords. 
For instance, there are none in the antipodally matched circle
for a surface of genus $k>1$. (The \emph{antipodally matched circle}
\index{matched circle!antipodally}%
is the one where the matched pairs of points are antipodal on the
circle.) On the other hand, for the surface of genus $k=1$, there is a
unique pointed matched circle, and it does have a special length three
chord, so the split pointed matched circle with any genus has
special length $3$ chords.

\begin{proposition}
  \label{prop:DDidUnique}
  Let ${M}$ be any type \DD\  bimodule over $\Alg(\PMC)\otimes\Alg(-\PMC)$,
  where $\PMC$ is any pointed matched circle with genus greater than one.
  Suppose ${M}$ satisfies the following properties:
  \begin{enumerate}
    \item Generators of $M$ are in one-to-one correspondence
      with pairs $I=(i,i_o)$ of
      complementary idempotents, in such a manner that
      if $\x(I)$ is the generator corresponding to $I=(i,i_o)$,
      then $I\cdot \x(I)=\x(I)$.
    \item The coefficients in the differential of 
      ${M}$ all lie in the diagonal subalgebra
      of $\Alg(\PMC)\otimes\Alg(-\PMC)$; i.e.,
      writing
      $$\partial \x(I)=\sum_{J} a(I,J)\cdot \x(J),$$
      where $a(I,J)\in\Alg(\PMC)\otimes\Alg(-\PMC)$, 
      then all of the $a(I,J)$ lie in the diagonal subalgebra
      of $\Alg(\PMC)\otimes\Alg(-\PMC)$.
    \item 
      \label{cond:NoZtorsion}
      $M$ is graded by a $\lambda$-free $G$-set $S$.
    \item The differential of $\x(I)$ contain all non-zero elements of
      the form $I\cdot {\widetilde a}(\xi)\cdot J\cdot \x(J)$ where
      $\xi$ is any length one chord.
  \end{enumerate}
  Then, ${M}$ is isomorphic to the bimodule $\DDmod(\Id_\PMC)$.
\end{proposition}

\begin{proof}
  First, we prove by induction on the length of the support of
  $\xi\in{\mathcal C}$ that all terms of the form $a(\xi)\otimes
  a_o(\xi)$ appear in the coefficient of the differential. 

  Assume for simplicity that our diagram does not have any special
  length three chords; we return to the general case later.  Our
  goal is to argue that, if $\xi$ is a chord of length $n>1$, then
  ${\widetilde a}(\xi)$ appears in the differential. We proceed by
  induction on the length of the support of $\xi$. Consider the
  part of the coefficient of $\partial^2\x(I)$ which is spanned
  by algebra elements whose supports
  coincides with the support of $\xi$.
  By the inductive hypothesis, the $\x(J)=\x(j,j_o)$ coefficient of
  $\partial^2 \x(I)$ with total support $\xi$ has the form
  $$\Biggl(\sum_{\{\xi_1,\xi_2\in{\mathcal C} \vert \xi_1\cup \xi_2=\xi\}}\!\!\!\!\!
    I\cdot {\widetilde a}(\xi_2) \cdot {\widetilde a}(\xi_1) \cdot J
  \Biggr) + \left(\sum I \cdot {\widetilde x}\cdot {\widetilde y}\cdot
    J\right) + d (I \cdot {\widetilde z}\cdot J),$$
  where ${\widetilde
    z}$ is the $\x(J)$ component of $\partial \x(I)$ with support
  equal to $\xi$, ${\widetilde x}$ and ${\widetilde y}$ are
  other basic elements in the diagonal algebra where at least one of
  $\widetilde{x}$ or $\widetilde{y}$ has at least one break (and
  $\supp({\widetilde x})+\supp({\widetilde y})=\supp(\xi)$). 

  According to
  Lemma~\ref{lem:MeetAlongBoundary}, terms appearing in the first sum
  cannot cancel with other terms appearing in the first sum or with
  terms appearing in the second sum;
  thus, they must cancel with terms in the third. 
  Moreover, as in the proof of Proposition~\ref{prop:DDsquaredZero}, every
  term in $d(I\widetilde{a}(\xi)J)$ occurs in the
  first sum; in particular, since $d(I\widetilde{a}(\xi)J)$ is
  nontrivial by Lemma~\ref{lem:special-3-special}, the first sum is
  nontrivial.  This forces
  $d(I\cdot {\widetilde z}\cdot J)$ to contain the terms in $d(I\cdot
  {\widetilde a}(\xi) \cdot J)$ with non-zero multiplicity. 
  It follows that $I\cdot {\widetilde z}\cdot J$ must contain
  $I\cdot a(\xi)\cdot J$ with non-zero multiplicity: 
  as in the proof of Lemma~\ref{lem:MeetAlongBoundary},
  the non-zero terms in the first sum have three non-horizontal strands in 
  them (two on one side and one on the other), so if they appear 
  in the differential of a homogeneous element, then that element
  must have exactly two moving strands in it, i.e., it must be of the
  form $I\cdot {\widetilde a}(\xi)\cdot J$.

  Having verified that the differential contains
  $\sum_{\xi\in {\mathcal C}} {\widetilde a}(\xi)$, we must verify
  that it contains no other terms.  This follows from grading reasons.
  Having established that if $I\cdot {\widetilde a}(\xi)\cdot J\neq
  0$ then $I\cdot {\widetilde a}(\xi)\cdot \x(J)$ appears in
  $\partial \x(I)$, we can conclude that
  \begin{equation}
    \label{eq:GradedEquality}
    \lambda^{-1}\grb(\x(I))=\grb(I\cdot {\widetilde a}(\xi)\cdot J)\cdot \grb(\x(J)),
  \end{equation}
  for any chord $\xi$.
  Let $a$ be any basic generator of the diagonal subalgebra, and
  suppose that $a\x(J)$ occurs in $\bdy\x(I)$. Then, in particular,
  $\grb(a)\grb(\x(J))=\lambda^{-1}\grb(\x(I))$.
  According to Lemma~\ref{lem:FactorDiagonalSubalgebra}, there
  is a sequence of chords $\{\xi_i\}_{i=1}^n$ with the property that
  $a=I\cdot \prod_{i=1}^n {\widetilde a}(\xi_i)\cdot J$. 
  By Equation~\eqref{eq:GradedEquality},
  $$\lambda^{-n}\cdot \grb(\x(I)) = \grb(a)\cdot \grb(\x(J));$$
  it follows that $n=1$ (thanks in part to Condition~\ref{cond:NoZtorsion}), and
  hence $a$ had to be a chord-like element.

\begin{figure}
\begin{center}
\input{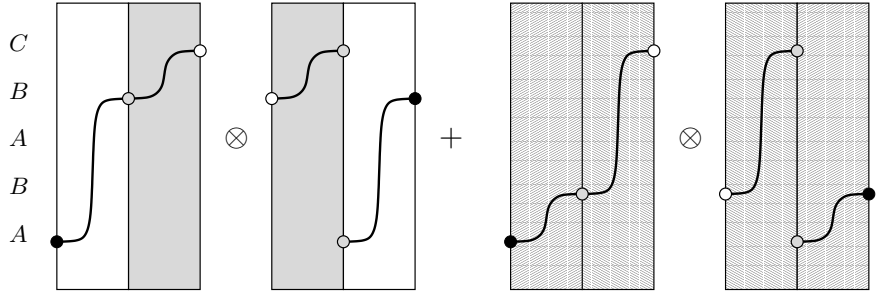}
\end{center}
\caption {{\bf Special length three chords.}
\label{fig:SpecialLengthThree}
The matching is indicated by the letters on the rows. The algebra element in the white box
represents the algebra element associated to a special length three chord; we multiply
it by the gray-shaded length one chord, to get a term which evidently appears in the differential
of the algebra element associated to the length four chord. The other term is illustrated to the
right (in the two hatched boxes). It decomposes as a product of length three and length one chords,
but now the length three chord is not special.}
\end{figure}

  This completes the proof of the proposition for pointed matched
  circles without special length three chords. When there are special length
  three chords $\xi$, we must show ${\widetilde a}(\xi)$ also appears
  in the differential.  This is seen by considering a length
  four chord $\eta$ which contains the given special length three
  chord $\xi$. Such a chord can be found since the genus $k$ is bigger
  than~$1$. 
  Now we have that
  $$d {\widetilde a}(\eta) = {\widetilde a}(\xi_1)\cdot {\widetilde a}(\xi_3) + 
  {\widetilde a}(\xi_3')\cdot {\widetilde a}(\xi_1'),$$ where $\xi_i$
  and $\xi_i'$ have length $i$, and exactly one of $\xi_3$ or $\xi_3'$
  $\xi$, while the other of $\xi_3$ or $\xi_3'$ is not
  special. Suppose for concreteness that $\xi=\xi_3$.
  Our inductive hypothesis ensures that ${\widetilde a}(\xi_3')$
  appears in the differential, and hence the term ${\widetilde
    a}(\xi_3')\cdot {\widetilde a}(\xi_1')$, appears in $\partial^2$,
  and must cancel against something.  According to
  Lemma~\ref{lem:MeetAlongBoundary}, the only term it can cancel is
  $d{\widetilde a}(\eta)$. Hence, we have established that
  ${\widetilde a}(\eta)$ appears with non-zero multiplicity in
  $\partial$. Thus, since $d{\widetilde a}(\eta)$ appears in
  $\partial^2$, we see that ${\widetilde a}(\xi)\cdot {\widetilde
    a}(\xi_3)$ appears in $\partial^2$. According to
  Lemma~\ref{lem:MeetAlongBoundary}, the only way this can cancel is
  if ${\widetilde a}(\xi_3)$---the algebra element associated to our
  special length three chord---also appears with non-zero
  multiplicity in the differential. Repeating this argument for each
  special length $3$ chord, we conclude that all length three chords
  appear in the differential. Now, we can proceed with the same
  induction as before.
\end{proof}

\begin{proof}[Proof of Theorem~\ref{thm:DDforIdentity}]
Consider the standard genus $2k$ Heegaard diagram for the identity
map of $S$, as pictured in Figure~\ref{fig:Genus2Identity}.
We verify that $\CFDDa$ of this Heegaard diagram satisfies the hypotheses
of Proposition~\ref{prop:DDidUnique}.

\begin{figure}
\begin{center}
\input{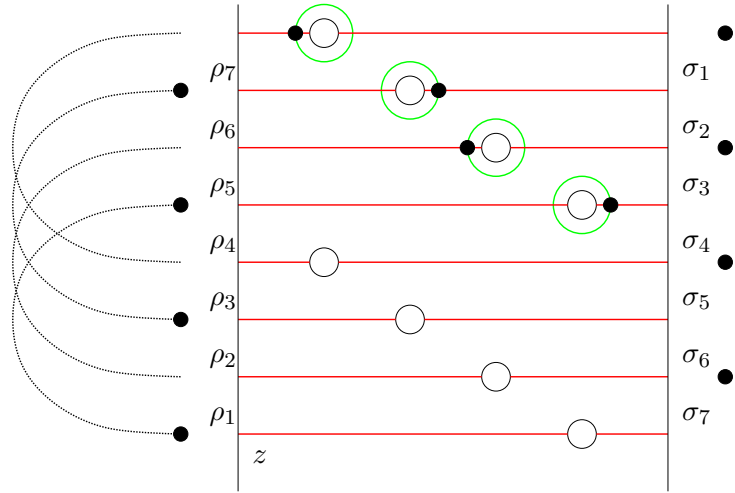}
\end{center}
\caption {{\bf Heegaard diagram for the identity map.}
\label{fig:Genus2Identity}
This is a Heegaard diagram for the identity cobordism of the genus two diagram with antipodal 
matching, as indicated by the arcs to the left of the diagram. To the left and the right of the
diagram, we have also
indicated a pair of complementary idempotents, along with its unique extension into the
diagram as a generator for the complex.}
\end{figure}

The generators $\Gen(\HD)$ of the
$\Alg(-\PMC)\oplus\Alg(\PMC)$-module $\CFDDa(\HD)$ are in one-to-one
correspondence with pairs of complementary idempotents, as follows. 
There are no $\alpha$-circles, and each $\beta_i$ meets exactly two
$\alpha$-arcs, $\alpha_i^L$ and $\alpha_i^R$, where we label the
$\alpha$-curves so that if $\alpha_i^L$
is the $\alpha$-arc which meets $\partial_L\HD$ in $\{p_i, q_i\}$ for $\PMC_L$ then $\alpha_i^{R}$ is
the $\alpha$-arc which meets $\partial_R\HD$ in $\{r(p_i),r(q_i)\}$.
The map $\x \mapsto I_{D,L}(\x)\times I_{D,R}(\x)$ sets
up the one-to-one correspondence.

Next, we claim that the coefficients of the boundary operator
$$\delta^1\co X\rightarrow (\Alg(\PMC)\oplus\Alg(-\PMC))\otimes X$$
take values in the diagonal subalgebra. We have already
checked this on the level of idempotents. To see that coefficients lie
in the diagonal subalgebra, notice that if $q^L_i$ is a position
between two consecutive places $p_i$ and $p_{i+1}$ on $\PMC_L$, and
$q^R_i$ is a position between the corresponding places $r(p_i)$ and
$r(p_{i+1})$ on $\PMC_R=-\PMC_L$, then $q_i^L$ and $q_i^R$ can be
connected by an arc in ${\overline\Sigma}$ which does not cross any of
the $\alpha$- or $\beta$-circles. Thus, if $B\in\pi_2(\x,\y)$ is any
homotopy class then the local multiplicity of $B$ at $q_i^L$
coincides with the local multiplicity of $B$ at $q_i^R$. It follows
that the coefficients for $\CFDDa(\Id_{\PMC})$
lie in the diagonal subalgebra.

A periodic domain in $\HD$ is uniquely determined by its local
multiplicities at the boundary of the Heegaard diagram. From this,
and the definition of the grading set for bimodules, it
follows that the grading set for the identity bimodule is $\lambda$-free.

If $\xi$ is a length one chord, the domains contributing to the
$J\cdot {\widetilde a}(\xi)\cdot I$ component of $\partial I$ are all
octagons; see Figure~\ref{fig:Genus2Flow}. These have holomorphic
representatives for any conformal structure, cf.~Lemma~\ref{lem:polygons-represent}.

\begin{figure}
\begin{center}
\input{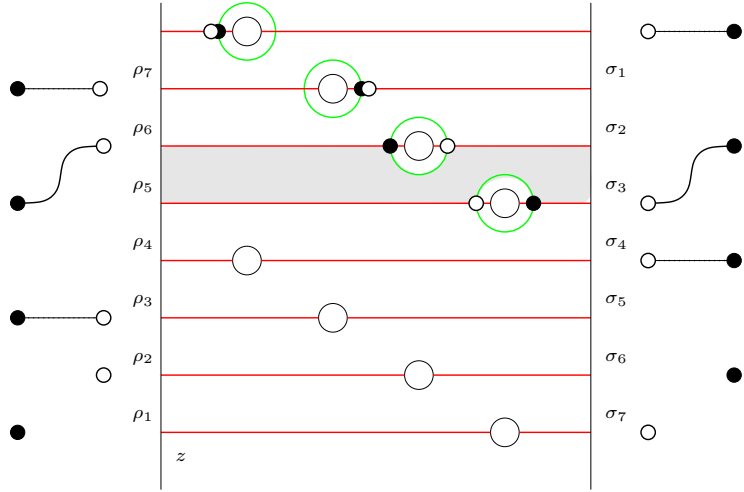}
\end{center}
\caption {{\bf A differential.}
\label{fig:Genus2Flow}
Here is another picture of the Heegaard diagram for the identity
cobordism on the antipodally
matched circle with genus two.  The shaded rectangle region represents a term
in the differential of the black idempotent, which has the form
$\rho_5\otimes \sigma_3$ times the white idempotent.}
\end{figure}

If the genus $k$ is bigger than $1$ then we have
verified that ${M}=\CFDDa(\Id_\PMC)$ satisfies all the hypotheses of
Proposition~\ref{prop:DDidUnique}. So, that proposition implies that
$\CFDDa(\Id_\PMC)$ is isomorphic to $\DDmod(\Id_\PMC)$, as desired.

The case that the genus $k=1$ is established as follows.  Embed the
given genus $1$ diagram for $\PMC$ inside a genus two diagram for
$\PMC\#\PMC_0$ (where here $\PMC_0$ is another genus one
diagram). This is shown in Figure~\ref{fig:EmbedGenusOne}. The diagram
defines a type $D$ structure over $\Alg(\PMC\#\PMC_0)\otimes
\Alg(-\PMC\#\PMC_0)$.  Setting all the algebra elements with
non-trivial support in $\PMC_0$ to zero, and restricting to elements with
some fixed idempotent $I(\SetS_0)$ in $\PMC_0$ and its complementary
idempotent in $-\PMC_0$, we obtain an induced module over
$\Alg(\PMC)\otimes\Alg(-\PMC)$. (In the figure, this corresponds to
setting $\rho_4=\rho_5=\rho_6=\rho_7=
\sigma_4=\sigma_5=\sigma_6=\sigma_7=0$ and restricting to generators
whose coordinates in the portion corresponding to $\PMC_0$ and
$-\PMC_0$ are at the displayed black dots.) The holomorphic
curve counts in this portion of the diagram correspond exactly to
holomorphic curve counts for the smaller genus one diagram gotten by
excising the portion corresponding to
$\Alg(\PMC_0)\otimes\Alg(-\PMC_0)$. (Algebraically, this corresponds to the
statement that the type \DD\ identity bimodule for the genus one
diagram coincides with the induced module
$Q_*\DDmod(\Id_{\PMC\#\PMC_0})$, compare Definition~\ref{def:StableModule}.)
Thus, the genus $1$ case follows from the genus $2$ case.

\begin{figure} 
  \begin{center}
    \input{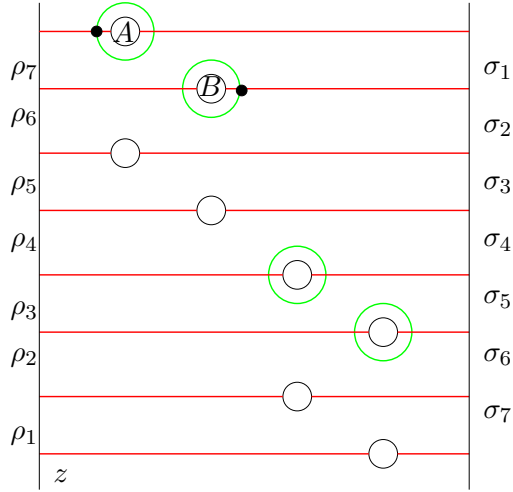}
  \end{center}
  \caption {{\bf Split genus two diagram.}
    \label{fig:EmbedGenusOne}
    A genus one diagram involving $\rho_i$ for $i=1,\dots, 3$
    and $\sigma_j$ $j=5,\dots,7$ is stabilized by adding another genus
    one diagram and, for example, the fixed generator.}
\end{figure}

The fact that the homotopy equivalences $\DDmod(\Id_\PMC)\simeq
\CFDDa(\Id_\PMC)$ are canonical, i.e., unique up to homotopy, follows
from~\cite[Corollary~\ref*{LOT2:Cor:MCG-Equivalences} and
Lemma~\ref*{LOT2:lem:bimodule-rigidity}]{LOT2}.
\end{proof}

\begin{remark}
  The grading set for the identity \DD\ bimodule
  can be explicitly determined from the Heegaard diagram, and is given
  as follows.
  Recall from Section~\ref{sec:coeff-algebra} that the map
  $$R\co \bigGroup(-\PMC)\to \bigGroup(\PMC)^{\op}$$
  defined by $R(k,\alpha)=(k,r_*(\alpha))$ is a group
  isomorphism.
  Using this, the grading set for $\CFDDa(\Id)$ is $\bigGroup(\PMC)$
  with the structure of a left
  $\bigGroup(\PMC)\times_{\ZZ}\bigGroup(-\PMC)$-set given by the rule
  \[(g_1\times_{\ZZ} g_2)* h := g_1\cdot h \cdot R(g_2)\]
  (where the operation $\cdot$ on the right-hand-side refers
  to multiplication in $\bigGroup(\PMC)$); for the proof,
  see~\cite[Lemma~\ref*{LOT2:lem:GradingIdentityDD}]{LOT2}.
  
  It follows that, in the language of 
  Sections~\ref{sec:coeff-algebra} and~\ref{sec:coeff-bimod},
  the diagonal subalgebra is the coefficient algebra of
  $\CFDDa(\Id)$.
\end{remark}


\newcommand\dischord{\flat}
\section{Bimodules for arc-slides}
\label{sec:Arc-Slides}
\setlength{\tabextrasep}{2pt}
Recall that Theorem~\ref{thm:DDforArc-Slides} provides a model for the
type \DD\  bimodule for an arc-slide. The aim of the present section is
to prove that theorem. As discussed in
Section~\ref{subsec:IdentificationArc-Slides}, the proof proceeds in
two steps:
\begin{enumerate}
\item Proving that the type \DD\ module associated to the standard
  Heegaard diagram for an
  arc-slide (Definition~\ref{def:StandardDiagram}) is a stable
  arc-slide bimodule (Definitions~\ref{def:Arc-SlideBimodule}
  and~\ref{def:StableModule}). This is proved in
  Section~\ref{sec:heeg-diagr-arc-slides}. The proof relies only on
  coarse properties of the Heegaard diagram (identification of
  generators, combinatorics of domains, and the existence of a few of
  the in principle many holomorphic curves which need to be computed).
\item Proving the uniqueness theorem for arc-slide bimodules
  (Proposition~\ref{prop:UniqueArc-Slide}). There
  are two combinatorially different cases of arc-slides: 
  {\em under-slides} and {\em over-slides} (see Definition~\ref{def:Over-Slide}).
  The argument is easier in the first case, where uniqueness
  is true in a slightly stronger form. In both cases, however, 
  the proof is broken
  down as follows:
  \begin{enumerate}
  \item The gradings on the modules restrict what terms can occur
    in the differential. The terms that are in correct gradings to
    occur are called \emph{near-chords} (for under-slides or
    over-slides---Definitions~\ref{def:UnderNearChord}
    and~\ref{def:OverNearChord}, respectively). The key properties of
    gradings required for this argument are collected in
    Section~\ref{sec:grad-near-diag}, and the restrictions are
    obtained in Lemmas~\ref{lem:under-slide-grading}
    and~\ref{lem:over-slide-grading}
    (Sections~\ref{subsec:Under-Slides}
    and~\ref{subsec:Over-Slides}). We use a comparison with gradings in the standard Heegaard diagram from
    Section~\ref{sec:heeg-diagr-arc-slides}.
  \item In the case of under-slides the equation $\bdy^2=0$ and the
    hypothesis that short near-chords appear in the differential implies that all near-chords appear; this is discussed in
    Section~\ref{subsec:Under-Slides}.
  \item In the case of over-slides there is some indeterminacy in
    which near-chords appear. Near-chords which may or may not appear
    are called \emph{indeterminate} (Definition~\ref{def:Indeterminate}); which indeterminate near-chords appear
    is determined by a so-called \emph{basic choice} (Definition~\ref{def:BasicChoice}). This is discussed in
    Section~\ref{subsec:Over-Slides}.
  \end{enumerate}
\end{enumerate}
Before turning to these proofs, we introduce a little more notation,
in Section~\ref{sec:arc-slide-notation}.

\subsection{More arc-slide notation and terminology}\label{sec:arc-slide-notation}
Let $m$ be an arc-slide taking a pointed matched circle $\PMC$ to
another pointed matched circle $\PMC'$.  Here, as in the introduction, $\PMC'$ is
obtained from $\PMC$ by sliding one of the feet
$b_1$ of an arc $B$ over another arc $C$; see
Figure~\ref{fig:ArcslideMatching}. The foot $b_1$ is connected to one
of the feet $c_1$ of $C$ by an arc $\sigma$ in $\PMC$; in $\PMC'$,
$b_1$ is replaced by the new foot $b_1'$ of $B'$, which is connected
by an arc $\sigma'$ to the foot $c_2$ of $C$ in $\PMC'$. 
\glsit{$\PMC'$}\glsit{$b_1$}\glsit{$B$}\glsit{$C$}\glsit{$\sigma$, $\sigma'$}\glsit{$b_1'$}\glsit{$B'$}%

\begin{convention}
  \label{label:conv}
  We focus on the case that $c_1$ is above $c_2$ in the $\PMC$
  matching, with respect to the orientation of $Z$; the case that
  $c_1$ is below $c_2$ is symmetric.
  With respect to a Heegaard diagram $\HD$
  for $m$, $\bdy\HD=-\PMC\amalg \PMC'$, so if we draw $\HD$ in the
  plane with handles attached, with $\PMC$ on the left and $\PMC'$ on
  the right, then $c_1$ is above $c_2$ in the plane as well; see
  Figure~\ref{fig:Genus2HandleSlide}.
\end{convention}

When the matching does not satisfy the assumption from Convention~\ref{label:conv}, we will switch the roles of $\PMC$ and $\PMC'$; 
see for example the remarks at the end of Definition~\ref{def:UnderNearChord}; see also Remark~\ref{rem:ConventionAgain}. 

\begin{definition}
  \label{def:Over-Slide}\index{over-slide}\index{arc-slide!over-slide}\index{under-slide}\index{arc-slide!under-slide}%
  An arc-slide is called an {\em over-slide} if $b_1$ is contained in
  the same component of $\PMC\setminus \{c_1, c_2\}$ as $z$.  (Note
  that this condition is symmetric in the roles of $\PMC$ and $\PMC'$.)
  Otherwise it is called an {\em under-slide}.
\end{definition}

The arc-slide on the left in Figure~\ref{fig:ArcslideMatching} is
an under-slide, while the one on the right is an over-slide.

We find it convenient to think of arc-slide bimodules as left-left
$\Alg(\PMC)\Hyph\Alg(-\PMC')$ bimodules (analogously to what was done
in Section~\ref{sec:DDforIdentity}) rather than left-right
$\Alg(\PMC)\Hyph\Alg(\PMC')$ bimodules (the point of view taken in
the introduction). Correspondingly, we can reformulate
Property~(\ref{AS:NearDiagonalSubalgebra}) for arc-slide bimodules
in terms of a subalgebra of $\Alg(\PMC)\otimes\Alg(-\PMC')$:

\begin{definition}\label{def:near-diag-alg}
  \index{near-diagonal subalgebra}\index{diagonal subalgebra|see{near-diagonal subalgebra}}\index{subalgebra!near-diagonal}%
  The {\em near-diagonal subalgebra} of $\Alg(\PMC)\otimes\Alg(-\PMC')$
  is the algebra generated by pairs of algebra elements 
  of the form $(j\cdot a\cdot i)\otimes (j' \cdot a'\cdot i')$
  where each of $(j,j')$ and $(i,i')$ is a pair of near-complementary
  idempotents, and also $\supp_R(a)=\supp_R(a')$.
\end{definition}

With this definition, Property~(\ref{AS:NearDiagonalSubalgebra})
of Definition~\ref{def:Arc-SlideBimodule} can be
reformulated as stating that the $\Alg(\PMC)\otimes
\Alg(-\PMC')$-module $N$ is induced from a module over the
near-diagonal subalgebra. 

Some definitions for $\Alg(\PMC)$ extend in obvious ways to the
near-diagonal subalgebra; for instance, a \emph{basic generator} of the
near-diagonal subalgebra is an element $a\otimes a'$ of the near-diagonal
subalgebra so that $a$ and $a'$ are basic generators of $\Alg(\PMC)$
and $\Alg(-\PMC)$.
\index{basic generator!of near-diagonal subalgebra}%

The restrictions placed by $\partial^2=0$ are strongest when we
restrict attention to the part of $\Alg(\PMC)$ with the property that
both the idempotents and the complementary idempotents have at least
two occupied positions; this is the portion with weight $-g+1<i<g-1$.
Since we are working with stable bimodules (in Propositions~\ref{prop:CalculateUnder-Slide}
and~\ref{prop:CalculateOver-Slide}), we can always stabilize (in the
sense of Definition~\ref{def:StableModule}) so that we are working
in this portion of the algebra. 
This will be the main way that we use
stability (see for example the proofs of Lemmas~\ref{lem:UFewerExist} and~\ref{lem:OExtraFewerExist};
we possibly stabilize more in the proof of Lemma~\ref{lem:IsBasicChoice}).

\subsection{Heegaard diagrams for arc-slides}
\label{sec:heeg-diagr-arc-slides}
\index{arc-slide!graph associated to}\index{graph associated to arc-slide}%
It is convenient to represent arc-slides by graphs embedded in an
annulus, as follows. Thinking of $\PMC$ and $\PMC'$ as two different
markings on the same circle $Z$, consider the annulus $[0,1]\times
Z$. This annulus has marked points on its boundary corresponding to
the positions in $\PMC$ (in $0\times Z$) and $\PMC'$ (in $1\times Z$), and
a special horizontal arc corresponding to the basepoint $z\in
\PMC$. Each position $p\in\PMC$, other than $b_1$, determines a
horizontal segment $[0,1]\times\{p\}$, connecting $p$ to its
corresponding point $p'\in\PMC'$. These horizontal segments, except for
the two horizontal segments corresponding to $c_1$ and $c_2$, are
edges for the graph. The horizontal segments for $c_1$ and $c_2$ are
both subdivided, by points $p_1$ and $p_2$ respectively, and we draw
two additional
edges, one connecting $p_1$ to $\{0\}\times b_1$, and another
connecting $p_2$ to $\{1\}\times b_1'$. The pictures in
Figure~\ref{fig:ArcslideMatching} can be thought of as illustrations
of these graphs (where the annuli have been cut along the horizontal
arcs corresponding to the basepoint $z$).

The graph for an arc-slide can
be turned into a Heegaard diagram $\HD(m)$, as
follows. 

\begin{definition}
\label{def:StandardDiagram}
\index{standard Heegaard diagram for arc-slide}\index{arc-slide!standard Heegaard diagram for}\glsit{$\HD(m)$}%
Let $m\colon\PMC\to\PMC'$ be an arc-slide. The {\em standard Heegaard
  diagram for the arc-slide $\HD(m)$} is the Heegaard diagram obtained
as follows.  Start from the graph associated to $m$, as
defined above. Attach a one-handle with feet at the two trivalent
points of the graph, in effect surgering out the two disks containing
trivalent points, and replacing them with an annulus equipped with
three arcs running along it.  These three arcs naturally extend to an arc
connecting $c_1$ and $c_2$, an arc connecting $c_1'$ and $c_2'$, and
an arc connecting $b_1$ and $b_1'$. Add a one-handle for each
matched pair $\{p,q\}$ in $\PMC$ other than $\{c_1,c_2\}$ and
$\{b_1,b_2\}$ with feet at the mid-points of the edges
corresponding to $p$ and $q$. The edges corresponding to $p$ and $q$ are surgered to get a pair
of arcs, one connecting $p$ to $q$ and the other connecting $p'$ to
$q'$. Performing one more such handle addition, one of whose feet is
on the edge corresponding to $b_2$ and the other at the edge
connecting $b_1$ to $b_1'$, we obtain the desired Heegaard surface.
The surgered arcs are the $\alpha$-arcs, and $\beta$-circles are
chosen to be meridians of the attached one-handles.
\end{definition}

\begin{figure}
    \begin{center}
      \input{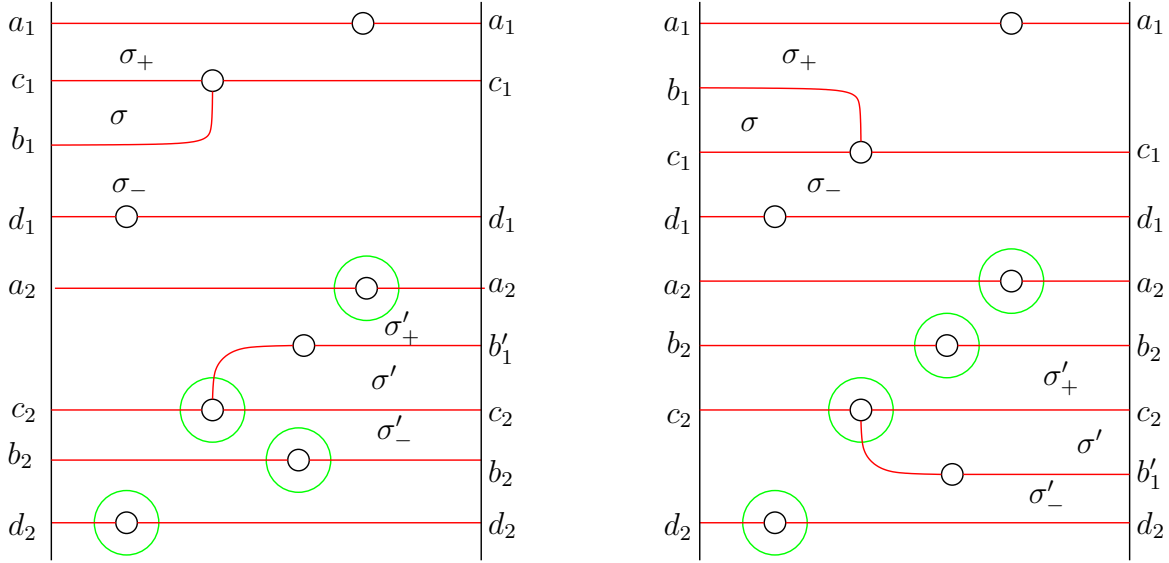}
    \end{center}
    \caption {{\bf Heegaard diagram for an arc-slide.}
      \label{fig:Genus2HandleSlide}
      Heegaard diagrams for the
      arc-slides from Figure~\ref{fig:ArcslideMatching}.
    The one on the left represents an under-slide, and the one on the
    right is an over-slide. Both satisfy Convention~\ref{label:conv}.
  In both cases, the basepoint $z$ in the pointed matched circle
  separates $d_2$ and $a_1$, and
  the picture is obtained by cutting the Heegaard diagram along the corresponding arc ${\mathbf z}$. }
\end{figure}

Figure~\ref{fig:Genus2HandleSlide} illustrates the result of this procedure.

Recall that if $m\colon \PMC\longrightarrow \PMC'$, then there is an
associated strongly-based mapping class $\PunctF(m)\co
\PunctF(\PMC)\to\PunctF(\PMC')$.
In~\cite[Definition~\ref*{LOT2:def:ConstructHeegaardDiagram}]{LOT2} we
  constructed a bordered Heegaard diagram from each strongly-based
  diffeomorphism. The diagram $\HD(m)$ is the Heegaard diagram
associated to the diffeomorphism $\PunctF(m)$.

\begin{figure}
    \begin{center}
      \input{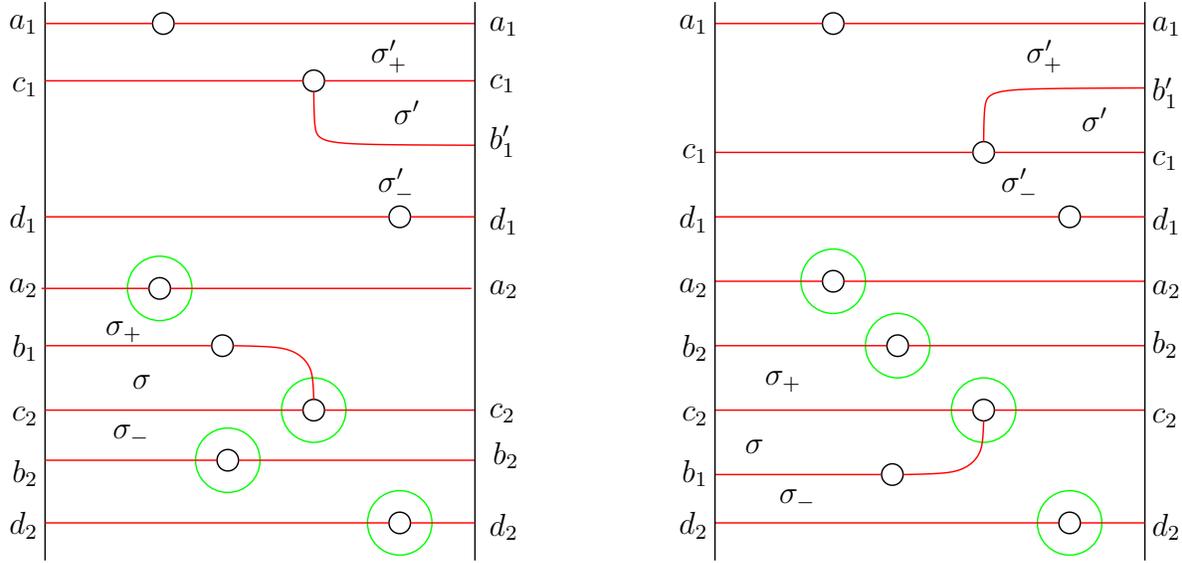}
    \end{center}
    \caption {{\bf Heegaard diagram for an arc-slide not satisfying Convention~\ref{label:conv}.}
      \label{fig:ReflectedDiagrams}
      This is a reflection of  Figure~\ref{fig:Genus2HandleSlide}.}
\end{figure}

\glsit{$\sigma$, $\sigma'$}\glsit{$\sigma_+,\sigma_,\sigma'_+,\sigma'_-$}%
For an arc-slide $m\co \PMC\to\PMC'$, 
the chord $\sigma$ in $\PMC$ lies on the boundary of a unique region
(component of $\Sigma\setminus(\alphas\cup\betas)$) in $\HD(m)$;
abusing notation, we will denote this region $\sigma$, as
well. Similarly, the chord $\sigma'$ in $\PMC'$ lies on the boundary
of a region $\sigma'$ in $\HD(m)$. 
Name the regions just
above and below $\sigma$ by $\sigma_+$ and $\sigma_-$, and the regions
just above and below $\sigma'$ by $\sigma'_+$ and $\sigma'_-$. 
All other regions look the same, and are strips across the diagram.
(For this notation, we are implicitly using the hypothesis from
Convention~\ref{label:conv}; for the other case, compare Figure~\ref{fig:ReflectedDiagrams}.) 

\begin{figure}
  \centering
  \includegraphics[scale=.5]{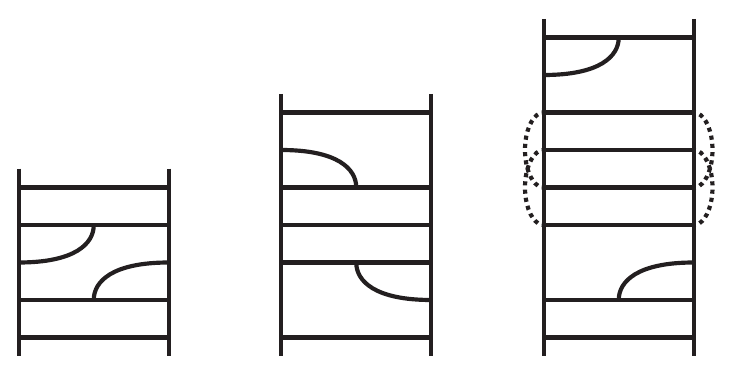}
  \caption{\textbf{Degenerate arc-slides and very special length $3$ chords.} Left: a degenerate under-slide (Definition~\ref{def:degen-slide}). Center: a degenerate over-slide (Definition~\ref{def:degen-slide}). Right: an under-slide containing a very special length $3$ chord (proof of Lemma~\ref{lem:UGenericExists}).}
  \label{fig:degen-slides}
\end{figure}

\begin{definition}\label{def:degen-slide}
  We call an arc-slide \emph{degenerate} if there is only one position
  between $c_1$ and $c_2$. (In the under-slide case satisfying
  Convention~\ref{label:conv}, this is equivalent to
  $\sigma_-=\sigma_+'$.) See Figure~\ref{fig:degen-slides}.
\end{definition}

For degenerate under-slides we need to allow one more kind of short
near-chord. We recall the definition of short near-chords
(Definition~\ref{def:RestrictedSupport}), extended to include this
case:
\begin{definition}\label{def:short-near-chord}
  A \emph{short near-chord} is a non-zero algebra element of the form
  $(i\cdot a \cdot j)\otimes (j'\cdot b'\cdot i')$ with the following
  four properties:
  \begin{enumerate}
    \item the pairs $(i\otimes i')$ and $(j\otimes j')$ are near-complementary idempotents;
    \item $\supp_R(a)=\supp_R(b)$;
    \item the support of at least one of $a$ or $b$ is non-zero; and
    \item the lengths of the (unrestricted) support of $a$ and the (unrestricted) support of $b$ are both no greater than $1$.
  \end{enumerate}

  In the degenerate case of under-slides with $\sigma_-=\sigma_+'$,
  we also call elements of the
  form $(i\cdot a(\sigma_-\cup\sigma) \cdot j)\otimes (j'\cdot
  a'_o(\sigma_-)\cdot i')$ and $(i\cdot a(\sigma_-) \cdot j)\otimes
  (j'\cdot a'_o(\sigma'\cup \sigma_-)\cdot i')$ short near-chords.
  Similarly, if $\sigma_+=\sigma'_-$, then we call
  $(i\cdot
  a(\sigma\cup\sigma_+) \cdot j)\otimes (j'\cdot a'_o(\sigma_+)\cdot
  i')$ and $(i\cdot a(\sigma_+) \cdot j)\otimes (j'\cdot
  a'_o(\sigma_+\cup \sigma')\cdot i')$ a short near-chord.
\end{definition}

\begin{proposition}
  \label{prop:CFDDisArcslideBimodule}
  If $m\colon \PMC\to\PMC'$ is an arc-slide and $\HD=\HD(m)$ is its
  associated standard Heegaard diagram (in the sense of
  Definition~\ref{def:StandardDiagram}), then the type \DD\
  bimodule $\CFDDa(\HD)$ is a stable arc-slide bimodule
  (in the sense of Definition~\ref{def:Arc-SlideBimodule} and 
  \ref{def:StableModule}).
\end{proposition}

\begin{proof}
  The generators correspond to near-complementary idempotents, as can
  be seen by adapting the corresponding fact for the identity type
  \DD\ bimodule (see the proof of Theorem~\ref{thm:DDforIdentity}).
  Generators of type $Y$ occur because now there is a $\beta$-circle
  which intersects three (rather than two) $\alpha$-arcs: two of those
  intersection points were of the type already encountered in the type
  \DD\ identity bimodule; the third, however, represents an
  intersection point of the $\alpha$-arc for the $B$-matched pair with
  the $\beta$-circle for the $C$ and $C'$ matched pair. See Figure~\ref{fig:GenArcslides}.

  \begin{figure}
    \begin{center}
      \input{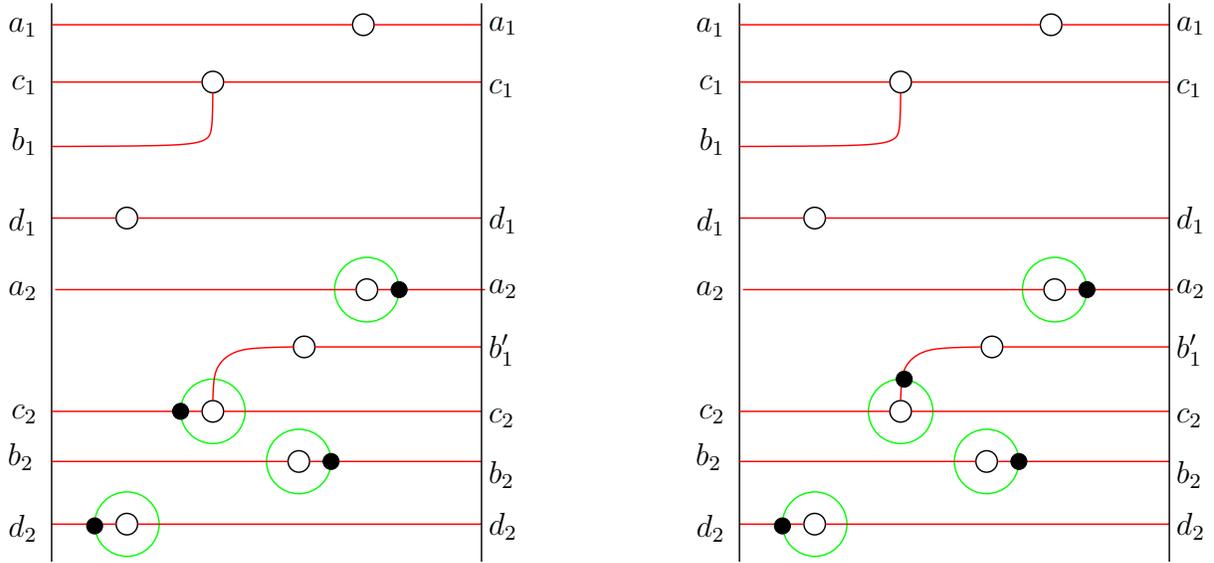}
    \end{center}
    \caption {{\bf Generators for arc-slide bimodule.}
      \label{fig:GenArcslides}
      At the left is a generator of type $X$ (indeed, $X_C$, for the notation used in the proof of Proposition~\ref{prop:near-diag-grading});
      at the right is a generator of type $Y$.}
  \end{figure}

  Since every region touches the boundary, the multiplicities at the
  boundary together with a choice of initial generator~$\x$ determine a domain $\Domain \in
  \pi_2(\x, \y)$.  Furthermore, local multiplicities at the boundary are
  possible if and only if the restricted supports agree on the two sides.  (Compare
  the proof of Theorem~\ref{thm:DDforIdentity}.)  Thus the coefficients
  of the differential on $\CFDDa(\HD(m))$
  lie in the near-diagonal subalgebra.
  Moreover, for each basic generator $I\cdot a\cdot
  J$ of the near-diagonal subalgebra there is a unique domain
  $B\in\pi_2(\x(I),\x(J))$ which could
  contribute $a\x(J)$ to $\bdy \x(I)$.
  \glsit{$\x(I)$}%

  In particular, as in the proof of
  Theorem~\ref{thm:DDforIdentity}, the grading set of $\CFDDa(\HD)$ is
  $\lambda$-free, as periodic domains are determined by their local
  multiplicities near the boundary.

  In the non-degenerate case, the differential on $\CFDDa(\HD)$
  contains all short near-chords (Condition~\ref{AS:ShortChords}), as
  they are represented by polygons (compare
  Lemma~\ref{lem:polygons-represent}). In the degenerate case, the
  domains corresponding to the additional near-chords $(i\cdot
  a(\sigma_-\cup\sigma) \cdot j)\otimes (j'\cdot a'_o(\sigma_-)\cdot
  i')$ and $(i\cdot a(\sigma_-) \cdot j)\otimes (j'\cdot
  a'_o(\sigma'\cup \sigma_-)\cdot i')$ are annuli. The combinatorics
  of these annuli (one $270^\circ$ corner, at which the $\alpha$- and
  $\beta$-curves go out to the other boundary component) are such that
  they always have algebraically $1$ holomorphic representative; see,
  for instance,~\cite[Lemma 9.10]{Rasmussen03:Knots} or the $\rho_3$
  case of the proof of~\cite[Proposition 10.6]{LOT2}.%

  Thus, we have verified that $\CFDDa(\HD)$ is an arc-slide bimodule
  in the sense of Definition~\ref{def:Arc-SlideBimodule}.
  To see it is
  stable in the sense of Definition~\ref{def:StableModule}, embed the
  Heegaard diagram $\HD$ for $m\co\PMC\rightarrow \PMC'$ in the
  Heegaard diagram $\HD_\#$ for the $\PMC_0$-stabilized arc-slide
  $m_\#$. Setting to zero algebra elements whose support intersects
  $\PMC_0$, and restricting the generator in the new region, we obtain
  a module representing $Q_*(\CFDDa(m_\#))$. If we choose
  almost-complex structures for $\HD$ and $\HD_\#$ compatibly then the
  holomorphic curve
  counts involved in $Q_*(\CFDDa(m_\#))$ coincide with the curve
  counts for $\CFDDa(\HD)$. 
\end{proof}

\begin{remark}
  The grading sets for $\CFDDa(\HD)$ are determined explicitly in
  Section~\ref{sec:mcg-grading}.
\end{remark}

\subsection{Gradings on the near-diagonal subalgebra}
\label{sec:grad-near-diag}

The near-diagonal subalgebra has a $\ZZ$-grading,
which can be used to exclude the appearance of many of its elements from 
the differential in an arc-slide bimodule.  This comes from thinking
of it as a coefficient algebra:  we will see
(Lemma~\ref{lem:coeff-near-diag} below) that for any arc-slide
bimodule~$N$ (including a standard one, $\CFDDa(\HD(m))$), $\Coeff(N)$
contains the near-diagonal subalgebra, which is therefore $\ZZ$-graded.
Furthermore, we will see (Proposition~\ref{prop:grading-unique}) that this
$\ZZ$-grading is independent of the choice of $N$.
Thus we can compute the $\ZZ$-grading on the near-diagonal subalgebra by
looking at the grading set~$S_\std$ for $\CFDDa(\HD(m))$; we do this
explicitly in Proposition~\ref{prop:near-diag-grading}.

By Lemma~\ref{lem:we-dont-care}, 
we are free here to choose the grading group to be $G(\PMC)$ or $G'(\PMC)$.
For definiteness, we work with the latter, i.e., $\Alg(\PMC)$ is thought
of as graded by $G'(\PMC)$, and the grading set of our module $N$ is a $G'(\PMC)$-set.
We will be working with both the grading set $S'_\std$ for
$\CFDDa(\HD(m))$, and the grading set~$S'_N$ for an arbitrary arc-slide
bimodule~$N$. Denote the
gradings with values in $S'_\std$ and $S'_N$ by $\grb_\std$ and $\grb_N$,
respectively. Also note that the generators of~$N$ and the generators
for $\CFDDa(\HD(m))$ are naturally identified,
since both are identified with idempotents for the near-diagonal subalgebra.
\glsit{$S'_\std$, $S'_N$, $\gr'_\std$, $gr'_N$}%

\begin{definition}
  \index{chain of algebra elements}%
  For a pair of basic idempotents $I$, $J$ in a based algebra~$A$ (as
  in Definition~\ref{def:based-alg}), a
  \emph{chain of algebra elements}~$\Chain$ connecting $I$ and $J$ is
  a sequence $(a_1, \epsilon_1)$, $(a_2, \epsilon_2)$, \dots, $(a_n,
  \epsilon_n)$ of pairs of algebra elements~$a_i$ and
  signs~$\epsilon_i \in \{-1, +1\}$ so that there is a sequence of
  basic idempotents
  \[
  I = I_0, I_1, \dots, I_{n-1}, I_n = J
  \]
  with
  \[
  a_i =
  \begin{cases}
    I_{i-1}\cdot a_i\cdot I_{i} & \text{if }\epsilon_i = +1 \\
    I_{i}\cdot a_i\cdot I_{i-1} & \text{if }\epsilon_i = -1
  \end{cases}
  \]
  for each $i$.
  (Note that the $I_i$ can be recovered from the $a_i$.)

  The \emph{inverse} of $\Chain$ is
  \index{inverse!of chain of algebra elements}%
  \[
  \Chain^{-1}=((a_n,-\epsilon_n),\dots,(a_1,-\epsilon_1)),
  \]
  which is a chain connecting $J$ to $I$.

  If $A$ is graded, the \emph{grading} of~$\Chain$ is
  \index{chain of algebra elements!support of}%
  \index{chain of algebra elements!grading of}%
  \index{support!of chain of algebra elements}\index{grading!of chain of algebra elements}%
  \begin{align*}
    \grb(\Chain) &= \prod_{i=1}^n (\lambda \grb(a_i))^{\epsilon_i}.\\
    \intertext{Similarly, if $A$ is $\Alg(\PMC)$ or the near-diagonal
      subalgebra, the \emph{support} of~$\Chain$ is}
    \supp(\Chain)&=\sum_{i=1}^n \epsilon_i\supp(a_i).
  \end{align*}
  \glsit{$\grb(\Chain)$}\glsit{$\supp(\Chain)$}%
\end{definition}

\index{extremal weight}\index{weight!extremal}%
We say that an idempotent $I=(i\otimes i')$ for the near-diagonal
subalgebra has {\em extremal weight} if the weight of $i$ is $\pm k$
(and hence the weight of $i'$ is $\mp k$). Note that there are exactly
two idempotents of extremal weight. For definiteness, we will say that
the weight of $I=(i\otimes i')$ is the weight of $i$.

\begin{lemma}\label{lem:alg-chain}
  For any two idempotents $I$ and $J$ for the near-diagonal
  subalgebra, with the same non-extremal weight, there is a chain $\Chain$ of
  short near-chords
  connecting $I$ and $J$.  Furthermore, if
  $\x$ and $\y$ are the corresponding generators for $\HD(m)$ and
  $\Domain \in \pi_2(\x, \y)$, we can choose the chain $\Chain$ so that
  \[
  \supp(\Chain) = -r_*(\bdy^\bdy(\Domain)).
  \]
\end{lemma}
(Recall that $\bdy\HD=-\PMC\amalg\PMC'$, while the near-diagonal
subalgebra lives inside $\Alg(\PMC)\otimes\Alg(-\PMC')$. This is the
reason for the presence of the map $r_*$ from
Section~\ref{sec:coeff-bimod}.)
\begin{proof}
  For definiteness, we will discuss the case of under-slides; the case
  of over-slides is similar. Let $d_1$ be the point in $\CircPts$ just
  below $b_1$ and $d_2$ the point matched with $d_1$. (For a
  degenerate arc-slide, $d_1=c_2$.)

  For the first part, connect any pair of idempotents by a chain of
  short near-chords by swapping adjacent pairs.
  More precisely, 
  we may assume that $I$ and $J$ are both idempotents
  of type $X$, by choosing $a_1$ and/or $a_n$ to be the short near-chord $\sigma\otimes 1$
  if necessary (and $\epsilon_1=1$ and/or $\epsilon_n=-1$).
  Define a multigraph\footnote{A \emph{multigraph} is a graph which may have multiple edges connecting the same pair of vertices.} $\Gamma$ as follows:
  \begin{itemize}
  \item $\Gamma$ has one vertex $A_i$ for each matched pair in $\PMC$
    (or equivalently $\PMC')$.
  \item $\Gamma$ has an edge connecting $A_i$ and $A_j$ each time a
    foot of $A_i$ is adjacent to a foot of $A_j$ in
    $\CircPts\setminus\{b_1\}$. In particular, the point $d_1$ is
    adjacent to $c_1$, and $b_1$ is not viewed as adjacent to
    anything.

    Let $(x_1,y_1)$ denote the edge connecting the matched pair
    $\{x_1,x_2\}$ to the matched pair $\{y_1,y_2\}$ coming from the
    fact that $x_1$ is adjacent to $y_1$.
  \item Add one more edge from the matched pair $\{d_1,d_2\}$ to the
    matched pair $\{b_1,b_2\}$. Call this edge $(d_1,b_1)$.
  \end{itemize}
  See Figure~\ref{fig:adj-graph}.

  \begin{figure}
    \centering
    \includegraphics{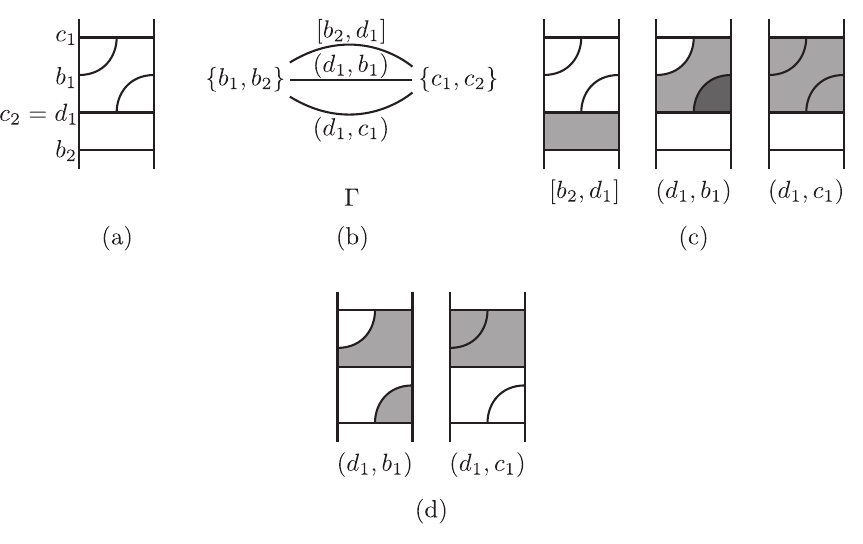}
    \caption{\textbf{The adjacency graph for an arc-slide.}  (a) A
      genus-one arc-slide. (b) The corresponding graph. (c) The
      domains corresponding to the edges in the graph. The
      darkly-shaded region is covered with multiplicity $2$. (d) The
      domains corresponding to $(d_1,b_1)$ and $(d_1,c_1)$ in the
      non-degenerate case.}
    \label{fig:adj-graph}
  \end{figure}
  
  The type $X$ idempotents of weight~$n$ correspond to subsets $S$ of
  the vertices of~$\Gamma$ with $|S|=n+k$. Most of the edges in
  $\Gamma$ correspond to short near-chords in an obvious way: except
  for the edges $(d_1,c_1)$ and $(d_1,b_1)$ the edges in $\Gamma$ come
  from length-1 intervals in $(Z,\CircPts)$. We will see that the edge
  $(d_1,b_1)$ corresponds to the pair of short near-chords
  $(\sigma_-\otimes\sigma_-,1\otimes \sigma')$ and the edge
  $(d_1,c_1)$ corresponds to the pair of short chords
  $(\sigma_-\otimes\sigma_-,\sigma\otimes 1)$ (in the non-degenerate
  case).

  Consider the symmetric group~$G$ on the vertices of~$\Gamma$. Since
  type $X$ idempotents correspond to subsets of the vertices of
  $\Gamma$, the group $G$ acts transitively on the set of type $X$ idempotents. 
  Since $\Gamma$ is connected, $G$ is generated by those
  transpositions that interchange vertices connected by an edge.  Let
  $\gamma$ be an edge of $\Gamma$, connecting some $A_i$ to some
  $A_j$, and $\tau_\gamma\in G$ the corresponding transposition. Let
  $I$ be a indecomposable idempotent.
  \begin{itemize}
  \item If $\gamma\not\in\{(d_1,b_1),(d_1,c_1)\}$ and exactly one of $A_i$
    or $A_j$ is occupied in $I$ then the action of $\tau_\gamma$ on
    $I$ can be achieved by multiplying by the corresponding short
    near-chord $a(\xi) \otimes a'_o(\xi)$ (with $\epsilon_i = \pm 1$).
  \item If $\gamma=(d_1,b_1)$ (so
    $\{A_i,A_j\}=\{\{d_1,d_2\},\{b_1,b_2\}\}$) and exactly one of
    $A_i$ or $A_j$ is occupied in $I$ then the action of $\tau_\gamma$
    on $I$ can be achieved by multiplying by either
    $(\sigma_-\otimes\sigma_-,\pm1),(1\otimes \sigma',\pm1)$ or
    $(1\otimes \sigma',\pm1),(\sigma_-\otimes\sigma_-,\pm1)$
    (depending on the occupancy of $\{c_1,c_2\}$ in $I$). (For a
    degenerate arc-slide, $\sigma_-\otimes\sigma_-$ is replaced by
    $(\sigma_-\cup\sigma)\otimes \sigma_-$.)
  \item If $\gamma=(d_1,c_1)$ (so
    $\{A_i,A_j\}=\{\{d_1,d_2\},\{c_1,c_2\}\}$) and exactly one of
    $A_i$ or $A_j$ is occupied in $I$ then the action of $\tau_\gamma$
    on $I$ can be achieved by multiplying by either
    $(\sigma_-\otimes\sigma_-,\pm1),(\sigma\otimes 1,\pm1)$ or
    $(\sigma\otimes 1,\pm1),(\sigma_-\otimes\sigma_-,\pm1)$, depending
    on the occupancy of $\{b_1,b_2\}$ in~$I$. (For a degenerate arc-slide, $\sigma_-\otimes\sigma_-$ is replaced by $\sigma_-\otimes (\sigma'\cup\sigma_-)$.)
  \item If neither $A_1$ nor $A_2$ is occupied in $I$ then the action
    of $\gamma$ on $I$ is trivial.
  \item If both $A_1$ and $A_2$ are occupied in $I$ then the action of
    $\gamma$ on $I$ is trivial.
  \end{itemize}
  Any two $(n+k)$-element sets of vertices of $\Gamma$ are related by an
  element of the symmetric group, so the first part of the claim
  follows.

  For the second part, it is enough to see that we can obtain the
  boundary of any periodic domain as the support of a chain. Then, since we
  can connect $I$ and $J$ by some chain $\Chain_0$, to get a chain
  representing any other domain connecting $I$ and $J$ we merely concatenate a chain
  representing the appropriate periodic domain. Moreover:
  \begin{itemize}
  \item We may work with any convenient starting and ending
    idempotent: given a periodic domain $\Domain$, a chain $\Chain_1$ connecting
    $I_0$ to $I_0$ with $\supp(\Chain) = -r_*(\bdy^\bdy \Domain)$, and another idempotent
    $J$, choose a chain $\OthChain$ connecting $I_0$ to $J$. Then $\OthChain^{-1}\Chain_1 \OthChain$ is
    a chain connecting $J$ to $J$ with $-r_*(\bdy^\bdy \Domain)=\supp(\OthChain^{-1}\Chain_1 \OthChain)=\supp(\Chain_1)$.
  \item Given a basis $\{\Domain_i\}$ for the space of periodic domains,
    it suffices to find a chain representing each basis element.
  \end{itemize}

  We will now show how cycles in the graph $\Gamma$ give periodic
  domains.  Pick a basis of $H_1(\Gamma)$ consisting of simple cycles
  (edge loops with no repeated edges).  For each such basis element,
  we can find a chain of short near-chords as follows.  Start in a
  type $X$ idempotent where there is at least one occupied and at
  least one unoccupied vertex in the cycle.  (This is possible since
  we are not in the extremal weight and the cycles in $\Gamma$ have at
  least $2$ vertices.)  Then swap any consecutive (occupied,
  unoccupied) vertex pair where the unoccupied vertex is clockwise
  from the occupied one.  Repeat until each occupied vertex has moved
  clockwise to the next occupied slot, or equivalently until we have
  swapped once on each edge. This sequence of swappings gives a
  periodic domain. The $2k$ independent cycles give a basis for the
  space of periodic domains.
\end{proof}

\begin{lemma}\label{lem:chain-grading}
  For any arc-slide bimodule $N$ for~$m$,
  near-complementary idempotents $I$ and $J$, and chain~$\Chain$ of
  short near-chords connecting $I$ and $J$,
  \[
  \grb_N(I) = \grb(\Chain)\grb_N(J).
  \]
\end{lemma}

\begin{proof}
  If $a_i$ is a short near-chord
  with $I_i \cdot a_i \cdot I_{i+1} = a_i$ then $a_i I_{i+1}$ appears in
  $\bdy I_i$ in~$N$ (by hypothesis of an arc-slide bimodule)
  and so
  \[
  \grb_N(I_i) = \lambda \grb(a_i) \grb_N(I_{i+1}).
  \]
  Similarly, if $I_{i+1}\cdot a_i \cdot I_{i} = a_i$ then $a_i \x_i$
  appears in $\bdy I_{i+1}$ and
  \[
  \grb_N(I_i) = (\lambda \grb(a_i))^{-1} \grb_N(I_{i+1}).
  \]
  The result follows by induction.
\end{proof}

Note that Lemma~\ref{lem:chain-grading} applies to any
arc-slide bimodule.  In particular, $\CFDDa(\HD(m))$ is an arc-slide
bimodule, and so it applies there.

\begin{lemma}\label{lem:coeff-near-diag}
  For any arc-slide bimodule~$N$, the near-diagonal subalgebra is
  contained in $\Coeff(N)$.
\end{lemma}

In principle, $\Coeff(N)$ could be larger than the near-diagonal
subalgebra, but Definition~\ref{def:Arc-SlideBimodule} guarantees that
any elements not in the near-diagonal subalgebra do not appear in the
differential.

\begin{proof}
  Let $a$ be an element of the near-diagonal subalgebra, with initial
  idempotent $I$ and final idempotent~$J$.  Then, as each region in
  $\HD(m)$ touches the boundary, there is a unique domain $\Domain$
  connecting the generators $\x(I)$, $\x(J)$ of $\HD(m)$ corresponding to
  $I$, $J$ so that $\supp(a) = -r_*(\bdy^\bdy \Domain)$.  Let $\Chain$ be a chain of
  algebra elements connecting $\x(I)$ to $\x(J)$ with support
  $\supp(\Chain)=-r_*(\bdy^\bdy \Domain)$, whose existence is guaranteed by
  Lemma~\ref{lem:alg-chain}.  Then
  \begin{align*}
    \supp(a) &= -r_*(\bdy^\bdy \Domain) = \supp(\Chain)\\
    \grb(a) &= \lambda^m \grb(\Chain)
  \end{align*}
  for some integer $m$.
  Thus by Lemma~\ref{lem:chain-grading},
  \[
  \grb_N(\x(I)) = \grb(\Chain) \grb_N(\x(J)) = \lambda^{-m} \grb(a) \grb_N(\x(J)),
  \]
  which says that $(\x(I),a,\x(J)) \in \Coeff(N)$, as desired.
\end{proof}

\begin{proposition}\label{prop:grading-unique}
  Let $N$ be an arc-slide bimodule for~$m$ in the sense of
  Definition~\ref{def:Arc-SlideBimodule}.
  Then there is a $G'$-set map $f \co S'_{\std}(m) \to S'_N$ so that for
  each generator~$\x(I)$ for $N$,
  \[
  f(\grb_{\std}(\x_{\std}(I))) = \grb_N(\x(I)),
  \]
  where $\x_{\std}(I)$ is the corresponding generator for
  $\CFDDa(\HD(m))$ with idempotent $I$.
  Furthermore, the $\ZZ$-gradings on the near-diagonal subalgebra from
  the two bimodules agree.
\end{proposition}
\begin{proof}
  Pick a base generator $\x_0$.  We first compare the
  stabilizer $\Stab_{\std}(\x_0)$ of $\grb_\std(\x_0)$ in $S'_\std(m)$
  and the stabilizer $\Stab_N(\x_0)$ of $\grb_N(\x_0)$ in~$S'_N$.
  For any $\Domain \in
  \pi_2(\x_0, \x_0)$, there is a chain~$\Chain$ of short near-chords with
  $-r_*(\bdy^\bdy \Domain) = \supp(\Chain)$ by
  Lemma~\ref{lem:alg-chain}, and so by Lemma~\ref{lem:chain-grading},
  \[
  \grb(\x_0) = \grb(\Chain)\grb(\x_0),
  \]
  where $\grb(\x_0)$ denotes either $\grb_\std(\x_0)$ or $\grb_N(\x_0)$.
  Thus $\grb(\Chain)$ is in both
  $\Stab_{\std}(\x_0)$ and $\Stab_N(\x_0)$.  From $\HD(m)$ we see
  that $\grb(\Chain)\grb(\x_0)$ must be $R(g'(\Domain))$: the homological components of
  $R(g'(\Domain))$ and of $\grb(\Chain)$ agree, and since the
  grading set on $\HD(m)$ is $\lambda$-free, there can be at most one
  such element in $\Stab_\std(\x_0)$.

  By hypothesis, $\Stab_{\std}(\x_0)$ is generated by elements of the
  form $R(g'(\Domain))$, and so we have $\Stab_{\std}(\x_0) \subset
  \Stab_N(\x_0)$.
  We can therefore define the map $f\co S'_{\std}(m) \to S'_N$ in a canonical way.

  Now, for any other generator $\x$, connect $\x$ to $\x_0$ by a chain~$\Chain$
  of short near-chords.
  Then
  \begin{equation}
    \label{eq:f-compat-gr}
  \begin{aligned}
  f(\grb_\std(\x))
    &= f(\grb(\Chain) \grb_\std(\x_0))\\
    &= \grb(\Chain) f(\grb_\std(\x_0))\\
    &= \grb(\Chain) \grb_N(\x_0)\\
    &= \grb_N(\x),
  \end{aligned}
  \end{equation}
  as desired.  We used Lemma~\ref{lem:chain-grading} twice
  (once in each grading set), as well as the fact that $f$ is a
  $G'(\HD)$-set map.

  Now consider any element of the near-diagonal subalgebra, which we can
  think of as a triple $(\x, a, \y)$ in the coefficient algebra of
  either $\CFDDa(\HD(m))$ or~$N$.  We have
  \begin{align*}
    \lambda^{\gr_\std(\x,a,\y)} \grb_\std(\x) &= \grb(a) \grb_\std(\y)\\
    \lambda^{\gr_\std(\x,a,\y)} f(\grb_\std(\x)) &= \grb(a) f(\grb_\std(\y))\\
    \lambda^{\gr_\std(\x,a,\y)} \grb_N(\x) &= \grb(a) \grb_N(\y),
  \end{align*}
  by (in order) the definition of $\grb_\std$ on
  $\Coeff(\CFDDa(\HD(m)))$
  (Lemma~\ref{lem:grading-coeff-implies}), the fact that $f$ is a $G'(\HD)$-set map,
  and Equation~\eqref{eq:f-compat-gr}.
  Since by assumption $S_N'$ is $\lambda$-free, the last equation
  implies that $\gr_N(\x,a,\y) = \gr_\std(\x,a,\y)$.
\end{proof}

Recall from Definition~\ref{def:SubComplementary} that the
near-complementary idempotents are either the analogues of
complementary idempotents, called type~$X$, or are sub-complementary
idempotents, called type~$Y$. Idempotents of type $X$ are further
divided into idempotents where there are horizontal strands at $C$
in~$\PMC$, which we call \emph{type~$\lsub{C}X$}, and idempotents where
there are horizontal strands at $C$ in~$\PMC'$, which we call \emph{type~$X_C$}.
\index{idempotent!type $\lsub{C}X$}\index{idempotent!type $X_C$}\index{type $\lsub{C}X$ idempotent}\index{type $X_C$ idempotent}%
\glsit{$\lsub{C}X$,$X_C$}%

Recall that $\sigma$, $\sigma'$, $\sigma_+$, $\sigma_-$, $\sigma'_+$
and $\sigma'_-$ denote regions in the standard Heegaard diagram
$\HD(m)$ for $m$; see Section~\ref{sec:heeg-diagr-arc-slides}.  For a
region $\Region$ in $\HD(m)$, $n_\Region$ is the
\glsit{$n_\Region$}%
multiplicity with which $\Region$ appears in a domain compatible
with~$a$.  Since every region of $\HD(m)$ touches the
boundary, this is actually a function of $\supp(a)$. These
multiplicities of $a$ are well-defined only if $a$ is in the
near-diagonal subalgebra, not for general $a$ in
$\Alg(\PMC)\otimes\Alg(-\PMC')$: if $a$ is not in the near-diagonal
subalgebra then there is no domain~$\Domain$ compatible with~$a$.

\begin{proposition}
  \label{prop:near-diag-grading}
  Consider an arc-slide satisfying Convention~\ref{label:conv}.
  In the $\ZZ$-grading on the near-diagonal subalgebra for an
  arc-slide of $b_1$ over $c_1$,
  the grading
  of a basic generator $a$ with initial idempotent $I$ and final
  idempotent $J$ is
  \begin{equation}\label{eq:gr-basic-near-diag}
  \iota(a) + c(I, \supp(a)) + c(J, \supp(a))
  \end{equation}
  where $c(I, \supp(a))$ is a correction term given by
  \index{grading!correction term}\index{correction term}%
  \glsit{$c(I, \supp(a))$}%
  \begin{equation}
    \label{eq:UnderslideGradingCorrection}
  c(I, \supp(a)) =
  \begin{cases}
    \OneQuart(n_{\sigma'_+} - n_{\sigma'}) &
      \text{$I$ of type $\lsub{C}X$}\\
    \OneQuart(-n_{\sigma} + n_{\sigma_-}) &
      \text{$I$ of type $X_C$}\\
    \OneQuart(n_{\sigma_+} - n_{\sigma} - n_{\sigma'} + n_{\sigma'_-}) &
      \text{$I$ of type $Y$}
  \end{cases}
  \end{equation}
  for an under-slide and by
  \begin{equation}
    \label{eq:OverslideGradingCorrection}
  c(I, \supp(a)) =
  \begin{cases}
    \OneQuart(- n_{\sigma'} + n_{\sigma'_-}) &
      \text{$I$ of type $\lsub{C}X$}\\
    \OneQuart(n_{\sigma_+}-n_{\sigma}) &
      \text{$I$ of type $X_C$}\\
    \OneQuart(- n_{\sigma} +n_{\sigma_-} + n_{\sigma'_+} - n_{\sigma'} ) &
      \text{$I$ of type $Y$}
  \end{cases}
  \end{equation}
  for an over-slide.
\end{proposition}

(In the degenerate case for under-slides in which $\sigma_-=\sigma_+'$, the formulas in Proposition~\ref{prop:near-diag-grading} hold without change.)

\begin{proof}
  By Proposition~\ref{prop:grading-unique}, we can compute the grading
  either in the coefficient algebra of the bimodule~$N$ or from the
  Heegaard diagram $\HD(m)$.  Inside $\HD(m)$, we can apply
  Lemma~\ref{lem:grading-coeff-DD}, as follows.  The $\alpha$- and
  $\beta$-curves divide the Heegaard surface $\Sigma$ into regions.
  Most of these regions~$\Region$, other
  than the six regions neighboring $\sigma$ or $\sigma'$, are strips
  across the diagram; these are all octagons, so have $e(\Region) = -1$, and
  have $n_\x(\Region) = n_\y(\Region) = \OneHalf$, for a total
  contribution of~$0$.  To analyze the six special
  regions more conveniently, divide the correction into two pieces, associated to
  the initial and final generators: for $\Domain\in\pi_2(\x,\y)$
  with $\bdy^\bdy\Domain=\supp(a)$,
  \[
  \gr(\x,a,\y) = \iota(a) + (-\OneHalf e(\Domain) - n_\x(\Domain))
     + (-\OneHalf e(\Domain) - n_\y(\Domain)).
  \]
  Then, for instance, if $b_1$ is below $c_1$, $\sigma_+$ is an
  octagon, and we find
  \[
  -\OneHalf e(\sigma_+) - n_I(\sigma_+) =
  \OneHalf +
  \begin{cases}
    -\OneHalf & \text{$I$ of type $X$}\\
    -\OneQuart & \text{$I$ of type $Y$.}
  \end{cases}
  \]
  This gives the contribution of $n_{\sigma_+}$ to $c(I, \supp(a))$ in
  this case.
  The other contributions are similar.
\end{proof}

\subsection{Under-slides}
\label{subsec:Under-Slides}

We explicitly describe the differential for the arc-slide bimodule
for under-slides in Proposition~\ref{prop:CalculateUnder-Slide}. We
first define near-chords in this case.

\begin{figure}
    \begin{center}
      \input{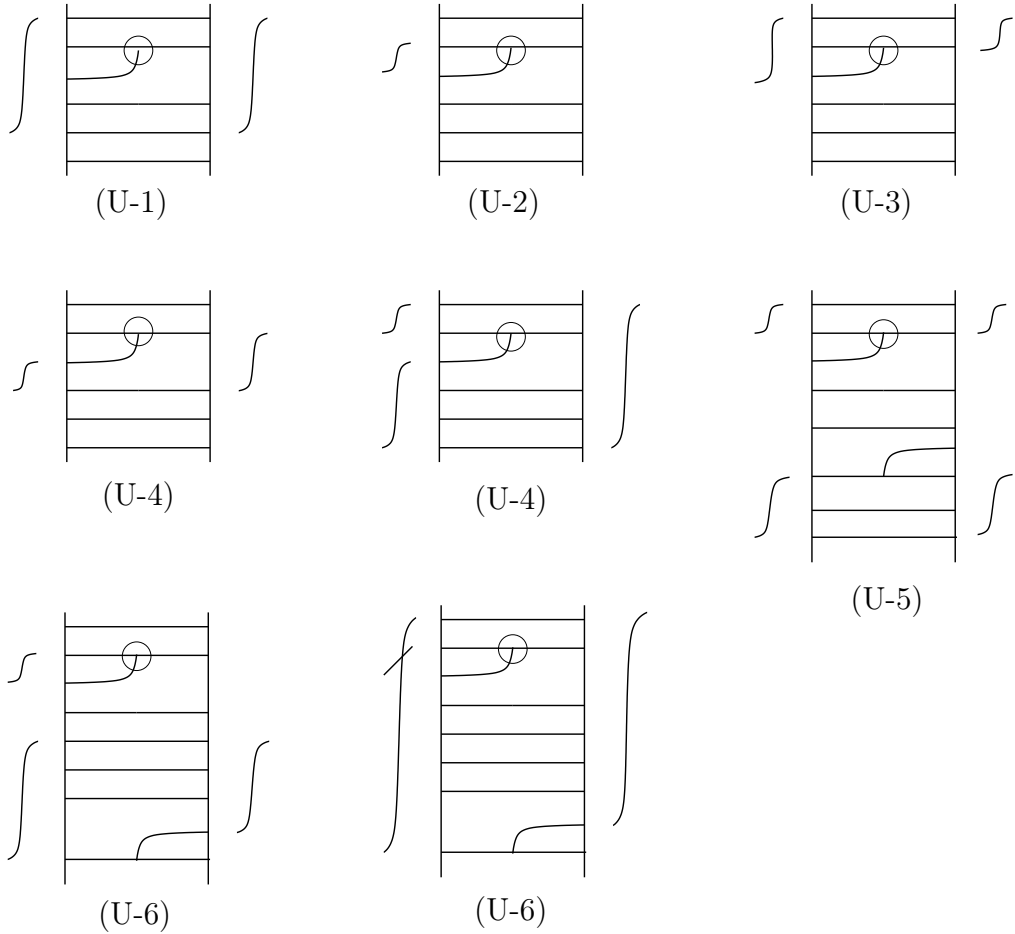}
    \end{center}
    \caption {{\bf Near-chords for under-slides.}
      \label{fig:UnderNearChords} We have illustrated
    here examples of all the types of near-chords for under-slides,
    Definition~\ref{def:UnderNearChord}. 
    (Note that there are two near-chords of Type~(\ref{typeU:FewerSigma})
    and of Type~(\ref{typeU:XX})
    since there are two distinct subtypes of these.)}
\end{figure}

\begin{definition}
  \label{def:UnderNearChord}
  A (non-zero) basic algebra element $x$ in the near-diagonal subalgebra for an
  under-slide satisfying Convention~\ref{label:conv} is
  called a \emph{near-chord} if it satisfies any of the following
  conditions:
  \begin{enumerate}[label=(U-\arabic*),ref=U-\arabic*]
    \item 
      \label{typeU:Generic}
      It has the form
      $x=I\cdot (a(\xi)\otimes a'_o(\xi))\cdot J$,
      where $I$ and $J$ are near-complementary idempotents and $\xi$ is some chord in $\PMC$ neither of whose endpoints is $b_1$
      (so that it can be interpreted, as it is in the above expression,
      as a chord also in $\PMC'$); furthermore, $\xi$ is required to be different from the chord
      $[c_2,c_1]$.
    \item
      \label{typeU:Sigma}
      It has the form $x=I\cdot (a(\sigma)\otimes
      1)\cdot J$ or $I\cdot (1\otimes a_o'(\sigma'))\cdot J$, where $I$
      and $J$ are near-complementary idempotents.
   \item
      \label{typeU:ExtraSigma}
        There is a chord $\xi$ for $\PMC$ so that the interior of
        $\xi$ is disjoint from $\sigma$
        and the support of $\xi\cup\sigma$ is connected,
        and $x$ has the form
        $x=I\cdot (a(\xi\cup \sigma)\otimes a'_o(\xi))\cdot J$;
        or there is a chord $\xi$ for $\PMC'$
        so that the interior of $\xi$ is disjoint from $\sigma'$,
        and the support of $\xi\cup\sigma'$ is connected, and
        $x$ has the form
        $x=I\cdot (a(\xi)\otimes a'_o(\xi\cup\sigma'))\cdot J$.
    \item 
      \label{typeU:FewerSigma}
      It has the form $x=I\cdot (a(\xi\setminus
      \sigma)\otimes a'_o(\xi))\cdot J$ where $\sigma\subset \xi$; or
      $x=I\cdot (a(\xi)\otimes a'_o(\xi\setminus \sigma'))\cdot J$ where $\sigma'\subset
      \xi$. (The element~$x$ can have two or three moving strands.)
    \item 
      \label{typeU:YY}
      It has the form
      $x=I\cdot (a(\xi\cup\eta)\otimes a_o'(\xi\cup \eta))\cdot J$,
      where:
      \begin{itemize}
      \item $\xi$ and $\eta$ are chords.
      \item Neither $b_1$ nor $b_2$ appears in the boundary of $\xi$.
      \item $c_1$ appears in the boundary of $\xi$ and $c_2$ appears in the boundary of $\eta$ 
        with the opposite orientation.
      \item The points $b_1$ and $b_1'$ do not appear in
        the support of $\xi$ (or, consequently, $\eta$). 
      \end{itemize}
    \item 
      \label{typeU:XX}
      It has the form $x=I\cdot (a(\xi\cup\sigma)\otimes a_o'(\xi\setminus\sigma'))\cdot J$ where
      $\xi$ is a chord other than $[c_1,c_2]$  such that
      $\sigma'\subset \xi$ but $\sigma'$ is not contained in the
      interior of $\xi$ (so $\xi$ has $c_2$ as one endpoint);
      or
      $x=I\cdot (a(\xi\setminus \sigma)\otimes
      a_o'(\xi\cup\sigma'))\cdot J$, where $\xi$ is a chord
      such that $\sigma\subset \xi$ but $\sigma$ is not
      contained in the interior of $\xi$ (so $\xi$ has $c_1$ as an
      endpoint).  (There are two subcases:
      either all local multiplicities are $0$ or $1$, or there is some
      local multiplicity of two, which occurs when $\sigma$ or
      $\sigma'$ is contained in $\xi$.)
      \end{enumerate}
      Near-chords for under-slides are illustrated in
      Figure~\ref{fig:UnderNearChords}.  In this subsection we will
      almost always call these elements simply ``near-chords.'' (The
      reader is forewarned that there is a different set of
      near-chords for over-slides,
      Definition~\ref{def:OverNearChord}, used in
      Section~\ref{subsec:Over-Slides}.) 

      When Convention~\ref{label:conv} does not hold,
      in the definition of ``near-chords'', we switch the roles of the two tensor factors.
\end{definition}

Note that near-chords of Type~(\ref{typeU:Sigma}) are short, so they appear in the differential of any arc-slide bimodule,
by hypothesis.

The aim of the present subsection is to verify the following special
case of Proposition~\ref{prop:UniqueArc-Slide}:

\begin{proposition}
  \label{prop:CalculateUnder-Slide}
  If $N$ is a stable arc-slide bimodule for an under-slide
  then the differential on $N$ contains precisely the near-chords.
\end{proposition}

We return to the proof at the end of the subsection, after some
preliminary results.  We first establish that the near-chords listed
above are all the elements of the near-diagonal subalgebra of
grading~$-1$.

\begin{lemma}\label{lem:under-slide-grading}
  In the near-diagonal subalgebra of an under-slide $m\co\PMC\to\PMC'$,
  there are no elements of positive grading;
  the elements of
  grading~$0$ are the idempotents; and the basic elements of
  grading~$-1$ are the near-chords for under-slides.
\end{lemma}

\begin{proof}
  Assume that the arc-slide satisfies Convention~\ref{label:conv}.
  Let $a=a_L\otimes a_R$ be a basic generator of the near-diagonal subalgebra of
  grading $\ge -1$.  Then $a$
  uniquely determines a pair of generators $\x,\y\in\Gen(\HD(m))$ and a
  non-negative domain~$\Domain\in \pi_2(\x, \y)$ in $\HD(m)$, with $\supp(a) = \bdy^\bdy \Domain$.
  By Proposition~\ref{prop:near-diag-grading}, the grading of $a$ is
  determined by $\Domain$, via formula~\eqref{eq:gr-basic-near-diag}.
  We will proceed by
  cases on the types of $\x$ and $\y$ (whether they are of type $\lsub{C}X$, $X_C$,
  or~$Y$) and the multiplicities of~$\Domain$ in the 6 regions
  $\sigma_{+}$, $\sigma$, $\sigma_-$, $\sigma'_+$, $\sigma'$,
  and~$\sigma'_-$.

  Recall from Lemma~\ref{lem:NegGradings} that for any $a_0 \in \Alg(\PMC)$, we have
  \[
  \iota(a_0) \le -m/2
  \]
  where $m$ is the minimal number of intervals needed to get the
  multiplicities in $\supp(a_0)$.  Thus for an element~$a$ of the
  near-diagonal subalgebra,
  \[
  \iota(a) \le -(m_L + m_R)/2,
  \]
  where $m_L$ and $m_R$ are the number of intervals needed to cover
  $\supp(a_L)$ and $\supp(a_R)$, respectively.  Note that
  $m_L=m_R$ unless $n_\sigma(a)$ or $n_{\sigma'}(a)$ is non-zero.

  If all six multiplicities $n_{\sigma_+}, \dots, n_{\sigma'_-}$
  are~$0$ then $\Domain$ consists of a union of horizontal strips.  Since
  the correction term is zero in this case, $\Domain$ can have at most one
  connected set of horizontal strips, as in Type (\ref{typeU:Generic}), or
  no strips at all, in which case $a$ is an idempotent.  In the case analysis that
  follows we will assume that not all multiplicities are zero.

  In general, there are other constraints on these
  multiplicities:
  \begin{itemize}
  \item The multiplicity difference across any $\alpha$-arc is at most~$1$.
  \item The multiplicity differences are constrained by the
    idempotents:
    \begin{align}
      \label{eq:under-mult-l}
      n_{\sigma_+} - n_\sigma + n_{\sigma'_+} - n_{\sigma'_-}
        &= \text{$-1$, $0$, or $+1$}\\
      \label{eq:under-mult-r}
      n_{\sigma_+} - n_{\sigma_-} + n_{\sigma'} - n_{\sigma'_-}
        &= \text{$-1$, $0$, or $+1$}
    \end{align}
    where the right hand side is determined by what happens to the occupancy of
    the $C$ idempotent on the left in \eqref{eq:under-mult-l} and on the right
    in \eqref{eq:under-mult-r}. This will be spelled out case by case below.
  \end{itemize}
  There are also some more constraints that depend more closely on the
  idempotents. For instance, if $C$ is not occupied on the left in
  either the starting or ending idempotents, then $n_{\sigma_+} =
  n_{\sigma}$ and $n_{\sigma'_+} = n_{\sigma'_-}$, as no strand can
  start or end at the endpoints of~$C$.

  \textbf{Case $\lsub{C}X \to \lsub{C}X$.}
  In this case we have
  \begin{align*}
    n_{\sigma_+} - n_\sigma + n_{\sigma'_+} - n_{\sigma'_-} &= 0\\
    n_{\sigma_+} - n_{\sigma_-} = n_{\sigma'} - n_{\sigma'_-} &= 0.
  \end{align*}
  (The second set of equations come from the fact that the pair~$C$ is not
  occupied on the right in either the initial or final idempotent, so
  there can be no strand starting or ending there.)
  According to Proposition~\ref{prop:near-diag-grading}, the
  correction to the grading is given by $c(I, \supp(a)) + c(J, \supp(a))
  = \frac{1}{2}(n_{\sigma'_+} - n_{\sigma'})$.

  The linear equations tell us that the
  multiplicities are of the following forms:

  \begin{center}
  \begin{tabular}{McMcMcMcMcMcMc}
    \toprule
    \text{Corr}&n_{\sigma_+}&n_{\sigma}&n_{\sigma_-}&n_{\sigma'_+}&n_{\sigma'}&n_{\sigma'_-}\\
    \midrule
    \epsilon/2 & m & m+\epsilon & m & l+\epsilon & l & l\\
    \bottomrule
  \end{tabular}
  \end{center}
  Here $\epsilon\in \{-1,0,1\}$ (as the difference in multiplicities
  is at most one), and the correction to the grading is $\epsilon/2$,
  as indicated.

  If $\epsilon=-1$, we would need to have $\iota(a) = -\OneHalf$,
  which is not possible with the given multiplicities.

  If $\epsilon=0$, we have complete horizontal
  strips, giving a near-chord of
  Type (\ref{typeU:Generic}) or an idempotent (if there are no strips at all).

  If $\epsilon=+1$, the left side of~$\Domain$ will need at
  least two intervals to cover it, leaving only one interval for the
  right.  This implies that $l=0$ and $m\in\{0,1\}$.  In both cases,
  this gives a near-chord of Type (\ref{typeU:XX}).

  In
  summary, the possibilities are:

  \begin{center}
  \begin{tabular}{cMcMcMcMcMcMcMcMc}
    \toprule
    Type&\text{Corr}&\text{Grading}&n_{\sigma_+}&n_{\sigma}&n_{\sigma_-}&n_{\sigma'_+}&n_{\sigma'}&n_{\sigma'_-}\\
    \midrule
    \ref{typeU:Generic} & 0         & -1 & 1 & 1 & 1 & 0 & 0 & 0\\
    \ref{typeU:Generic} & 0         & -1 & 0 & 0 & 0 & 1 & 1 & 1\\
    \ref{typeU:Generic} & 0         & -1 & 1 & 1 & 1 & 1 & 1 & 1\\
    \ref{typeU:XX} & +\OneHalf & -1 & 0 & 1 & 0 & 1 & 0 & 0\\
    \ref{typeU:XX} & +\OneHalf & -1 & 1 & 2 & 1 & 1 & 0 & 0\\
    \bottomrule
  \end{tabular}
  \end{center}

  \textbf{Case $\lsub{C}X \to X_C$.}
  In this case we have
  \begin{align*}
    n_{\sigma_+} - n_\sigma + n_{\sigma'_+} - n_{\sigma'_-} &= 1\\
    n_{\sigma_+} - n_{\sigma_-} + n_{\sigma'} - n_{\sigma'_-} &= 1.
  \end{align*}
  The
  correction to the grading is given by 
  \[c(I, \supp(a)) + c(J, \supp(a))
  = \OneQuart(- n_{\sigma} + n_{\sigma_-} + n_{\sigma'_+} -
  n_{\sigma'}) =0.\]
  The linear equations tell us that the multiplicities are given by
  \begin{center}
  \begin{tabular}{McMcMcMcMcMcMc}
    \toprule
    \text{Corr}&n_{\sigma_+}&n_{\sigma}&n_{\sigma_-}&n_{\sigma'_+}&n_{\sigma'}&n_{\sigma'_-}\\
    \midrule
    0 & m & m+\epsilon-1 & m+\delta-1 & l+\epsilon & l+\delta & l \\
    \bottomrule
  \end{tabular}
  \end{center}
  where $\epsilon,\delta\in\{0,1\}$.

  The only solutions to
  these equations which yield a connected domain on both sides
  are the following:

  \begin{center}
  \begin{tabular}{cMcMcMcMcMcMcMcMc}
    \toprule
    Type&\text{Corr}&\text{Grading}&n_{\sigma_+}&n_{\sigma}&n_{\sigma_-}&n_{\sigma'_+}&n_{\sigma'}&n_{\sigma'_-}\\
    \midrule
    \ref{typeU:Generic} & 0 & -1 & 1 & 0 & 0 & 0 & 0 & 0\\
    \ref{typeU:Generic} & 0 & -1 & 1 & 1 & 1 & 1 & 1 & 0\\
    \ref{typeU:Generic} & 0 & -1 & 0 & 0 & 0 & 1 & 1 & 0\\
    \bottomrule
  \end{tabular}
  \end{center}

  \textbf{Case $\lsub{C}X \to Y$.}
  In this case we have
  \begin{align*}
    n_{\sigma_+} - n_\sigma + n_{\sigma'_+} - n_{\sigma'_-} &= 0\\
    n_{\sigma_+} - n_{\sigma_-} + n_{\sigma'} - n_{\sigma'_-} &= 1.
  \end{align*}
  The
  correction to the grading is given by 
  \[c(I, \supp(a)) + c(J, \supp(a))
  = \OneQuart(n_{\sigma_+} - n_{\sigma} + n_{\sigma'_+} - 2
    n_{\sigma'} + n_{\sigma'_-}) = \OneHalf(-n_{\sigma'} + n_{\sigma'_-}).
  \]
  The linear equations tell us that the multiplicities are

  \begin{center}
  \begin{tabular}{McMcMcMcMcMcMc}
    \toprule
    \text{Corr}&n_{\sigma_+}&n_{\sigma}&n_{\sigma_-}&n_{\sigma'_+}&n_{\sigma'}&n_{\sigma'_-}\\
    \midrule
    -\delta/2 & m & m+\epsilon & m+\delta-1 & l+\epsilon & l+\delta & l \\
    \bottomrule
  \end{tabular}
  \end{center}
  with $\delta \in \{0,1\}$ and $\epsilon \in \{-1,0,1\}$.

  The solutions to
  these equations which can have grading $\ge -1$ are

  \begin{center}
  \begin{tabular}{cMcMcMcMcMcMcMcMc}
    \toprule
    Type&\text{Corr}&\text{Grading}&n_{\sigma_+}&n_{\sigma}&n_{\sigma_-}&n_{\sigma'_+}&n_{\sigma'}&n_{\sigma'_-}\\
    \midrule
    \ref{typeU:Sigma} &-\OneHalf &-1 & 0 & 0 & 0 & 0 & 1 & 0\\
    \ref{typeU:ExtraSigma} & 0        &-1 & 1 & 1 & 0 & 0 & 0 & 0\\
    \bottomrule
  \end{tabular}
  \end{center}

  \textbf{Case $X_C \to \lsub{C}X$.}  
  This is related to the case $\lsub{C}X \to X_C$ by rotating the
  diagram $180^\circ$.  Again, the solutions are all of Type (\ref{typeU:Generic}).

  \textbf{Case $X_C \to X_C$.}
  This is related to the case $\lsub{C}X \to \lsub{C}X$ by rotating the diagram $180^\circ$.
  The solutions are idempotents or near-chords of Type (\ref{typeU:Generic}) or (\ref{typeU:XX}).

  \textbf{Case $X_C \to Y$.}
  This is related to the case $\lsub{C}X \to Y$ by rotating the diagram $180^\circ$.  The
  solutions are of Type (\ref{typeU:Sigma}) or (\ref{typeU:ExtraSigma}).

  \textbf{Case $Y \to \lsub{C}X$.}
  In this case we have
  \begin{align*}
    n_{\sigma_+} - n_\sigma + n_{\sigma'_+} - n_{\sigma'_-} &= 0\\
    n_{\sigma_+} - n_{\sigma_-} + n_{\sigma'} - n_{\sigma'_-} &= -1.
  \end{align*}
  The
  correction to the grading is given by
  \[c(I, \supp(a)) + c(J, \supp(a))
  = \OneQuart(n_{\sigma_+} - n_{\sigma} + n_{\sigma'_+} - 2
    n_{\sigma'} + n_{\sigma'_-}) = \OneHalf(-n_{\sigma'} + n_{\sigma'_-}).
  \]
  The linear equations tell us that the multiplicities are

  \begin{center}
  \begin{tabular}{McMcMcMcMcMcMc}
    \toprule
    \text{Corr}&n_{\sigma_+}&n_{\sigma}&n_{\sigma_-}&n_{\sigma'_+}&n_{\sigma'}&n_{\sigma'_-}\\
    \midrule
    -\delta/2 & m & m+\epsilon & m+\delta+1 & l+\epsilon & l+\delta & l \\
    \bottomrule
  \end{tabular}
  \end{center}
  with $\delta \in \{-1,0\}$ and $\epsilon \in \{-1,0,1\}$.

  The solutions which can have grading $\ge -1$ are

  \begin{center}
  \begin{tabular}{cMcMcMcMcMcMcMcMc}
    \toprule
    Type&\text{Corr}&\text{Grading}&n_{\sigma_+}&n_{\sigma}&n_{\sigma_-}&n_{\sigma'_+}&n_{\sigma'}&n_{\sigma'_-}\\
    \midrule
    \ref{typeU:FewerSigma} & 0          & -1 & 0 & 1 & 1 & 1 & 0 & 0\\
    \ref{typeU:FewerSigma} & 0          & -1 & 0 & 0 & 1 & 0 & 0 & 0\\
    \ref{typeU:FewerSigma} & 0          & -1 & 0 & 0 & 1 & 1 & 1 & 1\\
    \ref{typeU:FewerSigma} & \OneHalf & -1 & 0 & 0 & 0 & 1 & 0 & 1\\
    \ref{typeU:FewerSigma} & \OneHalf & -1 & 1 & 1 & 1 & 1 & 0 & 1\\
    \bottomrule
  \end{tabular}
  \end{center}

  \textbf{Case $Y \to X_C$.}
  This is related to the case $Y \to \lsub{C}X$ by rotating the diagram $180^\circ$.  The
  solutions are of Type (\ref{typeU:Generic}) and (\ref{typeU:FewerSigma}).

  \textbf{Case $Y \to Y$.}
  In this case we have
  \begin{align*}
    n_{\sigma_+} - n_\sigma + n_{\sigma'_+} - n_{\sigma'_-} &= 0\\
    n_{\sigma_+} - n_{\sigma_-} + n_{\sigma'} - n_{\sigma'_-} &= 0\\
    n_{\sigma} - n_{\sigma_-} = n_{\sigma'_+} - n_{\sigma'} &= 0.
  \end{align*}
  (The last equations come from the fact that the $B$
  strand is not
  occupied in either the initial or final idempotent on either the
  left or right.)
  The
  correction to the grading is given by 
  \[c(I, \supp(a)) + c(J, \supp(a))
  = \OneHalf(n_{\sigma_+} - n_{\sigma}
    - n_{\sigma'} + n_{\sigma'_-}) = n_{\sigma_+} - n_\sigma.
  \]
  The linear equations tell us that the multiplicities are

  \begin{center}
  \begin{tabular}{McMcMcMcMcMcMc}
    \toprule
    \text{Corr}&n_{\sigma_+}&n_{\sigma}&n_{\sigma_-}&n_{\sigma'_+}&n_{\sigma'}&n_{\sigma'_-}\\
    \midrule
    -\epsilon & m & m+\epsilon & m+\epsilon & l+\epsilon & l+\epsilon & l \\
    \bottomrule
  \end{tabular}
  \end{center}
  with $\epsilon \in \{-1,0,1\}$.

  The solutions which can have grading $\ge -1$ are

  \begin{center}
  \begin{tabular}{cMcMcMcMcMcMcMcMc}
    \toprule
    Type&\text{Corr}&\text{Grading}&n_{\sigma_+}&n_{\sigma}&n_{\sigma_-}&n_{\sigma'_+}&n_{\sigma'}&n_{\sigma'_-}\\
    \midrule
    \ref{typeU:Generic} & 0 & -1 & 1 & 1 & 1 & 0 & 0 & 0\\
    \ref{typeU:Generic} & 0 & -1 & 0 & 0 & 0 & 1 & 1 & 1\\
    \ref{typeU:Generic} & 0 & -1 & 1 & 1 & 1 & 1 & 1 & 1\\
    \ref{typeU:YY} & 1 & -1 & 1 & 0 & 0 & 0 & 0 & 1\\
    \bottomrule
  \end{tabular}
  \end{center}

  We have shown that only the idempotents have degree~$0$ and that the
  elements of degree~$-1$ are all near-chords.  Checking through the
  various cases of near-chords verifies that they all appear in one of
  the cases above.  This completes the proof of
  Lemma~\ref{lem:under-slide-grading}, under the hypotheses of Convention~\ref{label:conv}.

  When Convention~\ref{label:conv} does not hold, the above discussion
  nonetheless applies, by reflecting the Heegaard diagram
  vertically. For example, the analogues of the
  formulas from Proposition~\ref{prop:near-diag-grading}
  (Equations~\eqref{eq:UnderslideGradingCorrection}
  and~\eqref{eq:OverslideGradingCorrection}) hold, where we relabel
  regions after reflection, as indicated in
  Figure~\ref{fig:ReflectedDiagrams}. For the new formulas, we swap the roles
  of $\sigma_+$, $\sigma_-$, and $\sigma$ with $\sigma_+'$, $\sigma_-'$, and $\sigma'$, respectively; so, for instance,
  the new correction formula for underslides (replacing Equation~\eqref{eq:UnderslideGradingCorrection}) reads:
  \begin{equation}
    \label{eq:UnderslideGradingCorrectionNoConv}
  c(I, \supp(a)) =
  \begin{cases}
    \OneQuart(n_{\sigma_+} - n_{\sigma}) &
      \text{$I$ of type $\lsub{C}X$}\\
    \OneQuart(-n_{\sigma'} + n_{\sigma_-'}) &
      \text{$I$ of type $X_C$}\\
    \OneQuart(n_{\sigma_+'} - n_{\sigma'} - n_{\sigma} + n_{\sigma_-}) &
      \text{$I$ of type $Y$.}
  \end{cases}
  \end{equation}

  With this understood, the above case analysis holds as stated, in the case where Convention~\ref{label:conv} is not met.
\end{proof}

\begin{remark}
  Instead of considering the grading on the near-diagonal subalgebra,
  we could instead prove an analogue of
  Lemma~\ref{lem:FactorDiagonalSubalgebra}: every element of the
  near-diagonal subalgebra can be factored into near-chords.
\end{remark}

\begin{remark}
  \label{rem:ConventionAgain}
In the proofs of the following sequence of lemmas, we will implicitly assume Convention~\ref{label:conv}.
This is not essential for the statements, and the proofs adapt easily to the other case (sometimes at the cost of switching the
orders of products, and replacing ``initial idempotents'' by ``terminal idempotents'').
\end{remark}
  
We break Proposition~\ref{prop:CalculateUnder-Slide} into a series of
lemmas. We focus first on the non-degenerate case (in which there are
at least two positions between $c_1$ and $c_2$; see
Definition~\ref{def:degen-slide}), and return to the degenerate case
in Lemma~\ref{lem:DegenerateUnderslide}.

\begin{lemma}
  \label{lem:UGenericExists}
  Let $N$ be a stable arc-slide bimodule for a non-degenerate
  under-slide.  Then the differential in $N$ contains all near-chords
  of Type~(\ref{typeU:Generic}).
\end{lemma}

\begin{proof}
  Consider the intervals connecting $b_2$ (which is matched with $b_1$)
  to the basepoint $z$.  There are two such intervals; we call the one
  which does not contain $b_1$ the {\em outside region}.
  We say that a pointed matched circle $\PMC$ is {\em big enough} for
  the arc-slide if the number of matched pairs in the outside region
  exceeds the number of positions in the complement of the 
  outside region. By stabilizing, it suffices to consider the case of pointed
  matched circles which are big enough for the arc-slide. 

  A {\em very special length $3$ chord} is a special length $3$ chord (Definition~\ref{def:SpecialLengthThree})
  contained between $c_1$ and $c_2$, which is adjacent to $b_1$ in
  $\PMC$ and to $b_1'$ in $\PMC'$. (See
  Figure~\ref{fig:degen-slides}.)  Very special length $3$ chords, when
  they exist, need special attention.

  {\bf Claim 1}: {\em Suppose $x$ is a near-chord of
    Type~(\ref{typeU:Generic}) which is not a very special length $3$
    chord, and suppose moreover that neither $b_1$ nor $b_1'$ is
    contained in the interior of (the support of) $x$. Then $x$
    appears in the differential.}
  
  Claim 1 follows from the same argument used to prove
  Proposition~\ref{prop:DDidUnique}, by induction on the length
  of the support. Note that each special length $3$ chord which is
  not very special is adjacent to a length $1$ chord which appears in
  the differential by hypothesis. Thus, these chords all appear in the
  differential, by the same principle which established the existence
  of special length $3$ chords in \DD\ identity bimodule (see the end
  of the proof of Proposition~\ref{prop:DDidUnique}).
  
\smallskip
  Consider now a near-chord $x=I\cdot(a(\xi)\otimes a'_o(\xi'))\cdot
  J$ of Type~(\ref{typeU:Generic}) whose restricted support has length
  one.   This near-chord appears in the differential by hypothesis except
  in the special cases when $\xi$ contains either $b_1$ or $b_1'$ in
  its interior.

  In the case that $\xi$ contains $b_1$ or $b'_1$ in its interior
  there are two subcases, according to whether or not $dx=0$.  If
  $dx\neq 0$ it is straightforward to see that $dx$
  factors as a product $y\cdot z$ of two short near-chords, and hence
  $dx$ appears in the expression for $\partial^2$. (See the first row
  of Figure~\ref{fig:UnderU1}.)  This product has no alternative
  factorization, and hence $x$ must appear in the differential. 

   \begin{figure}
     \begin{center}
       \includegraphics[scale=.5]{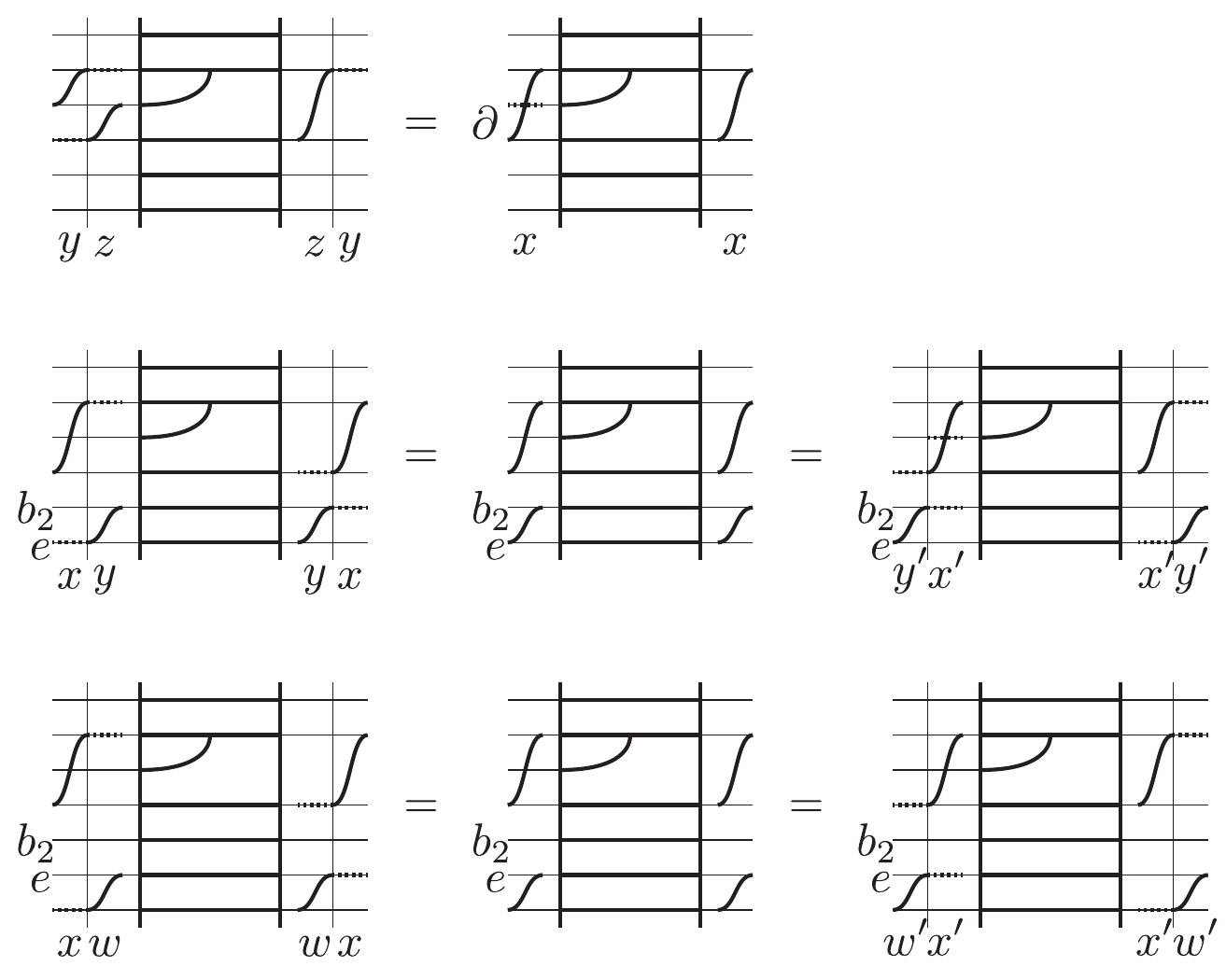}
     \end{center}
     \caption {\textbf{Existence of near-chords of
         Type~(\ref{typeU:Generic}) with length-$1$ restricted support.}
       \label{fig:UnderU1}
     }
   \end{figure}

  By induction on the length of the support (again, see the proof of Proposition~\ref{prop:DDidUnique}), we can conclude the following:

  {\bf Claim 2}: {\em
    Suppose $x$ is any near-chord of Type~(\ref{typeU:Generic}) with the following properties:
    \begin{itemize}
    \item the support of $x$ contains exactly one of $b_1$ or $b_1'$, and
    \item the position at $b_1$ or $b_1'$ (whichever is contained in the support of $x$) is occupied
      in both the initial and terminal idempotents for $x$ (on the $\PMC$ or the  $\PMC'$ side, respectively).
    \end{itemize}
    Then $x$ appears in the differential.} 

  (Note that any near-chord satisfying the criteria of Claim 2 has $dx\neq 0$.)
  
  Next we consider near-chords $x$ of Type~(\ref{typeU:Generic}) whose restricted support has length $1$ and with $dx=0$.  We consider two subcases: either
  $b_2$ lies below $b_1$ or above it.  For the time being, we make the following:

  {\bf Assumption 1}: Suppose that $b_2$ lies below
  $b_1$.  

  Let $e$ denote the position immediately below $b_2$. Clearly, $e\neq
  c_2$ (otherwise $\PMC'$ would not be a valid pointed matched
  circle); by Assumption 1, $e\neq c_1$.  In addition, $e$ cannot be
  matched with the position directly below~$b_1$ (since otherwise
  $\PMC$ would not be a valid pointed matched circle).

  {\bf Claim 3}: {\em Under Assumption 1, if $x$ is a
    Type~(\ref{typeU:Generic}) near-chord whose restricted support has
    length $1$, $b_1$ is in the support of $x$, and $e$ is contained
    in the left (hence also the right)
    idempotent of $x$ on the $\PMC$ side, then $x$ appears in the
    differential.}

  There is a short near-chord $y$ with support
  $[e,b_2]$
  and $x\cdot y \neq 0$. Now, there are Type~(\ref{typeU:Generic})
  near-chords $x'$ and $y'$ with $x\cdot y=y'\cdot x'$ with
  $\supp(x)=\supp(x')$, $\supp(y)=\supp(y')$ and $dx'\neq 0$.  (See
  the second row in Figure~\ref{fig:UnderU1}.) In fact, this is the
  only other way to factor $x\cdot y$ in the near-diagonal subalgebra.
  We already showed that $x'$ occurs in the differential, and $y'$
  occurs in the differential by hypothesis; it follows that $x$ occurs
  in the differential as well. Thus, we have established Claim~$3$.

  We now generalize Claim~$3$.

  {\bf Claim 4}:
  \emph{Under Assumption 1, if $x$ is a Type~(\ref{typeU:Generic}) near-chord whose
    restricted support has length $1$ and $b_1$ is contained in the 
    support of $x$, then $x$ appears in the differential.}
  
  If $e$ does not appear in the left idempotent of $x$ on the $\PMC$ side,
  since our pointed matched circle is big enough, we can find
  some near-chord $w$ of Type~(\ref{typeU:Generic}) supported entirely
  in the outside region, terminating at $e$, 
  and so that $x \cdot w\neq 0$.  Again, the product $x\cdot w$ does
  not appear in the differential of any algebra element, and it has a
  unique alternative factorization as $w'\cdot x'$, with
  $\supp(w)=\supp(w')$ and $\supp(x)=\supp(x')$. In Claim 3, we
  established that $x'$ appears in the differential (since $e$ is in
  the left idempotent of $x'$ on the $\PMC$ side). So, to show that $x$
  appears in the differential, it suffices to show that $w'$ appears
  in the differential.  Note that $b_1$ is not contained in the
  support of $w$. If $b_1'$ is also not in the support of $w$, then
  $w'$ appears in the differential by Claim 1; if $b_1'$ is contained
  in the support of $w$, then, since $dx$ was assumed to vanish, $b_1$
  is not contained in the left idempotent of $x$ (and hence of $w'$) on the
  $\PMC$ side, so $b_1'$ is contained in the left idempotent of $w'$ on the
  $\PMC'$ side; and hence $w'$ appears in the differential by Claim 2.
  We conclude that $x$ appears in the differential, proving
  Claim~$4$. (See the third
  line of Figure~\ref{fig:UnderU1}.)

  Next we consider the case when Assumption 1 does not hold, i.e.:
  
  {\bf Assumption 1$\boldsymbol{'}$}: {\em Suppose that $b_2$ lies above $b_1$.}

  Let $f$ denote the position immediately above $b_2$. The same argument
  used to establish Claim~3 (with some products reversed) shows the
  following:

  {\bf Claim 3$\boldsymbol{'}$}: {\em Under Assumption 1$'$, if $x$ is a
    Type~(\ref{typeU:Generic}) near-chord whose restricted support has
    length $1$, $b_1$ is in the support of $x$, and $f$ is contained
    in the left idempotent of $x$ on the $\PMC$ side, then $x$ appears in the
    differential.}

  Using this, the argument used to deduce Claim~4 can be modified to give the following:

  {\bf Claim 4$\boldsymbol{'}$}: {\em Under Assumption 1$'$, if $x$ is a Type~(\ref{typeU:Generic})
  near-chord whose restricted support has length $1$ and $b_1$ is contained in the support of $x$, then $x$ appears in the differential.}

  {\bf Claim 5}: \emph{If $x$ is any Type~(\ref{typeU:Generic}) near-chord whose
  restricted support has length $1$, then $x$ appears in the
  differential of $N$.}

  The only cases of Claim~5 not covered by hypothesis are
  those where $b_1$ or $b_1'$ are contained in the interior of $x$;
  but in this case Claim 4 or Claim 4$'$ applies (perhaps after reversing the
  roles of $\PMC$ and $\PMC'$).

  \textbf{Claim 6}: \emph{If $x$ is a near-chord of Type~(\ref{typeU:Generic}) which is not a
    special length~$3$ chord, then $x$ appears in the differential.}

  If the restricted support of~$x$ has length~$1$, this is covered by
  Claim~$5$.
  Otherwise, we claim that $dx$ contains terms of the form $y\cdot t$ with the following
  properties:
  \begin{itemize}
  \item Each of $y$ and $t$ is of Type~(\ref{typeU:Generic}).
  \item The product $y\cdot t$ has no alternative factorization.
  \item The product $y\cdot t$ does not appear in the differential
    of any other basic generator of the near-diagonal sub-algebra.
  \end{itemize}
  To ensure that the factorization is into two near-chords of
  Type~(\ref{typeU:Generic}), we break at a position in $x$ other than
  $b_1$, $b_2$, $c_1$, and $c_2$, which in turn can be done since we
  assumed the arc-slide is non-degenerate.  Therefore an induction on the
  length of the restricted support establishes Claim~$6$.

  Finally, since any special length $3$ chord $x$ is adjacent to a
  near-chord of Type~(\ref{typeU:Generic}) with length $1$ restricted
  support, which in turn appears in the differential according to
  Claim 5 established above (even in the case of a very special length
  three chord), the argument from Proposition~\ref{prop:DDidUnique}
  ensures the existence of $x$ in the differential, as well.%
\end{proof}

By hypothesis, the differential on $N$ contains all near-chords of
Type~(\ref{typeU:Sigma}).

\begin{lemma}
  \label{lem:UExtraExist}
  Let $N$ be a stable arc-slide bimodule for a non-degenerate under-slide. Then the
  differential on $N$ 
  contains all near-chords of
  Types~(\ref{typeU:ExtraSigma}) and~(\ref{typeU:XX}).
\end{lemma}

\begin{figure}
    \begin{center}
      \includegraphics[scale=.5]{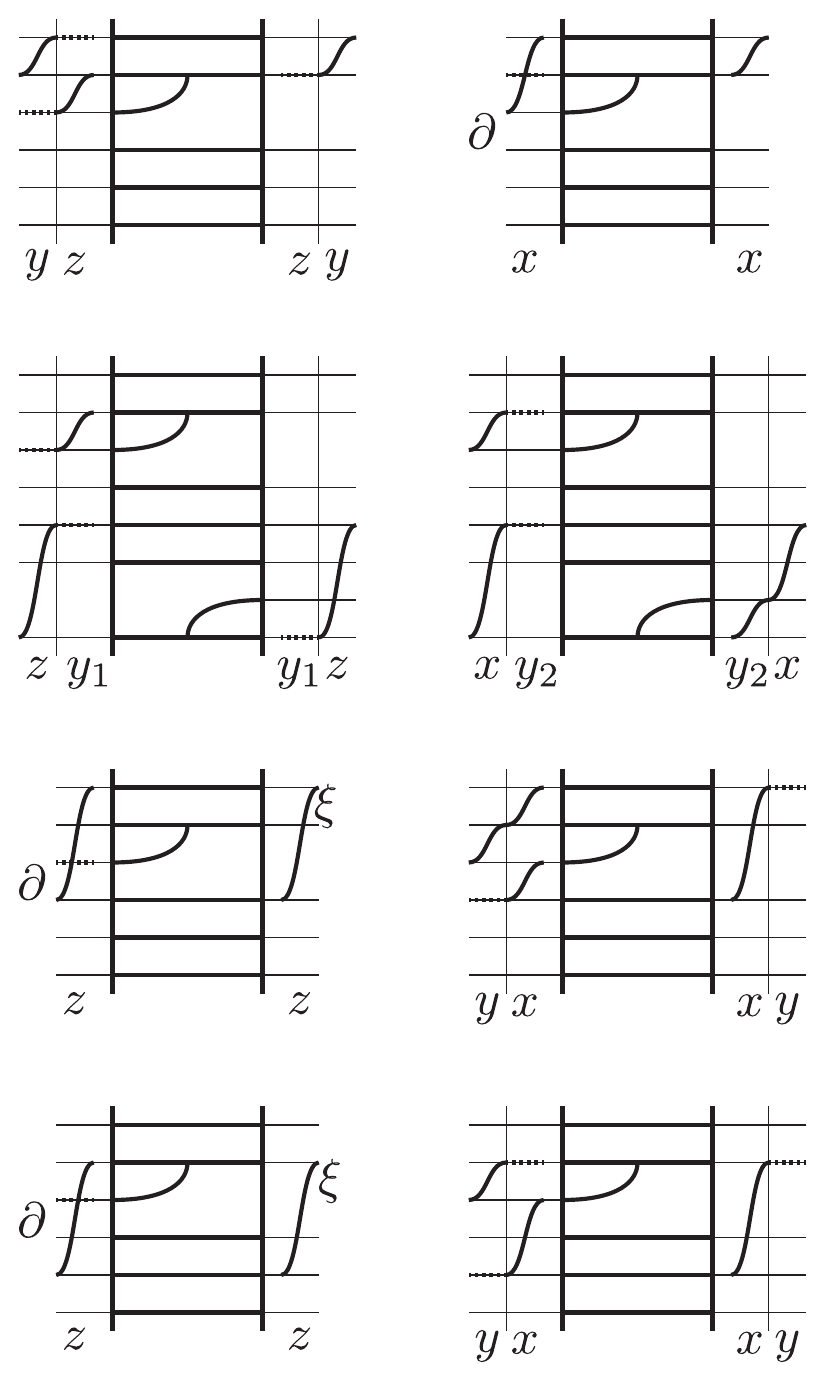}
    \end{center}
    \caption {\textbf{Existence of near-chords of
        Type~(\ref{typeU:ExtraSigma}), near-chords of
        Type~(\ref{typeU:XX}), and certain near-chords of
        Type~(\ref{typeU:FewerSigma}).}
      \label{fig:UnderEasy}
      In the left column, we illustrate terms we 
      already know contribute to $\partial^2$. On the
      right, we have the only possible alternative terms 
      which could cancel (some) terms on the left. Recall that we draw
      $a\cdot b$ with $a$ on the outside and $b$ on the inside;
      compare Figure~\ref{fig:BreakInProduct}. Here and later, some of
      the horizontal lines in the algebra elements have been suppressed.}
\end{figure}

\begin{proof}
  Let $x$ be a near-chord of Type~(\ref{typeU:ExtraSigma}).
  As illustrated in the first line in Figure~\ref{fig:UnderEasy},
  we can find near-chords $y$ and $z$ so that:
  \begin{itemize}
  \item $y$ is of Type~(\ref{typeU:Generic}),
  \item $z$ is of Type~(\ref{typeU:Sigma}),
  \item $y\cdot z$ has no alternative factorizations as a product of
    two homogeneous elements in the near-diagonal subalgebra, with
    non-trivial support, and
  \item $y\cdot z$ appears in $dx$ (and not in $dw$ for any basic
    generator $w\neq x$).
  \end{itemize}
  It follows that $x$ must appear in the differential, so that
  $\partial^2=0$. 

  Next, let $x$ be a near-chord of Type~(\ref{typeU:XX}).
  There are near-chords $y_1$ and $y_2$
  of Type~(\ref{typeU:Sigma}) and a near-chord $z$ of
  Type~(\ref{typeU:Generic})
  with the property that
  $z\cdot y_1=x\cdot y_2$ has exactly these two factorizations
  (with non-trivial support in the near-diagonal subalgebra), and
  $z\cdot y_1$ does not appear in the differential
  of any other algebra element. Since $z$ and $y_1$ appear in the
  differential (according to Lemmas~\ref{lem:UGenericExists}; and the 
  fact that all near-chords of Type~(\ref{typeU:Sigma}) are short), 
  it follows that $x$ must, as well. This is illustrated
  in the second line of Figure~\ref{fig:UnderEasy}.
\end{proof}

\begin{lemma}
  \label{lem:UFewerExist}
  Let $N$ be a stable arc-slide bimodule  for a non-degenerate under-slide. Then the differential on $N$ contains
  all near-chords of Type~(\ref{typeU:FewerSigma}).
\end{lemma}

\begin{proof}
  We will show that all near-chords of Type~(\ref{typeU:FewerSigma}) appear in the differential
  of any idempotent which has at least two occupied positions on both the $\PMC$ and the $\PMC'$ sides.
  This assumption on the idempotent can be made freely, in view of the
  stability hypothesis for $N$. 

  There
  are two types of near-chords of Type~(\ref{typeU:FewerSigma}): those
  for which the support has three moving strands, and those for which
  it has only two. In other words, after rotating $180^\circ$ if
  necessary, near-chords of Type~(\ref{typeU:FewerSigma}) have the
  form $I\cdot \bigl(a(\xi\setminus \sigma)\otimes
  a_o(\xi)\bigr)\cdot J$. 
  For the first type, $\sigma$ is a subset of the interior of $\xi$; 
  for the second type, $\sigma$ is a subset of $\xi$, but one of its
  boundary points is on the boundary of $\xi$.

  We handle first the case where the near-chord $x$ contains three
  strands.  In this case, we can find a near-chord $z$ of
  Type~(\ref{typeU:Generic}) and a near-chord $y$ of
  Type~(\ref{typeU:Sigma}) with the properties that $dz$ contains
  $y\cdot x$, as in the third line in Figure~\ref{fig:UnderEasy}.
  We claim that $y\cdot x$ has no alternative factorizations as
  a product of two near-chords with non-trivial support. More
  specifically, $y\cdot x$ has three moving strands, like~$x$.  
  The only other products of basic elements of the near-diagonal
  sub-algebra with the same support as $y\cdot x$ and three moving
  strands are products
  $y'\cdot x'$ where $y'$ has Type~(\ref{typeU:ExtraSigma}) and $x'$
  has Type~(\ref{typeU:FewerSigma}). However, a closer
  look at the idempotents shows that $y\cdot x$ does not equal
  $y'\cdot x'$.  Specifically, suppose with out loss of generality
  that two of the moving strands of $x$ are on the $\PMC$-side. Then,
  the final idempotent of $y\cdot x$ contains $C$ on the $\PMC'$
   side ($x$ is of type $Y \to X_C$), whereas the final idempotent of
  $y'\cdot x'$ contains $C$ on the $\PMC$ side ($x'$ is of type
  $Y \to \lsub{C}X$).

  It follows that $\partial^2=0$ forces $x$ to appear in the
  differential.

  \begin{figure}
  \begin{center}
      \includegraphics[scale=.5]{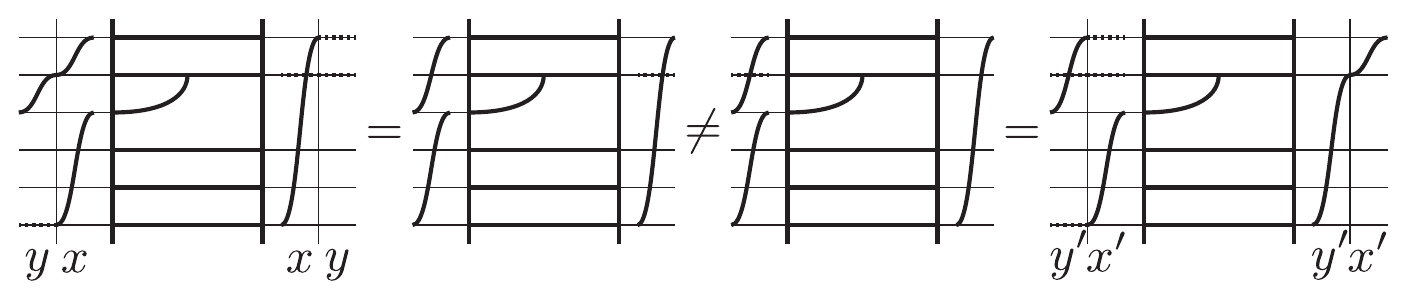}
    \end{center}
    \caption {{\bf Not an alternative factorization.}
      \label{fig:NotAFactorizationU}
      One of the cases in Figure~\ref{fig:UnderEasy} (third line, right hand
      column) might
      appear to have an alternate factorization; however, a more
      careful look at idempotents (as indicated) shows that this
      alternative factorization does not exist.}
  \end{figure}

  With one exception, the same argument also applies to near-chords $x$ of
  Type~(\ref{typeU:FewerSigma}) with only two moving strands, as
  illustrated in the fourth line in Figure~\ref{fig:UnderEasy}.
  The exception is when the restricted support of $\xi$ is
  $[c_2,c_1]\times [c_2,c_1]$.  In this case, there is no near-chord $z$ of
  Type~(\ref{typeU:Generic}) as required for the argument above, and we use a
  different argument.

  \begin{figure}
    \begin{center}
      \includegraphics[scale=.5]{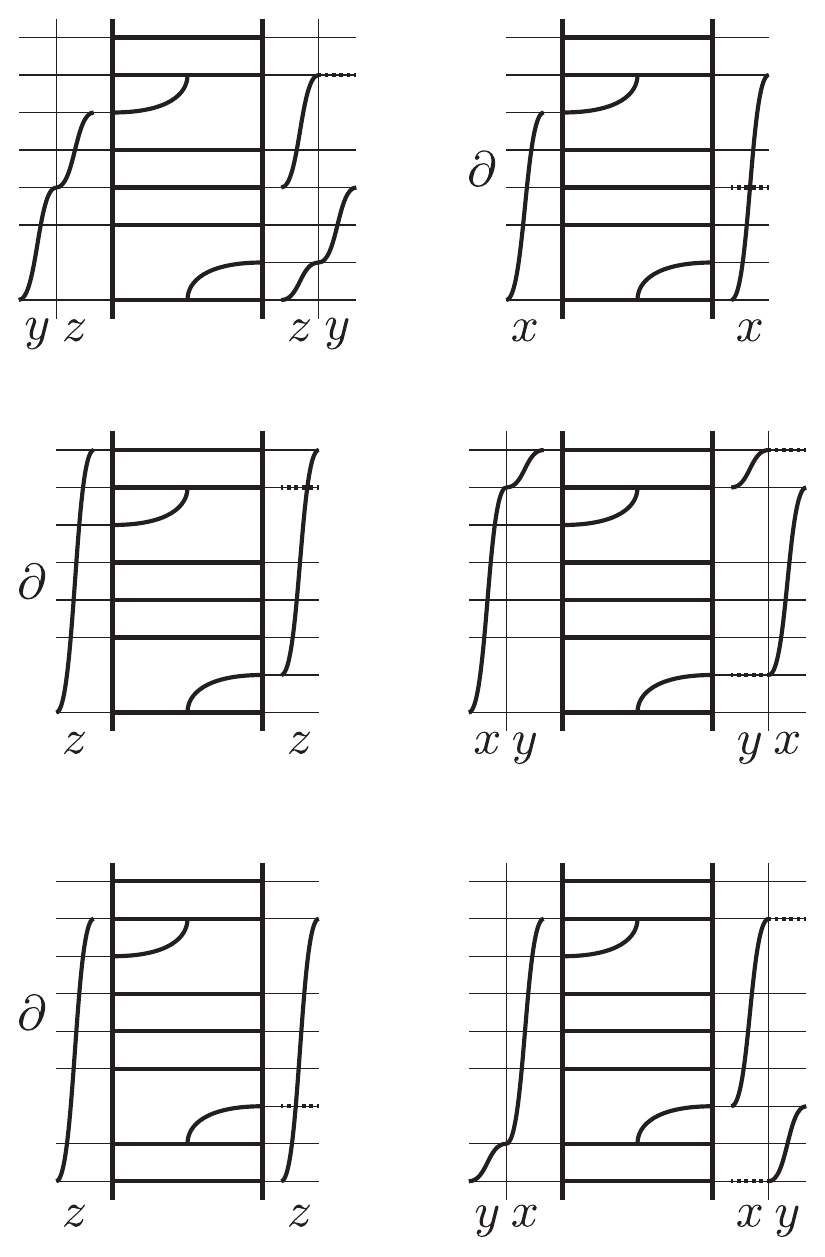}
    \end{center}
    \caption {\textbf{Existence of remaining near-chords of
        Type~(\ref{typeU:FewerSigma}).}  \label{fig:UnderFewer} We
      continue with the conventions from Figure~\ref{fig:UnderEasy}.
      This establishes the remaining cases of near-chords of
      Type~(\ref{typeU:FewerSigma}) which have not already been
      covered in Figure~\ref{fig:UnderEasy}.}
  \end{figure}
  There are three subcases of Type~(\ref{typeU:FewerSigma}) with
  restricted support $[c_2,c_1] \times [c_2,c_1]$, according to the
  placement of the initial
  idempotent. Specifically, if we rotate so that the support of the
  near-chord is $[c_2,c_1] \times [b_1',c_1]$, then the three
  (not necessarily distinct)
  subcases are:
  \begin{enumerate}
    \item The initial idempotent on the $\PMC'$ side
      contains a position in the open interval $(c_2,c_1)$. 
    \item The initial idempotent on the $\PMC'$
      side contains a position strictly above $c_1$.
    \item The initial idempotent on the $\PMC'$ side
      contains a position strictly below $c_2$.
  \end{enumerate}
  The three cases are illustrated in Figure~\ref{fig:UnderFewer}.
  In the above verification, we are using the hypothesis on our idempotent
  that in both $\PMC$ and $\PMC'$ there are at least two unoccupied positions.

  In the first of these cases, we can find near-chords $y$ and $z$
  with the following properties:
  \begin{itemize}
  \item $y$ is of Type~(\ref{typeU:FewerSigma}), but of the type which
  we have already verified appear in the differential.
  \item $z$ is of Type~(\ref{typeU:XX})
    (and hence appears in the differential, by Lemma~\ref{lem:UExtraExist}).
  \item The product $y\cdot z$ has no alternative factorizations into
    homogeneous elements in the near-diagonal subalgebra.
  \item The term $y\cdot z$ appears in the differential $dx$ of our
    near-chord $x$
    (and $y\cdot z$ does not appear in $dx'$ for any other homogeneous element of the
    near-diagonal subalgebra).
  \end{itemize}
  It follows from the above properties that
  $x$ appears in the differential.

  In the second case, we find a near-chord $y$ of
  Type~(\ref{typeU:Generic}) for which $x\cdot y$ has no alternative
  factorization, but $x\cdot y$ appears in the differential of
  another near-chord $z$ of Type~(\ref{typeU:FewerSigma}), which we have
  already verified occurs in the differential, and $x\cdot y$ does not
  appear in the differential of any other basic algebra element.
  
  In the third case, we find a near-chord $y$ of
  Type~(\ref{typeU:ExtraSigma}) for which $y\cdot x$ has no alternative
  factorization, and $y\cdot x$ appears in $dz$ for a near-chord
  $z$ of Type~(\ref{typeU:Generic}) which we have already verified
  occurs in the differential, and $y\cdot x$ does not appear in the
  differential of any other basic algebra element.
\end{proof}

\begin{lemma}
  \label{lem:UYYExist}
  Let $N$ be a stable arc-slide bimodule for a non-degenerate under-slide. Then, the
  differential in $N$ contains all near-chords of
  Type~(\ref{typeU:YY}).
\end{lemma}

\begin{figure}
\begin{center}
  \includegraphics[scale=.5]{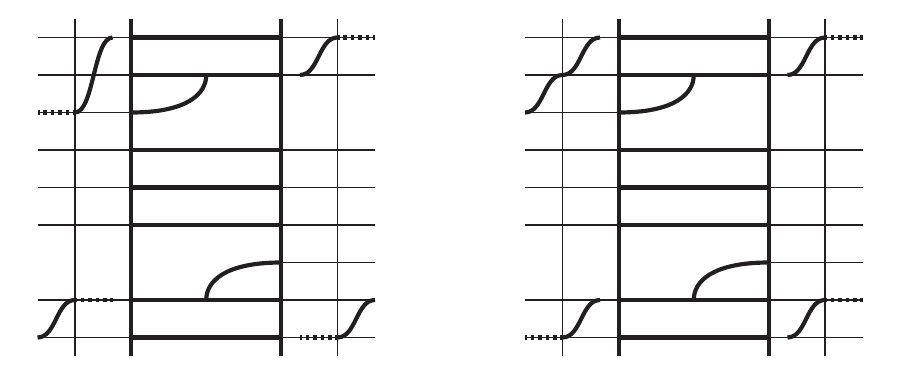}
\end{center}
\caption {{\bf Proof of Lemma~\ref{lem:UYYExist}.}
  \label{fig:UnderYY}
Verification that near-chords of Type~(\ref{typeU:YY}) occur in the differential.}
\end{figure}

\begin{proof}
  This follows from the
  observation that a near-chord of Type~(\ref{typeU:Sigma}) times one of
  Type~(\ref{typeU:YY}) has a unique alternative factorization as a
  near-chord of Type~(\ref{typeU:Generic}) times a near-chord of
  Type~(\ref{typeU:ExtraSigma}) (and it does not appear in the
  differential of any algebra element); together with
  Lemmas~\ref{lem:UGenericExists} and~\ref{lem:UExtraExist}. (See
  Figure~\ref{fig:UnderYY}.)
\end{proof}

We return now to the degenerate case.

\begin{lemma}
  \label{lem:DegenerateUnderslide}%
  Let $N$ be a stable arc-slide bimodule for a degenerate
  under-slide. Then the differential on $N$ contains all near-chords.
\end{lemma}

\begin{proof}
  We prove this in a sequence of claims:
  
  {\bf Claim 1}: {\em The differential contains all near-chords of Type~(\ref{typeU:Generic}) supported outside $[c_1,c_2]$.}

  This follows from a straightforward induction on the length of the support, as in Proposition~\ref{prop:DDidUnique}.
  
  {\bf Claim 2}:  {\em The differential contains all near-chords of Type~(\ref{typeU:ExtraSigma}).} 

  Apply the argument from
  Lemma~\ref{lem:UExtraExist} and Claim 1.

  {\bf Claim 3}: {\em The differential contains all near-chords $x$ of
    Type~(\ref{typeU:Generic}) with support $[e,c_2]$
    such that $b_1$ is contained in the 
    idempotent of $x$ on the $\PMC$ side, for any $e$ above $c_1$.}

  \begin{figure}
    \centering
    \includegraphics[scale=.5]{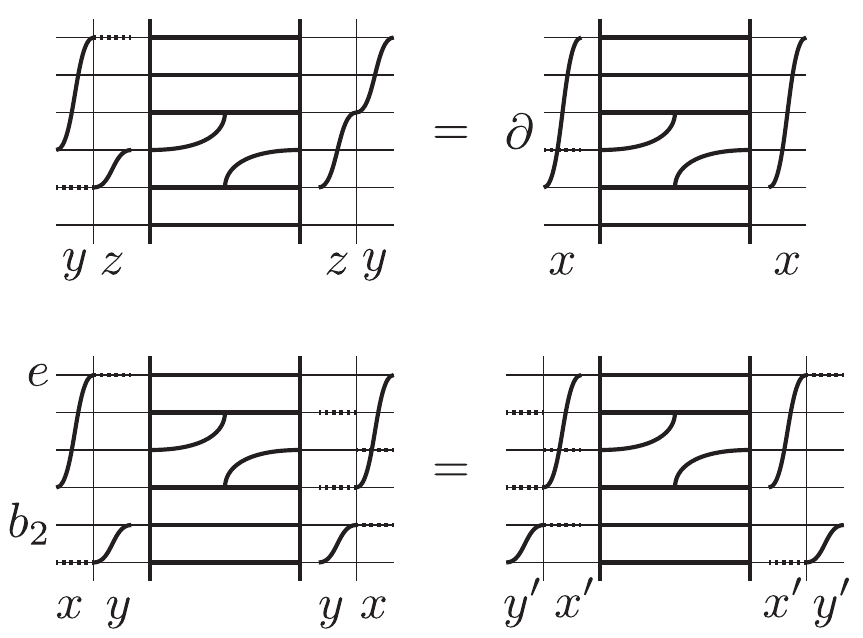}
    \caption{\textbf{Existence of Type~(\ref{typeU:Generic}) chords in
        the degenerate case.} Top line: Claim 3 of the proof of
      Lemma~\ref{lem:DegenerateUnderslide}. Bottom line: Claim 4 of
      the proof of Lemma~\ref{lem:DegenerateUnderslide}.}
    \label{fig:degen-chords}
  \end{figure}

  In this case, (a term in) $dx$ factors as $dx=y\cdot z$, where $y$ is of
  Type~(\ref{typeU:ExtraSigma}), while $z$ is the short near-chord
  added in the degenerate case (see
  Definition~\ref{def:short-near-chord}); see Figure~\ref{fig:degen-chords}.
  Thus, we know that $y$ and $z$ appear in the differential. Again,
  since there is no alternative factorization of this element, we
  conclude that $x$, too, must appear in the differential.

  {\bf Claim 4}: {\em The differential contains all near-chords $x$ of
    Type~(\ref{typeU:Generic}) with support $[e,c_2]$ for any $e$ above $c_1$.} 

  It remains to consider the case where $b_1$ is not in the initial
  idempotent $x$.  After stabilizing, it suffices to consider the case
  that the diagram is big enough, in the sense described in the proof
  of Lemma~\ref{lem:UGenericExists}.  
  Assume that $b_2$ is below
  $c_2$; the case that $b_2$ is above $c_1$ is similar. In this case,
  there is some chord $y$ of Type~(\ref{typeU:Generic}) in the outside
  region (in the sense of the proof of Lemma~\ref{lem:UGenericExists})
  which terminates at $b_2$, so that $x\cdot y$ does not vanish.  Now,
  there is a unique alternative factorization $x\cdot y=y'\cdot x'$,
  where $\supp(x)=\supp(x')$ and $\supp(y)=\supp(y')$.  The existence
  of $x'$ and $y'$ in the differential are ensured by Claims~3 and~1,
  respectively. The product $x \cdot y$ does not
  appear in the differential of any algebra element, and hence $x$
  must appear in the differential. (Again, see Figure~\ref{fig:degen-chords}.)

  {\bf Claim 5}:  {\em The differential contains all near-chords $x$ of
    Type~(\ref{typeU:Generic}) with support $[c_1,e]$ for any
    $e$ below $c_2$.} 

  The proof of Claim 4 applies {\em mutatis mutandis}.
  
  {\bf Claim 6}: {\em The differential contains all near-chord $x$ of
    Type~(\ref{typeU:Generic}).}
    
  It remains to consider cases where $[c_1,c_2]$ is contained in the
  interior of the support. For such a near-chord $x$, the differential
  has a term with a unique factorization as $y\cdot z$, where $y$ 
  $z$ exist because of Claims~4 and~5.
  
  Having established the existence of all near-chords of
  Type~(\ref{typeU:Generic}), the proofs of Lemmas~\ref{lem:UExtraExist},~\ref{lem:UFewerExist} 
  and~\ref{lem:UYYExist} apply to establish the existence of the
  remaining near-chords.
\end{proof}

\begin{proof} [Proof of Proposition~\ref{prop:CalculateUnder-Slide}]
  Lemmas~\ref{lem:grading-coeff-implies} and
  Lemma~\ref{lem:under-slide-grading} imply that the
  near-chords are the only elements in the near-diagonal subalgebra whose
  gradings are compatible with appearing in the differential.
  On the other hand,
  by Lemmas~\ref{lem:UGenericExists}, \ref{lem:UExtraExist},
  \ref{lem:UFewerExist}, \ref{lem:UYYExist} (all in the non-degenerate
  case) and Lemma~\ref{lem:DegenerateUnderslide} (in the degenerate case), we know that each
  near-chord appears in the differential.
\end{proof}

\subsection{Over-slides}
\label{subsec:Over-Slides}

We now turn to bimodules for over-slides. The appropriate notion of
near-chords in this case is given in
Definition~\ref{def:OverNearChord}. The first $6$ types correspond to
the types of near-chords for under-slides (except that one of the kinds
of (\ref{typeU:XX}) near-chords has no analogue for over-slides). There
are two kinds of near-chords for over-slides which have no under-slide
analogues, types (\ref{typeO:dC1C2}) and (\ref{typeO:Disconnected}).

\begin{figure}
    \begin{center}
      \input{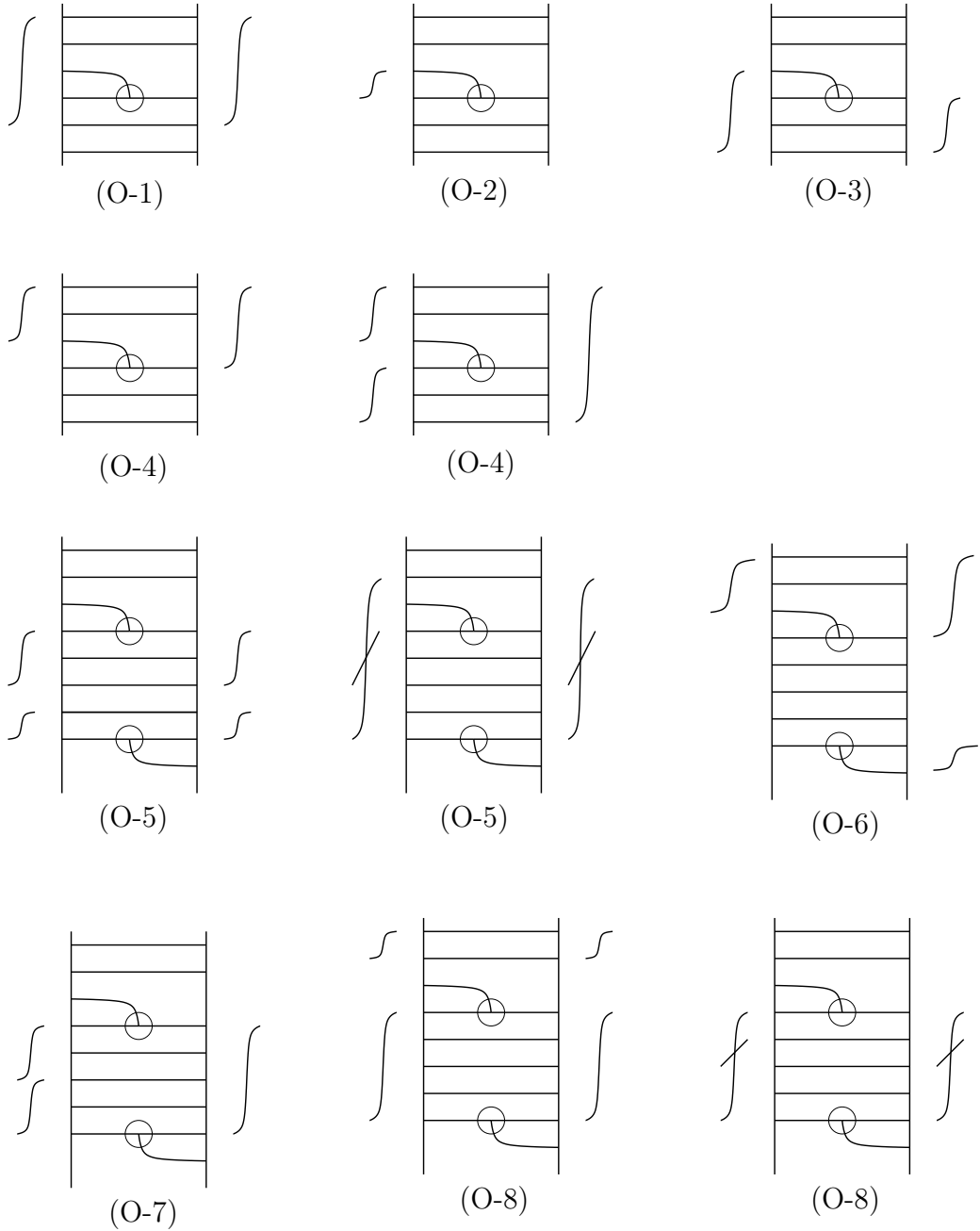}
    \end{center}
    \caption {{\bf Near-chords for over-slides.}
      \label{fig:NearChords} We have illustrated
    here examples of all the types of near-chords for over-slides
    appearing in Definition~\ref{def:OverNearChord}.
    (Note that there are two illustrations for Type~(\ref{typeO:FewerSigma})
    and Type~(\ref{typeO:Disconnected}.))}
\end{figure}

\begin{definition}
  \label{def:OverNearChord}
  \index{near-chord!for over-slides}%
  A (non-zero) basic algebra element $x$ in the near-diagonal
  subalgebra for an over-slide satisfying Convention~\ref{label:conv} is
  called a \emph{near-chord} (for the over-slide) if it satisfies
  any of the following $8$ conditions:
  \begin{enumerate}[label=(O-\arabic*),ref=O-\arabic*]
    \item 
      \label{typeO:Generic}
      It has the form
      $x=I\cdot (a(\xi)\otimes a'_o(\xi))\cdot J$,
      where $\xi$ is some chord in $\PMC$ neither of whose endpoints is $b_1$
      (so that it can be interpreted, as it is in the above expression,
      as a chord in $\PMC'$); furthermore, $\xi$ is required to be different from the chord
      $[c_2,c_1]$.
    \item
      \label{typeO:Sigma}
      It has the form
        $x=I\cdot(a(\sigma)\otimes 1)\cdot J$ or $I\cdot(1\otimes
        a_o'(\sigma'))\cdot J$, where $I$ and $J$ are near-complementary idempotents.
   \item
      \label{typeO:ExtraSigma}
        There is a chord $\xi$ with the 
        property that the interior of $\xi$ is disjoint from $\sigma$
        and the support of $\xi\cup\sigma$ is connected,
        and
        $x=I\cdot(a(\xi\cup \sigma)\otimes a'_o(\xi))\cdot J$;
        or the interior of $r(\xi)$ is disjoint from $\sigma'$ and
        $r(\xi)\cup\sigma'$ is connected, and
        $x=I\cdot(a(\xi)\otimes a'_o(r(\xi)\cup\sigma'))\cdot J$.
    \item 
      \label{typeO:FewerSigma}
      It has the form
      $x=I\cdot(a(\xi\setminus \sigma)\otimes a'_o(\xi))\cdot J$
      where $\sigma\subset \xi$
      or $x=I\cdot(a(\xi)\otimes a'_o(\xi\setminus \sigma'))\cdot J$
      where $\sigma'\subset \xi$. (Note that $\xi\setminus\sigma$
      or $\xi\setminus \sigma'$ can be disconnected in this case).
    \item 
      \label{typeO:YY}
      $x=I\cdot(a(\xi\cup\eta)\otimes a_o'(r(\xi)\cup r(\eta)))\cdot J$,
      where here:
      \begin{itemize}
      \item $\xi$ and $\eta$ are disjoint chords, or one is contained in the other,
      \item Neither $b_1$ nor $b_2$ appear in the boundary of $\xi$,
      \item $c_1$ appears in the boundary of $\xi$, and
      \item $c_2$ appears in the boundary of $\eta$ 
        with the opposite orientation.
      \end{itemize}
    \item
      \label{typeO:XX}
      The non-zero element $x$ has the form $x=I\cdot (a(\xi\cup\sigma)\otimes a_o'(\xi\setminus\sigma'))\cdot J$ where
      $\sigma'\subset \xi$ but $\sigma'$ is not contained in the interior of $\xi$,
      and $\xi\cap\sigma=\emptyset$;
      or
      $x=I\cdot (a(\xi\setminus \sigma)\otimes
      a_o'(\xi\cup\sigma'))\cdot J$, where the $\sigma\subset\xi$ but
      $\sigma$ is not
      contained in the interior of $\xi$, and $\xi\cap\sigma'=\emptyset$.
    \item 
      \label{typeO:dC1C2}
      $x$ has support $[c_2,c_1]\times [c_2,c_1]$ and
      exactly three moving strands.
    \item
      \label{typeO:Disconnected}
      $x$ factors as a product of 
      $I\cdot(a([c_2,c_1])\otimes a'_o([c_2,c_1]))\cdot I$
      and $I\cdot(a(\xi)\otimes a_o(\xi))\cdot J$,
      where $\xi$ is disjoint from $[c_2,c_1]$ or $\xi$
      is properly contained inside $[c_2,c_1]$.
      \end{enumerate}
    Near-chords for over-slides (satisfying Convention~\ref{label:conv}) are illustrated in Figure~\ref{fig:NearChords}.

      When Convention~\ref{label:conv} does not hold for the over-slide,
      as before we switch the roles of the two tensor factors in the definition of near-chords.
\end{definition}

\begin{figure}
    \begin{center}
      \input{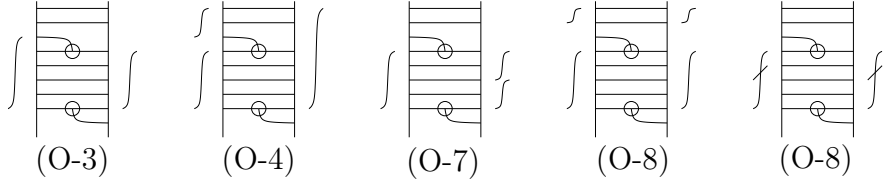}
    \end{center}
    \caption {{\bf Indeterminate near-chords.}
      \label{fig:Indeterminates} Examples of the different kinds of indeterminate near-chords.}
\end{figure}

The calculation of the arc-slide bimodule for an over-slide is not
quite as straightforward as for under-slides: we cannot say that
all near-chords appear in the differential, and indeed the bimodule is
determined uniquely only up to isomorphism.  The near-chords
which might or might not appear in the differential for a given
arc-slide bimodule are the following:

\begin{definition}
  \label{def:Indeterminate}
  \index{indeterminate near-chord}\index{near-chord!indeterminate}%
  A near-chord for an over-slide is called \emph{indeterminate} if it is
  of one of the following types:
  \begin{itemize}
    \item It is a near-chord of Type~(\ref{typeO:ExtraSigma}), and its
      restricted support is $[c_2,c_1]$.
    \item It is a near-chord of Type~(\ref{typeO:FewerSigma}), and
      the boundary of its restricted support contains $c_1$ or $c_2$, 
      and the restricted support contains the 
      interval $[c_2,c_1]$.
    \item It is a near-chord of Type~(\ref{typeO:dC1C2}). 
    \item It is a near-chord of Type~(\ref{typeO:Disconnected}).
  \end{itemize}
  These cases are illustrated in Figure~\ref{fig:Indeterminates}.
\end{definition}

A given indeterminate near-chord might or might not appear in the differential
of an arc-slide bimodule $N$. Exactly which ones appear are governed
by the following:

\begin{definition}
  For arc-slides satisfying Convention~\ref{label:conv},
  \label{def:BasicChoice}\index{basic choice}%
  a {\em basic choice} is a collection ${\mathcal B}$ of indeterminate
  near-chords of Type~(\ref{typeO:ExtraSigma}), satisfying the
  following condition: if $x$ and $x'$ are two distinct indeterminate
  near-chords of Type~(\ref{typeO:ExtraSigma}) with the same initial idempotent (i.e., there is some $I$
  with $I\cdot x=x$ and $I\cdot x'=x'$), then exactly one of $x$ or
  $x'$ is in ${\mathcal B}$.
 
  For arc-slides not satisfying Convention~\ref{label:conv},
  a basic choice is analogous, but with the terminal idempotent now playing the crucial role;
  i.e.,
  if $x$ and $x'$ are two distinct indeterminate
  near-chords of Type~(\ref{typeO:ExtraSigma}) with the same terminal idempotent (i.e., there is some $I$
  with $x=x\cdot I$ and $x'=x'\cdot I$), then exactly one of $x$ or
  $x'$ is in ${\mathcal B}$. 
\end{definition}

\begin{lemma}
  \label{lem:IsBasicChoice}
  If $N$ is a stable arc-slide bimodule, then the set of indeterminate
  near-chords of Type~(\ref{typeO:ExtraSigma}) which are contained in the
  differential on $N$ forms a basic choice,
  in the sense of Definition~\ref{def:BasicChoice}.
\end{lemma}

We defer the proof of Lemma~\ref{lem:IsBasicChoice} until
page~\pageref{page:IsBasicChoicepfpage}.

\begin{definition}
  \label{def:BasicChoiceOfN}
  \index{arc-slide bimodule!compatible with basic choice}%
  \index{basic choice!bimodule compatible with}%
  Let $N$ be a stable arc-slide bimodule. If ${\mathcal B}$ denotes
  the basic choice of indeterminate near-chords which appear as
  coefficients in the boundary operator for $N$ then we say that $N$ is
  {\em compatible with the basic choice~${\mathcal B}$};
  we also say that {\em ${\mathcal B}$ is the basic choice of $N$}.
\end{definition}

\begin{proposition}
  \label{prop:CalculateOver-Slide}
  If $N$ is a stable arc-slide bimodule which is compatible with a basic
  choice ${\mathcal B}$ then only near-chords can appear in the
  differential on $N$, and precisely which ones do appear are
  uniquely determined by the basic choice ${\mathcal B}$.  If
  $N_1$ and $N_2$ are two arc-slide bimodules which are compatible
  with basic choices ${\mathcal B}_1$ and~${\mathcal B}_2$ then there
  is an isomorphism between $N_1$ and~$N_2$.
\end{proposition}

Before proving Proposition~\ref{prop:CalculateOver-Slide} we
establish some preliminary results.

As in the case of under-slides, we study the elements of the
near-diagonal subalgebra of grading~$\ge -1$.  This time, there are
some elements of grading~$0$, which are responsible for the
indeterminacy. We give them a name:

\begin{definition}\label{def:dischord}
  \index{dischord}%
  A basic generator of the near-diagonal subalgebra with support
  $[c_2,c_1]\times[c_2,c_1]$ and exactly two moving strands (one in
  $\PMC$ and one in $\PMC'$) is called a \emph{dischord}.
\end{definition}

\begin{lemma}\label{lem:over-slide-grading}
  In the near-diagonal subalgebra of an over-slide $m\co\PMC\to\PMC'$,
  there are no elements of positive grading;
  the basic elements of
  grading~$0$ are the idempotents and the dischords
  (Definition~\ref{def:dischord}); and the basic elements of
  grading~$-1$ are the near-chords for over-slides.
\end{lemma}
\begin{proof}
  The proof is similar to the under-slide case (Lemma~\ref{lem:under-slide-grading}).
  Again, we assume Convention~\ref{label:conv}, reducing to this case via reflection if it is needed.
  As there, let $a$ be a basic generator in the near-diagonal subalgebra of
  grading $\ge -1$, and let $\Domain$ be the corresponding domain in the
  standard Heegaard diagram $\HD(m)$ for $m$.
  If all six multiplicities $n_{\sigma_+}, \dots, n_{\sigma'_-}$
  are~$0$ then $\Domain$ consists of a union of horizontal strips, and $a$
  is of Type (\ref{typeO:Generic}) or an idempotent.

  In general, the constraints are
  \begin{itemize}
  \item The multiplicity difference along any line is at most~$1$.
  \item The multiplicity differences from the idempotents are
    \begin{align}
      \label{eq:over-mult-l}
      n_{\sigma} - n_{\sigma_-} + n_{\sigma'_+} - n_{\sigma'_-}
        &= \text{$-1$, $0$, or $+1$}\\
      \label{eq:over-mult-r}
      n_{\sigma_+} - n_{\sigma_-} + n_{\sigma'_+} - n_{\sigma'}
        &= \text{$-1$, $0$, or $+1$}
    \end{align}
    where the right hand side is determined by what happens to the occupancy of
    the $C$ idempotent on the left in \eqref{eq:over-mult-l} and on the right
    in \eqref{eq:over-mult-r}.
  \end{itemize}

  \textbf{Case $\lsub{C}X \to \lsub{C}X$.}
  In this case we have
  \begin{align*}
    n_{\sigma} - n_{\sigma_-} + n_{\sigma'_+} - n_{\sigma'_-} &= 0\\
    n_{\sigma_+} - n_{\sigma_-} = n_{\sigma'_+} - n_{\sigma'} &= 0.
  \end{align*}
  (The second set of equations come from the fact that the strand~$C$ is not
  occupied on the right in either the initial or final idempotent, so
  there can be no strand starting or ending there.)
  According to Proposition~\ref{prop:near-diag-grading}, the
  correction to the grading is given by $c(I, \supp(a)) + c(J, \supp(a))
  = \frac{1}{2}(-n_{\sigma'} + n_{\sigma'_-})$.

  The linear equations tell us that the
  multiplicities are of the following forms:

  \begin{center}
  \begin{tabular}{McMcMcMcMcMcMc}
    \toprule
    \text{Corr}&n_{\sigma_+}&n_{\sigma}&n_{\sigma_-}&n_{\sigma'_+}&n_{\sigma'}&n_{\sigma'_-}\\
    \midrule
    \epsilon/2 & m & m+\epsilon & m & l & l & l+\epsilon\\
    \bottomrule
  \end{tabular}
  \end{center}
  Here, $\epsilon\in \{-1,0,1\}$ (as the difference in multiplicities
  is at most one).

  If $\epsilon=-1$, we would need to have $M(\grb(a)) = -\OneHalf$,
  which is not possible with the given multiplicities.

  If $\epsilon=0$, we have complete horizontal
  strips, giving near-chords of
  Type (\ref{typeO:Generic}) or, in the case of no strips at all, an idempotent.

  If $\epsilon=+1$, the left side of~$\Domain$ will need at
  least two intervals to cover it, leaving only one interval for the
  right.  This implies that $m=l=0$, which gives a domain of Type (\ref{typeO:XX}).

  In summary, the possibilities are:

  \begin{center}
  \begin{tabular}{cMcMcMcMcMcMcMcMc}
    \toprule
    Type&\text{Corr}&\text{Grading}&n_{\sigma_+}&n_{\sigma}&n_{\sigma_-}&n_{\sigma'_+}&n_{\sigma'}&n_{\sigma'_-}\\
    \midrule
    \ref{typeO:Generic} & 0         & -1 & 1 & 1 & 1 & 0 & 0 & 0\\
    \ref{typeO:Generic} & 0         & -1 & 0 & 0 & 0 & 1 & 1 & 1\\
    \ref{typeO:Generic} & 0         & -1 & 1 & 1 & 1 & 1 & 1 & 1\\
    \ref{typeO:XX} & +\OneHalf & -1 & 0 & 1 & 0 & 0 & 0 & 1\\
    \bottomrule
  \end{tabular}
  \end{center}

  \textbf{Case $\lsub{C}X \to X_C$.}
  In this case we have
  \begin{align*}
    n_{\sigma} - n_{\sigma_-} + n_{\sigma'_+} - n_{\sigma'_-} &= 1\\
    n_{\sigma_+} - n_{\sigma_-} + n_{\sigma'_+} - n_{\sigma'} &= 1.
  \end{align*}
  The
  correction to the grading is given by 
  \[c(I, \supp(a)) + c(J, \supp(a))
  = \OneQuart(n_{\sigma_+} - n_{\sigma} - n_{\sigma'} +
  n_{\sigma'_-}) =0.\]
  The linear equations tell us that the multiplicities are given by
  \begin{center}
  \begin{tabular}{McMcMcMcMcMcMc}
    \toprule
    \text{Corr}&n_{\sigma_+}&n_{\sigma}&n_{\sigma_-}&n_{\sigma'_+}&n_{\sigma'}&n_{\sigma'_-}\\
    \midrule
    0&m+\delta & m+\epsilon & m & l & l+\delta-1 & l+\epsilon-1 \\
    \bottomrule
  \end{tabular}
  \end{center}
  where $\epsilon,\delta\in\{0,1\}$.

  The only solutions to
  these equations which yield a connected domain on both sides
  (as required by the gradings)
  are the following:

  \begin{center}
  \begin{tabular}{cMcMcMcMcMcMcMcMc}
    \toprule
    Type&\text{Corr}&\text{Grading}&n_{\sigma_+}&n_{\sigma}&n_{\sigma_-}&n_{\sigma'_+}&n_{\sigma'}&n_{\sigma'_-}\\
    \midrule
    \ref{typeO:Generic} & 0 & -1 & 1 & 1 & 0 & 0 & 0 & 0\\ 
    \ref{typeO:Generic} & 0 & -1 & 1 & 1 & 1 & 1 & 0 & 0\\ 
    \ref{typeO:Generic} & 0 & -1 & 0 & 0 & 0 & 1 & 0 & 0\\ 
    \bottomrule
  \end{tabular}
  \end{center}

  \textbf{Case $\lsub{C}X \to Y$.}
  In this case we have
  \begin{align*}
    n_{\sigma} - n_{\sigma_-} + n_{\sigma'_+} - n_{\sigma'_-} &= 0\\
    n_{\sigma_+} - n_{\sigma_-} + n_{\sigma'_+} - n_{\sigma'} &= 1.
  \end{align*}
  The
  correction to the grading is given by 
  \[c(I, \supp(a)) + c(J, \supp(a))
  = \OneQuart(-n_{\sigma} + n_{\sigma_-} + n_{\sigma'_+} - 2
    n_{\sigma'} + n_{\sigma'_-}) = \OneHalf(n_{\sigma'_+} - n_{\sigma'}).
  \]
  The linear equations tell us that the multiplicities are

  \begin{center}
  \begin{tabular}{McMcMcMcMcMcMc}
    \toprule
    \text{Corr}&n_{\sigma_+}&n_{\sigma}&n_{\sigma_-}&n_{\sigma'_+}&n_{\sigma'}&n_{\sigma'_-}\\
    \midrule
    (1-\delta)/2 & m+\delta & m+\epsilon & m & l & l+\delta-1 & l+\epsilon \\
    \bottomrule
  \end{tabular}
  \end{center}
  with $\delta \in \{0,1\}$ and $\epsilon \in \{-1,0,1\}$.

  The solutions to
  these equations which can have grading $\ge -1$ are

  \begin{center}
  \begin{tabular}{cMcMcMcMcMcMcMcMc}
    \toprule
    Type&\text{Corr}&\text{Grading}&n_{\sigma_+}&n_{\sigma}&n_{\sigma_-}&n_{\sigma'_+}&n_{\sigma'}&n_{\sigma'_-}\\
    \midrule
    \ref{typeO:FewerSigma} & 0        & -1 & 1 & 0 & 0 & 0 & 0 & 0 \\ 
    \ref{typeO:FewerSigma} & \OneHalf & -1 & 0 & 0 & 0 & 1 & 0 & 1 \\ 
    \ref{typeO:FewerSigma} & \OneHalf & -1 & 1 & 1 & 1 & 1 & 0 & 1 \\ 
    \ref{typeO:FewerSigma} & \OneHalf & -1 & 1 & 0 & 1 & 1 & 0 & 0 \\ 
    \bottomrule
  \end{tabular}
  \end{center}
  (The last case is an indeterminate near-chord.)

  \textbf{Case $X_C \to \lsub{C}X$.}  
  This is related to the case $\lsub{C}X \to X_C$ by rotating the
  diagram $180^\circ$.  Again, the solutions are all of Type (\ref{typeO:Generic}).

  \textbf{Case $X_C \to X_C$.}
  This is related to the case $\lsub{C}X \to \lsub{C}X$ by rotating
  the diagram $180^\circ$.
  The solutions are idempotents and near-chords of Types (\ref{typeO:Generic}) and (\ref{typeO:XX}).

  \textbf{Case $X_C \to Y$.}
  This is related to the case $\lsub{C}X \to Y$ by rotating the
  diagram $180^\circ$.  The
  solutions are of Types (\ref{typeO:Generic}) and (\ref{typeO:FewerSigma}).

  \textbf{Case $Y \to \lsub{C}X$.}
  In this case we have
  \begin{align*}
    n_{\sigma} - n_{\sigma_-} + n_{\sigma'_+} - n_{\sigma'_-} &= 0\\
    n_{\sigma_+} - n_{\sigma_-} + n_{\sigma'_+} - n_{\sigma'} &= -1.
  \end{align*}
  The
  correction to the grading is given by
    \[c(I, \supp(a)) + c(J, \supp(a))
  = \OneQuart(-n_{\sigma} + n_{\sigma_-} + n_{\sigma'_+} - 2
    n_{\sigma'} + n_{\sigma'_-}) = \OneHalf(n_{\sigma'_+} - n_{\sigma'}).
  \]
  The linear equations tell us that the multiplicities are

  \begin{center}
  \begin{tabular}{McMcMcMcMcMcMc}
    \toprule
    \text{Corr}&n_{\sigma_+}&n_{\sigma}&n_{\sigma_-}&n_{\sigma'_+}&n_{\sigma'}&n_{\sigma'_-}\\
    \midrule
    -\delta/2 & m+\delta-1 & m+\epsilon & m & l & l+\delta & l+\epsilon \\
    \bottomrule
  \end{tabular}
  \end{center}
  with $\delta \in \{0,1\}$ and $\epsilon \in \{-1,0,1\}$.

  The solutions which can have grading $\ge -1$ are

  \begin{center}
  \begin{tabular}{cMcMcMcMcMcMcMcMc}
    \toprule
    Type&\text{Corr}&\text{Grading}&n_{\sigma_+}&n_{\sigma}&n_{\sigma_-}&n_{\sigma'_+}&n_{\sigma'}&n_{\sigma'_-}\\
    \midrule
    \ref{typeO:Sigma} & -\OneHalf & -1 & 0 & 0 & 0 & 0 & 1 & 0 \\ 
    \ref{typeO:ExtraSigma} & 0         & -1 & 0 & 1 & 1 & 0 & 0 & 0 \\ 
    \ref{typeO:ExtraSigma} & 0         & -1 & 0 & 1 & 1 & 1 & 1 & 1 \\ 
    \ref{typeO:ExtraSigma} & 0         & -1 & 0 & 0 & 1 & 1 & 1 & 0 \\ 
    \bottomrule
  \end{tabular}
  \end{center}
  (The last case is an indeterminate near-chord.)

  \textbf{Case $Y \to X_C$.}
  This is related to the case $Y \to \lsub{C}X$ by rotating the
  diagram $180^\circ$.  The
  solutions are of Type (\ref{typeO:Sigma}) and (\ref{typeO:ExtraSigma}).

  \textbf{Case $Y \to Y$.}
  In this case we have
  \begin{align*}
    n_{\sigma} - n_{\sigma_-} + n_{\sigma'_+} - n_{\sigma'_-} &= 0\\
    n_{\sigma_+} - n_{\sigma_-} + n_{\sigma'_+} - n_{\sigma'} &= 0\\
    n_{\sigma_+} - n_{\sigma} = n_{\sigma'} - n_{\sigma'_-} &= 0.
  \end{align*}
  (As in the under-slide case, the last equations come from the fact
  that the $B$ strand is not
  occupied in either idempotent on either side.)
  The
  correction to the grading is given by 
  \[c(I, \supp(a)) + c(J, \supp(a))
  = \OneHalf(-n_{\sigma} + n_{\sigma_-}
    + n_{\sigma'_+} - n_{\sigma'}) = -n_{\sigma} + n_{\sigma_-}.
  \]
  The linear equations tell us that the multiplicities are

  \begin{center}
  \begin{tabular}{McMcMcMcMcMcMc}
    \toprule
    \text{Corr}&n_{\sigma_+}&n_{\sigma}&n_{\sigma_-}&n_{\sigma'_+}&n_{\sigma'}&n_{\sigma'_-}\\
    \midrule
    -\epsilon & m+\epsilon & m+\epsilon & m & l & l+\epsilon & l+\epsilon \\
    \bottomrule
  \end{tabular}
  \end{center}
  with $\epsilon \in \{-1,0,1\}$.

  The solutions which can have grading $\ge -1$ are

  \begin{center}
  \begin{tabular}{cMcMcMcMcMcMcMcMc}
    \toprule
    Type&\text{Corr}&\text{Grading}&n_{\sigma_+}&n_{\sigma}&n_{\sigma_-}&n_{\sigma'_+}&n_{\sigma'}&n_{\sigma'_-}\\
    \midrule
    \ref{typeO:Generic} & 0 & -1  & 1 & 1 & 1 & 0 & 0 & 0 \\ 
    \ref{typeO:Generic} & 0  & -1 & 0 & 0 & 0 & 1 & 1 & 1 \\ 
    \ref{typeO:Generic} & 0  & -1 & 1 & 1 & 1 & 1 & 1 & 1 \\ 
    Dischord & -1 & 0 & 0 & 0 & 1 & 1 & 0 & 0\\
    \ref{typeO:YY}, \ref{typeO:dC1C2}, \ref{typeO:Disconnected} & -1 & -1 & 0 & 0 & 1 & 1 & 0 & 0 \\
    \ref{typeO:YY} & -1 & -1 & 1 & 1 & 2 & 1 & 0 & 0 \\ 
    \ref{typeO:YY} & -1 & -1 & 0 & 0 & 1 & 2 & 1 & 1 \\ 
    \bottomrule
  \end{tabular}
  \end{center}

Here there is a grading~$0$ solution, the dischord (with support
$[c_2,c_1] \times [c_2,c_1]$), which can be
modified without changing the multiplicities
near $\sigma$ or $\sigma'$ in a variety of ways: introducing a break
in the support (Type (\ref{typeO:YY})), introducing a break in the chord on one
side without changing the support (Type (\ref{typeO:dC1C2})), or adding a new chord
somewhere else, either overlapping with the existing support or not
(Type (\ref{typeO:Disconnected})).

This is the end of the case analysis.  It is straight-forward (if
somewhat tedious) to verify that every near-chord for over-slides
appears in the list of grading $-1$ elements. The idempotents and
dischords were exactly the grading $0$ elements which occurred in the
case analysis, and no positive grading elements occurred in the case
analysis.  This concludes the proof.
\end{proof}

\begin{lemma}
  \label{lem:ONoOthers}
  If $N$ is an arc-slide bimodule then the only algebra elements which
  appear in the differential are near-chords.
\end{lemma}

\begin{proof}
  This is an immediate consequence of the definition of the grading on
  the coefficient algebra and Lemma~\ref{lem:over-slide-grading}.
\end{proof}

\begin{lemma}
  \label{lem:OGenericExists}
  Let $N$  be a stable arc-slide bimodule for an over-slide. Then the differential
  on $N$ contains all near-chords of Type~(\ref{typeO:Generic}).
\end{lemma}

\begin{proof}
  The proof follows along the lines of Lemma~\ref{lem:UGenericExists}.

  We define the {\em outside region} and {\em big enough} as
  in the proof of that lemma; and by stability, we restrict attention
  to the case where the pointed matched circle is big enough.
  By stabilizing, we can further assume that any special length $3$ chord
  which is adjacent to one of $b_1$ or $b_1'$ is not adjacent to the
  basepoint. 

  Now the arguments proving Lemma~\ref{lem:UGenericExists} give the following
  sequence of claims:

  {\bf Claim 1:} {\em If $x$ is a near-chord of
    Type~(\ref{typeO:Generic}) and the support of $x$ is disjoint from
    $b_1$ and $b_1'$ then $x$ appears in the differential.}

  The case of special length $3$ chords are handled as they were in
  the proof of Proposition~\ref{prop:DDidUnique}. (There is no analogue of very special length $3$ chords in the over-slide case.)
 
  {\bf Claim 2}: {\em
    Suppose $x$ is any near-chord of Type~(\ref{typeO:Generic}) with the following properties:
    \begin{itemize}
    \item the support of $x$ contains exactly one of $b_1$ or $b_1'$, and
    \item the position at $b_1$ or $b_1'$ (whichever is contained in the support of $x$) is occupied
      in both the initial and terminal idempotents for $x$ (on the $\PMC$ or the  $\PMC'$ side, respectively).
    \end{itemize}
    Then $x$ appears in the differential.} 
  
  {\bf Claim 3:} {\em If $x$ is a Type~(\ref{typeO:Generic}) chord with
    restricted support of length $1$, then $x$ appears in the
    differential.}

  This is the analogue of Claim 5 from the proof of
  Lemma~\ref{lem:UGenericExists}, and its proof is similar
  (e.g., establishing first the case where $dx\neq 0$).
  
  We turn next to the inductive proof that longer near-chords of Type~(\ref{typeO:Generic}) appear.
  We call a near-chord $x$ of of Type~(\ref{typeO:Generic}) {\em simplifiable} if 
  $dx$ contains terms of the form $y\cdot t$ with the following
  properties:
  \begin{itemize}
  \item Each of $y$ and $t$ is of Type~(\ref{typeO:Generic}).
  \item The product $y\cdot t$ has no alternative factorization.
  \item The product $y\cdot t$ does not appear in the differential
    of any other basic generator of the near-diagonal sub-algebra.
  \end{itemize}
  
  The proof has one extra complication for degenerate over-slides (Definition~\ref{def:degen-slide}), so we separate out that case.
  
  {\bf Non-degenerate case.}  In this case, all near-chords of
  Type~(\ref{typeO:Generic}) either have restricted support of length
  $1$ (so they appear in the differential by Claim 3); they
  are special length $3$ chords (which can be shown to appear as they were in
  the proof of Proposition~\ref{prop:DDidUnique}); or they are
  simplifiable, in which case we can apply induction on the length of
  the restricted support of the near-chord, to conclude that all
  near-chords of Type~(\ref{typeO:Generic}) appear in the differential.

  \textbf{Degenerate case.} There are two sub-cases, depending on whether the unique position between $c_1$ and $c_2$ is $b_2$.

  {\bf Sub-case: $b_2$ lies between $c_1$ and $c_2$.}  In this case,
  there is one more distinguished type of near chords: chords $x$ of
  Type~(\ref{typeO:Generic}) such that:
  \begin{itemize}
  \item $x$ has length $5$ on both $\PMC$ and $\PMC'$,
  \item $x$ contains all of $c_1$, $c_2$, and $b_2$ in the interior of
    its support (and hence has restricted length $3$) and
  \item the initial idempotent of $x$ has type $Y$ (and hence so does
    the terminal idempotent of $x$).
  \end{itemize}
  We call such a near-chord {\em extra special}.  Under the present
  assumptions, there are three kinds of non-simplifiable chords of
  Type~(\ref{typeO:Generic}): those which have length $1$ restricted
  support (which appear in the differential by Claim 3), those which
  are special length $3$ chords (which appear in the differential by
  the usual arguments), and those which are extra special. (Chords
  like extra special ones but with idempotent of type $X$ do not need
  a special argument, but see the first line of
  Figure~\ref{fig:extra-special}.)

  \begin{figure}
    \centering
    \includegraphics[scale=.5]{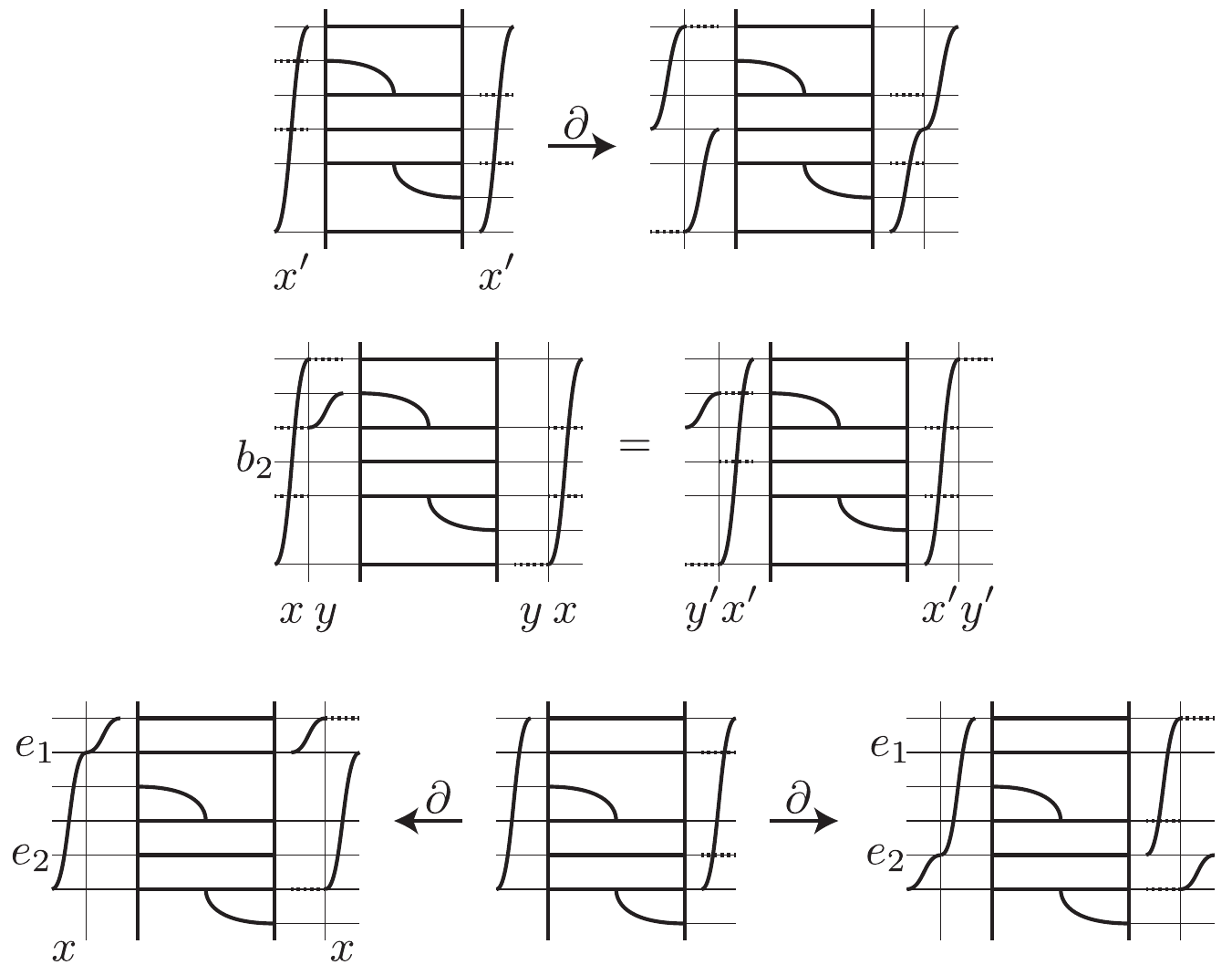}
    \caption{\textbf{Existence of extra-special chords $x$.} The chord
      in the top line is not extra-special: the idempotent has type
      $X$. The second and third lines demonstrate the existence of
      extra-special chords in the cases that the position between $c_1$
      and $c_2$ is or is not $b_2$, respectively.}
    \label{fig:extra-special}
  \end{figure}

  We must argue that every extra-special near-chord $x$ appears in the differential. Note that there is a chord $y$
  of Type~(\ref{typeO:Sigma}) with $x\cdot y\neq 0$. Now, $x\cdot y$ has
  a unique alternate factorization, which has the form $x \cdot y =
  y'\cdot x'$ (with $\supp(x)=\supp(x')$ and $\supp(y)=\supp(y')$, so
  that $x'$ is of Type~(\ref{typeO:Generic}) and $y'$ of
  Type~(\ref{typeO:Sigma})). Moreover, the term in $dx'$ corresponding to the crossing at $b_2$ factors as a product of two (not extra special) near-chords of Type~(\ref{typeO:Generic}), so $x'$ appears in the differential. The near-chord $y'$ appears in the differential by hypothesis.
  Moreover, $x\cdot y$ does not appear in
  the differential of any other algebra element, so it follows that $x$ appears in the differential. See the second line of Figure~\ref{fig:extra-special}.

  Now induction on the length of the restricted support establishes the proposition in this case.

  {\bf Sub-case: some $e_2\neq b_2$ lies between $c_1$ and $c_2$.}  If
  the position $e_1$ matched with $e_2$ is adjacent to neither $b_1$
  nor $b_1'$ then the same argument as in the non-degenerate case
  applies.  Special attention is needed when the position $e_1$
  matched with $e_2$ is adjacent to either $b_1$ or $b_1'$.  Assume
  for definiteness that $e_1$ is above and adjacent to $b_1$; the
  other cases are similar.  In this case, a near chord $x$ of
  Type~(\ref{typeO:Generic}) is called {\em extra special} if its
  support is $[e_1,c_2]$.  For an extra special near chord $x$,
  $dx=0$.

  Under the current assumptions, any near-chord $x$ of
  Type~(\ref{typeO:Generic}) either has restricted support of length
  $1$ (in which case $x$ appears in the differential by Claim 3), is a
  special length $3$ chord (in which case $x$ appears in the
  differential by the argument from
  Proposition~\ref{prop:DDidUnique}), is an extra-special chord,
  or is simplifiable.  Extra special chords appear in the
  differential by the same argument used for special length $3$
  chords. (See the last line in Figure~\ref{fig:extra-special}.)
  So, again, induction on the length of the restricted support establishes the proposition.
\end{proof}

By hypothesis, the differential on $N$ contains all near-chords of
Type~(\ref{typeO:Sigma}).

\begin{lemma}
  \label{lem:OExtraFewerExist}
  Let $N$  be a stable arc-slide bimodule for an over-slide. Then the differential
  on $N$ contains all non-indeterminate near-chords of Types~(\ref{typeO:ExtraSigma})
  and~(\ref{typeO:FewerSigma}).
\end{lemma}

\begin{proof}
  We describe the proof; but indeed, most of it is encapsulated in
  Figure~\ref{fig:ExtraFewerExist}, and it runs parallel
  to the proof of the corresponding fact 
  in the under-slide case (Lemma~\ref{lem:UFewerExist}).
  \begin{figure}
    \begin{center}
      \includegraphics[scale=.5]{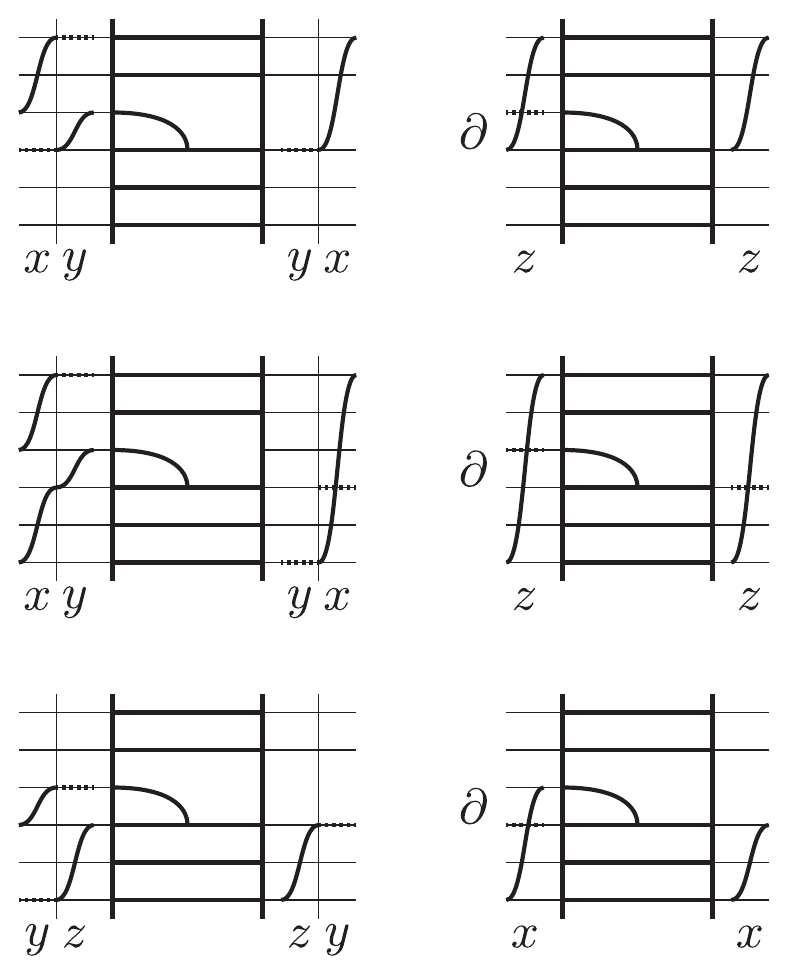}
    \end{center}
    \caption {{\bf The three cases of Lemma~\ref{lem:OExtraFewerExist}.}
      \label{fig:ExtraFewerExist}}
  \end{figure}

  For any near-chord $x$ of Type~(\ref{typeO:FewerSigma}),
  we can find a near-chord $y$ of Type~(\ref{typeO:Sigma}) with the
  property that $x\cdot y$ appears in the differential of a near-chord $z$
  of Type~(\ref{typeO:Generic}), as in the top two rows of
  Figure~\ref{fig:ExtraFewerExist}. 
  This product $x\cdot y$
  has no alternative factorization as a product of two near-chords
  with non-trivial support. (This is where we use that $x$ is not
  indeterminate.) More specifically, the case where $x$ has only two
  moving chords is obvious. In the case
  where $x$ has three moving chords, $x\cdot y$ has three moving
  chords as well.
  There is exactly one other
  factorization of $x\cdot y$, into 
  a near-chord $x'$ of Type~(\ref{typeO:FewerSigma}) and
  another $y'$ of Type~(\ref{typeO:ExtraSigma}) whose product
  $x'\cdot y'$ has the same three moving chords as $x\cdot y$.  (See
  Figure~\ref{fig:NotAFactorization}.) However, a closer look at the
  idempotents shows that $x\cdot y$ does not equal $x'\cdot y'$.
  Specifically, suppose, without loss of generality, that two of the
  moving strands of $x$ are on the $\PMC$-side.
  Then, since $x$ is not indeterminate, the initial idempotent of $x$ contains $C$ on the
  $\PMC'$-side ($x$ has type $X_C \to Y$).
  On the other hand, the initial idempotent of $x'$
  contains $C$ on the $\PMC$ side ($x'$ has type $\lsub{C}X \to Y$). (Again, see
  Figure~\ref{fig:NotAFactorization}.)

  Since, by Lemma~\ref{lem:OGenericExists}, $z$ appears in the
  differential, we conclude that $x$ must, as well.
  \begin{figure}
  \begin{center}
      \includegraphics[scale=.5]{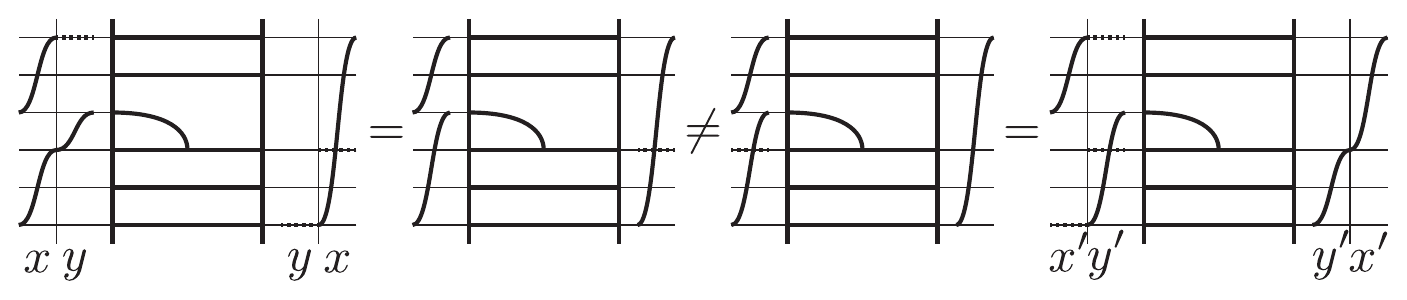}
    \end{center}
    \caption {{\bf Not an alternative factorization.}
      \label{fig:NotAFactorization}
      One of the cases in Figure~\ref{fig:ExtraFewerExist} might
      appear to have an alternate factorization; however, a more
      careful look at idempotents (as indicated) shows that this
      alternative factorization does not exist.}
  \end{figure}

  To prove the other half of the lemma, consider a near-chord $x$ of
  Type~(\ref{typeO:ExtraSigma}) which is
  not indeterminate. Then, the differential of $x$ contains an element
  which factors as $y\cdot z$, where $y$ is of Type~(\ref{typeO:Sigma})
  and $z$ is of Type~(\ref{typeO:Generic}), as in the
  third line of Figure~\ref{fig:ExtraFewerExist}. (This is where we
  use the fact that $x$ is not indeterminate.) That factorization
  is
  unique in the near-diagonal subalgebra, and $y\cdot z$ does not appear in the differential of
  any other algebra element. According to
  Lemma~\ref{lem:OGenericExists}, $y\cdot z$ appears in $\bdy^2$, so
  $x$ must appear in $\bdy$.
\end{proof}

\begin{lemma}
  \label{lem:OYYexist}
  Let $N$  be a stable arc-slide bimodule for an over-slide.
  Then the differential on $N$ contains all
  near-chords of Type~(\ref{typeO:YY}) and Type~(\ref{typeO:XX}).
\end{lemma}

\begin{proof}
  The proof is illustrated in  Figure~\ref{fig:YYexist}.

  Let $x$ be a near-chord of Type~(\ref{typeO:YY}). Post-multiply $x$
  by a near-chord $y$ of Type~(\ref{typeO:Sigma}). The resulting algebra
  element $x\cdot y$ does not appear in the differential of any other algebra
  element. Moreover, $x\cdot y$ has a unique alternative factorization
  as a  product of two elements of the near-diagonal subalgebra,
  $x'\cdot y'$, where $x'$ has Type~(\ref{typeO:ExtraSigma}),
  $y'$ has Type~(\ref{typeO:Generic}), and the boundary of $y'$ meets $C$.
  (In the case where the chords $\xi$ and $\eta$ are neither disjoint nor nested,
  this alternative factorization does not make sense. In fact, if $x$ is
  analogous to a Type~(\ref{typeO:YY}) near-chord, except that
  the chords $\xi$ and $\eta$ are neither disjoint nor nested, then 
  $x$ already admits a factorization into near-chords.)
  This is illustrated in the first row of Figure~\ref{fig:YYexist}.
  \begin{figure}
    \begin{center}
      \includegraphics[scale=.5]{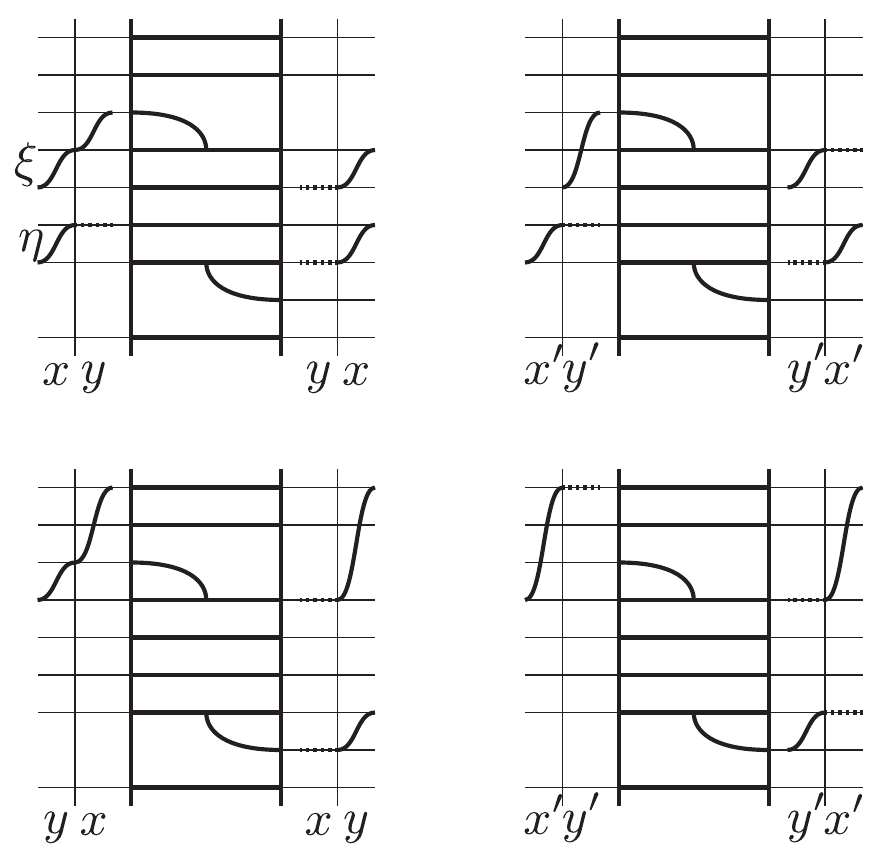}
    \end{center}
    \caption {{\bf Proof of Lemma~\ref{lem:OYYexist}.} The top line is
        for near-chords of Type~(\ref{typeO:YY}). 
        The bottom line is for near-chords of Type~(\ref{typeO:XX}). 
      }
      \label{fig:YYexist}
  \end{figure}

  If $x$ is of Type~(\ref{typeO:XX}), we pre-multiply $x$ by a
  near-chord $y$ of Type~(\ref{typeO:Sigma}).  The resulting algebra
  element $y\cdot x$ has an alternative factorization as $x'\cdot y'$
  where $y'$ is of Type~(\ref{typeO:Sigma}) and $x'$ is of
  Type~(\ref{typeO:Generic}).
  See the second row of Figure~\ref{fig:YYexist}.
\end{proof}

Next, we turn our attention to the indeterminate near-chords.

\begin{lemma}
  \label{lem:OIndeterminateFewerSigma}
  Let $N$ be a stable arc-slide bimodule for an over-slide.
  Which indeterminate near-chords of
  Type~(\ref{typeO:FewerSigma}) appear in the differential is
  uniquely determined by which indeterminate
  near-chords of Type~(\ref{typeO:ExtraSigma}) appear.
\end{lemma}

\begin{proof}
  Let $x$ be an indeterminate near-chord of
  Type~(\ref{typeO:FewerSigma}).  We can find a near-chord $y$ of
  Type~(\ref{typeO:Sigma}) with the
  property that $x\cdot y$ appears in the differential of a
  near-chord $z$ of Type~(\ref{typeO:Generic}).  The term $x\cdot y$ has a
  unique alternative factorization as $x'\cdot y'$, where $x'$ is of
  Type~(\ref{typeO:FewerSigma}) (but $x'$ is not indeterminate), and
  $y'$ is an indeterminate near-chord of Type~(\ref{typeO:ExtraSigma}).
  Since $dz$ appears in $\partial^2$, it follows that the term
  $y'$ appears in the differential if and only if $x$ does not. 
  See Figure~\ref{fig:IndeterminateFewerSigma} for an illustration.
  \begin{figure}
    \begin{center}
      \includegraphics[scale=.5]{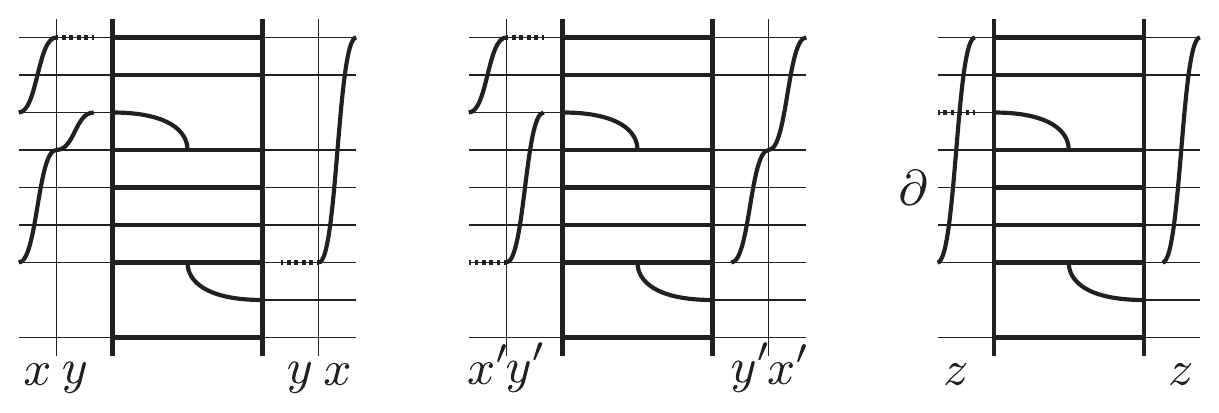}
    \end{center}
    \caption {{\bf The three types of terms in
        Lemma~\ref{lem:OIndeterminateFewerSigma}.} 
      \label{fig:IndeterminateFewerSigma}}
  \end{figure}
\end{proof}

Lemma~\ref{lem:OIndeterminateFewerSigma}, or rather its proof, can be used to
establish Lemma~\ref{lem:IsBasicChoice}.
To this end, note that an indeterminate near-chord of
Type~(\ref{typeO:ExtraSigma}) has the form $\dischord\cdot (\sigma\otimes
1)$ or $\dischord\cdot (1\otimes \sigma')$, where $\dischord$ is a dischord (in
the sense of Definition~\ref{def:dischord}). 

\begin{proof}[Proof of Lemma~\ref{lem:IsBasicChoice}.]
  \label{page:IsBasicChoicepfpage}
  Our aim is to show the following: for each dischord $\dischord$, the
  indeterminate near-chord $y=\dischord\cdot (\sigma\otimes 1)$ appears in the
  differential of a stable arc-slide bimodule if and only if $y'=\dischord\cdot
  (1\otimes \sigma')$ does not.
  (Again, we are using Convention~\ref{label:conv}, bearing in mind Remark~\ref{rem:ConventionAgain}.)

  After stabilizing, we work in a portion of the algebra where the number of
  occupied strands exceeds the number of positions between $c_1$ and $c_2$.
  We also stabilize so that there is some position above $b_1$.

  Consider first the case where the initial idempotent $I$ of $y$
  (i.e., the basic idempotent $I$ so that $I\cdot y=y$) contains 
  some occupied position in $\PMC$ which is above $b_1$.
  In this case, the proof of
  Lemma~\ref{lem:OIndeterminateFewerSigma} gives an indeterminate
  near-chord $x$ of Type~(\ref{typeO:FewerSigma}) which appears in
  the differential precisely if $y$ does not. (The initial idempotent
  of $y$ differs from the initial idempotent of $x$ in only one
  position.)  Now, $x\cdot (1\otimes \sigma')$ does not appear in the
  differential of any element, and it has a unique alternative
  factorization as $x'\cdot y'$, where $x'$ is a (not indeterminate) near-chord of
  Type~(\ref{typeO:FewerSigma}) which appears in the differential,
  by Lemma~\ref{lem:OExtraFewerExist}.  It follows that $y'$ appears
  in the differential if and only if $x$ does, which in turn appears
  in the differential if and only if $y$ does not. See
  Figure~\ref{fig:IsBasicChoice} for an illustration.
  \begin{figure}
    \begin{center}
      \includegraphics[scale=.5]{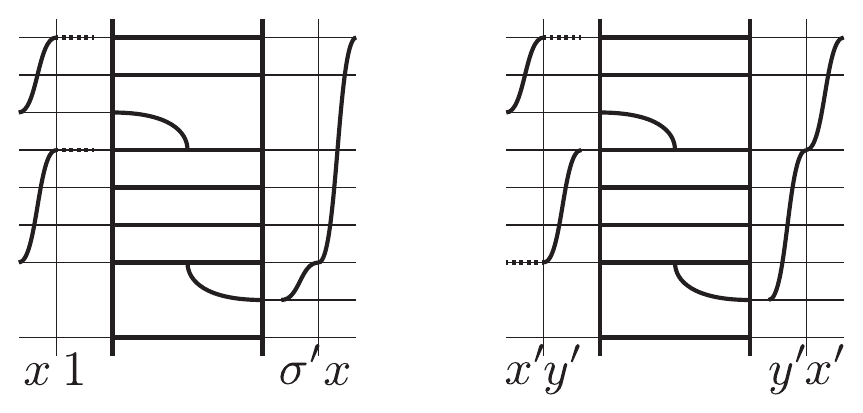}
    \end{center}
    \caption {{\bf Proof of Lemma~\ref{lem:IsBasicChoice}.}
      \label{fig:IsBasicChoice} The existence of the term 
      on the left in $\partial^2$ 
      is determined by the appearance in the differential of a
      corresponding near-chord $y$ of Type~(\ref{typeO:ExtraSigma})
      according to Lemma~\ref{lem:OIndeterminateFewerSigma}; the
      second factorization ($x'\cdot y'$) is determined by the existence of a
      different near-chord $y'$ of Type~(\ref{typeO:ExtraSigma}).}
  \end{figure}

  Suppose instead that in $I$ there is no occupied position above $b_1$. By our hypothesis
  on the total number of occupied positions, this forces there to be
  some occupied position $d_1$ below $c_2$.
  Thus, we can find a near-chord $x$ of Type~(\ref{typeO:Generic}) (connecting $d_1$ to some
  position above $b_1$) so that:
  \begin{itemize}
    \item $y\cdot x$ is non-zero, and has a unique alternate factorization as $x_2\cdot y_2$,
      where $x$ and $x_2$ have the same support, and $y$ and $y_2$ have the same support.
    \item $y'\cdot x$ is non-zero, and has a unique alternate factorization as $x_2\cdot y_2'$.
    \item $y\cdot x$ and $y'\cdot x$ do not appear in the differential of any other basic algebra element.
    \item the initial idempotent of $y_2$ (which is the same as the initial idempotent of $y_2'$) 
      contains some occupied
      position in $\PMC$ which is above $b_1$.
  \end{itemize}
  In particular, the pair $y_2$, $y'_2$ fit into the case of the
  lemma that we already established, so $y_2$ is contained in the differential if and only if $y_2'$
  is not. By $\partial^2=0$ we see that $y$ is contained in the differential if and only if $y_2$ is,
  and $y'$ is contained in the differential if and only if $y_2'$ is.
  It follows that $y$ is contained in the differential if and only if $y'$ is not, as claimed.
\end{proof}

Using the terminology of Definition~\ref{def:BasicChoiceOfN},
Lemma~\ref{lem:OIndeterminateFewerSigma} says that
that for a stable arc-slide bimodule, its basic choice ${\mathcal B}$
uniquely specifies which indeterminate near-chords of
Type~(\ref{typeO:FewerSigma}) appear in the differential. In the same vein, we have:

\begin{lemma}
  \label{lem:dC1C2}
  If $N$ is a stable arc-slide bimodule then
  its basic choice ${\mathcal B}$ uniquely specifies
  which indeterminate near-chords of Type~(\ref{typeO:dC1C2})
  appear in the differential on $N$.
\end{lemma}
\begin{proof}
  We consider terms in $\partial^2$ which have support, say,
  $[c_2,b_1]\times [c_2,c_1]$. Let $\alpha$ denote the
  sum of all terms of Type~(\ref{typeO:dC1C2}) which appear in the differential.
  The terms in $\partial^2$ coming
  from near-chords are of the following types:
  \begin{itemize}
  \item Terms of Type~(\ref{typeO:ExtraSigma}) times
    terms of Type~(\ref{typeO:Generic}).
  \item Terms of Type~(\ref{typeO:dC1C2})
    times terms of Type~(\ref{typeO:Sigma}).
  \item Differentials of indeterminate near-chords of Type~(\ref{typeO:ExtraSigma}).
  \end{itemize}
  See Figure~\ref{fig:dC1C2} for an illustration.
  \begin{figure}
    \begin{center}
      \includegraphics[scale=.5]{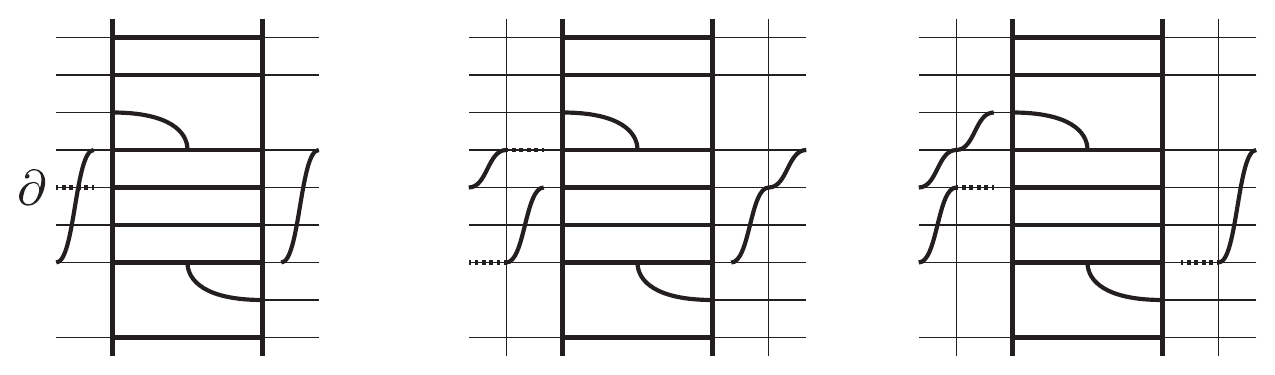}
    \end{center}
    \caption {{\bf The three types of terms in
        Lemma~\ref{lem:dC1C2}.} Left: a Type
      (\ref{typeO:ExtraSigma}) near-chord. Center: the product of a Type
      (\ref{typeO:ExtraSigma}) near-chord and a Type
      (\ref{typeO:Generic}) near-chord. Right: the product of a
      Type~\ref{typeO:dC1C2} near-chord and a Type~\ref{typeO:Sigma} near-chord.}
      \label{fig:dC1C2}
  \end{figure}
  Since any near-chord of Type~(\ref{typeO:dC1C2}) has non-trivial product with
  any near chord of Type~(\ref{typeO:Sigma}), it follows at once that the set of elements of
  Type~(\ref{typeO:dC1C2}) are determined by the terms of Type~(\ref{typeO:ExtraSigma}) which appear in the differential;
  and this latter is the basic choice.
\end{proof}

\begin{lemma}
  \label{lem:OIndeterminateDisconnected}
  If $N$ is a stable arc-slide bimodule then 
  its basic choice ${\mathcal B}$ uniquely specifies which 
  near-chords of Type~(\ref{typeO:Disconnected})
  appear in the differential on $N$.
\end{lemma}

\begin{proof}
  Suppose $x$ is a near-chord of
  Type~(\ref{typeO:Disconnected}).  We can find a near-chord~$y$ of
  Type~(\ref{typeO:Sigma}) so that $x\cdot y$ has exactly two
  alternate factorizations: $x\cdot y = w\cdot z = z'\cdot w'$, where
  $w$ and $w'$ are indeterminate of Type~(\ref{typeO:ExtraSigma})
  and $z$ and $z'$ are of Type~(\ref{typeO:Generic}).  Moreover,
  $x \cdot y$ does not appear in the differential of any other algebra
  element.  From this (and Lemma~\ref{lem:OGenericExists}), we see that $x$ appears in the differential if
  exactly one of $w$ or $w'$ appears in the differential. See
  Figure~\ref{fig:IndeterminateDisconnected} for an illustration.
  \begin{figure}
    \begin{center}
      \includegraphics[scale=.5]{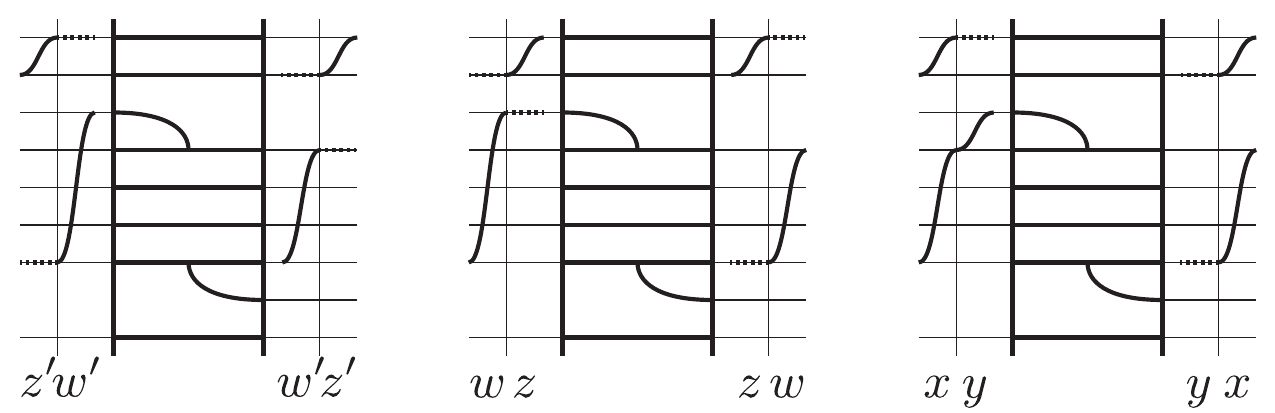}
    \end{center}
    \caption {{\bf The three types of terms in
        Lemma~\ref{lem:OIndeterminateDisconnected}.} Left: $z'\cdot
        w'$. Center: $w\cdot z$. Right: $x\cdot y$. 
      \label{fig:IndeterminateDisconnected}}
  \end{figure}
\end{proof}

\begin{lemma}
  \label{lem:NearChordsTimesDischord}
  Let $x$ be a near-chord and $y$ be a dischord. If $x\cdot y$ is
  non-zero then $x\cdot y$ is a near-chord as well, and similarly for
  $y\cdot x$.
\end{lemma}

\begin{proof}
  This is a simple case analysis on the type of $x$.

  If $x$ has Type~(\ref{typeO:Generic}) and $x\cdot y\neq 0$ then
  $x\cdot y$ has Type~(\ref{typeO:Disconnected}).
  If $x$ has Type~(\ref{typeO:Sigma}) and $x\cdot y\neq 0$ then
  $x\cdot y$ has Type~(\ref{typeO:ExtraSigma}).
  If $x$ has Type~(\ref{typeO:FewerSigma}) and $y\cdot x\neq 0$ then
  $y\cdot x$ has Type~(\ref{typeO:ExtraSigma}).
  All other cases vanish.

  (Alternately, the result follows from noting that the grading of $x\cdot y$
  in the near-diagonal subalgebra must be $-1$, and so by
  Lemma~\ref{lem:over-slide-grading} $x\cdot y$ must be a near-chord.)
\end{proof}

\begin{proof}[Proof of Proposition~\ref{prop:CalculateOver-Slide}]
  Let $N$ be a stable arc-slide bimodule and let ${\mathcal B}$ be its
  associated basic choice (whose existence is guaranteed by
  Lemma~\ref{lem:IsBasicChoice}).  
  By hypothesis, near-chords of Type~(\ref{typeO:Sigma}) exist in the differential; 
  combining this with
  Lemmas~\ref{lem:OGenericExists},
  \ref{lem:OExtraFewerExist}, and \ref{lem:OYYexist}, we conclude that that all the
  near-chords which are not indeterminate appear in the differential.
  According to Lemmas~\ref{lem:OIndeterminateFewerSigma},
  \ref{lem:dC1C2} and~\ref{lem:OIndeterminateDisconnected}, the basic
  choice ${\mathcal B}$ uniquely determines all the indeterminate
  near-chords which contribute to the differential.  According to
  Lemma~\ref{lem:ONoOthers} no other terms
  can contribute to the differential. In sum, the basic choice ${\mathcal B}$
  uniquely determines $N$.

  Next, let $N$ be compatible with a basic choice ${\mathcal B}$, and
  let  ${\mathcal B'}$ be a different basic choice. Then we can find a sum
  of dischords ${\mathcal Q}$ with the property that
  \begin{equation}
    \label{eq:FindGaugeTransformation}
    {\mathcal B}+{\mathcal B}' =
    {\mathcal Q}\cdot (a(\sigma)\otimes 1 + 1\otimes a(\sigma')),
  \end{equation}
  where in the above equation we do not distinguish between a
  basic choice ${\mathcal B}$ and its associated algebra element
  $$\sum_{b\in{\mathcal B}} b.$$
  We use ${\mathcal Q}$ to construct a new bimodule $N'$,
  with the same generators as $N$ and differential given by
  \begin{equation}\label{eq:Gauge}
  \partial_N'=(\Id+\cdot {\mathcal Q})\circ \partial_N\circ
  (\Id+\cdot {\mathcal Q}).
  \end{equation}
  Here, $\Id$ denotes the identity map
  $N\to N$ while $\cdot {\mathcal Q}$ denotes the map induced by
  $a\cdot \x(I)\mapsto a\cdot {\mathcal Q}\cdot \x(I)$ (where $\x(I)$
  is the generator of $N$ corresponding to the near-complementary
  idempotent $I$).  Since $(\Id+\cdot {\mathcal Q})^2=\Id$ (as
  ${\mathcal Q}^2=0$), it follows that $\partial_N'$ is a
  differential.  Note that $N'$ is also an arc-slide bimodule:
  Properties~(\ref{AS:Generators}) and
  (\ref{AS:NearDiagonalSubalgebra}) are clear;
  Property~(\ref{AS:Graded}) follows from the fact that the
  elements appearing in $\partial_N'$ are elements of the near-diagonal
  subalgebra with grading~$-1$, see
  Lemma~\ref{lem:NearChordsTimesDischord};
  Property~(\ref{AS:ShortChords}) continues
  to hold since the operation of replacing $\partial_N$ by $\partial_N'$ does
  not affect the short chords which appear in the differential.
  Similarly, $N'$ is also stable.
  Now, if $N$ is compatible with ${\mathcal B}$ then, according to
  Equation~\eqref{eq:FindGaugeTransformation}, $N'$ is compatible with
  ${\mathcal B'}$. Moreover, the map $f\co N \rightarrow N'$ induced
  by
  $$(\Id+\cdot {\mathcal Q})\colon N \to N'$$
  is an isomorphism of chain complexes.
\end{proof}

\subsection{Arc-slide bimodules}
\label{subsec:IdentificationArc-Slides}

We put together the results (Propositions~\ref{prop:CalculateUnder-Slide}
and~\ref{prop:CalculateOver-Slide}) 
from the previous sections to deduce Theorem~\ref{thm:DDforArc-Slides}.
First, we have the following:

\begin{proof}[Proof of Proposition~\ref{prop:UniqueArc-Slide}]
  For under-slides, this is
  Proposition~\ref{prop:CalculateUnder-Slide}.

  For over-slides, this is Proposition~\ref{prop:CalculateOver-Slide}.
  (Recall that the exact form of the differential is then determined by the basic
  choice, but by
  Proposition~\ref{prop:CalculateOver-Slide}, all of these choices
  give isomorphic bimodules.) The result follows.
\end{proof}

And next:

\begin{proof}[Proof of Theorem~\ref{thm:DDforArc-Slides}]
  This is immediate from the fact that $\CFDDa(\HD(m))$ is a stable
  arc-slide
  bimodule---Proposition~\ref{prop:CFDDisArcslideBimodule}---and the
  fact that all such bimodules are isomorphic---Proposition~\ref{prop:UniqueArc-Slide}. Moreover, it follows
  from~\cite[Corollary~\ref*{LOT2:Cor:MCG-Equivalences} and Lemma~\ref*{LOT2:lem:bimodule-rigidity}]{LOT2}
  that the homotopy equivalence $\DDmod(m)\cong
  (\Alg(\PMC)\otimes\Alg(-\PMC'))\DT\CFDDa(\PunctF(m))$ is unique up to homotopy.
\end{proof}

\begin{remark}
  It is worth noting that the above proof in fact performs all
  the holomorphic curve counts for $\CFDDa(\HD(m))$ when $m$ is an
  under-slide. This is not the case for over-slides. For
  instance, disconnected domains never contribute to a differential in
  $\CFDDa(\HD(m))$.  This can be used to give a criterion which
  ensures that certain near-chords of
  Type~(\ref{typeO:Disconnected}) do not contribute to the
  differential on $\CFDDa(\HD(m))$ (although one can arrange that they do
for other stable arc-slide bimodules). However, this point is irrelevant
  for our purposes: we are interested
  in $\CFDDa(\HD(m))$ only up to homotopy equivalence.
\end{remark}
\setlength{\tabextrasep}{0pt}


\section{The genus-one case}
\label{sec:GenusOne}

In this section, we illustrate the computations of the arc-slide
bimodules by spelling out the answers explicitly in the genus-one
case. We will focus on the part of the algebra with weight zero (i.e., one
moving strand), as $\Alg(T^2,-1)=\Field$ and $\Alg(T^2,1)\simeq
\Field$, and so the bimodules over these algebras are not very
interesting (free of rank $1$ with trivial differential, when viewed
as bimodules modules over $\Field$).

The (unique) pointed matched circle for the torus has $4$
matched points, $1,2,3,4$, with $1$ matched to $3$ and $2$ matched to
$4$. Let $\iota_0$ denote the idempotent in $\Alg(T^2,0)$
corresponding to $\{1,3\}$ and $\iota_1$ the idempotent corresponding
to $\{2,4\}$. Let $\rho_i$ denote the (short) chord from $i$ to $i+1$,
so
\[
\iota_0\rho_1\iota_1=\rho_1\qquad 
\iota_1\rho_2\iota_0=\rho_2\qquad 
\iota_0\rho_3\iota_1=\rho_3.
\]
Let
\[
\rho_{12}=\rho_1\rho_2 \qquad \rho_{23}=\rho_2\rho_3\qquad\rho_{123}=\rho_1\rho_2\rho_3.
\]

Of course, we are considering two copies of $\Alg(T^2)$. In the second
copy, we will denote the idempotents corresponding to $\iota_0$ and
$\iota_1$ by $j_0$ and $j_1$, and the chord corresponding to
$\rho_i$ by $\sigma_i$.

With this notation in hand, consider the \DD\ identity module, which
was already computed in this case
in~\cite[Proposition~\ref*{LOT2:prop:TypeDDTorus}]{LOT2}. The module
$\CFDDa(\Id)$ has two generators, $p=(\iota_0\otimes j_0)$ and
$q=(\iota_1\otimes j_1)$, and differential
\[
\partial \mathbf p = (\rho_1\sigma_3+ \rho_3\sigma_1 +
\rho_{123}\sigma_{123}) \otimes \mathbf q \qquad
\partial \mathbf q = (\rho_2\sigma_2)\otimes \mathbf p.
\]
This is in obvious agreement with
Theorem~\ref{thm:DDforIdentity}. Note that the term
$\rho_{123}\sigma_{123}$ comes from a special length $3$ chord, and is
not forced by $\bdy^2=0$.

\begin{figure}
  \centering
  \includegraphics{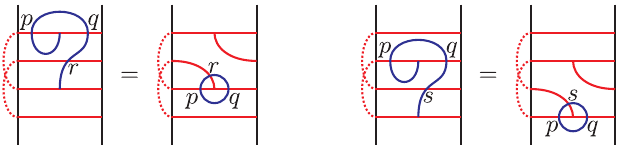}
  \caption{\textbf{Generators of the genus-one mapping class group.}
    This figure may be compared with~\cite[Figure~\ref*{LOT2:fig:DehnTwistsGenusOne}]{LOT2}.}
  \label{fig:genus-1-mcg-gens}
\end{figure}

Next we turn to the arc-slide bimodules. The genus $1$ mapping class
group is generated by Dehn twists $\tau_m$, $\tau_l$ around two curves
$m$, $l$ in $T^2$. We can view $\tau_m$ and $\tau_l$ as underslides;
see Figure~\ref{fig:genus-1-mcg-gens}.  We give the answers, and then
explain the terms.

The module $\CFDDa(\tau_m)$ has three generators,
$p=\iota_0\otimes j_0$, $q=\iota_1\otimes j_1$ and $r=\iota_1\otimes
j_0$. The differential is given by
\begin{align*}
  \bdy(p)&=(\stackrel{\text{\tiny{(U-1)}}}{\rho_1\sigma_3}+\stackrel{\text{\tiny{{(U-1)}}}}{\rho_{123}\sigma_{123}})\otimes
  q+(\stackrel{(\text{U-4})}{\rho_3\sigma_{12}})\otimes r\\
  \bdy(q)&=(\stackrel{(\text{U-4})}{\rho_{23}\sigma_2})\otimes r\\
  \bdy(r)&=(\stackrel{(\text{U-2})}{\rho_2})\otimes p + (\stackrel{(\text{U-2})}{\sigma_1})\otimes q.
\end{align*}
The type of the near-chord corresponding to each term is indicated
above it. There are no contributions from chords of types
(U-3), (U-5) or (U-6).

The terms $\rho_1\sigma_3$, $\rho_2$ and $\sigma_1$ are short
near-chords, and so appear by hypothesis / directly counting
holomorphic curves. The remaining terms are not forced by $\bdy^2=0$
in this diagram, but instead follow from the fact that our bimodule is
stable, together with $\bdy^2=0$ on a bigger diagram.

Similarly, the module $\CFDDa(\tau_l)$ has three generators, 
$p=\iota_0\otimes j_0$, $q=\iota_1\otimes j_1$ and $s=\iota_0\otimes
j_1$. The differential is given by
\begin{align*}
  \bdy(p)&=(\stackrel{\text{(U-1)}}{\rho_3\sigma_1}+\stackrel{\text{(U-1)}}{\rho_{123}\sigma_{123}})\otimes
  q+(\stackrel{\text{(U-4)}}{\rho_{12}\sigma_3})\otimes
  s\\
  \bdy(q)&=(\stackrel{\text{(U-4)}}{\rho_2\sigma_{23}})\otimes s\\
  \bdy(s)&=(\stackrel{\text{(U-2)}}{\sigma_2})p+(\stackrel{\text{(U-2)}}{\rho_1})\otimes
  q.
\end{align*}

(The modules $\CFDAa(\tau_m)$ and $\CFDAa(\tau_l)$ were computed
directly in~\cite[Section~\ref*{LOT2:subsec:DA-mcg}]{LOT2}. Tensoring
these modules with $\CFDDa(\Id)$ gives another computation of
$\CFDDa(\tau_m)$ and $\CFDDa(\tau_l)$.)

\begin{figure}
  \centering
  \includegraphics{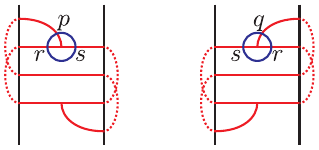}
  \caption{\textbf{Genus-one overslides.}}
  \label{fig:genus-1-over}
\end{figure}

There are also two overslides in genus $1$; they are shown in
Figure~\ref{fig:genus-1-over}. The bimodule $\CFDDa$ associated to the
overslide on the left of Figure~\ref{fig:genus-1-over} has three
generators, $p=\iota_0\otimes j_0$, $r=\iota_1\otimes j_0$ and
$s=\iota_0\otimes j_1$, and differential given by
\begin{align*}
  \bdy(p)&=(\stackrel{\text{(O-2)}}{\rho_3}
  +\textcolor{red}{\stackrel{\text{(O-3)}}{\rho_{123}\sigma_{12}}})\otimes
  r+(\stackrel{\text{(O-2)}}{\sigma_3}
  +\textcolor{red}{\stackrel{\text{(O-3)}}{\rho_{12}\sigma_{123}}})\otimes
  s\\
  \bdy(r)&=(\stackrel{\text{(O-1)}}{\rho_2\sigma_1})\otimes s\\
  \bdy(s)&=(\stackrel{\text{(O-1)}}{\rho_1\sigma_2})\otimes r,
\end{align*}
with the understanding that exactly one of the terms labeled (O-3) occurs.

The (O-1) and (O-2) terms are short near-chords, and so occur by
hypothesis / direct computation. Both (O-3) terms are indeterminate;
the fact that exactly one of them occurs is part of the definition of
a basic choice (Definition~\ref{def:BasicChoice}). (In terms of
holomorphic curves, the (O-3) near-chords do not follow from
$\bdy^2=0$ in this diagram, but do follow from $\bdy^2=0$ after
stabilizing; see the proof of Lemma~\ref{lem:IsBasicChoice}.)

We can switch between basic choices by making the change of variables
$p'=(1+\rho_{12}\sigma_{12}) \otimes p$ (compare Formula~\eqref{eq:Gauge}).

For completeness, the bimodule corresponding to the overslide on the
right has generators $q=\iota_1\otimes j_1$, $r=\iota_1\otimes j_0$ and
$s=\iota_0\otimes j_1$, and differential given by
\begin{align*}
  \bdy(q)&=0\\
  \bdy(r)&=(\stackrel{\text{(O-1)}}{\rho_2\sigma_3})\otimes s+(\stackrel{\text{(O-2)}}{\sigma_1}+\textcolor{red}{\stackrel{\text{(O-3)}}{\rho_{23}\sigma_{123}}})\otimes q\\
  \bdy(s)&=(\stackrel{\text{(O-1)}}{\rho_3\sigma_2})\otimes r
  +(\stackrel{\text{(O-2)}}{\rho_1}+\textcolor{red}{\stackrel{\text{(O-3)}}{\rho_{123}\sigma_{23}}})\otimes q,
\end{align*}
again with the understanding that exactly one of the terms labeled (O-3) contributes.
(Which one contributes is determined by the basic choice; bear in mind that this diagram is mirror to Convention~\ref{label:conv}.)


\section{Gradings on bimodules for arc-slides and mapping classes}
\label{sec:mcg-grading}
In this section, we discuss further the gradings on the type \DD\
modules associated to arc-slides. The main goal is to
compute explicitly the gradings on $\CFDDa(\PunctF(m))$, the bimodule
associated to an arc-slide $m$; this is done in
Section~\ref{sec:comp-gr}. As we explain in
Section~\ref{sec:CompleteProof}, this allows one to compute both the
(relative) Maslov gradings on $\HFa(Y)$ and the decomposition of $\HFa(Y)$
into $\SpinC$-structures. Section~\ref{sec:general-grads} is a brief
digression to compute the grading sets for general surface
homeomorphisms; Section~\ref{sec:general-grads} is not needed for the
rest of the paper, but answers a question which arises naturally.

\subsection{Gradings on arc-slide bimodules}\label{sec:comp-gr}
In Section~\ref{sec:grade-periodic-domain}, we finish computing the
gradings of periodic domains for the standard Heegaard diagrams for
arc-slides; this was begun in Section~\ref{sec:grad-near-diag}. In
Section~\ref{sec:grading-set-refined}, we give the grading set for
$\CFDDa(\PunctF(m))$ (i.e., the range of the grading function), with
respect to both the big and small grading groups. In
Section~\ref{sec:gr-of-gens} we compute the gradings of
generators of $\CFDDa(\PunctF(m))$, i.e., the grading function itself.

In this section, we will work with the bimodules $\CFDDa(\PunctF(m))$
associated to a Heegaard diagram,
rather than a general arc-slide bimodule $\DDmod(m)$. 
No generality is lost, according to the following:
\begin{proposition}
  Suppose that $\DDmod(m\co\PMC\to\PMC')$ is a stable arc-slide bimodule, and that
  the actions of $G'(\PMC)$ and $G'(-\PMC')$ on the grading set for
  $\DDmod(m)$ are free and transitive. Then the homotopy equivalence
  $\DDmod(m)\simeq \CFDDa(m)$ is a $G'$-set graded homotopy equivalence.
\end{proposition}
\begin{proof}
  This follows from Proposition~\ref{prop:grading-unique}.
\end{proof}
 
\subsubsection{Gradings of periodic domains}\label{sec:grade-periodic-domain}
Consider the standard Heegaard diagram $\HD(m)$ associated to an
arc-slide $m\co \PMC\to\PMC'$. Label the matched pairs in $\PMC$ by
$B_1,\dots,B_{2k}$, with $B_{2k}=B$ (in the notation of
Section~\ref{sec:Arc-Slides}).
\glsit{$B_1,\dots,B_{2k}$}%

\begin{figure}
  \centering
  \includegraphics{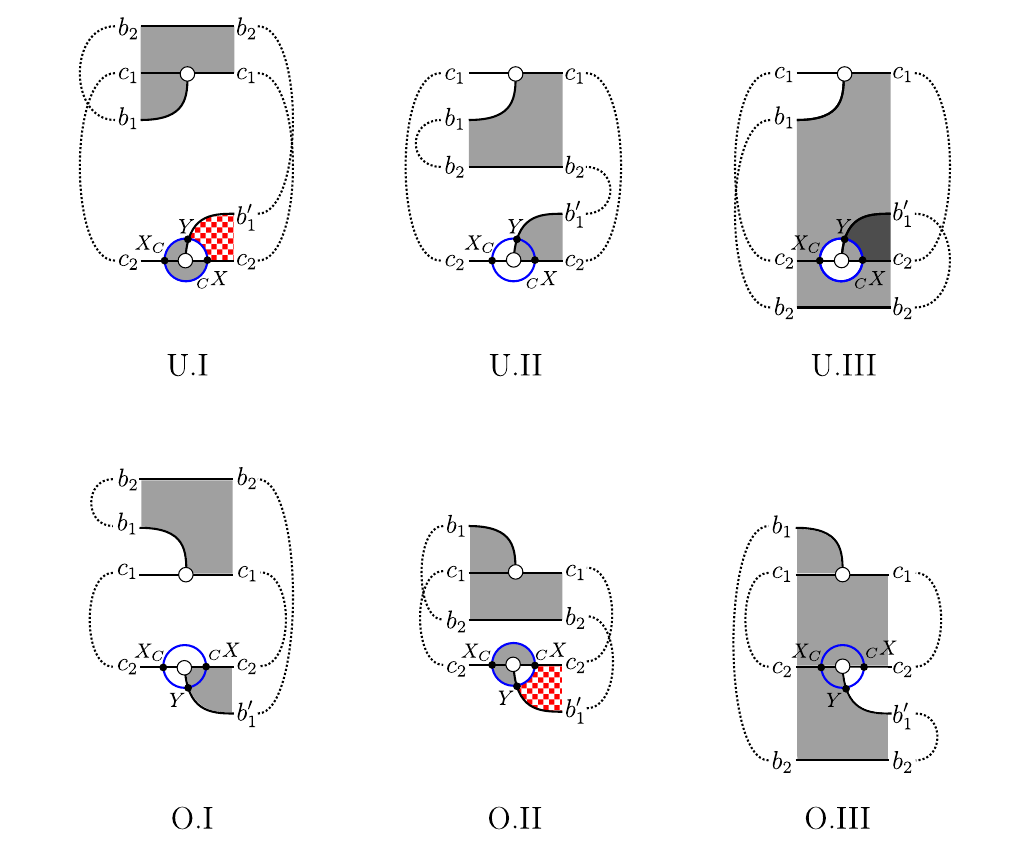}
  \caption{\textbf{Combinatorial cases of arc-slides.} The domain $P_{2k}$ is shaded in each diagram. Light gray shading indicates multiplicity $1$, dark gray shading indicates multiplicity $2$, and red, checkered shading indicates multiplicity $-1$.}
  \label{fig:comb-slide-cases}
\end{figure}

There are various combinatorial cases for the arc-slide $m$. 
For the remainder of Section~\ref{sec:comp-gr}, we will assume that
$c_1$ is above $c_2$ in $\PMC$, as in Section~\ref{sec:Arc-Slides};
the other case is obtained by reflecting the diagrams horizontally, or
equivalently by replacing $m$ by $m^{-1}$.
Then we have the following cases:
\begin{itemize}
\item Under-slides: $b_1$ is between $c_1$ and $c_2$. This case is further divided as follows:
  \begin{itemize}
  \item[(U.I)] $c_1$ is between $b_1$ and $b_2$.
  \item[(U.II)] $b_2$ is between $c_1$ and $c_2$.
  \item[(U.III)] $c_2$ is between $b_1$ and $b_2$.
  \end{itemize}
\item Over-slides: $b_1$ is above $c_1$. This case is further divided as follows:
  \begin{itemize}
  \item[(O.I)] $b_2$ is above $c_1$.
  \item[(O.II)] $b_2$ is between $c_1$ and $c_2$ (so $c_1$ is between $b_1$ and $b_2$).
  \item[(O.III)] $c_1$ and $c_2$ are both between $b_1$ and $b_2$.
  \end{itemize}
\end{itemize}
These cases are illustrated in Figure~\ref{fig:comb-slide-cases}.

We can find a basis for the space of periodic domains on $\HD(m)$,
$\pi_2(\x,\x)$ (for any generator $\x$), given by elements
$P_1,\dots,P_{2k}$ with the following properties. For $i<2k$, $\bdy
P_i\cap \PMC$ is the interval between the two points in $B_i$, and $\bdy P_i\cap \PMC'$ is the interval between the two points in $B'_i$. The domain $P_{2k}$ is such that:
\begin{itemize}
\item[(U.I)] $\bdy P_{2k}\cap \PMC$ consists of $[b_1,b_2]$, while
  $\bdy P_{2k}\cap \PMC'$ consists of
  $r([c_1,b_2])-r([c_2,b'_1])$. (The region $\sigma'$ is covered with
  multiplicity $-1$.)
\item[(U.II)] $\bdy P_{2k}\cap \PMC$ consists of $[b_2,b_1]$, while
  $\bdy P_{2k}\cap \PMC'$ consists of $r([b_2,c_1])+r([c_2,b'_1])$.
\item[(U.III)] $\bdy P_{2k}\cap \PMC$ consists of $[b_2,b_1]$, while
  $\bdy P_{2k}\cap \PMC'$ consists of
  $r([b_2,c_1])+r([c_2,b'_1])$. (The region $\sigma'$ is covered with
  multiplicity two.)
\item[(O.I)] $\bdy P_{2k}\cap \PMC$ consists of $[b_1,b_2]$, while
  $\bdy P_{2k}\cap \PMC'$ consists of $r([b'_1,c_2])+r([c_1,b_2])$.
\item[(O.II)] $\bdy P_{2k}\cap \PMC$ consists of $[b_2,b_1]$, while
  $\bdy P_{2k}\cap \PMC'$ consists of $r([b_2,c_1])-r([b'_1,c_2])$. (The
  region $\sigma'$ is covered with multiplicity $-1$.)
\item[(O.III)] $\bdy P_{2k}\cap \PMC$ consists of $[b_2,b_1]$, while
  $\bdy P_{2k}\cap \PMC'$ consists of $r([b_2,b'_1])+r([c_2,c_1])$.
\end{itemize}

See Figure~\ref{fig:comb-slide-cases} for an illustration.

As discussed in Section~\ref{sec:coeff-algebra}, for each generator
$\x$ and each $P_i\in\pi_2(\x,\x)$ there is an element
\[
g'(P_i)=(-e(P_i)-2n_\x(P_i),\bdy^\bdy P_i)\in G'(\PMC)\times_\ZZ G'(-\PMC').
\]
\glsit{$g'$}%
We compute these elements $g'(P_i)$:
\begin{lemma}\label{lem:comp-per-dom-gr}
  For the domain $P_i\in\pi_2(\x(I),\x(I))$ with corresponding
  matched pair~$B_i$, the Maslov component of
  $g'(P_i)$ is
  \begin{itemize}
  \item $0$ if $B_i\neq B,C$;
  \item $0$ if $B_i=C$ and $I$ has type $\lsub{C}X$ or $X_C$;
  \item $1$ if $B_i=C$, $m$ is an over-slide, and $I$ has type $Y$;
  \item $-1$ if $B_i=C$, $m$ is an under-slide, and $I$ has type $Y$; and
  \item given by the tables below for $P_{2k}$.
  \end{itemize}

\noindent
\hfill
\parbox[t]{.3\textwidth}{
If $I$ has type $X_C$:

\vspace{2pt}
\begin{tabular}{cMc}
  \toprule
  Handleslide & \text{Maslov}\\
  type        & \text{component}\\
  \midrule
  U.I & -1/2\\
  U.II & \hphantom{-}1/2\\
  U.III & \hphantom{-}1/2\\
  O.I & \hphantom{-}1/2\\
  O.II & -1/2\\
  O.III & -1/2\\
  \bottomrule
\end{tabular}
}\hfill
\parbox[t]{.3\textwidth}{
If $I$ has type $\lsub{C}X$:

\vspace{2pt}
\begin{tabular}{cMc}
  \toprule
  Handleslide & \text{Maslov}\\
  type        & \text{component}\\
  \midrule
  U.I & \hphantom{-}1/2\\
  U.II & -1/2\\
  U.III & -1/2\\
  O.I & -1/2\\
  O.II & \hphantom{-}1/2\\
  O.III & \hphantom{-}1/2\\
  \bottomrule
\end{tabular}
}\hfill
\parbox[t]{.3\textwidth}{
If $I$ has type $Y$:

\vspace{2pt}
\begin{tabular}{cMc}
  \toprule
  Handleslide & \text{Maslov}\\
  type        & \text{component}\\
  \midrule
  U.I & \hphantom{-}1/2\\
  U.II & -1/2\\
  U.III & -1/2\\
  O.I & -1/2\\
  O.II & \hphantom{-}1/2\\
  O.III & \hphantom{-}1/2\\
  \bottomrule
\end{tabular}
}\hfill
\end{lemma}

\begin{proof}
  This follows by inspecting the periodic domains in
  Figure~\ref{fig:comb-slide-cases} and applying the calculations from
  Proposition~\ref{prop:near-diag-grading}: in the notation from that
  proposition, the Maslov component of $g'(P_i)$ is $2c(I(\x),\supp(a))$.
\end{proof}


\subsubsection{The grading set and refined grading
  set}\label{sec:grading-set-refined}
\index{grading!set!for arc-slide}%
Following \cite[Section~\ref*{LOT2:sec:cf-gradings}]{LOT2}, as a $G'$-set graded module, the
module $\CFDDa(\HD(m))$ is graded by
\[
G'(\PMC)\times_\ZZ G'(-\PMC')/\langle R(g'(P_i))\rangle,
\]
\glsit{{$\protect\langle R(g'(P_i))\protect\rangle$}}%
where $\langle g'(P_i)\rangle$ denotes the subgroup of
$G'(-\PMC)\times_\ZZ G'(\PMC')$ generated the gradings of the
periodic domains, and $R$ is the orientation-reversal map from
Sections~\ref{sec:coeff-algebra} and~\ref{sec:coeff-bimod}.
The Maslov components of these gradings are computed
in Lemma~\ref{lem:comp-per-dom-gr}; the homology class component of
$R(g'(P_i))$ is simply $r_*(\bdy^\bdy P_i)$.

In the notation of Section~\ref{sec:grade-periodic-domain},
since $P_i$ is a periodic domain, $R(g'(P_i))$ lies in the
smaller grading group $G(\PMC)\times_\ZZ
G(-\PMC')$~\cite[Lemma~\ref*{LOT2:lem:Refinable}]{LOT2}.
As a $G$-set graded module, the module $\CFDDa(\HD(m))$ is graded by 
\[
S=G(\PMC)\times_\ZZ G(-\PMC')/\langle R(g'(P_i))\rangle.
\]

\glsit{$S$}\glsit{$H$}%
Let $H=\langle R(g'(P_i))\rangle$.
The subgroup $H$ is a complement to $G(\PMC)$ (respectively
$G(-\PMC')$), in the sense that $G(\PMC)$ (respectively $G(-\PMC')$)
intersects each coset of $H$ exactly once.
So, each of $G(\PMC)$ and $G(-\PMC')$ act freely and
transitively on $S$, i.e., $S$ is a \emph{principal left-right
  $G(\PMC)$-$G(\PMC')$-set}.
\index{principal $G(\PMC)$-$G(\PMC')$-set}%
(Here we interpret the left action of
$G(-\PMC')$ as a right action of $G(\PMC')$.)  We
spend the next few paragraphs re-interpreting $S$ as a map.

\begin{definition}\label{def:outer-aut}
  \glsit{$\Isom(G_1,G_2)$}\glsit{$\Out(G_1,G_2)$}%
  Let $G_1$ and $G_2$ be isomorphic groups. Let $\Isom(G_1,G_2)$ denote the set of isomorphisms from $G_1$ to $G_2$. Then $G_1$ acts on $\Isom(G_1,G_2)$ as follows: for $\phi\in \Isom(G_1,G_2)$ and $g\in G_1$ define
  \[
  (g\cdot \phi)(h) = \phi(g^{-1}hg).
  \]
  Let $\Out(G_1,G_2)=G_1\backslash \Isom(G_1,G_2)$. An element of
  $\Out(G_1,G_2)$ is an \emph{outer isomorphism from $G_1$ to $G_2$.}
  \index{outer isomorphism}%
\end{definition}
In Definition~\ref{def:outer-aut}, one could alternatively take the 
quotient $G_2 \backslash \Isom(G_1,G_2)$, with the action given by
$(g'\cdot\phi) (h)=g'\phi(h)(g')^{-1}$. This gives the same
equivalence relation on $\Isom(G_1,G_2)$.

Composition of maps induces a map $\Out(G_1,G_2)\times \Out(G_2,G_3)\to \Out(G_1,G_3)$. The following lemma is straightforward:

\begin{lemma}\label{lem:G-torseur-is-map}
  For any (isomorphic) groups $G_1$ and $G_2$,
  there is a canonical bijection between isomorphism classes of principal
  left-right $G_1$-$G_2$-sets~$S$ and outer isomorphisms
  $\Out(G_1,G_2)$. Given a left-right $G_1$-$G_2$ set $S$,
  let $\phi_S\in \Out(G_1,G_2)$ denote its corresponding outer isomorphism
  under this canonical bijection.
  If $S$ is a left-right principal $G_1$-$G_2$-set, and $T$ is a left-right principal
  $G_2$-$G_3$ set, we can form the orbit space $S\times_{G_2}T$. This
  is a left-right principal $G_1$-$G_3$-set. Moreover, 
  $\phi_{S {\times_{G_2}} T} = \phi_T \circ \phi_S$.
\end{lemma}

\begin{proof}
  For a given principal left-right $G_1$-$G_2$-set~$S$, pick any $x \in S$,
  and define $\phi_S$ by
  \[
  g\cdot x = x \cdot \phi_S(g).
  \]
  (This defines $\phi_S$ uniquely up to conjugation because the actions of both $G_1$
  and $G_2$ on~$S$
  are simply transitive.)  The map $\phi_S$ is a homomorphism:
  \begin{align*}
    x\cdot\phi_S(gh)&=(gh) \cdot x\\
    &= (g \cdot x) \cdot \phi_S(h)\\
      &= (x \cdot \phi_S(g))\cdot \phi_S(h) \\
      &= x \cdot \phi_S(g)\phi_S(h).
  \end{align*}
  $\phi_S$ is also clearly bijective.

  If we choose a different element $x'\in S$ to define $\phi_S$, the
  isomorphism changes by an inner automorphism, as follows.  Let
  $\phi'_S$ be the automorphism defined with respect to $x'$.  Choose
  $g'_0$ so that $x' = g_0 x$.  Then
  \begin{align*}
    x'\cdot \phi'_S(g) &= g\cdot x'\\
    &=g\cdot g_0 x\\
    &=x\cdot \phi_S(g g_0)\\
    &=(g_0)^{-1} x' \phi_S(g g_0)\\
    &=x' \phi_S((g_0)^{-1} g g_0).
  \end{align*}
  Thus the element $\phi_S\in\Out(G_1,G_2)$ is well-defined.

  We claim the assignment $S\mapsto \phi_S$ gives a bijection between the set of isomorphism classes of left-right transitive $G_1$-$G_2$ sets and the set $\Out(G_1,G_2)$.
  To invert this bijection, suppose that $\phi\colon G_1\to G_2$ is a
  group isomorphism and consider the associated principal $G_1$-$G_2$-set 
  $S_\phi$ defined to be the quotient of $G_1\times G_2$ by the equivalence relation
  $(g_1\cdot g,g_2)\sim (g_1,\phi(g)\cdot g_2)$ where $g\in G_1$. This quotient clearly retains its left-right $G_1$-$G_2$-set structure.
  We leave it as an exercise for the reader to verify the following facts:
  \begin{itemize}
  \item The quotient space $S_\phi$ defined above is indeed a transitive $G_1$-$G_2$ set.
  \item Given two isomorphisms $\phi$ and $\phi'$ from $G_1$ to $G_2$ representing the same outer isomorphism,
    there is a canonical isomorphism of left-right $G_1$-$G_2$-sets $S_\phi\cong S_{\phi'}$.
  \item Given any transitive $G_1$-$G_2$ set $S$, there is a canonical isomorphism $S_{\phi_S}\cong S$. 
  \item Given an isomorphism $\phi\colon G_1\to G_2$, $\phi_{S_\phi}=\phi$, as elements of $\Out(G_1,G_2)$.
  \end{itemize}

  Finally, let $S_{12}$ be a principal $G_1$-$G_2$-set and $S_{23}$ a principal $G_2$-$G_3$-set,
  with corresponding automorphisms $\phi_{12}$ and $\phi_{23}$, defined with
  respect to $x_{12} \in S_{12}$ and $x_{23} \in S_{23}$.  Then
  \begin{align*}
    (x_{12} \times x_{23})\cdot \phi_{S_{12} \times_{G_2} S_{23}}(g)&=
    g\cdot (x_{12}\times x_{23})\\
    &=(x_{12}\cdot \phi_{12}(g))\times x_{23}\\
    &=x_{12}\times (\phi_{12}(g)\cdot x_{23})\\
    &=(x_{12}\times x_{23})\cdot \phi_{23}(\phi_{12}(g)).
  \end{align*}
  This verifies the last part of the claim.
\end{proof}

\glsit{$\phi=\phi_S=\phi_m$}%
Our next goal is to describe explicitly
the map $\phi=\phi_S=\phi_m\co G(\PMC)\to G(\PMC')$
associated to the grading set for the arc-slide $m$. For each matched pair $B_i$ in $\PMC$
there is an associated homology class $h(B_i)\in H_1(F(\PMC))$, namely the
\glsit{$h(B_i)$}%
sum of the core of the handle attached to the pair $B_i\subset
Z$ and the segment in $Z$ between the points in $B_i$. We
orient $h(B_i)$ so that the induced orientation of $h(B_i)\cap Z$
agrees with the orientation of $Z$. The
group $G(\PMC)$ is generated by the element $\lambda$ and the $2k$
elements
\[
\gpElt(B_i)=(-1/2;h(B_i)),
\]
\glsit{$\gpElt(B_i)$}%
where we write elements of $G(\PMC)$ as pairs $(m;h)$ where
\glsit{$(m;h)$}\glsit{$\lambda$}%
$m\in\frac{1}{2}\ZZ$ and $h\in H_1(F(\PMC))$ subject to a congruency
condition (compare Section~\ref{sec:gradings-alg} and~\cite[Section
3.1]{LOT2}); in this notation, $\lambda=(1;0)$. Recall that for
$B_i\neq B$ (i.e., $i\neq 2k$), there is a corresponding matched pair
$B_i$ in $\PMC'$, while $B$ corresponds to a matched pair $B'$. We use
$\othh(B_i)$ (respectively $\othgpElt(B_i)$) to denote the homology
class in $F(\PMC')$ (respectively group element in $G(\PMC')$)
associated to the matched pair $B_i$.
\begin{lemma}\label{lem:gr-map-is} The map $\phi$ is given (up to inner isomorphism) by
  $\phi(\lambda)=\lambda$,
  $\phi(\gpElt(B_i))=\othgpElt(B_i)$ if $B_i\neq B$, and
  \[
  \phi(\gpElt(B))=
  \begin{cases}
    \lambda^{-1}\othgpElt(B')\othgpElt(C)^{-1}=(0;\othh(B')-\othh(C)) & \text{in case U.I}\\
    \lambda^{-1}\othgpElt(B')^{-1}\othgpElt(C)=(-1;-\othh(B')+\othh(C)) & \text{in case U.II}\\
     \lambda^{-1}\othgpElt(B')\othgpElt(C)=(-1;\othh(B')+\othh(C)) & \text{in case U.III}\\
    \lambda^{-1}\othgpElt(B')\othgpElt(C)^{-1}=(-1;\othh(B')-\othh(C)) & \text{in case O.I}\\
    \lambda^{+1}\othgpElt(B')^{-1}\othgpElt(C)=(0;-\othh(B')+\othh(C)) & \text{in case O.II}\\
    \lambda^{+1}\othgpElt(B')\othgpElt(C)=(0;\othh(B')+\othh(C)) & \text{in case O.III}
  \end{cases}
  \]
\end{lemma}
\begin{proof}
  Choose a base idempotent of type $\lsub{C}X$. By Lemma~\ref{lem:comp-per-dom-gr}, the grading set for $\HD(m)$ is given by $G(\PMC)\times_\ZZ G(-\PMC')/H$ where $H$ is generated by $(0;h(B_i),-h(B'_i))$ for $i\neq 2k$, and one additional element:
  
  \begin{center}
    \begin{tabular}{rc}
      \toprule
      Handleslide type & $g(P_{2k})$\\ \midrule
      U.I &$(-1/2;h(B),-h'(B')+h'(C'))$\\
      U.II &$(1/2; h(B), h(B')-h'(C'))$\\
      U.III &$(1/2; h(B), -h(B')-h(C'))$\\
      O.I & $(1/2;h(B),-h'(B')+h'(C'))$\\
      O.II & $(-1/2;h(B), h'(B')-h'(C'))$\\
      O.III & $(-1/2;h(B), -h'(B')-h'(C'))$\\
      \bottomrule
    \end{tabular}
  \end{center}
  (Note that we have reversed all of the signs from Lemma~\ref{lem:comp-per-dom-gr}; the domains shown there have boundary $-h(B)$ in $H_1(F(\PMC))$.)
  By Lemma~\ref{lem:G-torseur-is-map}, this data is reinterpreted as a map. Using 
$(0;0,0)\in G(\PMC)\times_\ZZ G(-\PMC')/H$ as the element $x$ in proof of Lemma~\ref{lem:G-torseur-is-map}, 
for the case U.II, say, we have
  \begin{align*}
  (1/2;h(B))\cdot 
  (0;0,0)\cdot 
  (0;\othh(B')-\othh(C)) 
  &=(0;0,0)\\
  (-1/2; h(B)) \cdot (0;0,0) &= (0; 0,0) \cdot (-1;-\othh(B')+\othh(C)),
  \end{align*}
  which corresponds to the statement that $\phi((-1/2;h(B)))=(-1;-\othh(B')+\othh(C))$.
\end{proof}

\subsubsection{Gradings of generators and algebra elements}\label{sec:gr-of-gens}
Now that we understand the grading sets for an arc-slide $m\co
\PMC\to \PMC'$, the next task is to understand the gradings of
generators of $\CFDDa(\HD(m))$ and of algebra elements of $\Alg(\PMC)$ and
$\Alg(\PMC')$.

\glsit{$I_0$}%
Fix a base idempotent $I_0$ in the near-diagonal subalgebra of
$\Alg(\PMC,i)\otimes \Alg(-\PMC',-i)$. We choose $I_0$ to have type~$X$.
Recall that any such idempotent corresponds to a unique generator
$\x_0=\x(I_0)\in\Gen(\HD(m))$. Choose for each other idempotent $J$ an
element $Q_J\in\pi_2(\x_0,\x(J))$.
\glsit{$Q_J$}%
The unrefined ($G'$-set) grading of generators is given by
\begin{align*}
\gr'(\x(J))&=R(g'(Q_J))\cdot\langle R(g'(P_i))\rangle\\
&=(-e(Q_J)-n_{\x(I)}(Q_J)-n_{\x(J)}(Q_J); r_* \bdy^{\bdy_L}
Q_J, r_* \bdy^{\bdy_R} Q_J)\cdot\langle R(g'(P_i))\rangle,
\end{align*}
where $\bdy^{\bdy_L}Q_J$ and $\bdy^{\bdy_R}Q_J$ denote the
\glsit{$\bdy^{\bdy_L}$, $\bdy^{\bdy_R}$}%
intersections of $Q_J$ with the left and right boundaries of
$\HD(m)$, and we view $G'(\PMC)\times_\ZZ G'(\PMC)$ as a subset of
$\OneHalf\ZZ\times H_1(Z,\mathbf{a})\times H_1(Z',\mathbf{a}')$. These
cosets are independent of the choices of the~$Q_J$.

For definiteness, for idempotents $J$ of type $X$, take $Q_J$ to be a
union of horizontal strips in the graph for the arc-slide
(Figure~\ref{fig:ArcSlide}), so that
\begin{align*}
R(g'(Q_J))&=(0;r_*\bdy^{\bdy_L}Q_J,r_*\bdy^{\bdy_R}Q_J)\\
\gr'(\x(J))
&=(0; r_*\bdy^{\bdy_L} Q_J, r_* \bdy^{\bdy_R} Q_J)\cdot\langle R(g'(P_i))\rangle.
\end{align*}

For an idempotent $J$ of type~$Y$, there is an associated idempotent
\glsit{$J_X$}%
$J_X$ of type $X_C$ obtained by replacing the $C$ in the left of $J$
by a $B$. Then $Q_{J_X}*\sigma^{\epsilon}$ is a domain connecting
$I_0$ to $J$, where $\epsilon$ is $\pm 1$, depending on the geometry
of the arc-slide: $\epsilon$ is $+1$ if $b_1$ is below $c_1$, and $-1$
if $b_1$ is above $c_1$.  We will choose $Q_J$ to be
$Q_{J_X}*\sigma^\epsilon$.  We then have
\begin{align*}
  R(g'(Q_J)) &=  R(g'(\sigma))^\epsilon R(g'(Q_{J_X})) \\
    &= (-\epsilon/2; \epsilon r_*[\sigma] +
       r_* \bdy^{\bdy_L}Q_J,r_* \bdy^{\bdy_R}Q_J) \\
  \gr'(\x(J)) &= R(g'(Q_J)) \cdot \langle R(g'(P_i)) \rangle.
\end{align*}

Next we turn to the refined gradings.
To specify grading refinement data for $\Alg(\PMC)$ and $\Alg(-\PMC')$
(as in~\cite[Section~\ref*{LOT2:subsec:SmallGroup}]{LOT2}),
recall that for each idempotent $J$ of type $X$ there are corresponding
\glsit{$j_L$, $j_R$}%
idempotents $j_L$ and $j_R$ of $\Alg(\PMC)$ and $\Alg(-\PMC')$,
respectively, and that every idempotent of $\Alg(\PMC)$ (respectively
$\Alg(-\PMC')$) arises as $j_L$ (respectively $j_R$) for a unique
idempotent $J$ of type $X$.  

We will use the $Q_J \in \pi_2(\x(I_0),
\x(J))$ for idempotents $J$ of type $X$
 to define grading refinement data.
Specifically, define grading refinement data $\Xi$ for $\Alg(\PMC)$ and
$\Alg(-\PMC')$ by setting, for idempotents $J$ of type~$X$,
\begin{align*}
  \Xi(j_L)&=(\epsilon(\bdy^{\bdy_L}Q_J);-r_*\bdy^{\bdy_L}Q_J)\\
  \Xi(j_R)&=(\epsilon(\bdy^{\bdy_R}Q_J);-r_*\bdy^{\bdy_R}Q_J),
\end{align*}
where $\epsilon\co H_1(Z,\mathbf{a})\to \frac{1}{2}\ZZ$ is as in Section~\ref{sec:gradings-alg}. 
(The minus signs arise because we are giving grading refinement data on
$\Alg(\PMC)$, not on $\Alg(\bdy_L \HD(m)) = \Alg(-\PMC)=\Alg(\PMC)^\op$,
and similarly for $\Alg(-\PMC')$.)

With respect to this choice, if $J$ is an idempotent of type $X$ then
\[
  \gr(\x(J))=\Xi(j_L)\cdot\Xi(j_R)\cdot\gr'(\x(J))
    =(0;0)\cdot \langle R(g'(P_i))\rangle;
\]
\glsit{$\gr$}%
compare \cite[Section~\ref*{LOT1:sec:refined-gradings-2}]{LOT1} and \cite[Section~\ref*{LOT2:sec:cf-gradings}]{LOT2}.

For idempotents $J$ of type $Y$, we again use
$Q_{J_X}*\sigma^{\epsilon} \in \pi_2(I_0, J)$ to define the grading
refinement data.
Then
\[
  \gr(\x(J))=\Xi(J)R(g'(\sigma))^\epsilon \Xi(J_X)^{-1}
    \cdot\langle R(g'(P_i))\rangle.
\]
We can furthermore choose the domains $Q_J$ for all $J$ of type~$X$ to
have equal multiplicities at $\sigma$, $\sigma_+$, and $\sigma_-$.
(This can be achieved by adding or subtracting the periodic domain
corresponding to~$C$.)  Then $\Xi(j_R)$, $\Xi((J_X)_R)$, $\Xi(j_L)$,
and $\Xi((J_X)_L)$ all
commute with $R(g'(\sigma))^\epsilon$, and we can write
\[
  \gr(\x(J))=R(g'(\sigma))^\epsilon \Xi(j_L)\Xi((J_X)_L)^{-1}
    \cdot\langle R(g'(P_i))\rangle.
\]

With respect to the refined grading, the grading of an algebra element
$a\in \Alg(\PMC)$ (say)
with $j\cdot a\cdot k=a$ (for primitive idempotents $j$ and $k$) is
given by
\[
\gr_\Xi(a)=\Xi(j)\grb(a)\Xi(k)^{-1}.
\]
See~\cite[Definition~\ref*{LOT2:def:refined-grading}]{LOT2}.
\glsit{$\gr_\Xi$}%

\subsection{The refined grading set for general surface homeomorphisms}\label{sec:general-grads}
Much of the discussion of gradings for arc-slides extends to arbitrary
surface homeomorphisms.  In particular, given a strongly-based mapping
class $\psi\co F(\PMC_1)\to F(\PMC_2)$ there is an associated outer
isomorphism of grading groups $\phi_\psi \co G(\PMC_1)\to
G(\PMC_2)$, gotten from the grading set of any Heegaard diagram for
$\psi$. These assemble to give an action of the mapping class
group by outer automorphisms on $G(\PMC)$, as follows:

\begin{proposition}\label{prop:MCG-acts-on-G}
  For any pointed matched circle $\PMC$, there is a
  homomorphism
  \[
  \phi\co\MCG_0(F(\PMC)) \to \Out(G(\PMC)),
  \]
  extending the standard action on the homology
  \[
  \MCG_0(F(\PMC)) \to \Sp(2k, \ZZ) \subset \Aut(H_1(F(\PMC))).
  \]
  
  More generally, let $\{\PMC^k_i\}$ be the set of all pointed
  matched circles for surfaces of genus~$k$. Let $\Out(\{G(\PMC^k_i)\})$
  be the groupoid with object set $\{\PMC^k_i\}$ and morphisms
  $\Hom(\PMC_i,\PMC_j)=\Out(G(\PMC_i),G(\PMC_j))$ (see
  Definition~\ref{def:outer-aut}). Then there is
  a homomorphism
  \[
  \phi\co\MCG_0(k)\to \Out(\{G(\PMC^k_i)\})
  \]
  \glsit{$\Out(\{G(\PMC^k_i)\})$}%
  extending the standard action on $H_1$.
\end{proposition}
\begin{proof}
  Let $\psi\in\MCG_0(F(\PMC),F(\PMC'))$. Then the mapping cylinder $Y_\psi$ of
  $\psi$ has grading set $S_\psi$
  \glsit{$S_\psi$}%
  which is isomorphic to $G(\PMC)$
  (respectively $G(\PMC')$) as a left
  $G(\PMC)$-set (respectively right $G(\PMC')$-set), which by
  Lemma~\ref{lem:G-torseur-is-map} gives a map $\phi_\psi \in
  \Out(G(\PMC), G(\PMC'))$.
  If $\psi'$ is an element of $\MCG_0(F(\PMC'),F(\PMC''))$
  then $S_{\psi'\circ\psi}=S_\psi\times_{G(\PMC)}S_{\psi'}$. Thus,
  Lemma~\ref{lem:G-torseur-is-map} applied to the $S_\psi$ gives a map
  $\MCG_0(k)\to \Out(\{G(\PMC^k_i)\})$.

  It is immediate from the form of $S_\psi$ that the action on
  $\Out(\{G(\PMC^k_i)\})$ projects to the standard action of the
  mapping class group on $H_1$.
\end{proof}

The main goal of this subsection is to compute $\phi_\psi$ for general
mapping classes $\psi$.

We start with some lemmas regarding the structure of $G(\PMC)$.

Recall from Section~\ref{sec:grading-set-refined} that the pointed
matched circle $\PMC$ gives a canonical basis $\{h(B_i)\}$ for
$H_1(F(\PMC))$, gotten by joining the cores of the handles attached
to $\PMC$ with corresponding segments in $\PMC$. These can be upgraded to a generating
set of $G(\PMC)$ as follows:
\begin{lemma}
  Let $\{B_i\}_{i=1}^{2k}$ denote the set of matched pairs for $\PMC$,
  and $\gpElt(B_i)=(-1/2,h(B_i))$ denote the corresponding group
  element of $G(\PMC)$.  Then $G(\PMC)$ is generated by the elements
  $\{\gpElt(B_i)\}_{i=1}^{2k}$ and $\lambda$, subject to the
  relations:
  \begin{itemize}
  \item $\lambda$ is central.
  \item $\gpElt(B_i)\cdot \gpElt(B_j) = \lambda^{2 h(B_i)\cdot h(B_j)}\gpElt(B_j)\cdot
    \gpElt(B_i)$, where $\cdot$ denotes the intersection product in $H_1(F(\PMC))$.
  \end{itemize}
\end{lemma}
\begin{proof}
  It is immediate from the definition of $\epsilon$ that 
  $\{\gpElt(B_i)\}_{i=1}^{2k}$ lie in $G(\PMC)$.
  Together with $\lambda$, the elements generate an index two subgroup of 
  $\OneHalf\ZZ\times H_1(F(\PMC))$ (with its twisted
  multiplication), which must be $G(\PMC)$. The relations follow from
  the description of
  $G(\PMC)$ as a central extension of
  $H_1(F(\PMC))$~\cite[Section~\ref*{LOT1:sec:refined-grading}]{LOT1}.%
\end{proof}

\begin{lemma}\label{lem:inner}
  For any pointed matched circle $\PMC$ and any matched pair $B_i$ in
  $\PMC$, the map $\psi\co
  G(\PMC)\to G(\PMC)$ given by 
  \begin{align*}
    \psi(\lambda)&=\lambda\\
    \psi(\gpElt(B_i))&=\lambda^{2}\gpElt(B_i)\\
    \psi(\gpElt(B_j))&=\gpElt(B_j)\qquad \text{if }j\neq i.
  \end{align*}
  is an inner automorphism.
\end{lemma}
\begin{proof}
  Consider first the case when $\PMC$ is the split pointed matched circle, for
  which the argument is more concrete. For each $B_i$ there is a
  corresponding $B_{\tau(i)}$ so that $h(B_i)\cdot h(B_{\tau(i)})=\pm
  1$ and $h(B_i)\cdot h(B_j)=0$ for $j\neq \tau(i)$, where $\cdot$
  denotes the intersection product on $H_1(F(\PMC))$. The inner
  automorphism of $G(\PMC)$ given by conjugating by
  $\gpElt(B_{\tau(i)})$ is given by
   \begin{align*}
     \gpElt(B_i)&\mapsto \lambda^{\mp2}\gpElt(B_i)\\
     \gpElt(B_j)&\mapsto \gpElt(B_j)\qquad\text{if }j\neq i,
   \end{align*}
   as desired.

   For general pointed matched circles, instead of conjugating by
   $\gpElt(B_{\sigma(i)})$, let $V\subset H_1(F(\PMC))$ denote the
   subspace generated by all $h(B_j)$ for $j\neq i$. Choose a
   primitive vector $v\in H_1(F(\PMC))$ symplectically orthogonal to
   $V$ (so $v \cdot h(B_j) = 0$ for $j \ne i$ and $v \cdot h(B_i) =
   \pm 1$), and conjugate by an element of $G$ of the form $(m,\pm v)$
   (for some $m\in\frac{1}{2}\ZZ$).
\end{proof}

\glsit{$\GmodTwo(\PMC)$}%
Let $\GmodTwo(\PMC)=(\ZZ/2)\times H_1(F(\PMC))$, which we view as
the trivial $\ZZ/2$-central extension of $H_1(F(\PMC))$. 
Since the cocycle defining the $\ZZ$-central extension $G(\PMC)$ of
$H_1(F(\PMC))$ takes
values in $2\ZZ\subset\ZZ$, there is a homomorphism $\Upsilon\co
G(\PMC)\to \GmodTwo(\PMC)$ extending the projection map $\ZZ\to \ZZ/2$.
Explicitly,
we may take $\Upsilon$ to be defined by $\Upsilon(\gpElt(B_i))=(0,h(B_i))$ and $\Upsilon(\lambda)=(1,0)$.
\glsit{$\Upsilon$}%

\begin{lemma}\label{lem:mu-defd-mod-2}
  Let $\PMC_1$ and $\PMC_2$ be pointed matched circles. Let $\phi\co
  G(\PMC_1)\to G(\PMC_2)$ be an isomorphism so that
  $\phi(\lambda)=\lambda$. Then:
  \begin{enumerate}
  \item\label{item:mu-defd-1} The isomorphism $\phi$ descends to an isomorphism $\overline{\phi}\co
    \GmodTwo(\PMC_1)\to \GmodTwo(\PMC_2)$.
  \item\label{item:mu-defd-2} As an outer isomorphism with $\phi(\lambda)=\lambda$, $\phi$
    is determined by the induced map~$\overline{\phi}$.
  \end{enumerate}
\end{lemma}
\begin{proof}
  For Part~(\ref{item:mu-defd-1}), since the kernel of $\Upsilon$ is the subgroup generated by $\lambda^2$, $\phi$ descends to a homomorphism $\overline{\phi}\co \GmodTwo(\PMC_1)\to \GmodTwo(\PMC_2)$. Since $\phi$ is an isomorphism and $\phi(\lambda)=\lambda$, the elements $[\phi(0,h(B_i)]$ form a basis for $H_1(F(\PMC_2))$. It follows that the map $\overline{\phi}$ is an isomorphism.

  Part~(\ref{item:mu-defd-2}) follows from Lemma~\ref{lem:inner}.
\end{proof}

\begin{remark}
  In particular, Lemma~\ref{lem:mu-defd-mod-2} gives an injective
  homomorphism $\Out_\lambda(G(\PMC))\to \Aut(\GmodTwo(\PMC))$,
  where $\Out_\lambda(G(\PMC))=\{\phi\in\Out(G(\PMC))\mid
  \phi(\lambda)=\lambda\}$. The image of this homomorphism is all
  automorphisms of $\GmodTwo$ which cover symplectomorphisms of
  $H_1(F(\PMC))$ (with respect to the intersection form). 
\end{remark}

\begin{definition}
  For any class $a\in H_1(F(\PMC))$, write $a=\sum_i a_i h(B_i)$, and
  define the \emph{$\ell^1$-norm} of $a$ to be
  \[
  \|a\|_1=\sum_{i=1}^{2k}|a_i|.
  \]
\end{definition}

Given a map $\psi\co F(\PMC_1)\to F(\PMC_2)$, define a map
$\xi_\psi\co \GmodTwo(\PMC_1)\to \GmodTwo(\PMC_2)$
\glsit{$\xi_\psi$}%
by
\begin{equation}
  \label{eq:DefOfXi}
\xi_\psi(m,a)=(m+\|a\|_1+\|\psi_*(a)\|_1,\psi_*(a)).
\end{equation}
\begin{lemma}
  The map $\xi_\psi$ is a group homomorphism.
\end{lemma}
\begin{proof}
  This is a direct computation:
  \begin{align*}
    \xi_\psi(m,a)\xi_\psi(m',a')&=(m+\|a\|_1+\|\psi_*(a)\|_1,\psi_*(a))(m'+\|a'\|_1+\|\psi_*(a')\|_1,\psi_*(a'))\\
    &=(m+m'+\|a\|_1+\|a'\|_1+\|\psi_*(a)\|_1+\|\psi_*(a')\|_1,\psi_*(a)+\psi_*(a'))\\
    &=(m+m'+\|a+a'\|_1+\|\psi_*(a+a')\|_1,\psi_*(a+a'))\\
    &=\xi_\psi(m+m',a+a')=\xi_\psi((m,a)(m',a'))
  \end{align*}
  (The third equality uses the fact that $\psi_*$ is linear, and the
  modulo 2 reduction of $\|\cdot\|_1$ is also linear.)
\end{proof}

\begin{proposition}\label{prop:mu-is-ell-1}
  Let $\psi\co F(\PMC_1)\to F(\PMC_2)$ be a strongly-based mapping
  class. Then the outer isomorphism $\phi_\psi\co G(\PMC_1)\to
  G(\PMC_2)$ is characterized by $\overline{\phi_\psi}=\xi_\psi$
  (as defined in Equation~\eqref{eq:DefOfXi}).
  \end{proposition}
\begin{proof}
  The proof is in two steps:
  \begin{enumerate}
  \item\label{item:mu-ell-1} Verify that if $\psi=\PunctF(m)$ is the
    diffeomorphism induced by an arc-slide then
    $\xi_\psi=\overline{\phi_\psi}$.
  \item\label{item:mu-ell-3} Verify that if $\psi_1\co F(\PMC_1)\to
    F(\PMC_2)$ and $\psi_2\co F(\PMC_2)\to
    F(\PMC_3)$ then
    $\xi_{\psi_1}=\overline{\phi_{\psi_1}}$ and
    $\xi_{\psi_2}=\overline{\phi_{\psi_2}}$ imply
    $\xi_{\psi_2\circ\psi_1}=\overline{\phi_{\psi_2\circ\psi_1}}$.
  \end{enumerate}
  Since the arc-slides generate the strongly-based mapping class
  group, it follows from steps~\ref{item:mu-ell-1}
  and~\ref{item:mu-ell-3} that $\overline{\phi_\psi}=\xi_\psi$ for any
  mapping class $\psi$.

  We start with part~(\ref{item:mu-ell-1}).  Since both
  $\overline{\phi_\psi}$ and $\xi_\psi$ are group homomorphisms, it
  suffices to verify that
  $\overline{\phi_\psi}(\lambda)=\xi_\psi(\lambda)$ and
  $\overline{\phi_\psi}(\gpElt(B_i))=\xi_\psi(\gpElt(B_i))$. The
  statement about $\lambda$ holds trivially. For $i\neq 2k$,
  \begin{align*}
    \phi_\psi(\gpElt(B_i))&=\othgpElt(B_i)\\
    \overline{\phi_\psi}(\othgpElt(B_i))&=(0,\othh(B_i))\\
    \xi_\psi(\gpElt(B_i))&=(0,\othh(B_i))
  \end{align*}
  as desired. For $i=2k$, it follows from the formulas in
  Lemma~\ref{lem:gr-map-is} that
  \[
  \overline{\phi_\psi}(\gpElt(B_i))=(1,\psi_*(h(B_i)))
  \]
  where $\psi_*(h(B_i))$ has $\ell^1$-norm $2$. Similarly,
  \[
  \xi_\psi(\gpElt(B_i))=(0+1+2,\psi_*(h(B_i)))=(1,\psi_*(h(B_i)))=\overline{\phi_\psi}(\gpElt(B_i)),
  \]
  as desired.

  Recall that the cases (U.I)--(O.III) are only half of the combinatorial
  configurations for arc-slides; the other half are their
  inverses. But $\phi_{\psi^{-1}}=\phi_\psi^{-1}$ and
  $\xi_{\psi^{-1}}=\xi_{\psi}^{-1}$, so it follows that $\overline{\phi_\psi}$
  and $\xi_\psi$ agree on the remaining six types of arc-slides, as well.

  For part~\ref{item:mu-ell-3}, since we already know (by
  Proposition~\ref{prop:MCG-acts-on-G}) that
  $\phi_\bullet$ is a groupoid homomorphism $\MCG_0(k)\to
  \Out\{G(\PMC)\}$, and so $\overline{\phi_\bullet}$ is a groupoid
  homomorphism as well, it suffices to verify that $\xi_\bullet$ is a
  groupoid homomorphism. That is, we must check that
  \begin{align*}
    \xi_{\psi_2\circ\psi_1}&=\xi_{\psi_2}\circ\xi_{\psi_1}\\
    \xi_{\psi^{-1}}&=\xi_\psi^{-1}.
  \end{align*}
  The second property is obvious (and we already used it, above). For
  the first property,
  \begin{align*}
    \xi_{\psi_2}\circ\xi_{\psi_1}(m,a)&=\xi_{\psi_2}(m+\|a\|_1+\|(\psi_1)_*(a)\|_1,(\psi_1)_*(a))\\
    &=(m+\|a\|_1+\|(\psi_1)_*(a)\|_1+\|(\psi_1)_*(a)\|_1\\
    &\qquad\qquad+\|(\psi_2)_*\circ(\psi_1)_*(a)\|_1,(\psi_2)_*\circ(\psi_1)_*(a))\\
    &=(m+\|a\|_1+\|(\psi_2\circ\psi_1)_*(a)\|_1,(\psi_2\circ\psi_1)_*(a))\\
    &=\xi_{\psi_2\circ\psi_1}(m,a),
  \end{align*}
  as desired.
\end{proof}

\begin{remark}
  Since the $G'$-set grading $\CFDDa(\psi)$ is not transitive for the
  left and right actions, it does not correspond to a map $G'(\PMC)\to
  G'(\PMC')$. Note also that the $G$-set grading induces the $G'$-set
  grading (cf.~\cite[Lemma~\ref*{LOT2:lem:RefineGrading}]{LOT2}), so
  Proposition~\ref{prop:mu-is-ell-1} determines the $G'$-set which grades
  $\CFDDa(\psi)$.
\end{remark}


\section{Assembling the pieces to compute \textalt{$\HFa$}{HF\textasciicircum}.}
\label{sec:CompleteProof}
Using the computations from Section~\ref{sec:Arc-Slides} and the pairing theorems (Theorems~\ref{thm:PairingTheorem}
and~\ref{thm:PairingTheoremModBimod}),
 we finish the proofs of Theorems~\ref{thm:Handlebodies}
and~\ref{thm:HFa}. Then, in Section~\ref{sec:complete:gradings}, we use the computations of gradings from Section~\ref{sec:mcg-grading} to compute the decomposition of $\HFa(Y)$ according to $\SpinC$-structures, and the relative Maslov grading inside each $\SpinC$-structure.

First, we turn to the calculation of $\CFDa$ for a handlebody. We
start with the standard ``$0$-framed'' handlebody:
\begin{proof}[Proof of Proposition~\ref{prop:CFD-of-0-fr-hdlbdy}]
Recall from Section~\ref{subsec:Handlebodies} that $\HB^g$ denotes the
$0$-framed handlebody of genus $g$, the boundary sum of $g$ copies of
the $0$-framed solid torus.
A bordered Heegaard diagram $\HD^g_0$
\glsit{$\HD^g_0$}%
for $\HB^g$ can be constructed as the
boundary sum of $g$ bordered Heegaard diagrams for the $0$-framed
solid torus. We can draw this diagram $\HD_0^g$ on a genus $g$ surface with $2g$
alpha-arcs
$\{\alpha_1^c,\alpha_2^c,\dots,\alpha^c_{2g-1},\alpha^c_{2g}\}$ and %
$g$ beta-circles, so that $\beta_i$ meets $\alpha_{2i}^c$
transversely in a single point (and is disjoint from all other
$\alpha$-arcs).
The $g=2$ case $\HD^2_0$ is illustrated in
Figure~\ref{fig:GenusTwoBorderedDiagram}.

The type $D$ module $\CFDa(\HD^g_0)$ has a single generator $\x$. Its associated
idempotent $I_{D}(\x)$ is the idempotent where all the
odd-numbered strands are occupied.  The differential is specified by
$$\partial \x =
\left(\sum_{i=1}^g a(\rho_{ 4i-3})\cdot a(\rho_{4i-2})\right)\cdot
\x$$
where $\rho_j$ denotes the $j\th$ short chord in $-\bdy\HD_0^g$.
To see this, note that there are $g$ connected components of
$\Sigma\setminus(\alphas\cup\betas)$, each of which contributes one
term in the sum. Although these disk do not quite fit the conditions
of Definition~\ref{def:PolygonConnects}, it is easy to see that they
each contribute as indicated. Moreover, there are no other
contributions to the differential. 

The module $\CFDa(\HD^g_0)$ we have just described is exactly the
module $\Dmod(\HB^g)$ of Definition~\ref{subsec:Handlebodies}.  
\end{proof}

\begin{proof}[Proof of Theorem~\ref{thm:Handlebodies}]
  This follows immediately from
  Proposition~\ref{prop:CFD-of-0-fr-hdlbdy} and
  Theorems~\ref{thm:DDforIdentity} and~\ref{thm:DDforArc-Slides}, by
  inductively applying the pairing theorem
  Corollary~\ref{cor:PairingTheoremReparam} and the relation
  \begin{align*}
    \Mor_{D,C}(\lsub{D,C}X,\lsub{D}M\otimes_{\Field}&
    \Mor_{B,A}(\lsub{B,A}Y,\lsub{C,B}N\otimes_{\Field}
    \lsub{A}P))\\
    &\simeq 
    \Mor_{D,C,B,A}(\lsub{D,C}X\otimes_{\Field} \lsub{B,A}Y,
    \lsub{D}M\otimes_{\Field} \lsub{C,B}N\otimes_{\Field}
    \lsub{A}P).
    \qedhere
  \end{align*}
\end{proof}

\begin{proof}[Proof of Theorem~\ref{thm:HFa}]
  This follows immediately from Theorem~\ref{thm:Handlebodies} and
  Theorem~\ref{thm:PairingTheorem}.
\end{proof}

\subsection{Gradings}\label{sec:complete:gradings}
The one missing ingredient to compute the $\SpinC$ and Maslov gradings
on $\HFa(Y)$ is the computation of the gradings on
$\CFDa(\HD_0^g)$. Number the points in $\PMC_0^g=-\bdy\HD_0^g$ by
$1,\dots,4g$, from bottom to top. Write elements of $G'(\PMC)$ as
pairs $(m;a)$
\glsit{$(m;a)$}%
where $m\in\frac{1}{2}\ZZ$ and $a$ is a linear combination of
intervals $[i,j]$ in $\PMC$, subject to the congruence condition
$m\equiv \epsilon(a) \pmod 1$.

\begin{proposition}
  \label{prop:CFD-of-0-fr-hdlbdy-graded}
  As a $G'(\PMC_0^g)$-set graded module, $\CFDa(\HD_0^g)$ is graded by 
  $S'_0=G'(\PMC_0^g)/H$
  \glsit{$S'_0$}%
  where
  \[
  H=\langle \{(-1/2;-[4i+1,4i+3])\}_{i=0}^{g-1}\rangle.
  \]
  The generator $\x$ has grading $(0;0)$ in $G'(\PMC)$.
\end{proposition}
\begin{proof}
  Each generator of $H$ corresponds to a periodic domain $P$ consisting of a single component of $\Sigma\setminus(\alphas\cup\betas)$, with multiplicity $1$. Each of these domains $P$ has Euler measure $-1/2$ and point measure $1$.
\end{proof}

Now, to compute the gradings on $\CFa(Y)$, one follows a simple, seven-step process:
\begin{enumerate}
\item 
  Take a Heegaard decomposition $Y=\HB_0^g\cup_\psi \HB_0^g$ and
  factoring the gluing map into arc-slides, $\psi=\psi_1\circ\cdots\circ\psi_n$, where $\psi_i\co -F(\PMC_i)\to -F(\PMC_{i-1})$. By
  Theorem~\ref{thm:HFa}, $\CFa(Y)$ is homotopy equivalent to the
  complex
  \begin{equation}\label{eq:complete:mor-cx}
  \Mor\big(\DDmod(\Id_{\PMC_n})\otimes\cdots\otimes\DDmod(\Id_{\PMC_1}),\Dmod(\HB_0^g)\otimes\DDmod(\psi_n)\otimes\cdots\otimes\DDmod(\psi_1)\otimes\Dmod(\HB_0^g)\big),
  \end{equation}
  where the $\Mor$-complex is over
  $\Alg(\PMC_n)\otimes\cdots\otimes\Alg(\PMC_1)$.
\item Choose grading refinement data $\Xi_i$ for $\Alg(\PMC_i)$, as
  in~\cite[Section~\ref*{LOT2:subsec:SmallGroup}]{LOT2}. Each homogeneous element $a=j\cdot
  a\cdot j'\in \Alg(\PMC_i)$ (where $j$ and $j'$ are primitive
  idempotents) gets a grading $\gr_{\Xi}(a)=\Xi(j)^{-1}\gr'(a)\Xi(j')$.
  (Sometimes we will denote $\gr_\Xi(a)$ simply by $\gr(a)$.)

 (In Section~\ref{sec:grading-set-refined}, we explain how to use $\psi_i$ to choose convenient grading refinement data for $\PMC_i$ and $\PMC_{i-1}$. But the data chosen this way by $\psi_i$ and $\psi_{i+1}$ for $\PMC_i$ may not agree.)
\item As discussed in Section~\ref{sec:gr-of-gens}, each module $\DDmod(\psi_i)\simeq \CFDDa(\PunctF(\psi_i))$ is graded by a set $S'_i$ with a left action by $G'(\PMC_i)\times G'(-\PMC_{i-1})$. (These sets are computed in Section~\ref{sec:grading-set-refined}, and the gradings of generators are computed in Section~\ref{sec:gr-of-gens}.)

  Consequently, 
  \[ R=\Dmod(\HB_0^g)\otimes\DDmod(\psi_n)\otimes\cdots\otimes\DDmod(\psi_1)\otimes\Dmod(\HB_0^g)
  \]
  is graded by the $G'(-\PMC_n)\times G'(\PMC_n)\times G'(-\PMC_{n-1})\times\cdots\times G'(-\PMC_1)\times G'(\PMC_1)$-set 
  \[
  S'_R=S'_0\times S'_{n}\times\cdots\times S'_{1}\times S'_0.
  \]
  (Here, $\times$ means $\times_\ZZ$; we suppress the subscript in this section to keep the notation manageable.)
\item Using the grading refinement data $\Xi_i$, we can turn $S'_R$
  into a $G(-\PMC_n)\times G(\PMC_n)\times
  G(-\PMC_{n-1})\times\cdots\times G(-\PMC_1)\times G(\PMC_1)$-set
  $S_R$, as in~\cite[Lemma~\ref*{LOT2:lem:RefineGrading}]{LOT2}. In this set, the grading of a
  generator $y_0\otimes y_1\otimes\cdots\otimes y_n\otimes y_{n+1}$,
  where $j_i\cdot y_i\cdot j_{i+1}=y_i$ for primitive idempotents $j_1,\dots,j_n$, is given by 
  \begin{multline*}
    \gr(y_0\otimes y_1\otimes\cdots\otimes y_n\otimes
    y_{n+1})=\gr'(y_0)\Xi_1(j_1)\otimes
    \Xi_1(j_1)^{-1}\gr'(y_1)\Xi_2(j_2)\otimes\\
    \cdots\otimes
    \Xi_{n-1}(j_{n-1})^{-1}\gr'(y_n)\Xi_n(j_n)\otimes\Xi_n(j_n)^{-1}\gr'(y_{n+1}).
  \end{multline*}
\item The module $L=\DDmod(\Id_{\PMC_n})\otimes\cdots\otimes\DDmod(\Id_{\PMC_1})$ is graded by the $G(-\PMC_n)\times G(\PMC_n)\times G(-\PMC_{n-1})\times\cdots\times G(-\PMC_1)\times G(\PMC_1)$-set 
  \[
  S_L=G(\PMC_n)\times G(\PMC_{n-1})\times\cdots\times G(\PMC_1).
  \]
  Each complementary idempotent $I$ generating $\CFDDa(\Id_{\PMC_i})$
  has grading $(0;0)$. (This argument is a simpler version of the
  discussion in Section~\ref{sec:comp-gr}, and is also given
  in~\cite[Lemma~\ref*{LOT2:lem:GradingIdentityDD}]{LOT2}.)
\item Finally, as discussed in Section~\ref{sec:pairing-theorem}, the $\Mor$-complex~(\ref{eq:complete:mor-cx})  is graded by 
  \[
  S=S_L^{\op}\times_{G(\PMC_n)\times G(\PMC_n)^{\op}\times G(\PMC_{n-1})\times\cdots\times G(\PMC_1)^\op} S_R.
  \]
  The $\Mor$-complex is generated, as an $\Field$-vector space, by maps of the form
  \[
  f=x_1\otimes\dots\otimes x_n\mapsto y_0\otimes a_1\otimes y_1\otimes\cdots\otimes y_n\otimes a_n\otimes y_{n+1}
  \]
  where $x_1,\dots,x_n$ are generators of
  $\DDmod(\Id_{\PMC_n}),\dots,\DDmod(\Id_{\PMC_1})$; $y_0,y_n$ are generators of $\Dmod(\HB_0^g)$; $y_i$ is a generator for $\DDmod(\psi_{n+1-i})$ if $0<i<n+1$; and $a_i\in \Alg(\PMC_{n+1-i})$ is a basic generator (strands diagram). The grading of such a map $f$ is
  \begin{align*}
    \gr(f)&=(\gr(x_1)^{\op}\times\dots\times\gr(x_n)^{\op})\times(\gr(y_0)\gr(a_1)\times\gr(a'_1)\gr(y_1)\gr(a_2)\\
    &\qquad\qquad\times\dots\times\gr(a'_{n+1})
    \gr(y_{n+1}))\\
    &=\gr(y_0)\gr(a_1)\times\gr(x_1)^{\op}\times\gr(a'_1)\gr(y_1)\gr(a_2)\\
    &\qquad\qquad\times\dots\times\gr(x_n)^{\op}\times\gr(a'_{n+1})\gr(y_{n+1}))\\
    &=\gr(y_0)\gr(a_1)\times\gr(a'_1)\gr(y_1)\gr(a_2)\times\dots\times\gr(a'_{n+1})\gr(y_{n+1}))
  \end{align*}
  where the last equality uses the fact that each $x_i$ has grading
  $(0;0)\in G(\PMC_i)$. All $\times$'s are $\times_{G(\PMC_i)}$'s, for
  appropriate $i$'s.
\item The set $S$ retains an action by the central element $\lambda$.
  By the graded version of the pairing theorem, each $\SpinC$-structure on $\CFa(Y)$ corresponds to a $\lambda$-orbit. If $f$ and $g$ lie in the same $\SpinC$-structure then $\gr(f)=\lambda^i\gr(g)$, for some $i$, and the grading difference between $f$ and $g$ is $i$. (If $\lambda$ acts non-freely on that orbit then this grading difference is only well-defined modulo $\min\{n\mid \lambda^n\gr(f)=\gr(f)\}$; this $n$ is the divisibility of the first Chern class of the corresponding $\SpinC$-structure.)
\end{enumerate}

For practical computations, it is generally not necessary to refine the gradings. Working with the larger grading groups $G'$, let $S'_L=G'(\PMC_n)\times G'(\PMC_{n-1})\times\cdots\times G'(\PMC_1)$, a set with a left action of $G'(\PMC_n)\times G'(\PMC_n)^\op\times\cdots\times G'(\PMC_1)\times G'(\PMC_1)^\op$. Then, without worrying about grading refinements, elements of the $\Mor$ complex~(\ref{eq:complete:mor-cx}) are graded by
\[
S'=(S'_L)^{\op}\times_{G'(\PMC_n)\times G'(\PMC_n)^{\op}\times G'(\PMC_{n-1})\times\cdots\times G'(\PMC_1)^\op} S'_R.
\]
The grading of a generator $f$ is given by 
\[
\gr'(f)=(\gr'(x_1)^{\op}\times\dots\times\gr'(x_n)^{\op})\times(\gr'(y_0)\times\gr'(a_1)\times\dots\times\gr'(y_{n+1})).
\]
Two generators $f$ and $g$ lie in the same $\lambda$-orbit of $S'$ if and only if they represent the same $\SpinC$-structure. In this case, $\gr'(f)=\lambda^i\gr'(g)$ where $i$ is the grading difference between $f$ and $g$. The only part that breaks is that with respect to the unrefined grading there will be $\lambda$-orbits with no generators, which do not correspond to $\SpinC$-structures on $Y$.


\section{Elementary cobordisms and bordered invariants}
\label{sec:ElementaryCobordisms}

Although we have focused so far on computing $\HFa$ for a closed
three-manifold, the techniques of this paper can also be used to
calculate $\CFDa$ and $\CFDDa$ for bordered three-manifolds, as
well. (One can then compute the type $A$-, \DA- and \AAm-invariants,
too; see Section~\ref{sec:AModules}.)

Let $Y$ be a bordered three-manifold with boundary $\partial Y$
parameterized by $F(\PMC)$. The goal is to break up $Y$ into basic
pieces, calculate the bimodules associated to those pieces, and then
calculate $\CFDa(Y)$ by composing the individual bimodules. In
addition to arc-slides, we need one new kind of basic piece: an
elementary cobordism. Since we already know how to change the
parametrizations of the boundary by using arc-slides, we will only
need elementary cobordisms with a particular bordering. In
Section~\ref{sec:split-el-cob} we compute the invariants of these
elementary cobordisms. In Section~\ref{sec:compute-CFD} we discuss how to compute the invariant
$\CFDa$ of a bordered $3$-manifold with connected boundary. Because of
the particular form of the pairing theorem we have been using, there
are mild technical complications for the invariant $\CFDDa$ of a
$3$-manifold with two boundary components; we discuss (and overcome)
these in Section~\ref{sec:compute-CFDD}.

\subsection{The invariant of a split elementary cobordism}\label{sec:split-el-cob}

\begin{definition}
  Let $F_1$ be a connected surface of genus $g$ and $F_2$ be a surface
  of genus $g+1$.  Let $M_1=M_1(g)$
  \glsit{$M_1=M_1(g)$}%
  be the cobordism from $F_1$ to $F_2$
  gotten by adding a three-dimensional one-handle to $[0,1]\times F_1$
  at a pair of points in $\{1\}\times F_2$. Dually, let $M_2=M_2(g)$
  \glsit{$M_2=M_2(g)$}%
  be the
  cobordism from $F_2$ to $F_1$ gotten by adding a three-dimensional
  two-handle to $[0,1]\times F_2$ along a non-separating curve in
  $\{1\}\times F_2$.  The three-manifolds $M_1(g)$ and $M_2(g)$ are called
  {\em elementary cobordisms}.
  \index{elementary cobordism}%
\end{definition}

Up to diffeomorphism, an elementary cobordism is uniquely determined
by the genera of its boundary. Also, $M_1(g)$ is simply the
orientation-reverse of $M_2(g+1)$.

Next, we compute the invariant of an elementary cobordism with a
particular bordering.  Recall from Section~\ref{subsec:Handlebodies} that
$\PMC_1$ denotes the (unique) pointed matched circle for a genus $1$
surface. Given a pointed matched circle $\PMC$ for a
surface of genus $g$, the \emph{split bordering by $\PMC$}
\index{split bordering of elementary cobordism}%
\index{elementary cobordism!split bordering of}%
of $M_1(g)$ is the bordering by
$F(\PMC)$ and $F(\PMC\#\PMC_1)$ so that the circle in $F(\PMC_1)$
specified by $\{a,a'\}$ bounds a disk in $M_1(g)$. The \emph{split
  bordering by $\PMC$} of $M_2(g)$ is the orientation-reverse of the
split bordering of $M_1(g-1)$ by $-\PMC$.

Let $\PMC$ be a pointed matched circle, and let $\PMC^1$ denote the
pointed matched circle for a genus one surface. Let $\HB^1$ denote the
$0$-framed solid torus, as in Section~\ref{sec:0-fr-hb}. There is an
inclusion map
$$i\co \Alg(\PMC)\otimes\Alg(\PMC^1) \rightarrow \Alg(\PMC\#\PMC^1).$$
Now, $\Dmod(\Id_{\PMC})\otimes \Dmod(\HB^1)$ can be viewed as a type $D$
structure over $\Alg(-\PMC)\otimes\Alg(\PMC)\otimes\Alg(\PMC^1)$.
Promoting the homomorphism $i$ by the identity map on $\Alg(-\PMC)$, we obtain
a homomorphism
$$j\co\Alg(-\PMC)\otimes\Alg(\PMC)\otimes \Alg(\PMC^1)\rightarrow
\Alg(-\PMC)\otimes\Alg(\PMC\#\PMC^1).$$

\begin{proposition}
  \label{prop:CalculateElementaryCobordism}
  The induced module $j_*(\Dmod(\Id_{\PMC})\otimes \Dmod(\HB^1)))$,
  which is a (left-right) type \DD\ bimodule over 
  $\Alg(\PMC\#\PMC^1)$ and $\Alg(-\PMC)$, is isomorphic to the bimodule $\CFDDa$ of the
  elementary cobordism $M_1$, endowed with the split bordering by
  $\PMC$.
  Similarly, the induced bimodule $j_*(\Dmod(\Id_{-\PMC})\otimes
  \Dmod(-\HB^1))$, which is a type \DD\
  bimodule over $\Alg(\PMC)$ and $\Alg(\PMC\#\PMC^1)$, is isomorphic to
  the bimodule $\CFDDa$ of the elementary cobordism $M_2$,
  endowed with the split bordering by $\PMC$.
\end{proposition}

\begin{figure}
  \begin{center}
    \input{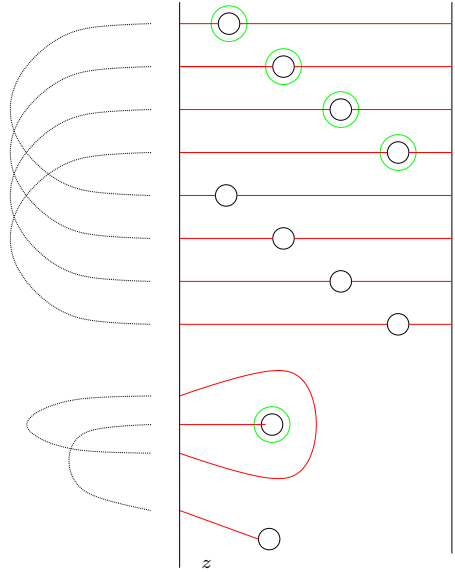}
  \end{center}
  \caption {{\bf An elementary cobordism.}
    \label{fig:ElementaryCobordism}
    A diagram for an elementary cobordism from a genus three surface
    to a genus two surface. The genus two surface here has the
    antipodal matching.}
\end{figure}

\begin{proof}
  Draw a Heegaard diagram for $(F(\PMC)\times[0,1])\#_\bdy
  \HB^1$ by forming the boundary
  sum of a Heegaard diagram $\HD(\Id_\PMC)$ for $\Id_\PMC$
  with the Heegaard diagram $\HD^1_0$ for the $0$-framed genus one
  handlebody, where the
  boundary sum is taken near $\mathbf{z}\cap (\bdy_L\HD(\Id_\PMC))\in
  \bdy \HD(\Id_\PMC)$ and $z\in \bdy \HD^1_0$.  See
  Figure~\ref{fig:ElementaryCobordism}.
  Holomorphic curves which contribute to the differential must have
  connected domains. For this reason, the holomorphic curves in
  the differential are either supported in the identity \DD\ bimodule
  region, or in the handlebody region. This implies the result for
  $M_1$. The result for $M_2$ is obtained by replacing $\PMC$ with
  $-\PMC$ and reflecting the picture horizontally. 
\end{proof}

\begin{remark}
  The above proof can be seen as a special case of a boundary
  connected sum formula for $\CFDa$; see~\cite{Zarev09:BorSut} for further 
  generalizations.
\end{remark}

\subsection{Computing \textalt{$\CFDa$}{CFD\textasciicircum} of a \textalt{$3$}{3}-manifold with
  connected boundary}\label{sec:compute-CFD}
By standard Morse theory, any connected three-manifold $Y$ with
connected boundary $F(\PMC)$ can be obtained from the $3$-ball by a
sequence of elementary cobordisms. That is,
\[
Y=\bD^3\cup_{\phi_1}M_{i_1}(g_1)\cup_{\phi_2}M_{i_2}(g_2)\cup_{\phi_3}\cdots\cup_{\phi_n}M_{i_n}(g_n),
\]
where each $i_j\in\{1,2\}$ and the genera $g_i$ are determined by the
sequence of $i_j$'s in the obvious way (i.e.,
$g_i=\sum_{j<i}2(3/2-i_j)$). 

Let $Y_k$ be the part of $Y$ obtained after attaching $k$ of the
elementary cobordisms, i.e.,
\[
Y_k=\bD^3\cup_{\phi_1}M_{i_1}(g_1)\cup_{\phi_2}M_{i_2}(g_2)\cup_{\phi_3}\cdots\cup_{\phi_n}M_{i_k}(g_k).
\]
We compute $\CFDa(Y_k)$ inductively as follows.  The gluing map
$\phi_{k}$ is a map $\bdy_L M_{i_k}(g_k)\to -\bdy Y_{k-1}$. Suppose we
are given $\CFDa(Y_{k-1})$ for some bordering $-F(\PMC^{g_{k-1}})\to
\bdy Y_{k-1}$, and want to compute $\CFDa(Y_k)$ with respect to
a bordering $-F(\PMC^{g_{k}})\to \bdy Y_k$. (Here, $\PMC^{g_{k-1}}$
and $\PMC^{g_k}$ are arbitrary pointed matched
circle for a surface of the right genus.) 

Choose the split bordering of $M_{i_k}(g_k)$ by $\PMC^{g_{k-1}}$
(say). Then $\phi$ (respectively $\phi'$) corresponds to some map
$\psi\co F(\PMC^{g_{k-1}})\to F(\PMC^{g_{k-1}})$ (respectively
$\psi'\co F(\PMC^{g_k})\to F(\PMC^{g_k,\prime})$). Factor $\psi$ and
$\psi'$ into arc-slides,
\[
\psi=m_1\circ\cdots\circ m_l\qquad\qquad \psi'=m'_1\circ\cdots\circ m'_{l'}.
\]
Then, by Theorem~\ref{thm:PairingTheoremModBimod},
\begin{equation}
  \begin{split}
    \CFDa(Y_k)\simeq\Mor(\CFDDa(\Id)\otimes\cdots\otimes\CFDDa(\Id),\CFDa(Y_{k-1})\otimes
    \CFDDa(m_l)\otimes\cdots\\\otimes\CFDDa(m_1) \otimes
    \CFDDa(M_{i_k}(g_k)) \otimes
    \CFDDa(m'_{l'})\otimes\cdots\otimes\CFDDa(m'_1)).
  \end{split}
\label{eq:compute-CFD}
\end{equation}
(As in Section~\ref{subsec:Handlebodies}, $\Mor$ denotes the chain
complex of bimodule homomorphisms. Also, $\CFDDa(M_{i_k}(g_k))$
denotes the bimodule computed with respect to the split bordering by
$\PMC^{g_{k-1}}$, the $\CFDDa(\Id)$'s in the formula are with respect
to the appropriate pointed matched circles.)

By now, we know how to compute all of the pieces of
Equation~\eqref{eq:compute-CFD}: the bimodule $\CFDDa(M_{i_k}(g_k))$
is computed in Proposition~\ref{prop:CalculateElementaryCobordism};
the bimodules $\CFDDa(m_i)$ are computed in
Theorem~\ref{thm:DDforArc-Slides}; and the bimodules $\CFDDa(\Id)$ are
computed in Theorem~\ref{thm:DDforIdentity}.

\begin{remark}
  In the inductive computation, the first step formally uses bordered
  Floer homology of a manifold with boundary $S^2$. In this degenerate
  case, the definitions from bordered Floer homology give
  $\Alg(\PMC)=\Field$, and the module $\CFDa(Y)$ is just
  $\CFa(Y\cup_\bdy \bD^3)$. In particular, $\CFDa(\bD^3)=\Field$.
\end{remark}

\begin{remark}
  Another way to compute $\CFDa(Y)$ is to decompose $Y$ as the union
  of a handlebody and a compression body. Computing the invariant of a
  compression body with a standard bordering is a simple extension of
  Proposition~\ref{prop:CalculateElementaryCobordism} (or, in fact,
  follows from that proposition). Section~\ref{subsec:Handlebodies}
  explains how to compute the invariant of a handlebody with an
  arbitrary framing. Theorem~\ref{thm:PairingTheoremModBimod} then
  says how to compute $\CFDa(Y)$.
\end{remark}

\subsection{Computing \textalt{$\CFDDa$}{CFDD\textasciicircum} of a \textalt{$3$}{3}-manifold with
  two boundary components}\label{sec:compute-CFDD}
Recall that associated to a strongly bordered $3$-manifold $Y$ with
two boundary components is a bordered $3$-manifold $\drY$
\glsit{$\drY$}%
with connected boundary, obtained by deleting a neighborhood of the framed
arc from $Y$. Suppose $Y$ is bordered by $-F(\PMC_1)$ and
$-F(\PMC_2)$. Then $\drY$ is bordered by $-F(\PMC_1\#\PMC_2)$. There is
a projection map $p\co \Alg(\PMC_1\#\PMC_2)\to \Alg(\PMC_1)\otimes
\Alg(\PMC_2)$
\glsit{$p$}%
which sets to zero any algebra element crossing between
$\PMC_1$ and $\PMC_2$. Then $\CFDDa(Y)$ is defined to be the induced
module $p_*(\CFDa(\drY))$.

Section~\ref{sec:compute-CFD} explains how to compute
$\CFDa(\drY)$. Thus, we now know how to compute $\CFDDa(Y)$.

\begin{remark}
  It is possible to give a more direct computation of $\CFDDa(Y)$,
  without using $\drY$. The version of the pairing theorem we have
  used so far is inconvenient for this. Suppose $Y_1$ and $Y_2$ are
  strongly bordered $3$-manifolds with two boundary components, and
  $\bdy_R(Y_1)=F(\PMC)=-\bdy_L(Y_2)$. Then 
  \[
  \Mor(\lsub{\Alg(\PMC)}\CFDDa(\Id_{\PMC})_{\Alg(\PMC)},\CFDDa(Y_1)\otimes\CFDDa(Y_2))
  \]
  is not $\CFDDa(Y_1\cup_\bdy Y_2)$, but rather
  $\CFDDa(\tau_{\bdy}^{-1}(Y_1\cup_\bdy Y_2))$, the result of gluing
  $Y_1$ to $Y_2$ and then decreasing the framing on the arc by one.

  To remedy this, one could use any of several other variants of the
  pairing theorem. One way forward will be apparent in
  Section~\ref{sec:AModules}.
\end{remark}


\section{Computing (with) type \textalt{$A$}{A} invariants}
\label{sec:AModules}
In~\cite{LOT1}, we consider two types of modules,
$\CFDa(Y)$ and $\CFAa(Y)$.  There are analogues $\CFDDa(Y)$,
$\CFDAa(Y)$ and $\CFAAa(Y)$ in the two boundary component case \cite{LOT2}. Until now, we
have focused exclusively on $\CFDa(Y)$ and $\CFDDa(Y)$.  There are several reasons for doing
this. Type $D$ modules are easier to compute, as they count fewer
holomorphic curves. They are algebraically simpler to describe,
because they are ordinary differential modules, rather than $\Ainf$
modules. Moreover, thanks to duality results from~\cite{LOT2}, it is
possible to formulate the theory purely in terms of $\CFDa$ and
$\CFDDa$; this formulation serves to shorten the exposition.

However, in some contexts it is useful to think of type $A$
modules. To this end, we recall how to extract the type $A$ modules, and type \DA\ and \AAm\
bimodules, from the type $D$ modules and \DD\ bimodules via the
duality
result,~\cite[Proposition~\ref*{LOT2:prop:DDAA-duality}]{LOT2}. (Alternatively,
one
could use any of several results from~\cite{LOTHomPair}.) This is
discussed in Section~\ref{sec:computing-type-a}.

One advantage of working with type $A$ modules is that one can work
with chain complexes with fewer generators; the cost is more
complicated algebra actions. This is a special case of a class of
results called homological perturbation theory. We review the relevant
instance in Section~\ref{sec:hpt}. Using this, we discuss
reconstructing the closed invariant using type $A$ modules in
Section~\ref{sec:reconstruct-via-A}. As an example of these
techniques, we compute a small (in fact, minimal) model for the type
\AAm\ module for the identity map of the torus in
Section~\ref{sec:AA-torus-id} (where here size is measured by the rank
of the underlying vector space). In Sections~\ref{sec:ok-computer}
and~\ref{sec:find-slides} we discuss some computer computations.

\subsection{Computing type \textalt{$A$}{A} invariants}\label{sec:computing-type-a}
The key step in computing the type $A$ modules and bimodules from the
type $D$ modules is understanding the type \AAm\ module for the
identity map. Fix a pointed matched circle $\PMC$. We proved
in~\cite[Proposition~\ref*{LOT2:prop:DDAA-duality}]{LOT2} that
\begin{equation}
  \label{eq:DefAA}
  \CFAAa(\Id_{\PMC})=
  \Mor_{\Alg(-\PMC)}(\lsub{\Alg(-\PMC),
  \Alg(\PMC)}\CFDDa(\Id_{\PMC}), \Alg(-\PMC)).
\end{equation}
\glsit{$\CFAAa(\Id_{\PMC})$}%
Here, $\Mor_{\Alg(-\PMC)}$ denotes the chain complex of
$\Alg(-\PMC)$-module maps between the bimodules.  This complex retains
commuting right actions by $\Alg(\PMC)$ and
$\Alg(-\PMC)$. Theorem~\ref{thm:DDforIdentity} calculates
$\CFDDa(\Id_{\PMC})$; this description can be combined with
Equation~\eqref{eq:DefAA} to give an explicit description of $\CFAAa(\Id)$.

For the purposes of this paper, the reader unfamiliar with the
definition of $\CFAAa$ in terms of holomorphic curves can safely take
Equation~\eqref{eq:DefAA} (as well as
Equations~\eqref{eq:comp-CFA},~\eqref{eq:comp-CFDA}
and~\eqref{eq:comp-CFAA}) as a definition.

Tensoring with the type \AAm\ bimodule for the identity map turns type
$D$ invariants into type $A$ invariants. More precisely, suppose that
$Y$ is a $\PMC$-bordered $3$-manifold. Then
\begin{equation}
\CFAa(Y)\simeq \CFAAa(\Id_\PMC)\DT\CFDa(Y).\label{eq:comp-CFA}
\end{equation}
Here, $\DT$,
\glsit{$\DT$}%
an operation taking in one
$\Ainf$-module and one type $D$ structure,
denotes a particular model for the derived tensor
product~\cite[Section~\ref*{LOT2:sec:tensor-products}]{LOT2}.

Similarly, if $Y'$ is a strongly bordered $3$-manifold with two boundary
components $F(\PMC_1)$ and $F(\PMC_2)$ then 
\begin{align}
  \lsub{F(\PMC_1)}\CFDAa(Y')_{F(\PMC_2)}&\simeq
  \CFAAa(\Id_{\PMC_1})\DT\CFDDa(Y')\label{eq:comp-CFDA}\\
  \lsub{F(\PMC_1),F(\PMC_2)}\CFAAa(Y')&\simeq
  \CFAAa(\Id_{\PMC_2})\DT \CFAAa(\Id_{\PMC_1})\DT\CFDDa(Y')\label{eq:comp-CFAA}.
\end{align}
\glsit{$\CFDAa(Y)$}\glsit{$\CFAAa(Y)$}%

In Section~\ref{sec:ElementaryCobordisms}, we explained how to compute
$\CFDa(Y)$ and $\CFDDa(Y')$. So, Equations~\eqref{eq:comp-CFA},
\eqref{eq:comp-CFDA} and~\eqref{eq:comp-CFAA} give algorithms for
computing the type $A$, \DA\ and \AAm\ invariants, as well.

We can describe this procedure for calculating the bimodules
in more detail, as follows. Theorem~\ref{thm:DDforArc-Slides}
computes the type \DD\ modules for arc-slides, and
Proposition~\ref{prop:CalculateElementaryCobordism} computes the type
\DD\ module for a split elementary cobordism. Using
Equation~\eqref{eq:comp-CFDA}, we can turn each of these bimodules
into a type \DA\ module. Any bordered $3$-manifold $Y$ can be factored
as
\[
Y=Y_1\sos{\bdy_R}{\cup}{\bdy_L}Y_2\sos{\bdy_R}{\cup}{\bdy_L}\cdots\sos{\bdy_R}{\cup}{\bdy_L}Y_k,
\]
where each $Y_i$ is either an arc-slide or a split elementary cobordism.
The pairing theorem~\cite[Theorem~\ref*{LOT2:thm:GenComposition}]{LOT2}
then gives
\begin{equation}
\CFDAa(Y)=\CFDAa(Y_1)\DT\CFDAa(Y_2)\DT\cdots\DT \CFDAa(Y_k).\label{eq:recon-CFDA}
\end{equation}
We can compute $\CFDDa(Y)$ or $\CFAAa(Y)$ by tensoring $\CFDAa(Y)$ with
$\CFDDa(\Id)$ or $\CFAAa(\Id)$, respectively.

As a special case, Equation~\eqref{eq:recon-CFDA} gives another
algorithm for reconstructing $\CFa(Y)$ for a closed $3$-manifold $Y$
via a Heegaard splitting of $Y$. (Actually, unpacking the definitions,
this is the same as the formula in Theorem~\ref{thm:HFa}.)

\subsection{Homological perturbation theory}\label{sec:hpt}
One key advantage of $\Ainf$-modules is that for any $\Ainf$-module
(and in particular, any \dg module) $M$, the homology $H_*(M)$ of $M$
carries an $\Ainf$-module structure which is (under mild assumptions)
homotopy equivalent to $M$. We state a version of this result
presently; for more on this, see for
example~\cite{Keller:OtherAinfAlg} (which discusses the case of
algebras, rather than modules).

Fix a ground ring (with trivial differential) $\Ground$ of
characteristic $2$. 
We consider $\Ainf$ algebras over $\Ground$, and
strictly unital $\Ainf$-modules. In particular, such $\Ainf$-modules
over $\Alg$ are honest differential modules over $\Ground$.
\index{strictly unital}%

\begin{lemma}
  \label{lem:HomologicalPerturbation}
  Let ${\mathcal M}$ be an $\Ainf$-module over an $\Ainf$-algebra
  $\Alg$ over
  $\Ground$, let $M$ denote its underlying chain complex over $\Ground$,
  and let $f\co N \rightarrow M$ be a homotopy equivalence of
  $\Ground$-modules. Then, we can find
  \begin{itemize}
  \item an $\Ainf$-module structure ${\mathcal N}$ on $N$ and
  \item an $\Ainf$ quasi-isomorphism
    $$F\co {\mathcal N} \rightarrow {\mathcal M}$$
    with the property that ${F}_1=f$.
  \end{itemize}
\end{lemma} 
(Of  course, if $\Ground$ is $\Field$, or more generally a direct sum of copies of
$\Field$ then we can replace the condition that $f$ be a homotopy
equivalence by the condition that $f$ be a quasi-isomorphism: for
these rings, any quasi-isomorphism is a homotopy equivalence.)
\begin{proof}
  We recall now some notation.
  Let $\Alg[1]$
  \glsit{$\Alg[1]$}%
denote the algebra $\Alg$ with grading shifted by $1$
  (or, in the group-graded context, shifted by $\lambda$; i.e.
  $\Alg[1]_g=\Alg_{\lambda^{-1}g}$). If $V$ is a $\Ground$-bimodule, let $\Tensor^*(V)$ denote the
  tensor algebra on $V$ (this includes the $0\th$
  tensor product), and let $\Tensor^+(V)\subset \Tensor^*(V)$ denote the ideal 
  generated by $V$.
  \glsit{$\Tensor^*(V)$}\glsit{$\Tensor^+(V)$}%
  Here (and below), tensor products will be taken over $\Ground$.
  In practice, we will apply this construction in the case where $V=\Alg[1]$,
  thought of as a $\Ground$-bimodule.

  Let $g\co M\to N$ be a homotopy inverse to $f$ and 
  $T\co M\to M$ be a homotopy between $f\circ g$ and $\Id_M$.
  
  The $\Ainf$-module structure on ${\mathcal N}$ is given as follows. We
  write the $\Ainf$-module structure on $\mathcal{M}$ as a map $m\co
  M\otimes \Tensor^*(A[1])\to M$.
  Comultiplication induces a map
  $$ \mu^*\co \Tensor^+(\Alg)\rightarrow
  \Tensor^*(\Tensor^+(\Alg[1])),$$
  defined as follows: for ${\overline a}=a_1\otimes\dots\otimes a_n\in \Tensor^+(\Alg)$,
  let
  \begin{multline*}
    \mu^*(a_1\otimes\dots\otimes a_n)= \\
\sum_{\mathclap{\{i_1,\dots, i_k | 1\leq i_1<\dots<i_k<n\}}}\ \ 
  (a_1\otimes\dots\otimes a_{i_1})\otimes (a_{i_1+1}\otimes\dots\otimes a_{i_2})
  \otimes\dots\otimes (a_{i_k+1}\otimes\dots\otimes a_n).
  \end{multline*}
  With these maps in hand, define the operations $m_i$ for $i>1$ on
  $\mathcal{N}$ by:

  \begin{equation}\label{eq:induced-mod-str}
  \mathcenter{
    {\begin{tikzpicture}[x=1cm,y=32pt]
        \node at (-1,-1) (empty) {};
        \node at (-1,7) (x) {$m^N(\x\otimes{\overline{a}})=$};
        \end{tikzpicture}}
    {\begin{tikzpicture}[x=1cm,y=32pt]
        \node at (-1,-2) (empty) {};
        \node at (-1,2) (terminal) {};
        \node at (-1,3) (g) {$g$};
        \node at (-1,4) (m1) {$m$};
        \node at (-1,5) (f) {$f$};
        \node at (0,5) (mu) {$\mu^*$};
        \node at (-1,6) (x) {${\mathbf x}$};
        \node at (0,6) (alg) {${\overline a}$};
        \draw[modarrow] (g) to (terminal);
        \draw[othmodarrow] (m1) to (g);
        \draw[modarrow] (x) to (f);
        \draw[othmodarrow] (f) to (m1);
        \draw[tensoralgarrow] (alg) to (mu);
        \draw[tensoralgarrow, bend left=15] (mu) to (m1); 
      \end{tikzpicture}}
    {\begin{tikzpicture}[x=1cm,y=32pt]
        \node at (-1,-1) (empty) {};
        \node at (-1,1) (terminal) {};
        \node at (-1,2) (g) {$g$};
        \node at (-1,3) (m2) {$m$};
        \node at (-1,4) (T1) {$T$};
        \node at (-1,5) (m1) {$m$};
        \node at (-1,6) (f) {$f$};
        \node at (0,6) (mu) {$\mu^*$};
        \node at (-2,7) (plus) {$+$};
        \node at (-1,7) (x) {${\mathbf x}$};
        \node at (0,7) (alg) {${\overline a}$};
        \draw[othmodarrow] (m2) to (g);
        \draw[modarrow] (g) to (terminal);
        \draw[othmodarrow] (T1) to (m2);
        \draw[othmodarrow] (m1) to (T1);
        \draw[modarrow] (x) to (f);
        \draw[othmodarrow] (f) to (m1);
        \draw[tensoralgarrow] (alg) to (mu);
        \draw[tensoralgarrow, bend left=15] (mu) to (m2); 
        \draw[tensoralgarrow, bend left=15] (mu) to (m1); 
      \end{tikzpicture}}
    {\begin{tikzpicture}[x=1cm,y=32pt]
        \node at (-1,-1) (terminal) {};
        \node at (-1,0) (g) {$g$};
        \node at (-1,1) (m3) {$m$};
        \node at (-1,2) (T2) {$T$};
        \node at (-1,3) (m2) {$m$};
        \node at (-1,4) (T1) {$T$};
        \node at (-1,5) (m1) {$m$};
        \node at (-1,6) (f) {$f$};
        \node at (0,6) (mu) {$\mu^*$};
        \node at (-1,7) (x) {${\mathbf x}$};
        \node at (-2,7) (plus) {$+$};
        \node at (1,7) (plus2) {$+$};
        \node at (2,7) (dots) {$\dots$};
        \node at (0,7) (alg) {${\overline a}$};
        \draw[othmodarrow] (m3) to (g);
        \draw[modarrow] (g) to (terminal);
        \draw[othmodarrow] (m3) to (g);
        \draw[othmodarrow] (T2) to (m3);
        \draw[othmodarrow] (m2) to (T2);
        \draw[othmodarrow] (T1) to (m2);
        \draw[othmodarrow] (m1) to (T1);
        \draw[modarrow] (x) to (f);
        \draw[othmodarrow] (f) to (m1);
        \draw[tensoralgarrow] (alg) to (mu);
        \draw[tensoralgarrow, bend left=15] (mu) to (m3); 
        \draw[tensoralgarrow, bend left=15] (mu) to (m2); 
        \draw[tensoralgarrow, bend left=15] (mu) to (m1); 
      \end{tikzpicture}}}
  \end{equation}
  Here, doubled arrows indicate elements of $\Tensor^+\Alg$; in
  particular, there is always at least one algebra element
  present. Dashed lines indicate elements of $N$, while solid lines
  indicate elements of $M$. It is a property of the comultiplication
  $\mu^*$ that for any given ${\overline a}$, there are only finitely
  many non-zero elements in this sum.

  We verify that this is indeed an $\Ainf$ module. As usual, $m^N$ induces
  an endomorphism ${\overline m}$ of 
  $N \otimes \Tensor^*(\Alg)$, and $\mu$ a differential $d$ on
  $\Tensor^*(\Alg)$. We must check that 
  \begin{equation}
    \label{eq:AInftyModule}
    {\overline m}^2(\x,{\overline a}) + {\overline m}(x,d{\overline
      a})=0,
  \end{equation} for any ${\overline a}=a_1\otimes\dots\otimes a_n$.
  We begin
  this verification by considering terms in the second sum
  ${\overline m}(x,d{\overline a})$. Applying the $\Ainf$ relation on
  $M$, we see that ${\overline m}(x,d{\overline a})$ can be
  interpreted as counting the same kinds of trees as in the definition
  of $m^N(\x\otimes {\overline a})$, except for one difference:
  whereas the $m$-labeled vertices in the definition of
  $m^N(\x\otimes {\overline a})$ all have incoming algebra elements
  (and no two of these $m$-labeled vertices are consecutive), the
  trees in ${\overline m}(x,d{\overline a})$ have a pair of consecutive
  $m$-labeled vertices (one of which may have no incoming algebra elements).
  Equivalently, we can think of these as counting trees obtained from
  trees counted in ${\overline m}(x,{\overline a})$, by
  applying one of the following operations:
  \begin{enumerate}[label=(T-\arabic*),ref=T-\arabic*]
  \item \label{type:BeforeT}
    Insert an $m$-labeled vertex with no incoming algebra elements
    immediately before some 
    $T$-labeled vertex.
  \item 
    \label{type:AfterT}
    Insert an $m$-labeled vertex with no incoming algebra elements
    immediately after some 
    $T$-labeled vertex.
  \item 
    \label{type:AfterF}
    Insert an $m$-labeled vertex with no incoming algebra elements
    immediately after the initial $f$-labeled vertex.
  \item 
    \label{type:BeforeG}
    Insert an $m$-labeled vertex with no incoming algebra elements
    before the terminal $g$-labeled vertex.
  \item
    \label{type:Inside}
    Split some $m$-labeled vertex with at least $2$ incoming algebra elements
    to a pair of consecutive $m$-labeled vertices, each of which has at least $1$
    incoming algebra element.
  \end{enumerate}
  Terms of Type~(\ref{type:BeforeT}) pair off with terms of
  Type~(\ref{type:AfterT}) (in view of the formula $d\circ T+T\circ d=\Id+f\circ
  g$)
  to produce a sum of two types of trees: one of these types match
  those of Type~(\ref{type:Inside}); the other type is gotten 
  by applying the following operation:
  \begin{enumerate}[resume*]
  \item 
    \label{type:gf}
    Replace a $T$-labeled vertex by a vertex labeled by $f\circ g$.
  \end{enumerate}
  Thus, 
  ${\overline m}(x,d{\overline
    a})$ counts trees of Types~(\ref{type:AfterF}),
  (\ref{type:BeforeG}) and
  (\ref{type:gf}). The fact that $f$ and $g$ are chain maps allow us to 
  move the differentials past the $f$ and $g$-labeled vertices in the trees
  of Types~(\ref{type:AfterF}) and (\ref{type:BeforeG}).

  But these terms are precisely the trees counted in ${\overline
    m}^2(\x,{\overline a})$; and hence Equation~\eqref{eq:AInftyModule} holds.

  We claim that $f$ extends to an $\Ainf$ homomorphism, which has the following graphical
  representation:
  \begin{equation}\label{eq:induced-quasi-iso}
  \mathcenter{
    {\begin{tikzpicture}[x=1cm,y=32pt]
        \node at (-1,-5) (empty) {};
        \node at (-1,3) (x) {$f(\x\otimes{\overline{a}})=$};
        \end{tikzpicture}}
    {\begin{tikzpicture}[x=1cm,y=32pt]
        \node at (-1,-4) (empty) {};
        \node at (-1,2) (terminal) {};
        \node at (-1,3) (f) {$f$};
        \node at (-1,4) (x) {${\mathbf x}$};
        \draw[modarrow] (x) to (f);
        \draw[algarrow] (f) to (terminal);
      \end{tikzpicture}}
    {\begin{tikzpicture}[x=1cm,y=32pt]
        \node at (-1,-2) (empty) {};
        \node at (-1,2) (terminal) {};
        \node at (-1,3) (T1) {$T$};
        \node at (-1,4) (m1) {$m$};
        \node at (-1,5) (f) {$f$};
        \node at (0,5) (mu) {$\mu^*$};
        \node at (-1,6) (x) {${\mathbf x}$};
        \node at (0,6) (alg) {${\overline a}$};
        \node at (-2,6) (plus) {$+$};
        \draw[algarrow] (T1) to (terminal);
        \draw[algarrow] (m1) to (T1);
        \draw[modarrow] (x) to (f);
        \draw[algarrow] (f) to (m1);
        \draw[tensoralgarrow] (alg) to (mu);
        \draw[tensoralgarrow, bend left=15] (mu) to (m1); 
      \end{tikzpicture}}
    {\begin{tikzpicture}[x=1cm,y=32pt]
        \node at (-1,-1) (empty) {};
        \node at (-1,1) (terminal) {};
        \node at (-1,2) (T2) {$T$};
        \node at (-1,3) (m2) {$m$};
        \node at (-1,4) (T1) {$T$};
        \node at (-1,5) (m1) {$m$};
        \node at (-1,6) (f) {$f$};
        \node at (0,6) (mu) {$\mu^*$};
        \node at (-2,7) (plus) {$+$};
        \node at (-1,7) (x) {${\mathbf x}$};
        \node at (0,7) (alg) {${\overline a}$};
        \draw[algarrow] (m2) to (T2);
        \draw[algarrow] (T2) to (terminal);
        \draw[algarrow] (T1) to (m2);
        \draw[algarrow] (m1) to (T1);
        \draw[modarrow] (x) to (f);
        \draw[algarrow] (f) to (m1);
        \draw[tensoralgarrow] (alg) to (mu);
        \draw[tensoralgarrow, bend left=15] (mu) to (m2); 
        \draw[tensoralgarrow, bend left=15] (mu) to (m1); 
      \end{tikzpicture}}
    {\begin{tikzpicture}[x=1cm,y=32pt]
        \node at (-1,-1) (terminal) {};
        \node at (-1,0) (T3) {$T$};
        \node at (-1,1) (m3) {$m$};
        \node at (-1,2) (T2) {$T$};
        \node at (-1,3) (m2) {$m$};
        \node at (-1,4) (T1) {$T$};
        \node at (-1,5) (m1) {$m$};
        \node at (-1,6) (f) {$f$};
        \node at (0,6) (mu) {$\mu^*$};
        \node at (-1,7) (x) {${\mathbf x}$};
        \node at (-2,7) (plus) {$+$};
        \node at (1,7) (plus2) {$+$};
        \node at (2,7) (dots) {$\dots$};
        \node at (0,7) (alg) {${\overline a}$};
        \draw[algarrow] (m3) to (T3);
        \draw[algarrow] (T3) to (terminal);
        \draw[algarrow] (m3) to (T3);
        \draw[algarrow] (T2) to (m3);
        \draw[algarrow] (m2) to (T2);
        \draw[algarrow] (T1) to (m2);
        \draw[algarrow] (m1) to (T1);
        \draw[modarrow] (x) to (f);
        \draw[algarrow] (f) to (m1);
        \draw[tensoralgarrow] (alg) to (mu);
        \draw[tensoralgarrow, bend left=15] (mu) to (m3); 
        \draw[tensoralgarrow, bend left=15] (mu) to (m2); 
        \draw[tensoralgarrow, bend left=15] (mu) to (m1); 
      \end{tikzpicture}}}
  \end{equation}
  We leave the verification that this extension is, indeed, an $\Ainf$ homomorphism as an exercise
  for the reader. The fact that $F_1=f$ is immediate, and hence $F$ is a quasi-isomorphism.
\end{proof}

\subsection{Reconstruction via type \textalt{$A$}{A} modules, and its advantages}\label{sec:reconstruct-via-A}
As discussed in Section~\ref{sec:computing-type-a},
Theorem~\ref{thm:HFa} can be reformulated as follows.  As before,
write $Y$ as a union of two $0$-framed handlebodies $\HB^g$,
identified via some gluing map $\phi$.  Once again,
decompose $\phi$ as a composition of arc-slides,
$\phi=m_1\circ\dots\circ m_k$. Next, calculate $\CFDAa(m_i)$ for each
arc-slide, using Theorem~\ref{thm:DDforArc-Slides} and
Equation~\eqref{eq:comp-CFDA}; and $\CFAa$ for the handlebody $\HB^g$
using Proposition~\ref{prop:CFD-of-0-fr-hdlbdy} and
Equation~\eqref{eq:comp-CFA} (or by simply writing down $\CFAa(\HB^g)$
directly by counting holomorphic curves).  Then, there is a
homotopy equivalence
\begin{equation}
  \label{eq:CalculateHFaAgain}
  \begin{split}
    \CFa(Y)&\simeq \CFAa(\HB^g)\DT\CFDAa(m_1)\DT\cdots \DT
    \CFDAa(m_k)\CFDa(\HB^g)\\
    &=\CFDa(\HB^g)\DT \CFAAa(\Id)\DT
    \CFDDa(m_1)\DT\CFAAa(\Id)\DT\cdots\\
    &\qquad\qquad\DT
    \CFDDa(m_k)\DT\CFAAa(\Id)\DT\CFDa(\HB^g).
  \end{split}
\end{equation}
(For clarity, we are viewing the bimodules as each having one left and
one right action.)

So far, we have not gained anything computationally. The point is
that, since the operation $\DT$ respects quasi-isomorphisms, we are
free to replace the bimodule $\CFAAa(\Id)$ with a smaller model than
the one given by Equation~\eqref{eq:DefAA}. 
Thinking of the type \AAm\ module associated to the standard diagram
for the identity map (Figure~\ref{fig:Genus2Identity}), 
it is clear that there is a model for $\CFAAa(\Id)$ whose rank
over $\Field$ is the number of idempotents in
$\Alg(\PMC)$.  One can arrive at a model of this
size by taking the homology of the chain complex $\CFAAa(\Id_{\PMC})$
(thought of as a vector space over $\Field$).  Although this homology $H_*(\CFAAa(\Id_{\PMC}))$
is no longer an honest $\Alg(\PMC)$ bimodule, by
Lemma~\ref{lem:HomologicalPerturbation} it does
retain the structure of an $\Ainf$-bimodule over
$\Alg(\PMC)$ (which is quasi-isomorphic to
$\CFAAa(\Id_{\PMC})$).

Replacing $\CFAAa(\Id_{\PMC})$ by its homology drastically reduces the
number of generators. For instance, we will see in
Section~\ref{sec:AA-torus-id} that in the torus boundary case, 
$\dim_{\Field}(\CFAAa(\Id_{\PMC},0))=30$,
while $\dim_{\Field}(H_*(\CFAAa(\Id_{\PMC},0)))=2$.
(The $0$ in the notation here refers to the restriction to 
the portion of the algebra with strands weight zero, i.e., with exactly one moving strand or 
two matched horizontal strands.)

In practice, to work with the fewest generators possible, one would
want to compute Equation~\eqref{eq:CalculateHFaAgain} by starting at
one end, and taking homology (using
Lemma~\ref{lem:HomologicalPerturbation}) after each tensor product.

Doing the replacement of $\CFAAa(\Id_{\PMC})$ with
$H_*(\CFAAa(\Id_{\PMC}))$, one has to take care about the boundedness
hypotheses needed for $\DT$ to be well-defined. If we take a bounded
model for $\CFDa(\HB^g)$ on the left of
Equation~\eqref{eq:CalculateHFaAgain}, and compute the tensor products
starting at the far right, the relevant boundedness hypotheses are
satisfied throughout, even if one takes homology after each tensor
product.

As a last computational point, we can use a smaller model for the
algebra $\Alg(\PMC)$: we can divide out by the differential ideal of algebra elements with
local multiplicity at least $2$ somewhere.  (See
\cite[Proposition~\ref*{LOT2:prop:SmallerModel}]{LOT2}.)

\subsection{Example: torus boundary}\label{sec:AA-torus-id}
\begin{figure}
\begin{center}
\input{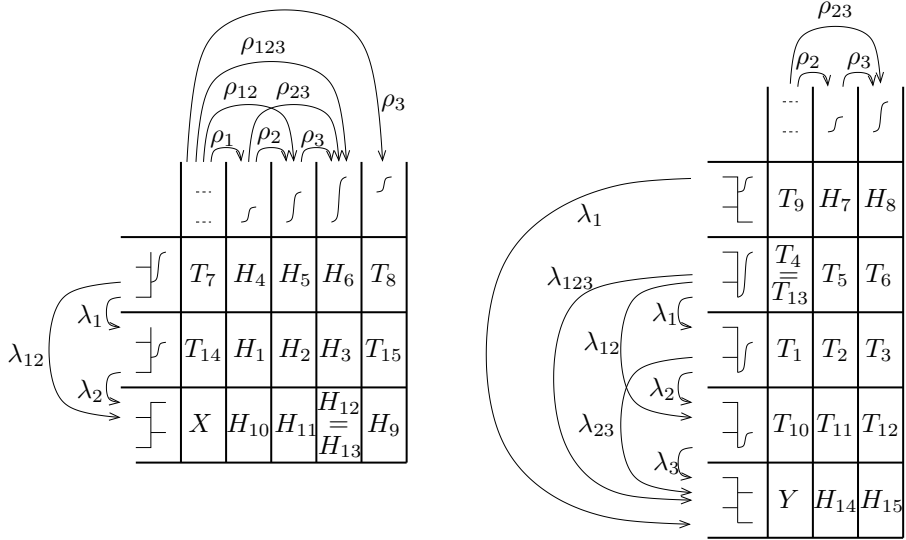}
\end{center}
\caption {{\bf Genus one identity.}
\label{fig:Genus1AA}
Dualizing the type \DD\  bimodule, in the case of genus one.}
\end{figure}

Let $\Id=\Id_{\PMC^1}$ denote the identity map of the torus.  In this
section we calculate the $AA$ identity bimodule for the torus with strands weight zero, 
$\CFAAa(\Id)=\CFAAa(\Id_{\PMC^1},0)$, using results of the previous
subsections. This module
was calculated before, by explicitly finding all the relevant
holomorphic curves,
in~\cite[Section~\ref*{LOT2:subsec:AAId1}]{LOT2}. The method here is
more algebraic, and generalizes in a straightforward way to arbitrary
genus; direct holomorphic curve counts in higher genus would be
complicated at best, and probably intractable. 

Let $\Alg=\Alg(\PMC^1,0)$ be the algebra associated to the torus with strands weight equal to zero, 
and $\Alg'=\Alg(-\PMC^1,0)$.
(The algebras $\Alg$ and $\Alg'$ are isomorphic, but
it will be clearer to treat them as distinct.)
Note that $\CFDDa(\Id)$ has eight generators as a left $\Alg$ module;
it also admits a left action by $\Alg'$, a ring whose generators we denote
$\lambda_i$ rather than $\rho_i$.

Similarly, $\Alg$ has eight generators as an $\Field$-vector space.

By Equation~\eqref{eq:DefAA}, $\CFAAa(\Id)$ is equivalent to the bimodule
of left $\Alg$-module maps from $\CFDDa(\Id)$ to $\Alg$. 
Such an $\Alg$-linear map from $\CFDDa(\Id)$ to $\Alg$ is determined by
the image of a basis for $\CFDDa(\Id)$; and it is zero unless elements are mapped to 
elements with compatible idempotents.

Using the description of $\CFDDa(\Id)$ from Theorem~\ref{thm:DDforIdentity},
there are $30$ generators of the bimodule
$\Hom_{\Alg}(\CFDDa(\Id),{\Alg})$.  It is straightforward to find
the $15$ differentials which connect various generators. The results of
this are illustrated in Figure~\ref{fig:Genus1AA}, with the following
convention. Generators are named $H_i$, $T_i$, $X$ or $Y$, with the
understanding that $H_i$ appears in the differential of $T_i$, while
the generators $X$ and $Y$ have no differentials either entering or leaving
them. The labels $H_{12}$ and $H_{13}$ refer to the same generator, as
do $T_4$ and $T_{13}$. In the table, rows
are indexed by generators of $\Alg'\otimes \CFDDa(\Id)$, and the
columns are indexed by generators of
$\Alg$. Each square corresponds to the morphism which takes the generator in that row
to the generator in that column (and all other generators to zero). The right action by $\Alg$ is indicated by
the arrows connecting the columns, while the right action by $\Alg'$ is indicated by the arrows connecting
the rows. (Note that left translation in $\Alg'$  dualizes to the stated right action.)

The homology of the complex is two-dimensional,
generated by the generators $X$ and $Y$.
Thus, $\CFAAa(\Id)$ is quasi-isomorphic to an $\Ainf$-bimodule with
just two generators, $X_0$ and $Y_0$ (as we already knew from the standard
Heegaard diagram for the identity map).
Indeed,
the quasi-isomorphism
$$f\co H(\CFAAa(\Id))\rightarrow \CFAAa(\Id)$$
defined by $f(X_0)=X$, $f(Y_0)=Y$ has a homotopy inverse
$$g\co \CFAAa(\Id) \rightarrow H(\CFAAa(\Id))$$
defined by
$g(X)=X_0$, $g(Y)=Y_0$, $g(H_i)=g(T_i)=0$. 
Note that $g\circ f=\Id_{H_*(\CFAAa(\Id))}$, while $f\circ g\simeq
\Id_{\CFAAa(\Id)}$ via the map
$$T\co \CFAAa(\Id) \rightarrow \CFAAa(\Id)$$
specified by
\begin{equation}\label{eq:DefOfT}
\begin{aligned}
T(X) &= 0  \\
T(Y) &= 0 \\
T(T_j) &= 0,&&1 \le j \le 15 \\
T(H_4) &= T_4 + T_{12} \\
T(H_i) &= T_i,&&1 \le i \le 15, i \notin \{4, 13\} 
\end{aligned}
\end{equation}
At first glance,
it might appear that $T$ is not defined on $H_{13}$; but in fact
$H_{13}=H_{12}$, so $T(H_{13})=T_{12}$.

\begin{figure}
\[
\begin{tikzpicture}[x=1.4cm,y=1.4cm]
  \matrix[column sep={1.4cm,between origins}, row sep = {1.4cm,between origins}] {
    \node[boxed] (X1) {$X$}; &
      \node (H10) {$H_{10}$}; &
        \node (T10) {$T_{10}$}; &
          \node[boxed] (Y1) {$Y$}; \\
    \\
    & & & \node[boxed] (X2) {$X$}; &
            \node (H1) {$H_1$}; &
              \node (T1) {$T_1$}; &
                \node[boxed] (Y2) {$Y$}; \\
    & \node[boxed] (Y3) {$Y$}; &
        \node (H14) {$H_{14}$}; &
          \node (T14) {$T_{14}$}; &
            \node (H2) {$H_2$}; &
              \node (T2) {$T_2$}; \\
    \node[boxed] (X3) {$X$}; &
      \node (H11) {$H_{11}$}; &
        \node (T11) {$T_{11}$}; & &
            \node (H3) {$H_3$}; &
              \node (T3) {$T_3$}; &
                \node[boxed] (X4) {$X$}; & &
                    \node[boxed] (Y4) {$Y$}; \\
    \node[boxed] (X5) {$X$}; &
      & \node (H12) {$H_{12} = H_{13}$}; &
          & \node (T12) {$T_{12}$}; &
              \node (H15) {$H_{15}$}; &
                \node (T15) {$T_{15}$}; &
                  \node (H9) {$H_9$}; &
                    \node (T9) {$T_9$}; &
                      \node (H8) {$H_8$}; &
                        \node (T8) {$T_8$}; \\
    & & \node[boxed] (Y5) {$Y$}; &
          & \node[boxed] (Y6) {$Y$}; &
              \node (T6) {$T_6$}; &
                \node (H6) {$H_6$}; &
                & \node (H7) {$H_7$}; \\
    & & \node (T13) {$T_{13} = T_4$}; &
          & \node (H4) {$H_4$}; &
              & \node (T7) {$T_7$}; &
                  & \node (H5) {$H_5$}; &
                      \node (T5) {$T_5$}; \\
    & & & & & & \node[boxed] (X6) {$X$}; \\
  };
  \draw[->] (X1) to node[above,cdlabel] {\rho_1} (H10);
    \draw[->,dashed] (H10) to (T10);
      \draw[->] (T10) to node [above,cdlabel] {\lambda_3} (Y1);
  \draw[->,dashed] (H1) to (T1);
    \draw[->] (T1) to node[above,cdlabel] {\lambda_{23}} (Y2);
  \draw[->] (Y3) to node [above,cdlabel] {\rho_2} (H14);
    \draw [->,dashed] (H14) to (T14);
      \draw [->] (T14) to node [left,cdlabel] {\lambda_2} (X2);
      \draw [->] (T14) to node [above left=-3pt,cdlabel] {\rho_1} (H1);
      \draw [->] (T14) to node [above,cdlabel] {\rho_{12}} (H2);
        \draw [->,dashed] (H2) to (T2);
        \draw [->] (T2) to[out=0,in=-90] (1.8,1.8) to[out=90,in=90]
           node[above right=-3pt,cdlabel,pos=0.2] {\lambda_{23}} (H14);
  \draw[->] (X3) to node[above,cdlabel] {\rho_{12}} (H11);
    \draw[->,dashed] (H11) to (T11);
      \draw[->] (T11) to node[left,cdlabel] {\lambda_3} (H14);
        \draw [->] (T14) to node[below left=-3pt,cdlabel] {\rho_{123}} (H3);
          \draw[->,dashed] (H3) to (T3);
  \draw[->] (X5) to node[above,cdlabel]{\rho_{123}} (H12);
    \draw[->,dashed] (H12) to (T12);
      \draw[->] (T12) to node[above,cdlabel]{\lambda_3} (H15);
        \draw[->] (T3) to node[right,cdlabel]{\lambda_{23}} (H15);
        \draw[->,dashed] (H15) to (T15);
          \draw[->] (X4) to node[above right=-3pt,cdlabel]{\rho_3} (H9);
          \draw[->] (T15) to node[below,cdlabel]{\lambda_2} (H9);
            \draw[->,dashed] (H9) to (T9);
              \draw[->] (T9) to node[left,cdlabel]{\lambda_1} (Y4);
              \draw[->] (T9) to node[above,cdlabel]{\rho_{23}} (H8);
                \draw[->,dashed] (H8) to (T8);
                \draw[->] (T8) to[out=90,in=0] (3.6,0.7) to[out=180,in=90] node[pos=0,cdlabel,above] {\lambda_{12}} (H9);
  \draw[->] (Y6) to node[above left=-3pt,cdlabel,pos=0.4] {\rho_{23}} (H15);
    \draw[->] (T6) to node [right,cdlabel] {\lambda_{123}} (H15);
    \draw[->,dashed] (H6) to (T6);
      \draw[->] (T9) to node [right,cdlabel] {\rho_2} (H7);
  \draw[->] (T13) to node[left,cdlabel] {\lambda_{123}} (Y5);
    \draw[->,dashed] (H4) to (T13);
    \draw[->,dashed,bend left=40] (H4) to (T12);
      \draw[->] (T7) to node[above,cdlabel] {\rho_1} (H4);
        \draw[->] (T7) to node[above,cdlabel] {\rho_{12}} (H5);
        \draw[->] (T7) to node[left,cdlabel] {\rho_{123}} (H6);
        \draw[->,dashed,bend right=10] (H7) to (T7);
          \draw[->,dashed] (H5) to (T5);
          \draw[->] (T5) to[out=0,in=-90] (6,0) to[out=90,in=4] node [above,cdlabel,pos=0.9] {\lambda_{123}} (0,3.7) to[out=184,in=120] (H14);
  \draw[->] (T7) to node[left,cdlabel] {\lambda_{12}} (X6);
\end{tikzpicture}
\]
\caption {{\bf Genus one identity type $AA$ bimodule.}
\label{fig:Genus1AAgraph}
A graphical representation of the $\Ainf$ operations
on $\CFAAa(\Id)$. $T$-labeled edges are shown dashed. All vertices
labeled $X$ are identified, as are all vertices labeled~$Y$. The
algebra-labeled edges can be immediately
read off from the algebra operations coming from Figure~\ref{fig:Genus1AA};
the $T$-labeled edges are determined by Equation~\eqref{eq:DefOfT}.}
\end{figure}
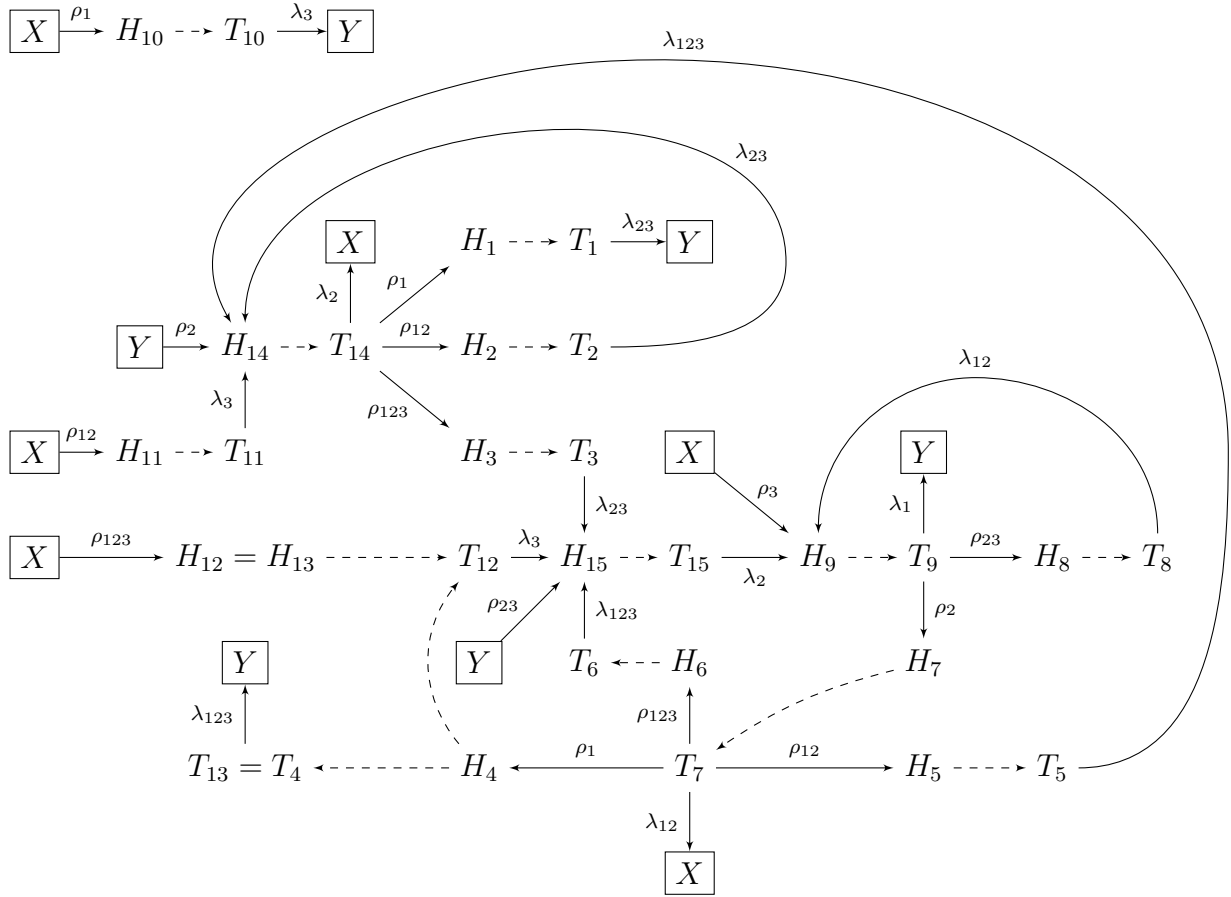

An explicit form of the $\Ainf$ structure on $H_*(\CFAAa(\Id))$ is
constructed in the proof of Lemma~\ref{lem:HomologicalPerturbation}.  We can
think of this graphically, as follows.  Consider the directed graph
whose nodes correspond to the $30$ generators of
$\Hom_{\Alg}(\CFDDa(\Id),{\Alg})$.  Draw an edge from $v_1$ to
$v_2$ labeled by a basis vector $a$ for the algebra if $v_2$ appears in
$m_2(v_1,a)$; include another kind of edge---a {\em $T$-labeled edge}---from $v_1$ to $v_2$ if $v_2$ appears in $T(v_1)$. All $\Ainf$
operations on $H_*(\CFAAa(\Id))$ correspond to paths with the following
properties:
\begin{itemize}
  \item The path starts and ends at vertices labeled $X$ or $Y$.
  \item The initial and final edges in the graph are labeled by
    algebra elements.
  \item the path alternates between edges labeled by algebra elements
    and $T$-labeled edges.
\end{itemize}
Each such path corresponds to a term in an $\Ainf$ operation, starting
at the initial vertex $v$, with coefficient $1$ in the terminal vertex
$w$; and the sequence of algebra elements give the sequence of operations.
More precisely, if $r_1,\dots r_m$ and $\ell_1,\dots, \ell_n$ are the
sequence of algebra elements in the order they are encountered, where
here $r_i$ correspond to those labeled by algebra elements gotten by
products of $\rho_j$, while $\ell_i$ correspond to those gotten as
products of the $\lambda_j$, then this path gives a contribution of
$w$ in the operation $m_{1,n,m}(v, (\ell_1\otimes \dots\otimes \ell_n),
(r_1\otimes \dots\otimes r_m))$. Note that paths of the above type in
this graph coincide with the smaller graph, where we include
algebra-labeled edges only in cases where the terminal point of the
edge is either the initial vertex of a $T$-labeled edge, or it is
one of $X$ or $Y$. This smaller graph is illustrated in
Figure~\ref{fig:Genus1AAgraph} (except that we have drawn several
vertices corresponding to $X$ and $Y$, for clarity).

For instance, we see at once that $m_{1,1,1}(X_0,\rho_3,\lambda_2)=Y_0$, 
via the path through $X$, $H_9$, $T_9$, and
$Y$. More generally, by traveling around the loop through
$H_9$, $T_9$, $H_8$, and $T_8$ we see that
$m(X_0,\rho_3,\overbrace{\rho_{23},\dots,\rho_{23}}^{n},
\overbrace{\lambda_{12},\dots,\lambda_{12}}^n,\lambda_2)=Y_0$
for any $n\geq 0$.
(Note that to draw this conclusion, we need to check that
there are no other paths in the graph which give rise to a cancelling
term.)

\subsection{A computer computation}\label{sec:ok-computer}
\begin{figure}
  \centering
  \includegraphics{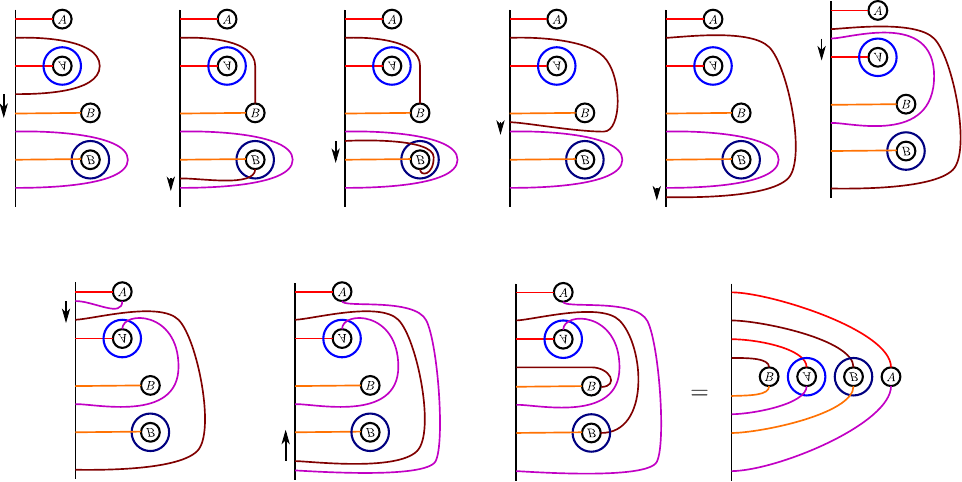}
  \caption{\textbf{Obtaining the genus $2$ self-gluing handlebody by arc-slides.}}
  \label{fig:self-gluing-hb}
\end{figure}

We conclude by describing a computer computation of $\HFa$ of the
Poincar\'e homology sphere, via an open book decomposition. For this,
we will use a particular class of handlebodies, useful for studying
open books:
\begin{definition}\label{def:sghb} Let $\PMC$ be a pointed matched circle.
  The \emph{self-gluing handlebody of $\PMC$}, denoted $\HB_{sg}(\PMC)$,
  \glsit{$\HB_{sg}(\PMC)$}\index{self-gluing handlebody}\index{handlebody!self-gluing}%
is the handlebody corresponding to a bordered Heegaard diagram $\HD_{sg}(\PMC)$ obtained from the standard Heegaard diagram for the identity map of $\PMC$ by deleting the arc $\mathbf{z}$ and placing a basepoint on one of the two sides.
\end{definition}
See the lower right of Figure~\ref{fig:self-gluing-hb} for the case
when $\PMC$ is the genus $1$ pointed matched circle.  Also, observe
that $\bdy\HB_{sg}(\PMC)=(-\PMC)\#\PMC$.

Of course, $\HB_{sg}(\PMC)$ can be obtained from the split handlebody of genus $g$ by doing a sequence of arc-slides. Suppose $\PMC$ is the pointed matched circle of genus $1$. Then we can get $\HB_{sg}(\PMC)$ from the split genus $2$ handlebody by performing $8$ handleslides. Numbering the points along the boundary from $1$ to $8$, we perform the handleslides: 5 over 4; 2 over 1; 3 over 2; 4 over 3; 2 over 1; 6 over 5; 7 over 6; 2 over 3. See Figure~\ref{fig:self-gluing-hb}.
In Section~\ref{sec:find-slides}, we explain how to find this sequence
of handleslides.

We can use a computer to calculate $\CFDa(\PMC)$, for instance by computing inductively:
\[
\CFDa(\HB_i)=\Mor(\CFDDa(-\psi_i),\CFDa(\HB_{i-1})),
\]
where $\psi_i$ is the $i\th$ arc-slide, $-\psi_i$ denotes the same map
but between orientation-reversed surfaces, $\HB_0$ is the
$0$-framed, split handlebody, and $\HB_i=\psi_i(\HB_{i-1})$.
This computation leads to type $D$ structures with the following number of generators before and after simplification (cancelling arrows labeled by idempotents):
\begin{center}
  \begin{tabular}{rcc}
    \toprule
    &Gens.\ before& Gens.\ after\\
    Diagram & simplifying & simplifying\\
    \midrule
    $\CFDa(\HB_1)$ & 34 & 2\\
     $\CFDa(\HB_2)$ & 56 & 2\\
    $\CFDa(\HB_3)$ & 57 & 1\\
     $\CFDa(\HB_4)$ & 31 & 3\\
    \bottomrule
  \end{tabular}
  \qquad
  \begin{tabular}{rcc}
    \toprule
    &Gens.\ before& Gens.\ after\\
    Diagram & simplifying & simplifying\\
    \midrule
     $\CFDa(\HB_5)$ & 33 & 1\\
     $\CFDa(\HB_6)$ & 50 & 2\\
     $\CFDa(\HB_7)$ & 130 & 2\\
     $\CFDa(\HB_8)$ & 134 & 4\\
    \bottomrule
  \end{tabular}
\end{center}
(These computations were done using Sage~\cite{sage}. Code for doing such
computations, together with documentation including this example, is available from \newline
\verb|http://math.columbia.edu/~lipshitz/research.html#Programming| .)

The result has a more conceptual description, via a mild
generalization of Theorem~\ref{thm:DDforIdentity}. First, let $\gamma$
be the arc in $(-\PMC)\#\PMC$ running from $-\PMC$ to $\PMC$.
There is a map $\Alg((-\PMC)\#\PMC)\to\Alg(-\PMC)\otimes\Alg(\PMC)$ gotten by setting to zero any basis element crossing $a$.
Call an element $a\in\Alg((-\PMC)\#\PMC)$
\begin{itemize}
\item \emph{symmetric of type I} if $a$ does not cross $\gamma$, and the image of $a$ in $\Alg(-\PMC)\otimes\Alg(\PMC)$ is a chord-like element of the diagonal subalgebra
\item \emph{symmetric of type II} if $a$ consists of a single chord $\xi$ from $4k-i$ to $4k+i$, where $k$ is the genus of $\PMC$.
\end{itemize}
Also, there is an inclusion map $\Alg(-\PMC)\otimes\Alg(\PMC)\to\Alg((-\PMC)\#\PMC)$.
\begin{theorem}\label{thm:sghb}
  Let $\HB_{sg}(\PMC)$ be the self-gluing handlebody of $\PMC$. Then $\CFDa(\HB_{sg}(\PMC))$ is generated by the images in $\Alg((-\PMC)\#\PMC)$ of all pairs of complementary idempotents. The differential is given by right multiplication by
\[
\sum_{a\text{ symmetric of type I or II}}a.
\]
\end{theorem}
\begin{proof}[Proof sketch.]
  This follows from a factorization lemma just like Lemma~\ref{lem:FactorDiagonalSubalgebra}, except that there is an additional type of chord: a chord covering the middle region of the diagram, going between $-\PMC$ and $\PMC$.
\end{proof}

\index{Poincar\'e homology sphere}%
The Poincar\'e homology sphere has a genus $1$, one boundary component open book decomposition with monodromy $(\tau_a\tau_b)^5$, where $\tau_a$ and $\tau_b$ are Dehn twists around a pair of dual curves in the punctured torus. Note that these Dehn twists can be viewed as arc-slides in the genus $1$ pointed matched circle, of point 2 over point 1 and point 3 over point 2, respectively.

\begin{figure}
  \centering
  \includegraphics{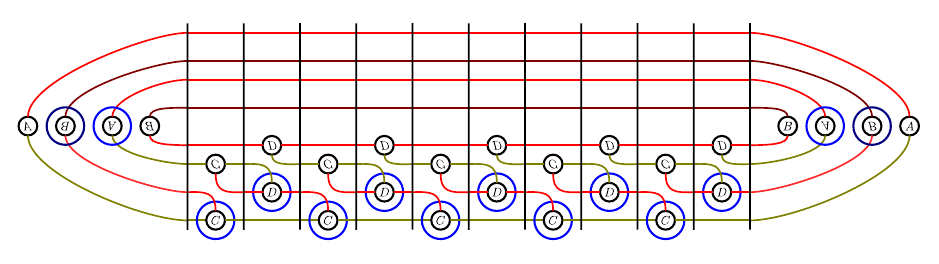}
  \caption{\textbf{Heegaard diagram for the Poincar\'e homology
      sphere.} This diagram is obtained by gluing together a diagram
    for $-\HB_{sg}$ and $\HB_{sg}$, the self-gluing handlebody and its
    reflection, with bordered Heegaard diagrams for ten arc-slides
    (Dehn twists) in between, and then destabilizing $30$ times.}
  \label{fig:for-Poincare}
\end{figure}

Extend these arc-slides by the identity map to arc-slides $A$ and $B$ of the genus $2$ split pointed matched circle. Then the Poincar\'e homology sphere is given by
\[
-\HB_{sg}(\PMC)\cup Y_{(AB)^5}\cup \HB_{sg}(\PMC)
\]
(see Figure~\ref{fig:for-Poincare}).

We can again compute this composition by computer. Let $\HB'_i=Y_{B(AB)^{i-1}}\cup\HB_{sg}$, $\HB''_i=Y_{(AB)^i}\cup\HB_{sg}$. Computing $\CFDa(\HB'_i)$ and $\CFDa(\HB''_i)$ inductively as
\begin{align*}
  \CFDa(\HB'_i)&=\Mor(\CFDDa(-Y_B),\CFDa(\HB''_{i-1}))\\
  \CFDa(\HB''_i)&=\Mor(\CFDDa(-Y_A),\CFDa(\HB'_i)),
\end{align*}
computer computation gives type $D$ structures with the following ranks before and after simplification:
\begin{center}
  \begin{tabular}{lcc}
    \toprule
   & Gens.\ before & Gens.\ after \\
   Diagram & simplifying & simplifying\\
   \midrule
   $\CFDa(\HB'_1)$ & 229 & 7\\
   $\CFDa(\HB'_2)$ & 250 & 6\\
   $\CFDa(\HB'_3)$ & 337 & 9\\
   $\CFDa(\HB'_4)$ & 445 & 11\\
   $\CFDa(\HB'_5)$ & 586 & 14\\
   \bottomrule
  \end{tabular}
  \qquad
  \begin{tabular}{lcc}
    \toprule
   & Gens.\ before & Gens.\ after \\
   Diagram & simplifying & simplifying\\
   \midrule
   $\CFDa(\HB''_1)$ & 317 & 5\\
   $\CFDa(\HB''_2)$ & 263 & 7\\
   $\CFDa(\HB''_3)$ & 374 & 10\\
   $\CFDa(\HB''_4)$ & 447 & 13\\
   $\CFDa(\HB''_5)$ & 567 & 15\\
   \bottomrule
  \end{tabular}
\end{center}

Finally, we have
\[
\CFa(Y)=\Mor(\CFDa(\HB_{sg}),\CFDa(\HB''_5))
\]
(where $Y$ is the Poincar\'e homology sphere). This is a complex with 405 generators, and $1$-dimensional homology.

\subsection{Finding handleslide sequences}
\label{sec:find-slides}

\begin{figure}
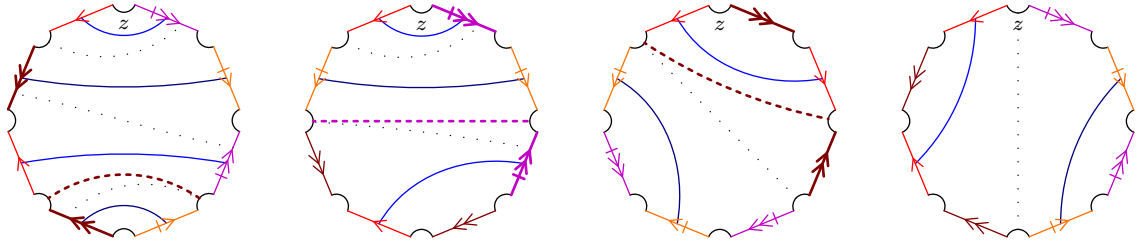

  \[
  \includegraphics{gluing-oct-0} \qquad
  \includegraphics{gluing-oct-1} \qquad
  \includegraphics{gluing-oct-2} \qquad
  \includegraphics{gluing-oct-3}
  \]
  \caption{\textbf{The source of the handleslides in
      Figure~\ref{fig:self-gluing-hb}.}
  We start with a Heegaard diagram for the self-gluing handlebody and
  end up with the split handlebody by a sequence of cuttings and
  re-gluings.  At each stage, we glue the edges drawn thick, and cut
  along the dashed line.  The dotted line is the desired separating curve.}
  \label{fig:explaining-slides}
\end{figure}

We now
explain how we found the sequence of handleslides used to construct
the self-gluing handlebody from the split handlebody, as in
Figure~\ref{fig:self-gluing-hb}. 
It is easiest to compute the
sequence of handleslides going in the other direction, starting from a
self-gluing handlebody and sliding until we obtain a split
handlebody.  The split handlebody has two distinguishing
characteristics:
\begin{itemize}
\item The two $\beta$-circles each intersect a single $\alpha$-arc once.
\item There is a separating loop $\gamma$ from the boundary to itself that does
  not intersect any $\alpha$-arcs.
\end{itemize}
We try to achieve these two features by successively performing
arc-slides among the $\alpha$-arcs, keeping the two $\beta$-circles
unchanged. To this end we choose also a separating loop $\gamma$. 
Each arc-slide is chosen to decrease the number of intersection points
between the $\alpha$-arcs, the $\beta$-circles, and $\gamma$.
This is illustrated in Figure~\ref{fig:explaining-slides}, where
the Heegaard
surface is shown split open along the $\alpha$-arcs,
rather than along the $\beta$-circles as in, for example,
Figure~\ref{fig:self-gluing-hb}. The curve $\gamma$ is indicated by a dotted arc.
In this representation, a sequence
of handleslides involving a single $\alpha$-arc $\alpha_i$ sliding over others
consists of
\begin{itemize}
\item gluing the sides of the diagram corresponding to $\alpha_i$ and
\item cutting open the resulting annulus along a new arc connecting
  the different sides of the annulus.
\end{itemize}
At each stage, we were able to choose the arc at the second step to reduce the
number of intersections of the $\alpha_i$ with the $\beta$-circles.


\bibliographystyle{hamsalpha}\bibliography{heegaardfloer}

\newcommand{\etalchar}[1]{$^{#1}$}
\providecommand{\noopsort}[1]{}
\providecommand{\bysame}{\leavevmode\hbox to3em{\hrulefill}\thinspace}
\providecommand{\href}[2]{#2}
\providecommand{\eprint}{\begingroup \urlstyle{rm}\Url}
\begin{thebibliography}{O{\relax Sz}04b}

\bibitem[ABP09]{AndersenBenePenner09:MCGroupoid}
J{\o}rgen~Ellegaard Andersen, Alex~James Bene, and R.~C. Penner, \emph{Groupoid
  extensions of mapping class representations for bordered surfaces}, Topology
  Appl. \textbf{156} (2009), no.~17, 2713--2725.

\bibitem[Ben10]{Bene08:ChordDiagrams}
Alex~James Bene, \emph{A chord diagrammatic presentation of the mapping class
  group of a once bordered surface}, Geom. Dedicata \textbf{144} (2010),
  171--190, \eprint{arXiv:0802.2747}.

\bibitem[Kel]{Keller:OtherAinfAlg}
Bernhard Keller, \emph{A brief introduction to {$A$}-infinity algebras},
  \url{http://people.math.jussieu.fr/~keller/publ/IntroAinfEdinb.pdf}.

\bibitem[LOT]{LOTCobordisms}
Robert Lipshitz, Peter~S. Ozsv{\'a}th, and Dylan~P. Thurston, \emph{Computing
  cobordism maps with bordered {F}loer homology}, in preparation.

\bibitem[LOT08]{LOT1}
\bysame, \emph{Bordered {H}eegaard {F}loer homology: {I}nvariance and pairing},
  2008, \eprint{arXiv:0810.0687v4}.

\bibitem[LOT11]{LOTHomPair}
Robert Lipshitz, Peter~S. Ozsv\'ath, and Dylan~P. Thurston, \emph{{H}eegaard
  {F}loer homology as morphism spaces}, Quantum Topology \textbf{2} (2011),
  no.~4, 384--449, \eprint{arXiv:1005.1248}.

\bibitem[LOT15]{LOT2}
Robert Lipshitz, Peter~S. Ozsv{\'a}th, and Dylan~P. Thurston, \emph{Bimodules
  in bordered {H}eegaard {F}loer homology}, 2015, \eprint{arXiv:1003.0598v3},
  To appear in {\em{Geometry \& Topology}}.

\bibitem[OS{\relax Sz}09]{OzsvathStipsiczSzabo:Nice}
Peter~S. Ozsv{\'a}th, Andr{\'a}s~I. Stipsicz, and Zolt{\'a}n {\relax
  Sz}ab{\'o}, \emph{Combinatorial {H}eegaard {F}loer homology and nice
  {H}eegaard diagrams}, 2009, \eprint{arXiv:0912.0830}.

\bibitem[O{\relax Sz}04a]{OS04:ThurstonNorm}
Peter~S. Ozsv{\'a}th and Zolt{\'a}n {\relax Sz}ab{\'o}, \emph{Holomorphic disks
  and genus bounds}, Geom. Topol. \textbf{8} (2004), 311--334.

\bibitem[O{\relax Sz}04b]{OS04:HolomorphicDisks}
\bysame, \emph{Holomorphic disks and topological invariants for closed
  three-manifolds}, Ann. of Math. (2) \textbf{159} (2004), no.~3, 1027--1158,
  \eprint{arXiv:math.SG/0101206}.

\bibitem[Ras03]{Rasmussen03:Knots}
Jacob Rasmussen, \emph{Floer homology and knot complements}, Ph.D. thesis,
  Harvard University, Cambridge, MA, 2003, \eprint{arXiv:math.GT/0306378}.

\bibitem[S{\etalchar{+}}13]{sage}
W.\thinspace{}A. Stein et~al., \emph{{S}age {M}athematics {S}oftware ({V}ersion
  5.9)}, The Sage Development Team, 2013, {\tt http://www.sagemath.org}.

\bibitem[SW10]{SarkarWang07:ComputingHFhat}
Sucharit Sarkar and Jiajun Wang, \emph{An algorithm for computing some
  {H}eegaard {F}loer homologies}, Ann. of Math. (2) \textbf{171} (2010), no.~2,
  1213--1236, \eprint{arXiv:math/0607777}.

\bibitem[Zar09]{Zarev09:BorSut}
Rumen Zarev, \emph{Bordered {F}loer homology for sutured manifolds}, 2009,
  \eprint{arXiv:0908.1106}.

\end{thebibliography}
\end{document}